\documentclass{ucthesis}
\usepackage{amssymb,amstext}
\usepackage{amsmath,amsbsy}
\usepackage{amsthm}
\usepackage{leftidx}
\usepackage[dvips]{graphicx}
\usepackage{mathrsfs}
\usepackage[color]{epsfig}
\usepackage{color}
\usepackage{epsf}
\usepackage{graphics}
\usepackage{psfrag}
\usepackage{epsf}
\newcommand*{\CopyCounter}[2]{%
  \expandafter\def\csname c@#2\endcsname{\csname c@#1\endcsname}%
  \expandafter\def\csname p@#2\endcsname{\csname p@#1\endcsname}%
  \expandafter\def\csname the#2\endcsname{\csname the#1\endcsname}}
\newcounter{Theorem}
\numberwithin{Theorem}{section}
\CopyCounter{Theorem}{Proposition}
\CopyCounter{Theorem}{Claim}
\CopyCounter{Theorem}{Lemma}
\CopyCounter{Theorem}{Corollary}
\CopyCounter{Theorem}{Conjecture}
\CopyCounter{Theorem}{Definition}
\CopyCounter{Theorem}{Example}
\CopyCounter{Theorem}{Remark}
\CopyCounter{Theorem}{Remarks}
\CopyCounter{Theorem}{Notation}
\theoremstyle{plain}
\newtheorem{thm}[Theorem]{Theorem}

\newtheorem{prop}[Proposition]{Proposition}

\newtheorem{claim}[Claim]{Claim}

\newtheorem{lemma}[Lemma]{Lemma}

\newtheorem{cor}[Corollary]{Corollary}

\theoremstyle{definition}
\newtheorem{definition}[Definition]{Definition}

\newtheorem{remark}[Remark]{Remark}
\newtheorem{remarks}[Remarks]{Remarks}
\newtheorem{notation}[Notation]{Notation}

\def\chaptername{Chapter}

\newcommand{\B}{\mathcal{B}}
\newcommand{\bigE}{\mathcal{E}}
\newcommand{\BofH}{\mathcal{B}(\mathcal{H})}
\newcommand{\br}[1]{\left({#1}\right)}

\newcommand{\cstar}{$\mathrm{C}\sp{*}$}
\newcommand{\CHI}{{\mbox{\large $\chi$}}}
\newcommand{\cj}{\mathcal{J}}
\newcommand{\curlyd}{\partial}
\newcommand{\dee}{\mathrm{d}}
\newcommand{\defeq}{\overset{\text{def}}{=}}
\newcommand{\EM}[2]{E^{M_{#1}}_{M_{#2}}}
\newcommand{\End}{\mathrm{End}}
\newcommand{\Finfinity}{{\mathbb F}_{\infty}}
\newcommand{\gothn}{\mathfrak{n}}
\newcommand{\gothm}{\mathfrak{m}}
\newcommand{\h}{\mathcal{H}}
\newcommand{\hatr}[1]{\left({#1}\right)\!\!\widehat{\phantom{A}}}
\newcommand{\hatt}{\widehat}
\newcommand{\Hom}{\mathrm{Hom}}
\newcommand{\id}{\mathrm{id}}
\newcommand{\ip}[1]{\left<{#1}\right>}
\newcommand{\isom}{\cong}
\newcommand{\II}{\mathrm{II}}
\newcommand{\IIone}{\mathrm{II}_1}
\newcommand{\IIinf}{\mathrm{II}_{\infty}}
\newcommand{\III}{\mathrm{III}}
\newcommand{\Ind}{\mathrm{Ind}}
\newcommand{\K}{\mathcal{K}}
\newcommand{\Ltwo}{{\mathrm L}^2}
\newcommand{\M}[2]{M_{#1}' \cap M_{#2}}
\newcommand{\Mi}[1]{M' \cap M_{#1} & \subset }
\newcommand{\modular}{\mathrm{mod}}

\newcommand{\Ni}[1]{N' \cap M_{#1} & \subset }
\newcommand{\ob}{\overline{b}}

\newcommand{\OB}{\overline{B}}

\newcommand{\ov}{\overline}
\newcommand{\PA}{\mathbb{P}}
\newcommand{\PAr}{\mathbb{P}^r}

\newcommand{\phat}{\widehat{\phantom{A}}}
\newcommand{\plushat}[1]{{#1}\hatt{_+}}   
\newcommand{\proofend}{\hfill $\Box$}
\newcommand{\Proof}{\noindent {\it Proof }}
\newcommand{\rig}{\mathrm{rigid}}
\newcommand{\rmspan}{\mathrm{span}}
\newcommand{\spdot}{\,\cdot\;}
\newcommand{\sph}{\mathrm{sph}}
\newcommand{\sqbr}[1]{\left[{#1}\right]}

\newcommand{\subf}[1]{\overset{#1}{\subset}}
\newcommand{\sumb}{\sum_{b \in B}}

\newcommand{\tensor}{\otimes}
\newcommand{\tN}{\underset{N}{\otimes}}
\newcommand{\tM}{\underset{M}{\otimes}}

\newcommand{\tlN}{\otimes_N}
\newcommand{\tlM}{\otimes_M}

\newcommand{\tensorlN}{\otimes_N}
\newcommand{\tensorlP}{\otimes_P}

\newcommand{\tensorN}{\underset{N}{\otimes}}
\newcommand{\tensorP}{\underset{P}{\otimes}}
\newcommand{\tensorM}{\underset{M}{\otimes}}
\newcommand{\tensorQ}{\underset{Q}{\otimes}}
\newcommand{\tensorSE}[1]{\underset{#1}{\otimes}}
\newcommand{\tensorS}[1]{\underset{#1}{\otimes}}

\newcommand{\trace}{\mathrm{tr}}
\newcommand{\Trace}{\mathrm{Tr}}

\newcommand{\wtilde}{\widetilde}

\degreeyear{2003}
\degreesemester{Spring}
\degree{Doctor of Philosophy}
\prevdegrees{
B.Sc.(Hons) \ (University of Canterbury) 1996
}

\chair{Professor Vaughan F.~R.~Jones}
\othermembers{
Professor Marc A.~Rieffel\\ 
Professor Martin B.~Halpern
}

\field{Mathematics}
\campus{BERKELEY}

\title{Subfactors, Planar Algebras and Rotations}

\author{Michael Burns}

\begin{document}

\maketitle
\copyrightpage
\abstract

Growing out of the initial connections between subfactors and knot
theory that gave rise to the Jones polynomial, Jones' axiomatization
of the standard invariant of an extremal finite index $\IIone$
subfactor as a spherical \cstar-planar algebra, presented
in~\cite{Jones1999}, is the most elegant and powerful description
available.

We make the natural extension of this axiomatization to the case of
finite index subfactors of arbitrary type.  We also provide the first
steps toward a limited planar structure in the infinite index case.
The central role of rotations, which provide the main non-trivial part
of the planar structure, is a recurring theme throughout this work.

In the finite index case the axioms of a \cstar-planar algebra need to
be weakened to disallow rotation of internal discs, giving rise to the
notion of a rigid \cstar-planar algebra.  We show that the standard
invariant of any finite index subfactor has a rigid \cstar-planar
algebra structure.  We then show that rotations can be re-introduced
with associated correction terms entirely controlled by the
Radon-Nikodym derivative of the two canonical states on the first
relative commutant, $N' \cap M$.

By deforming a rigid \cstar-planar algebra to obtain a spherical
\cstar-planar algebra and lifting the inverse construction to the
subfactor level we show that any rigid \cstar-planar algebra arises as
the standard invariant of a finite index $\IIone$ subfactor equipped
with a conditional expectation, which in general is not trace
preserving.  Jones' results thus extend completely to the general
finite index case.

We conclude by applying our machinery to the $\IIone$ case, shedding
new light on the rotations studied by Huang~\cite{Huang2000} and
touching briefly on the work of Popa~\cite{Popa2002}.

In the case of infinite index subfactors there are obstructions to
having a full planar algebra theory.  We constructing a periodic
rotation operator on the $\Ltwo$-spaces of the standard invariant of
an approximately extremal, infinite index $\IIone$ subfactor.  In the
finite index case we recover the usual rotation.  We also show that
the assumption of approximate extremality is necessary and sufficient
for rotations to exist on these $\Ltwo$-spaces.

The potential complexity of the standard invariant of an infinite
index subfactor is illustrated by the construction of a $\IIone$
subfactor with a type $\III$ central summand in the second relative
commutant, $N' \cap M_1$.  The restriction to $\Ltwo$-spaces does not
see this part of the standard invariant and Izumi, Longo and
Popa's~\cite{IzumiLongoPopa1998} examples of subfactors that are not
approximately extremal provide a further challenge to move beyond the
$\Ltwo$-spaces in the construction of rotation operators.  The present
construction is simply an initial step on the road to a planar
structure on the standard invariant of an infinite index subfactor.

\endabstract

\begin{frontmatter}

\tableofcontents

\begin{acknowledgements} 

I would like to thank my advisor, Vaughan Jones, a mathematician who
exemplifies the highest standards in every aspect of the profession,
from research and exposition to teaching and mentoring students.  His
enthusiasm for mathematics is infectious and his ability to produce
and expound elegant mathematics have served as an inspiration to me.
It has been a pleasure to work with him and I thank him for his
limitless patience and for having more faith in me than I sometimes
had in myself.

My colleagues and fellow students have been invaluable in their advice
and feedback.  In particular I would like to thank Dietmar Bisch, Zeph
Landau and Hsiang-Ping Huang, who were always eager to talk about
subfactors.  Michel Enock and Ryszard Nest's beautiful paper paved the
way for the results on infinite index subfactors.

To the many people in the Cal Hiking Club, CHAOS, who helped me avoid
working on my thesis I owe a debt of gratitude for keeping things in
perspective and for adding immeasurably the quality of my life.  In
particular Jeremy Schroeder and Mark Miller were steadfast partners
and great friends in many vertical adventures.

\end{acknowledgements}

\end{frontmatter}



\chapter{Introduction}

The study of subfactors was initiated by Jones' startling results on
the index for subfactors in~\cite{Jones1983}.  His work gave rise to a
powerful invariant of a subfactor known as the standard invariant.
This invariant has many equivalent descriptions, including Ocneanu's
paragroups, bimodule endomorphisms and 2-\cstar-tensor categories.
Popa's axiomatization of the standard invariant of a finite index
$II_1$ subfactor in terms of {\em standard
$\lambda$-lattices},~\cite{Popa1995}, was a major advance in the
field.  Jones'~\cite{Jones1999} planar algebra axiomatization for
extremal finite index $II_1$ subfactors builds on this to produce a
diagrammatic formulation in which the standard invariants ``seem to
have now found their most powerful and efficient formalism'' (quoting
Popa~\cite{Popa2002}).

The present work is concerned with extensions of the planar algebra
machinery to wider classes of subfactors than those considered in
Jones~\cite{Jones1999} and a recurring theme will be the properties of
rotation operators.  After these introductory remarks in chapter one,
the second chapter is concerned with extending Jones' subfactor-planar
algebra correspondence from extremal finite index $\IIone$ subfactors
to the general finite index case and proving some results with this
machinery.  Chapter three concerns infinite index subfactors and
defining rotations on their standard invariants as a step towards a
restricted 
planar structure on them.  There are obstructions to a full planar
algebra structure in the infinite index case.

Before embarking on a chapter by chapter summary we present a quick
overview of the area.  The reader is referred to Jones and
Sunder~\cite{JonesSunder1997} for basic material on subfactors, with
a focus on the finite index $\IIone$ case.

Let $N \subf{E} M$ be an inclusion of factors equipped with a normal
conditional expectation $E$ of finite index.  The Jones' basic
construction yields a factor $M_1$, generated by $M$ and the first
Jones' projection $e_1$, together with a normal conditional
expectation $E_M:M_1 \rightarrow M$.  Iterating this procedure we
obtain a tower of factors and conditional expectations, $N \subf{E} M
\subf{E_M} M_1 \subf{E_{M_1}} M_2 \subf{E_{M_2}}
\cdots$.  The standard invariant of $N \subf{E} M$ is the lattice of relative commutants obtained from this
tower:
\begin{equation*}
\begin{matrix}
\mathbb{C} = N' \cap N & \subset & N' \cap M & \subset & \Ni{1} & \Ni{2} 
	& \Ni{3} \cdots \\
 & & \cup & & \cup & & \cup & & \cup \\
 & & \mathbb{C} = M' \cap M & \subset & \Mi{1} & \Mi{2} & \Mi{3} \cdots
\end{matrix}
\end{equation*}
together with the conditional expectations.  The algebras $\M{i}{j}$
are all finite dimensional \cstar-algebras (multi-matrix algebras).
The standard invariant is a powerful invariant of the subfactor and it
is a complete invariant in the case of amenable $\II_1$ subfactors
(Popa~\cite{Popa1994}).

The case where $M$ (and hence $N$) is a $\IIone$ factor and $E$ is the
unique trace-preserving conditional expectation from $M$ onto $N$ is
the most extensively studied.  A special case of this is when the
subfactor is {\em extremal}, which is to say that the trace $\trace$
on $M$ and the unique trace $\trace'$ on $N'$ agree on $N' \cap M$.
In~\cite{Popa1995} Popa axiomatized the standard invariant of an
extremal $\IIone$ subfactor, proving that the standard invariant of an
extremal $\IIone$ subfactor forms an {\em extremal standard
$\lambda$-lattice} and conversely any extremal standard
$\lambda$-lattice arises in this way.  The general finite index
$\IIone$ case, a simple generalization of the proof in~\cite{Popa1995}
and known to Popa, first appears in print in~\cite{Popa2002}.

In~\cite{Jones1999} Jones characterizes the standard invariant of an
extremal $\IIone$ subfactor as a {\em spherical \cstar-planar
algebra}: loosely speaking, a sequence of finite dimensional vector
spaces $V_i$ with an action of the operad of (planar isotopy classes)
of planar tangles as multi-linear maps, consistent with composition of
tangles, equipped also with an involution $\phantom{l}^*$ and
satisfying certain positivity conditions and an additional spherical
isotopy invariance.

The proof that every spherical \cstar-planar algebra arises from an
extremal $\IIone$ subfactor is an application of Popa's standard
$\lambda$-lattice result.  The construction of the planar algebra of
an extremal $\IIone$ subfactor is the main result of
Jones~\cite{Jones1999} and a key ingredient is the periodicity of the
rotation operator which, after the result has been proved, can be
realized as the tangle below, illustrated in the case of $V_4=N' \cap
M_3$.

\newpage

\begin{figure}[htbp]
\begin{center}
\psfrag{*}{$*$}
\includegraphics{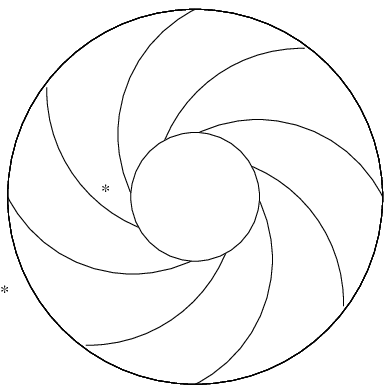}
\end{center}
\end{figure}

The planar algebra machinery has been a powerful tool in proving
results on the combinatorial structure of the standard invariant and
obtaining obstructions on the possible principal graphs and standard
invariants of subfactors.  See for example the work of Bisch and
Jones~\cite{BischJones1997, BischJones1997II, BischJones2000,
Jones2000, Jones2001}.

Extending the planar algebra formalism to the general $\IIone$ case
begins with proving the periodicity of the rotation operator.
In~\cite{Huang2000} Huang defines two rotation operators on the
standard invariant of a $\IIone$ subfactor and proves that each is
periodic.

\vspace{3mm}

The case of infinite index subfactors is first taken up by Herman and
Ocneanu in~\cite{HermanOcneanu1989}.  The results that they announced
were proved and expanded upon by Enock and Nest~\cite{EnockNest1996},
where the basic results for infinite index subfactors are laid down.
This paper also characterizes the subfactors arising as cross-products
by discrete type Kac algebras and compact type Kac algebras.

The tower and standard invariant can still be constructed in the case
of an infinite index inclusion, but the process becomes much more
technical.  Operator-valued weights must be used in place of
conditional expectations and the relative commutants need no longer be
finite dimensional.  Even in the simplest case of an infinite index
inclusion of $\IIone$ factors, only the odd Jones' projections exist,
implementing conditional expectations, and the other half of the maps
in the tower are operator-valued weights.  The standard invariant does
not form a planar algebra in this case.

\vspace{3mm}

In the following chapter we discuss the finite index case.  Background
material on general finite index subfactors is presented and a small
number of technical results developed before we turn our attention to
extending Jones' planar algebra results.  

Rigid planar algebras are defined in almost the same way as the planar
algebras of Jones except that we do not allow rotations of internal
discs.  After proving basic results about rigid \cstar-planar
algebras, including the central role of Radon-Nikodym derivatives, we
construct a rigid \cstar-planar algebra structure on the standard
invariant of any finite index subfactor.

Every rigid \cstar-planar algebra is shown to have a modular extension
in which discs can be rotated provided we insert correction terms
involving the Radon-Nikodym derivatives.  We then go on to deform the
modular extension so as to remove the correction terms, in essence by
incorporating them already in the action of the tangles.  We thus
obtain a spherical \cstar-planar algebra and thus an extremal $\IIone$
subfactor.

Finally we lift the inverse construction to the level of subfactors
and are able to prove that every rigid \cstar-planar algebra arises as
the standard invariant of a finite index subfactor, in fact a
subfactor of type $\IIone$, though the expectation may not be
trace-preserving.

We then apply the planar algebra machinery that we have developed to
the rotations studied by Huang~\cite{Huang2000} and illuminate some of
Popa's constructions in~\cite{Popa2002} in our context.

\vspace{3mm}

Chapter~\ref{chapter: infinite index} is concerned with infinite index
$\IIone$ subfactors.  After some initial material describing the basic
construction in this setting we take advantage of the additional
structure provided by the requirement that the first two factors in
the tower be of type $\IIone$.  We thus have half of the Jones
projections still at our disposal and are able to construct an
orthonormal basis for $M$ over $N$.  These tools allow us to extend
the notion of extremality to the infinite index case and show that it
has the usual properties.

Motivated by some of our work in the finite index case we can formally
define rotation operators on the $\Ltwo$-spaces of the standard
invariant.  The existence of these operators is then shown to be
equivalent to approximate extremality of the initial subfactor.

We conclude with an result indicating the potential complications that
arise once we move to finite index subfactors.  We construct of a
$\IIone$ subfactor with a type $\III$ central summand in the second
relative commutant, $N' \cap M_1$.

\chapter{Finite Index Subfactors of Arbitrary Type}
\label{general PA}


In this chapter we extend Jones' subfactor-planar algebra
correspondence in~\cite{Jones1999} from extremal finite index $\IIone$
subfactors to the general finite index case.  

We begin in sections \ref{sect: op-val weights} and~\ref{sect: Hilb
mod and rel tensor prod} with a review of the machinery of
operator-valued weights and relative tensor products as a prelude to
the discussion of index and the basic construction in
section~\ref{Kosaki background}.

A number of tedious, but necessary, computational results are proved
in Section~\ref{sect: tools}.  We connect the bimodule and relative
tensor product structure on the algebras $M_k$ with that on the
Hilbert modules $\Ltwo(M_k)$, before concluding with further results
on the multi-step basic construction.

Section~\ref{modular theory} is concerned with modular theory.  We
prove that relative commutants lie in the domains of the operators of
modular theory and that these operators leave the relative commutants
invariant.  We also show that the actions of these operators are
compatible with inclusions in the tower of higher relative commutants.

The rotation operator enters the picture in Section~\ref{section:
rotation} and we prove that it is quasi-periodic:
$
\left(\rho_k\right)^{k+1}=\Delta_k^{-1} .
$

\vspace{5pt}

In Section~\ref{sect: PA of fin index subfactor} we define a more
general notion of a planar algebra.  In a rigid planar algebra we
restrict attention to isotopies of tangles under which internal boxes
only undergo translations.  Every rigid \cstar-planar algebra has a
modular extension in which general planar isotopies are allowed, but
rotations of boxes change the action of the tangle.

There are two canonical states $\varphi$ and $\varphi'$ on a rigid
\cstar-planar algebra given by capping off boxes to the left or to the
right.  The Radon-Nikodym derivative $w$ or $\varphi'$ with respect to
$\varphi$ controls the effect of rotations in the following way.  Any
string along which the total angle changes after isotopy must be
modified by inserting a 1-box defined in terms of $w$ and the total
angle change.

In Section~\ref{section: subfactor to PA} we show that the standard
invariant of a finite index subfactor has a rigid planar algebra
structure.  Section~\ref{section: modular structure} describes the
construction of the modular extension to a rigid \cstar-planar algebra
and Section~\ref{section: rigid to spherical} contains the
construction of an associated spherical \cstar-planar algebra from a
rigid \cstar-planar algebra.  This recovers a result of Izumi that for
any finite index subfactor there is a $\IIone$ subfactor with the same
algebraic standard invariant.

We conclude in Section~\ref{section: PA to subfactors} by showing that
any rigid \cstar-planar algebra arises as the standard invariant of a
finite index subfactor.

\vspace{3mm}

Moving on to the specific case of a (not necessarily extremal) finite
index $\IIone$ subfactor, in Section~\ref{section: two rotations} we
establish the connection between Huang's two
rotations~\cite{Huang2000}.  We show that the two rotations are the
same if and only if the subfactor is extremal.  The method of proof
allows a very simple alternative proof of periodicity in the $\IIone$
case which will generalize to infinite index $\IIone$ subfactors in
Chapter~\ref{chapter: infinite index}.  We go on to prove some other
results for general finite index $\IIone$ subfactors using the planar
algebra machinery.  Section~\ref{section: 2-param family} sees a
two-parameter family of rotations defined on the standard invariant of
a finite index $\IIone$ subfactor, while Section~\ref{section:
additional results IIone} makes contact with the work of Popa
in~\cite{Popa2002}.


\section{General Finite Index Subfactors}
\label{sect: gen finite index subfactors}

While most of our work on finite index subfactors can proceed without
direct reference to the machinery of operator-valued weights and
relative tensor products, there are occasions when this material is
necessary.  We present a summary of technical results in
Sections~\ref{sect: op-val weights} and~\ref{sect: Hilb mod and rel
tensor prod}.  The reader who wishes to avoid this material can skip
these two sections and take as a starting point the results of Kosaki
quoted in Section~\ref{Kosaki background}.

In Chapter~\ref{chapter: infinite index} we will make heavy use of
operator-valued weights.


\subsection{Operator-valued weights}
\label{sect: op-val weights}

Here we summarize some key definitions and results from
Haagerup's foundational paper~\cite{Haagerup1979I}.

\begin{definition}
Let $M$ be a von Neumann algebra.  The extended positive part
$\hatt{M}_+$ of $M$ is the set of ``weights on the predual of $M$'',
namely $\hatt{M}_+$ is the set of maps $m: M_*^+\rightarrow [0,\infty
]$ such that $m$ is lower semi-continuous and
$m(\lambda\varphi+\mu\psi)=\lambda m(\varphi)+\mu m(\psi)$ for all
$\lambda,\mu\in [0,\infty ], \varphi,\psi\in M_*^+$.
\end{definition}

Note that $M_+$ embeds in $\hatt{M}_+$ by $x\mapsto m_x$ where
$m_x(\varphi)=\varphi(x)$.

Addition and positive scalar multiplication are defined on
$\hatt{M}_+$ in the obvious way.  For $a\in M$, $m\in \hatt{M}_+$
define $a^*ma$ by $(a^*ma)(\varphi)=m(\varphi(a^*\spdot a))$.  For $S
\subset \hatt{M}_+$ define $\sum_{m\in S} m$ pointwise.

\begin{prop}[Haagerup~\cite{Haagerup1979I} 1.2, 1.4, 1.5, 1.6, 1.9]
\label{prop: Haagerup ext pos part}
There are several alternative characterizations of $\hatt{M}_+$,
including:
\begin{description}
\item{(i)} Any pointwise limit of an
increasing sequence of bounded operators in $M_+$ is in $\hatt{M}_+$
and every element of $\hatt{M}_+$ arises this way.
\item{(ii)} $m\in \plushat{\B(\h)}$ is
in $\hatt{M}_+$ iff it is affiliated with $M$ ($u^*mu=m$ for all
unitary elements $u\in M'$).
\item{(iii)} Let $M$ be
represented on a Hilbert space $H$.  Let $p\in M$ and let $A$ be a
positive self-adjoint operator (possibly unbounded) on $p\h$
affiliated with $M$.  Define $m\in \plushat{\B(\h)}$ by
\[
m(\omega_\xi)=\begin{cases}
||A^{1/2}\xi||^2 & \xi \in D(A^{1/2}) \\
\infty		 & \text{otherwise}
\end{cases}
\]
where $\omega_\xi=\ip{\spdot\xi,\xi}$.  Then $m\in\hatt{M}_+$.  Every
element of $\hatt{M}_+$ arises this way.
\item{(iv)} Every $m\in\hatt{M}_+$ has a unique spectral resolution
\begin{align*}
m(\varphi)=\int_0^\infty \lambda \dee\varphi(e_\lambda) +\infty\varphi(p)
& \varphi\in M_*^+
\end{align*}
where $\{ e_\lambda\}_{\lambda\in[0,\infty)}$ is an increasing family
of projections in $M$, strongly continuous from the right and with
$p=1-\lim e_\lambda$.  In addition $e_0=0$ iff $m(\varphi)>0$ for all
nonzero $\varphi\in M_*^+$ and $p=0$ iff
$\{\varphi:m(\varphi)<\infty\}$ is dense in $M_*^+$.
\end{description}
\end{prop}

\begin{prop}[Haagerup~\cite{Haagerup1979I} 1.10]
\label{prop: Haagerup 1.10}
Any normal weight $\varphi$ on $M$ has a unique extension (also
denoted $\varphi$) to $\hatt{M}_+$ such that: (i) $\varphi(\lambda
m+\mu n)=\lambda \varphi(m)+\mu \varphi(n)$ for all $\lambda,\mu\in
[0,\infty ], m,n\in \hatt{M}_+$ and (ii) if $m_i \nearrow m$ then
$\varphi(m_i)\nearrow \varphi(m)$.
\end{prop}

\begin{definition}
Let $N \subset M$ be von Neumann algebras.  An {\em operator-valued
weight} from $M$ to $N$ is $T:M_+\rightarrow \hatt{N}_+$ satisfying
\begin{description}
\item{1.} $T(\lambda x+\mu y)=\lambda T(x)+\mu T(y)$ for all
$\lambda,\mu\in [0,\infty ], x,y\in M_+$.
\item{2.} $T(a^*xa)=a^*T(x)a$ for all $a\in N$, $x\in M_+$.
\end{description}
$T$ is normal if
\begin{description}
\item{3.} $x_i \nearrow x$ implies $T(x_i) \nearrow T(x)$.
\end{description}
\end{definition}

\begin{remarks}
\begin{description}
\item
\item{$\bullet$} A normal operator-valued weight $T:M_+\rightarrow \hatt{N}_+$
has a unique extension $T:\hatt{M}_+ \rightarrow \hatt{N}_+$ also
satisfying 1, 2 and 3 above.
\item{$\bullet$} $T$ is a conditional expectation iff $T(1)=1$.
\end{description}
\end{remarks}

\begin{definition}
As for ordinary weights define
\begin{align*}
\gothn_T &= \{ x\in M : ||T(x^*x)||<\infty \} \\
\gothm_T &= \gothn_T^* \gothn_T
          = \rmspan \{ x^*y | x,y \in \gothn_T \}
\end{align*}
[In the case when $T$ is a trace $\Trace$ these are the
Hilbert-Schmidt and Trace Class operators respectively.]
\end{definition}
Note that $\gothn_T$ is a left ideal, $\gothn_T$ and $\gothm_T$ are
$N-N$ bimodules and $T$ has a unique extension to a map
$\gothm_T\rightarrow N$.  For $x \in \gothm_T$, $a,b\in N$,
$T(axb)=aT(x)b$.

\begin{definition}
$T$ is {\em faithful} if $T(x^*x)=0$ implies $x=0$.  $T$ is {\em
semifinite} if $\gothn_T$ is $\sigma$-weakly dense in $M$.
n.f.s. will be used to denote ``normal faithful semifinite''.
\end{definition}

\begin{prop}[Haagerup~\cite{Haagerup1979I} 2.3]
Let $T$ be an operator-valued weight from $M$ to $N$ and let $\varphi$
be a weight on $N$.  If $T$ and $\varphi$ are normal (resp. n.f.s.)
then $\varphi\circ T$ is normal (resp. n.f.s).
\end{prop}

\begin{thm}[Haagerup~\cite{Haagerup1979I} 2.7]
\label{thm: ! n.f.s tr-pres op-val weight}
Let $N \subset M$ be semifinite von Neumann algebras with traces
$\Trace_N$ and $\Trace_M$.  Then there exists a unique
n.f.s. operator-valued weight $T:M_+\rightarrow \hatt{N}_+$ such that
$\Trace_M=\Trace_N\circ T$ ($\Trace_N$ on the right side of the
equality denotes the extension of $\Trace_N$ to $\hatt{N}_+$).
\end{thm}

\begin{remark}
In the proof of Theorem~\ref{thm: ! n.f.s tr-pres op-val weight}
Haagerup shows that, for $x \in M_+$, $T(x)$ is the unique element of
$\hatt{N}_+$ such that
\begin{align*}
&\Trace_M(y^{1/2}xy^{1/2})=\Trace_N(y^{1/2}T(x)y^{1/2})
& \text{for all } y \in M_+ .
\end{align*}
\end{remark}


\subsection{Hilbert $A$-modules and the relative tensor product}
\label{sect: Hilb mod and rel tensor prod}

The material in this section on general Hilbert $A$-modules is taken
from Sauvageot~\cite{Sauvageot1983} which also draws on
Connes~\cite{Connes1980}.  We make considerable use of the
Tomita-Takesaki Theory (see Takesaki~\cite{Takesaki2002II})

\begin{definition}
Let $A$ be a von Neumann algebra.  A left-Hilbert-$A$-module
$\leftidx{_A}{\K}{}$ is a nondegenerate normal representation of $A$
on a Hilbert space $\K$, while a right-Hilbert-$A$-module $\h_A$ is a
left-Hilbert-$A^{\text{op}}$-module.

Let $\varphi$ be an n.f.s. weight on $A$.  Given a
left-Hilbert-$A$-module $\leftidx{_A}{\K}{}$, the set
$D(\leftidx{_A}{\K}{},\varphi)$ (also denoted $D(\K,\varphi)$ or
$D(\leftidx{_A}{\K}{})$) of right-bounded vectors consists of those
vectors $\xi\in\K$ such that $\hatt{a} \mapsto a\xi$ extends to a
bounded operator $R(\xi)=R^{\varphi}(\xi):\Ltwo(A,\varphi)\rightarrow
\K$.  For a right-Hilbert-$A$-module ${\h}_A$, the set
$D({\h}_A,\varphi)=D({\h}_A)$ of left-bounded vectors consists of
those vectors $\xi\in\h$ such that $J\hatt{a} \mapsto \xi a^*$ extends
to a bounded linear operator
$L(\xi)=L^{\varphi}(\xi):\Ltwo(A,\varphi)\rightarrow \h$.
\end{definition}

\begin{remark}
The left-bounded vectors in $\h_A$ can equivalently be defined as
$D(\h,\varphi^{\text{op}})$, the set of right-bounded vectors for $\h$
considered as an $A^{\text{op}}$-module by defining $\pi(a)\xi=\xi a$.
In other words by requiring the boundedness of the maps
$R^{\varphi^{\text{op}}}(\xi):\Ltwo(A^{\text{op}},\varphi^{\text{op}})\rightarrow\h$
defined by $\hatt{a}\mapsto \pi(a) \xi=\xi a$.
\end{remark}

\subsubsection{Relative Tensor Product}

Given a left-Hilbert-$A$-module $\leftidx{_A}{\K}{}$ and right-bounded
vectors $\eta_1,\eta_2\in D(\leftidx{_A}{\K}{},\varphi)$ note that
$R(\eta_1)^*R(\eta_2)\in A' \cap \B(\Ltwo(A,\varphi))=JAJ$ and so defines an
element of $A$
\[
\leftidx{_A}{\ip{\eta_1,\eta_2}}{} = JR(\eta_1)^*R(\eta_2)J .
\]
Similarly, given a right-Hilbert-$A$-module $\h_A$ and left-bounded
vectors $\xi_1,\xi_2\in D(\h_A,\varphi)$ note that $\ip{\xi_1,\xi_2}_A \defeq
L(\xi_1)^*L(\xi_2)\in (JAJ)'=A$.

Although both of these pairings satisfy
$\ip{\zeta_1,\zeta_2}^*=\ip{\zeta_2,\zeta_1}$ and $\ip{\zeta,\zeta}\geq
0$, in general they are not $A$-valued inner products in the regular
sense as they are not $A$-linear in either component.  However, we
have the following result.

\begin{lemma}
Recall the definition of the modular automorphism group,
$\sigma^{\varphi}_t=\Delta_{\varphi}^t \spdot \Delta_{\varphi}^{-t}$
$(t\in\mathbb{R})$. $D(\leftidx{_A}{\K}{},\varphi)$ is stable under
elements of $D(\sigma^{\varphi}_{i/2})$ and
\[
\leftidx{_A}{\ip{a\eta_1,\eta_2}}{}
=\sigma^{\varphi}_{i/2}(a)\leftidx{_A}{\ip{\eta_1,\eta_2}}{} .
\]
$D(\h_A,\varphi)$ is stable under elements of
$D(\sigma^{\varphi}_{-i/2})$ and
\[
\ip{\xi_1,\xi_2 a}_A
=\ip{\xi_1,\xi_2}_A \sigma^{\varphi}_{i/2}(a) .
\]
\end{lemma}

\begin{lemma}[Sauvageot~\cite{Sauvageot1983} 1.5]
\begin{description}
\item
\item{(i)}
$\varphi\br{\leftidx{_A}{\ip{\eta_1,\eta_2}}{}}=\ip{\eta_1,\eta_2}$,
$\varphi\br{\ip{\xi_1,\xi_2}_A}=\ip{\xi_2,\xi_1}$.
\item{(ii)} $\leftidx{_A}{\ip{\eta_1,\eta_2}}{}, \ip{\xi_1,\xi_2}_A
\in \gothn_{\varphi}$ and
\begin{align*}
&\hatr{\leftidx{_A}{\ip{\eta_1,\eta_2}}{}}=JR(\eta_1)^*\eta_2,
&\hatr{\ip{\xi_1,\xi_2}_A}=L(\xi_2)^*\xi_1 .
\end{align*}
\item{(iii)} $\ip{\ip{\xi_2,\xi_1}_A\eta_1,\eta_2}
=\ip{\xi_2,\xi_1\leftidx{_A}{\ip{\eta_1,\eta_2}}{}}$.
\end{description}
\end{lemma}

\begin{definition}[Relative tensor product]
Given Hilbert $A$-modules $\h_A$ and $\leftidx{_A}{\mathcal{K}}{}$
define the relative tensor product $\h \tensor_{\varphi} \K$
(sometimes denoted $\h \tensor_A \K$ when the choice of $\varphi$ is
clear) to be the Hilbert space completion of the algebraic tensor
product $D(\h_A) \odot D(\leftidx{_A}{\mathcal{K}}{})$ equipped with
the inner product
\begin{align}
\ip{\xi_1 \odot \eta_1,\xi_2 \odot \eta_2}
&=\ip{\ip{\xi_2,\xi_1}_A\eta_1,\eta_2} \label{eq: rel tensor 1} \\
&=\ip{\xi_2,\xi_1\leftidx{_A}{\ip{\eta_1,\eta_2}}{}}
  \label{eq: rel tensor 2}
\end{align}
[First quotient by the space of length-zero vectors, then complete].
The image of $\xi \odot \eta$ in $\h \tensor_A \mathcal{K}$ is denoted
$\xi \tensor_A \eta$ or $\xi \tensor_{\varphi} \eta$.  If $\h$ is a $B$-$A$ bimodule and $\mathcal{K}$
an $A$-$C$ bimodule then $\h \tensor_A \mathcal{K}$ is naturally a
$B$-$C$ bimodule.
\end{definition}

\begin{remark}
$\h \tensor_A \mathcal{K}$ is also the completion of $D(\h_A) \odot_A
\mathcal{K}$ using (\ref{eq: rel tensor 1}) or the completion of $\h
\odot_A D(\leftidx{_A}{\mathcal{K}}{})$ using (\ref{eq: rel tensor
2}).
The relative tensor product is not $A$-middle-linear, but we have the
following result.
\end{remark}

\begin{lemma}
\label{lemma: almost middle linear}
For $\xi\in\h, \eta\in D(\K,\varphi), a\in D(\sigma^{\varphi}_{-i/2})$ we have
\[
\xi a \tensor_{\varphi} \eta 
= \xi \tensor_{\varphi} \sigma^{\varphi}_{-i/2}(a) \eta .
\]
\end{lemma}

\begin{notation}
\label{notation: left tensor}
\begin{description}
\item
\item{(1)} For $\eta_1,\eta_2\in D(\leftidx{_A}{\K}{},\varphi)$ define
\[
\theta^{\varphi}\br{\eta_1,\eta_2}
=R(\eta_1)R(\eta_2)^*
\in A' \cap \B(\K) .
\]

\item{(2)} For $\xi \in D(\h_A)$ let $L_\xi:\mathcal{K} \rightarrow \h
\tensor_A \mathcal{K}$ denote the map
$
L_\xi : \eta \rightarrow \xi \tensorSE{A} \eta .
$

For $\eta \in D(\leftidx{_A}{\mathcal{K}}{})$ let $R_\eta:\h
\rightarrow \h \tensor_A \mathcal{K}$ denote the map
$
R_\eta : \xi \rightarrow \xi \tensorSE{A} \eta .
$

By (\ref{eq: rel tensor 1}) and (\ref{eq: rel tensor 2}) $L_\xi$ and
$R_\eta$ are bounded and $L_\xi^*L_\xi=\ip{\xi,\xi}_A$, $R_\eta^*R_\eta=\leftidx{_A}{\ip{\xi,\xi}}{}$.
\end{description}
\end{notation}

\begin{remark}
In the case of a semifinite von Neumann algebra $A$ with trace
$\Trace$ we have $J\hatt{a}=\hatt{a^*}$ so $D(\h_A)=\{ \xi\in\h :
\hatt{a} \mapsto \xi a \text{ is bounded}\}$ and $L(\xi)\hatt{a}=\xi
a$.  The modular automorphisms are trivial and hence
$\ip{\spdot,\spdot}_A$ is right-$A$-linear and
$\leftidx{_A}{\ip{\spdot,\spdot}}{}$ is left-$A$-linear.
\end{remark}


\subsection{Background material from Kosaki}
\label{Kosaki background}

Let $N \subf{E} M$, sometimes denoted $(N \subset M,E)$, be an
inclusion of ($\sigma$-finite) factors with a normal conditional
expectation $E:M\rightarrow N$.  We recall some material from
Kosaki~\cite{Kosaki1986, Kosaki1998} on index and the basic
construction in this general setting. 


\subsubsection{Index}

Let $P(M,N)$ denote the set of all normal faithful semi-finite
(n.f.s.) operator-valued weights from $M$ to $N$.  Let $M$ be
represented on $\mathcal{H}$.  By Haagerup~\cite{Haagerup1979I} there
is a bijection between $P(M,N)$ and $P(N',M')$.  In~\cite{Kosaki1986}
Kosaki constructs this bijection in a canonical way, so that there is
a unique n.f.s. operator valued weight $E^{-1}:(N')_+ \rightarrow
\widehat{(M')}_+$ such that
\[
\frac{\dee(\varphi\circ E)}{\dee\psi}=\frac{\dee\varphi}{\dee(\psi\circ E^{-1})}
\]
for all n.f.s. weights $\varphi$ on $N$ and all n.f.s weights $\psi$
on $M'$, where the derivatives are Connes' spatial derivatives
(\cite{Connes1980}).  We have the following alternative
characterization of $E^{-1}$

\begin{lemma}[Kosaki~\cite{Kosaki1998}~3.4]
\label{lemma: Kosaki 1998 3.4}
Whenever both sides are defined, 
\[
E^{-1}\br{\theta^{\varphi}\br{\xi,\xi}}
=\theta^{\varphi\circ E}\br{\xi,\xi}.
\]

\end{lemma}

The index of $E$ is defined to be $\Ind(E)=E^{-1}(1)\in
\plushat{Z(M)} =[0,\infty]$ and is independent of the Hilbert
space $\mathcal{H}$ on which $M$ is represented.  We will use $\tau$
to denote $\Ind(E)^{-1}$ (note that Kosaki uses $\lambda$ rather than
$\tau$).


\subsubsection{Basic Construction}

We will assume henceforth that the index is finite, in which case
$E'=\tau E^{-1}:N' \rightarrow M'$ is a normal conditional
expectation.

The basic construction is performed as follows.  Take any faithful
normal state $\varphi$ on $N$ and extend it to $M$ by $\varphi\circ
E$.  Let $\mathcal{H}=\Ltwo(M,\varphi)$ and denote the inclusion map
of $M$ into $\Ltwo(M,\varphi)$ by either $x \mapsto \hatt{x}$ or $x
\mapsto \Lambda(x)$.  The inner product on $\Ltwo(M,\varphi)$ is
$\ip{\hatt{x},\hatt{y}}=\varphi(y^* x)$ for $x, y \in M$.  Define the
Jones' projection $e_1$ by
\[
e_1 \hatt{x} = \widehat{E(x)} .
\]
$e_1$ extends to a projection in $N'$ and one defines $M_1$ to be the
von Neumann algebra $<M, e_1 >$ generated by $M$ and $e_1$.  The usual
properties are satisfied

\begin{prop}[Kosaki~\cite{Kosaki1986} Lemma 3.2]
\label{prop: basic properties of basic constr}
\begin{description}
\item
\item{(i)} $e_1 x e_1 = E(x) e_1$ for all $x \in M$.
\item{(ii)} $N = M \cap \{e_1\}'$.
\item{(iii)} $J e_1 J = e_1$ where $J=J_0=J_{\varphi}$.
\item{(iv)} $M_1=J N' J$.
\item{(v)} $M_1=\rmspan \br{M e_1 M} = \rmspan \{a e_1 b : a, b \in M \}$.
\end{description}
\end{prop}
\noindent There is a canonical conditional expectation $E_M=E_M:M_1\rightarrow
M$ given by
\[
E_M(x)=J_0 E'(J_0 x J_0) J_0
\]
and one has $\Ind(E_M)=\Ind(E)=\tau^{-1}$ and $E_M(e_1)=\tau$.  In
addition, we have

\begin{lemma}[Pull-down Lemma]
\label{lemma: fin pull-down}

$M_1 e_1 = M e_1$.  For $z\in M_1$ we have $ze_1=xe_1$ where
$x=\tau^{-1} E_M(xe_1)$
\end{lemma}

Iterating the construction as usual we obtain a sequence of Jones'
projections $\{ e_i \}_{i\geq 1}$ and a tower of factors
\[
N \subset M \subset M_1 \subset M_2 \subset \ldots
\]
The state $\varphi$ is extended to the entire tower via $\varphi \circ
E \circ E_M \circ E_{M_1} \circ \cdots E_{M_k}$ and will simply be
denoted $\varphi$.  We will use $\hatt{\phantom{A}}$ or $\Lambda_k$ to
denote the inclusion of $M_k$ in $\Ltwo(M_k,\varphi)$ and $\pi_k$ to
denote the representation of $M_k$ on $\Ltwo(M_k,\varphi)$ by left
multiplication.  Reference to $\Lambda_k$ will be suppressed when it
is clear that we are considering an element of $M_k$ in
$\Ltwo(M_k,\varphi)$.  Reference to $\pi_k$ will often be suppressed
when the representation is clear.

We have the following additional properties:

\begin{prop}
\label{prop: basic properties of ei's}
\begin{description}
\item
\item{(i)} $E_{M_i}(e_{i+1})=\tau$;
\item{(ii)} $e_i e_{i\pm 1} e_i = \tau e_i$ and $[e_i,e_j]=0$ for
$|i-j| \geq 2$;
\item{(iii)} $e_i$ is in the centralizer of $\varphi$ on $M_j$ for all
$j \geq i$ (i.e. $\varphi(e_i a) = \varphi(a e_i)$ for all $a\in
M_j$);
\item{(iv)} hence $\varphi$ is a trace on $\{ 1,e_1, e_2, \ldots , e_j
\}''$ which forces the usual restrictions on the value of the index
originally found in Jones~\cite{Jones1983}.
\end{description}
\end{prop}

\begin{notation}
Following Jones~\cite{Jones1999} let
$\delta=\tau^{-1/2}=(\Ind(E))^{1/2}$, $E_k=\delta e_k$, \\ $v_k=E_k
E_{k-1} \ldots E_1$.  Note that:
\begin{align*}
E_k^2 & = \delta E_k ,\\
E_k E_{k\pm 1} E_k & = E_k ,\\
v_kv_k^* & = \delta E_k ,\\
v_k^*v_k & = \delta E_1 ,\\
v_kxv_k^* & = \delta E_N(x)E_k && \text{for } x \in M .
\end{align*}
\end{notation}


\subsubsection{The Multi-step Basic Construction}
\label{multi-step bc}

Use $\EM{k}{j}$ to denote the conditional expectation $E_{M_j}
E_{M_{j+1}} \cdots E_{M_{k-1}}: M_k \rightarrow M_j$.  When $k$ is
clear from context we will sometimes abuse notation and use $E_{M_j}$
to denote $\EM{k}{j}$.  We have the following result, originally
proved in the $\IIone$ case by Pimsner and Popa~\cite{PimsnerPopa1988}
and found as 4.3.6 of Jones and Sunder~\cite{JonesSunder1997}.  The
proof only involves properties of the Jones projections found in
Propositions~\ref{prop: basic properties of basic constr}
and~\ref{prop: basic properties of ei's} and thus holds in the general
finite index case.

\begin{thm}
\label{multi-step bc thm}
For all $-1 \leq i < j < k = 2j-i$ let $m=j-i$ and let $f_{[i,j]}$ be
the Jones projection for $(M_i \subset M_j, \EM{j}{i})$.  Let
$F_{[i,j]}=\Ind\br{\EM{j}{i}}^{1/2} f_{[i,j]} =\delta^m f_{[i,j]}$ and
let $<M_j,f_{[i,j]}>$ be the factor resulting from the basic
construction.  Then there is a unique isomorphism
$\pi^k_j:<M_j,f_{[i,j]}> \rightarrow M_k$ which is the identity on
$M_j$ and such that
\begin{align*}
\pi^k_j(f_{[i,j]})
	& = \delta^{m(m-1)} (e_{j+1}e_j \cdots e_{i+2})
	(e_{j+2} \cdots e_{i+3}) \cdots (e_k \cdots e_{j+1}) \\
	& \defeq e_{[i,j]}
\end{align*}
or, equivalently,
\begin{align*}
\pi^k_j(F_{[i,j]})
	& = (E_{j+1}E_j \cdots E_{i+2})
	(E_{j+2} \cdots E_{i+3}) \cdots
	(E_k \cdots E_{j+1}) \\
	& \defeq E_{[i,j]} .
\end{align*}
\end{thm}

Note that we could express this theorem as saying that for $j < k \leq
2j+1$ there is a unique representation $\pi^k_j$ of $M_k$ on
$\Ltwo(M_j,\varphi)$ such that $M_j$ acts as left multiplication and
$e_{[i,j]}$ acts as expectation onto $M_i$ ($i=2j-k$).  In fact we can
be more explicit in our description of $\pi_j^k$.

\begin{prop}[Bisch~~\cite{Bisch1995} 2.2]
\label{prop: fin multi-step}
For $z \in M_k$, $y \in M_j$ we have
\[
\pi^{k}_j(z) \hatt{y}
	= \tau^{j-k} \EM{k}{j} (zy e_{[i,j]})
	= \delta^{k-j} \EM{k}{j} (zy E_{[i,j]}) .
\]
\begin{proof}
For $k=2j+1$ the same proof as in Prop 2.2 of Bisch~\cite{Bisch1995}
yields
\[
\pi^{2j+1}_j(z) \hatt{y}
	= \tau^{-j-1} \EM{2j+1}{j} (zy e_{[-1,j]})
	= \delta^{j+1} \EM{2j+1}{j} (zy E_{[-1,j]}) .
\]
Applying this to the tower obtained from $M_{2j-k}\subset M_{2j-k+1}$
we get
\[
\pi^{k}_j(z) \hatt{y}
	= \tau^{j-k} \EM{k}{j} (zy e_{[i,j]})
	= \delta^{k-j} \EM{k}{j} (zy E_{[i,j]}) 
\]
for $z \in M_k$, $y \in M_j$. 
\end{proof}
\end{prop}

For more on the multi-step basic construction see Section~\ref{sect:
tools}, Proposition~\ref{rep M_k on M_j}, and 
Proposition~\ref{prop: EN and PP multi-steps are the same}.


\subsubsection{Local index; Finite dimensional relative commutants}

As in the $\IIone$ case there is a local index formula which implies
the finite dimensionality of the relative commutants arising from a
finite index subfactor.  Given $p \in N' \cap M$ define $E_p: pMp
\rightarrow Np$ by
\[
E_p(x)=E(p)^{-1}E(x)p .
\]
Then one has
\begin{itemize}
\item If $\sigma_t^E(p)=p$ for all $t \in \mathbb R$ then
\[
\Ind(E_p)=E(p) E^{-1}(p) = \Ind(E) E(p) E'(p)
\]
\item In general
\[
\Ind(E_p) \leq E(p) E^{-1}(p) = \Ind(E) E(p) E'(p)
\]
\end{itemize}
This proves the finite dimensionality of $N' \cap M$ exactly as in
Jones and Sunder~\cite{JonesSunder1997}~2.3.12 (originally in
Jones~\cite{Jones1983}).


\subsubsection{Basis; relative tensor product}
\label{basis}

As in the type $\IIone$ case there exists a (right-module) basis for
$M$ over $N$.  That is, there exists a finite set $B = \{ b_i \}_{i
\in I} \subset M$ such that $\sum_{b \in B} b e_1 b^* =1$.  In fact
there exists an orthonormal basis, one in which $E_N(b^*
\tilde{b})=\delta_{b,\tilde{b}} q_b$ where $q_b$ are projections in $N$.

It is worth noting that in the type $\III$ case this basis can be
chosen to have one element $u$ with $u e_1 u^*=1$ and $E_N(u^*u)=1$.

Also following the $\IIone$ case we have $M_{i+1}\isom M_i
\otimes_{M_{i-1}} M_i$ via $x E_{i+1} y \mapsto x \otimes_{M_{i-1}}
y$, where $x, y \in M_i$.  Hence there exists an isomorphism
$\theta=\theta_k: \otimes^{k+1}_N M \rightarrow M_k$ given by
\begin{eqnarray}
\theta\left(x_1 \tN x_2 \tN \cdots \tN x_{k+1}\right)
 & = & x_1 v_1 x_2 v_2 \cdots v_k x_{k+1} \label{theta1} \\
 & = & x_1 v_k^* x_2 v_{k-1}^* \cdots v_1^* x_{k+1} . \label{theta2} 
\end{eqnarray}

Note that as in Jones and Sunder~\cite{JonesSunder1997}~4.3.4 we have:

\newpage  

\begin{lemma}
\label{basislemma}
\begin{description}
\item
\item{(i)} If $B$ is a basis for $N \subf{E} M$ and $\tilde{B}$ is a 
basis for $M \subf{\tilde{E}} P$, then $\tilde{B}B=\{ \tilde{b}b :
b\in B, \tilde{b} \in \tilde{B}\}$ is a basis for $N
\subf{E\circ\tilde{E}} P$.
\item{(ii)} $Bv_i^*=\{ b v_i^* : b \in B \}$ is a basis for $M_i$ 
over $M_{i-1}$
\item{(iii)} $B_k=\{\theta(b_{i_1}\tN \cdots \tN b_{i_{k+1}}) 
: i_j \in I \}$ is a basis for $M_k$ over $N$.
\end{description}
\begin{proof}
\begin{description}
\item
\item{(i)} Simply note that for all $x \in P$
\[
x=\sum_{\tilde{b}} \tilde{b}\tilde{E}(\tilde{b}^*x) 
	=\sum_{\tilde{b}} \sum_b \tilde{b}bE(b^*\tilde{E}(\tilde{b}^*x)) 
	=\sum_{\tilde{b},b} \tilde{b}bE(\tilde{E}(b^*\tilde{b}^*x)) .
\]
\item{(ii)} $\sum bv_i^* e_{i+1} v_i b^*
 = \sum \delta^{-2} b v_{i+1}^* v_{i+1} b^*
 = \sum b e_1 b^*
 = 1$.
\item{(iii)} Using (ii) and iterating (i) we obtain the basis
$B v_k^* B v_{k-1}^* \cdots B v_1 B = B_k$.
\end{description}
\end{proof}
\end{lemma}


Finally, the basis can be used to implement the conditional
expectation from $N'$ onto $M'$.  This result is proved in the
$\IIone$ case in Bisch~\cite{Bisch1995}~2.7.

\begin{prop}
\label{E_M' via basis}
Let $B$ be a basis for $M$ over $N$.  Then $E':N' \rightarrow M'$ is
given by
\[
E'(x)=\tau \sum_{b \in B} b x b^* .
\]
\begin{proof}
Let $\Phi(x)=\sum b x b^*$.  It is equivalent to show that
$E^{-1}(x)=\Phi(x)$.  Note that if $\xi\in D(\h,\varphi)$ then for
$x\in M$,
\begin{align*}
||x\xi||
&=\left|\left|\sum_b bE(b^*x)\xi \right|\right| \\
&\leq \sup||b|| \sum_b ||E(b^*x)||_2 ||R^{\varphi}(\xi)|| \\
&\leq \sup||b|| \br{\sum_b \varphi\br{E(x^*b)E(b^*x)}}^{1/2}
  K^{1/2} ||R^{\varphi}(\xi)|| \\
&= \sup||b|| \varphi\br{E(x^*x)}^{1/2} K^{1/2} ||R^{\varphi}(\xi)||
\end{align*}
where $K$ is the cardinality of $B$.  Hence $R^{\varphi\circ E}(\xi)$
is bounded and so $\xi\in D(\h,\varphi\circ E)$.

By Connes~\cite{Connes1980}~Prop~3, $\rmspan
\{\theta^{\varphi}\br{\xi,\xi} : \xi\in D(\h,\varphi) \}$ is weakly
dense in $N'$.  As both $E^{-1}$ and $\Phi$ are weakly continuous, it
suffices to show that $\Phi\br{\theta^{\varphi}\br{\xi,\xi}}
=E^{-1}\br{\theta^{\varphi}\br{\xi,\xi}}$.

For $a\in M$ define $L(a):\Ltwo(N) \rightarrow \Ltwo(M)$ by
$L(a)\hatt{x}=\hatt{ax}$.  Then
$L(a)L(a)^*=ae_1a^*:\Ltwo(M)\rightarrow\Ltwo(M)$.  Also note that for
$x\in N$, $bR^{\varphi}(\xi)\hatt{x}=bx\xi=R^{\varphi\circ
E}(\xi)L(b)\hatt{x}$, so $bR^{\varphi}(\xi)=R^{\varphi\circ
E}(\xi)L(b)$.  Hence
\begin{align*}
\Phi\br{\theta^{\varphi}\br{\xi,\xi}} 
&= \sum_b bR^{\varphi}(\xi)R^{\varphi}(\xi)^*b^* \\
&= \sum_b R^{\varphi\circ E}(\xi)L(b)L(b)^*R^{\varphi\circ E}(\xi)^* \\
&= R^{\varphi\circ E}(\xi)R^{\varphi\circ E}(\xi)^* \\
&= \theta^{\varphi\circ E}\br{\xi,\xi} \\
&= E^{-1}\br{\theta^{\varphi}\br{\xi,\xi}}
\end{align*}
by Lemma~\ref{lemma: Kosaki 1998 3.4}.  Hence $\Phi=E^{-1}$.
\end{proof}
\end{prop}


\subsection{Computational tools}
\label{sect: tools}

Here we discuss the relationship between the bimodule structure and
relative tensor products of the algebras $M_k$ and those of the
Hilbert modules $\Ltwo(M_k)$.  We then look at the conditional
expectation in terms of the isomorphism theta taking $M_{k-1}$ to the
$k$-fold algebraic relative tensor product of $M$ over $N$ .  Finally
we take another look at the multi-step basic construction.

\begin{prop}
\label{prop: L^2 tensor and alg tensor}
There is an isomorphism of bimodules $\Ltwo(M_k)\isom \tlN^{k+1}
\Ltwo(M)$ given by $u_k:\theta\hatr{x_1 \tlN \cdots \tlN
x_{k+1}} \mapsto \hatt{x_1} \tlN \cdots \tlN \hatt{x_{k+1}}$.
\end{prop}

\begin{remark}
This is not immediately obvious.  The bimodule structure on
$\Ltwo(M_k)$ does not restrict to that on the algebra $M_k$.  The
right action of $N$ on $\Ltwo(M_k)$ is $J_k n^* J_k$, {\em not} right
multiplication $S_kn^*S_k$ which may be unbounded.  However, by
Lemma~\ref{lemma: almost middle linear},
\[
\hatt{x} \tN n\hatt{y}
= \hatt{x}\cdot \sigma_{i/2}(n) \tN \hatt{y}
= J(\Delta^{-1/2}n\Delta^{1/2})^*J\hatt{x} \tN \hatt{y}
= S n^* S\hatt{x} \tN \hatt{y}
=\hatt{xn} \tN \hatt{y} .
\]
\end{remark}

\begin{lemma}
$\Ltwo(M_1)\isom \Ltwo(M) \tlN \Ltwo(M)$ via $\theta\hatr{x
\tlN y}=\hatr{xE_1 y} \mapsto \hatt{x} \tN \hatt{y}$.
\end{lemma}

\begin{remark}
Note that, although $y$ is not necessarily in
$D(\leftidx{_N}{\Ltwo(M)}{})$, $x$ is in $D(\Ltwo(M)_N)$ because
$L(\hatt{x}):J_{-1}\hatt{n} \mapsto \hatt{x}\cdot n^* = J_0 n J_0
\hatt{x}$ is bounded.  To see this note that $J_0|_{\Ltwo(N)}=J_{-1}$
and
\begin{equation*}
J_0 n J_0\hatt{x}
=J_0 n J_0 x J_0 \Omega
=J_0 J_0 x J_0 n \Omega
=x J_0 \hatt{n}
=x J_{-1} \hatt{n} .
\end{equation*}
\end{remark}

\Proof {\it of Lemma.}

Both sides have dense span, so it suffices to check that the map
preserves the inner product:
\begin{align*}
\ip{\hatt{x'}\tN\hatt{y'},\hatt{x}\tN\hatt{y}}
&= \ip{\ip{\hatt{x},\hatt{x'}}_N\hatt{y'},\hatt{y}}
 = \ip{E_N(x^*x')\hatt{y'},\hatt{y}} \\
&= \varphi\br{y^*E_N(x^*x')y'}
 = \delta^2 \varphi\br{y^*E_N(x^*x')e_1y'} \\
&= \delta^2 \varphi\br{y^*e_1x^*x'e_1y'}
 = \ip{\hatr{x'E_1y'},\hatr{xE_1y}} .
\end{align*}


\proofend

\Proof {\it of Prop~\ref{prop: L^2 tensor and alg tensor}.}

The proposition is true for $k=0,1$.  Suppose the result is true for
some $k \geq 1$.  Applying this to $(M\subset M_1, E_M)$ we have
$\Ltwo(M_{k+1})\isom \tlM^{k+1} \Ltwo(M_1)$ via
\[
\hatr{A_1E_2A_2E_3E_2\cdots A_kE_{k+1}E_k\cdots E_2 A_{k+1}}
\mapsto
\hatt{A_1} \tM \cdots \tM \hatt{A_{k+1}} .
\]
Note that $\Ltwo(M)\tlM\Ltwo(M)\isom\Ltwo(M)$ via
$m:\hatt{x}\tlM \hatt{y}\mapsto \hatt{xy}$, so
$\tensor_M^{k+1}\Ltwo(M_1) \isom \tlN^{k+2}\Ltwo(M)$ via
$V=\br{\id \tlN \br{\tensor_N^k m} \tlN \id}\circ
\br{\tensor_M^{k+1} u_1}$.

Let $A_1=x_1E_1x_2, A_i=E_ix_{i+1}$, $2\leq i\leq k+1$.  Then
\[
\hatr{A_1E_2A_2E_3E_2\cdots A_kE_{k+1}E_k\cdots E_2 A_{k+1}}
=\theta\br{x_1 \tN \cdots \tN x_{k+2}}\phat ,
\]
and
\[
V\br{\hatt{A_1} \tM \cdots \tM \hatt{A_{k+1}}}
=V\br{\hatt{x_1E_1x_2} \tM \hatt{1E_1x_3} \tM \cdots \tM 
  \hatt{1E_1x_{k+2}}}
=\hatt{x_1} \tN \cdots \tN \hatt{x_{k+2}} .
\]

\proofend

\begin{prop}
\label{prop: properties of tensorN}
Let $(N \subset M,E)$ be a finite index subfactor.  Let $a_i\in M$
$(i\geq 0)$.Then
\[
E_{M_{k-1}}\br{ \theta\br{a_1 \tN \cdots \tN a_{k+1}} } \\
= \begin{cases}
   \delta^{-1}\theta\br{a_1 \tN \cdots \tN a_r a_{r+1} \tN 
    \cdots \tN a_{k+1}}
   &k=2r-1\\
   \theta\br{a_1 \tN \cdots \tN a_r E(a_{r+1}) \tN \cdots 
    \tN a_{k+1}}
   &k=2r\\
  \end{cases}
\]

\begin{proof}
Note that if $X_i=y_iE_1z_i$, then
\begin{align*}
\theta\br{X_1 \tM \cdots \tM X_k}
&= X_1 E_2 X_2 E_3E_2 \cdots X_{k-1} E_k \cdots E_3E_2 X_k \\
&= y_1 E_1 z_1y_2E_2E_1z_2y_3E_3E_2E_1\cdots z_{k-1}y_k E_k\cdots E_2E_1 z_k\\
&= \theta\br{y_1 \tN z_1y_2 \tN \cdots \tN z_{k-1}y_k \tN
   z_k} .
\end{align*}
In particular
\[
\theta\br{A_1 \tM \cdots \tM A_{k-1} \tM \overline{A}_k}
=\theta\br{a_1 \tN \cdots \tN a_{k+1}} ,
\]
where $A_i=a_iE_1$ and $\overline{A}_k=a_kE_1a_{k+1}$.

\vspace{3mm}

\noindent For $k=0$, $E_{M_{-1}}(a_1)=E(a_1)$.  Assume the result
holds for some $k\geq 0$.  Note that
$A_rA_{r+1}=a_rE_1a_{r+1}E_1=\delta a_rE(a_{r+1})E_1$ and
$A_rE_M((A_{r+1}))=a_r E_1 E_M(a_{r+1}E_1)=\delta^{-1} a_rE_1a_{r+1}$.
Hence, with the first of the two cases denoting $k=2r-1$
and the second $k=2r$,
\begin{align*}
E_{M_k}\br{\theta\br{a_1 \tN \cdots \tN a_{k+2}}} 
&= E_{M_k}\br{\theta\br{A_1 \tM \cdots \tM A_k \tM 
  \overline{A}_{k+1}}} \\
&= \begin{cases}
 \delta^{-1}\theta\br{A_1 \tM \cdots \tM A_rA_{r+1} \tM \cdots 
  \tM A_k \tM \overline{A}_{k+1}} \\
 \theta\br{A_1 \tM \cdots \tM A_r E_M(A_{r+1}) \tM \cdots 
  \tM A_k \tM \overline{A}_{k+1}}
 \end{cases} \\
&= \begin{cases}
 \theta\br{a_1 \tN \cdots \tN a_{r-1} \tN a_rE(a_{r+1}) 
  \tN a_{r+2} \tN \cdots \tN a_{k+2}} \\
 \delta^{-1}\theta\br{a_1 \tN \cdots \tN a_r \tN a_{r+1}
  a_{r+2}\tN a_{r+3} \tN \cdots a_{k+2}}
 \end{cases} .
\end{align*}
\end{proof}
\end{prop}

\begin{remarks}
We could have proved many other properties of $\theta$ with almost
identical arguments to those in Prop~\ref{prop: tensorP}, but the result
above is all that we require here.

We conclude with some further results on the multi-step basic
construction of Theorem~\ref{multi-step bc thm}.  We first clarify a
certain compatibility of the representations $\pi^k_j$.
\end{remarks}

\begin{prop}
\label{rep M_k on M_j}
Let $j \leq k \leq 2j$ and let $z \in M_k$.  Then
$\pi^k_j(z)=\pi^{k+1}_j(z)$.

\begin{proof}
Using the explicit formula for $\pi^{k+1}_j(z)$ from
Proposition~\ref{prop: fin multi-step} this is basically just a long
exercise in simplifying words in the $E_i$'s.  Here are the details.

Without loss of generality we may assume that $k=2j$ (just use
$M_{2j-k-1}\subset M_{2j-k}$ in place of $N\subset M$).  For $r \geq s$
let $V_{r,s}=E_r E_{r-1} \cdots E_s$ and for $r<s$ let $V_{r,s}=1$.
Note that
\[
E_{[a,b]}=V_{b+1,a+2} V_{b+2,a+3} \cdots V_{2b-a,b+1} .
\]
We will make two very simple observations and then prove the lemma.
First note that for $a \geq d > b$, $c \geq d+2$, $d \geq e$ we have
\begin{equation}
\label{eq: simp1}
V_{a,b} V_{c,d+2} V_{d,e} = V_{a,d} V_{d-2,b} V_{c,d+2} V_{d-1,e} .
\end{equation}
This may be paraphrased as follows.  Think of $V_{c,d+2}V_{d,e}$ as
$V_{c,e}$ with a missing term, or {\em gap}, at $d+1$.  We could
similarly talk or larger gaps with more terms missing.  Then the
equation above says that a gap in one $V$ propagates to the left into
the previous $V$ and leaves a bigger gap in the original $V$.  More
succinctly we could say ``gaps propagate left leaving bigger gaps''.
The proof is quite simple:
\begin{align*}
V_{a,b} V_{c,d+2} V_{d,e}
&= V_{a,b} V_{c,d+2} E_d V_{d-1,e} \\
&= V_{a,b} E_d V_{c,d+2} V_{d-1,e} \\
&= \br{E_aE_{a-1}\cdots E_dE_{d-1}E_dE_{d-2}E_{d-3}\cdots E_b}
    V_{c,d+2} V_{d-1,e} \\
&= \br{E_aE_{a-1}\cdots E_dE_{d-2}E_{d-3} \cdots E_b} V_{c,d+2} V_{d-1,e} \\
&= V_{a,d} V_{d-2,b} V_{c,d+2} V_{d-1,e} .
\end{align*}
Second note that for $a \geq b \geq d$, $c \geq b+2$ we have
\begin{equation}
\label{eq: simp2}
V_{a,b} V_{c,b+2} V_{b-1,d} = V_{b,d} V_{c,b+2} .
\end{equation}
which follows directly from the fact that the second and third terms
on the left commute.

With these preliminary results we can now begin the main proof.
\begin{align*}
\delta E_{M_{2j}}\br{E_{[-1,j]}}
&= \delta E_{M_{2j}}\br{V_{j+1,1} V_{j+2,2} \cdots V_{2j+1,j+1}} \\
&= V_{j+1,1} V_{j+2,2} \cdots V_{2j,j} \delta E_{M_{2j}}\br{V_{2j+1,j+1}} \\
&=V_{j+1,1}V_{j+2,2}\cdots V_{2j,j}\delta E_{M_{2j}}\br{E_{2j+1}}V_{2j,j+1}\\
&= \sqbr{V_{j+1,1} V_{j+2,2} \cdots V_{2j,j}} V_{2j,j+1} .
\end{align*}
Iterating (\ref{eq: simp1}) from right to left
\begin{eqnarray*}
\delta E_{M_{2j}}(E_{[-1,j]})
&=& V_{j+1,1}V_{j+2,2}\cdots V_{2j-2,j-2}V_{2j-1,j-1}[V_{2j,j} V_{2j,j+1}] \\
&=& V_{j+1,1} V_{j+2,2} \cdots V_{2j-2,j-2} V_{2j-1,j-1} [V_{2j,2j} 
     V_{2j-2,j} V_{2j-1,j+1}] \\
&=& V_{j+1,1} V_{j+2,2} \cdots V_{2j-2,j-2} [V_{2j-1,j-1} V_{2j,2j} 
     V_{2j-2,j}] V_{2j-1,j+1} \\
&=& V_{j+1,1} V_{j+2,2} \cdots V_{2j-2,j-2} [V_{2j-1,2j-2} V_{2j-4,j-1} 
     V_{2j,2j} V_{2j-3,j}] V_{2j-1,j+1} \\
&=& V_{j+1,1} V_{j+2,2} \cdots [V_{2j-2,j-2} V_{2j-1,2j-2} V_{2j-4,j-1}] 
     V_{2j,2j} V_{2j-3,j} V_{2j-1,j+1} \\
& & \vdots \\
&=& (V_{j+1,2}) (V_{j+2,4}) (V_{j+3,6} V_{3,3}) (V_{j+4,8} V_{5,4}) \cdots \\
& & \cdots (V_{2j-2,2j-4} V_{2j-7,j-2}) (V_{2j-1,2j-2} V_{2j-5,j-1}) 
     (V_{2j,2j} V_{2j-3,j}) V_{2j-1,j+1} .
\end{eqnarray*}
Iterating (\ref{eq: simp2}) from left to right,
\begin{eqnarray*}
\delta E_{M_{2j}}(E_{[-1,j]})
&=& V_{j+1,2} [V_{j+2,4} V_{j+3,6} V_{3,3}] (V_{j+4,8} V_{5,4}) \cdots \\
& & \cdots (V_{2j-2,2j-4} V_{2j-7,j-2}) (V_{2j-1,2j-2} V_{2j-5,j-1}) (V_{2j,2j} V_{2j-3,j}) V_{2j-1,j+1} \\
&=& V_{j+1,2} [V_{j+2,3} V_{j+3,6}] (V_{j+4,8} V_{5,4}) \cdots \\
& & \cdots (V_{2j-2,2j-4} V_{2j-7,j-2}) (V_{2j-1,2j-2} V_{2j-5,j-1}) (V_{2j,2j} V_{2j-3,j}) V_{2j-1,j+1} \\
&=& V_{j+1,2} V_{j+2,3} [V_{j+3,6} V_{j+4,8} V_{5,4}] \cdots \\
& & (V_{2j-2,2j-4} V_{2j-7,j-2}) (V_{2j-1,2j-2} V_{2j-5,j-1}) (V_{2j,2j} V_{2j-3,j}) V_{2j-1,j+1} \\
& & \vdots \\
&=& V_{j+1,2} V_{j+2,3} V_{j+3,4} \cdots V_{2j-1,j} V_{2j,j+1} \\
&=& E_{[0,j]}
\end{eqnarray*}
Hence, for $y \in M_j$,
\begin{align*}
\pi^{2j+1}_j(z) \hatt{y}
&= \delta^{j+1} \EM{2j+1}{j} (zy E_{[-1,j]})
 = \delta^{j} \EM{2j}{j} (zy \delta E_{M_{2j}}(E_{[-1,j]})) \\
&= \delta^{j} \EM{2j}{j} (zy E_{[0,j]})
 = \pi^{2j}_j(z) \hatt{y}
\end{align*}
\end{proof}
\end{prop}

\begin{notation}
From lemma \ref{rep M_k on M_j} we see that if $x \in M_k$, $k \leq
2j+1$, then $\pi^k_j(x)=\pi^l_j(x)$ for all $k \leq l \leq 2j+1$ and
hence we will use $\pi_j$ to denote this representation, with no reference
to the algebra $M_k$ that is acting.
\end{notation}


\begin{prop}
\label{prop: EN and PP multi-steps are the same}
For $R\in M_{2j+1}$, $\pi_{j+t}(R)=\pi_j(R)
\tlN \br{\id_{\Ltwo(M)}}^{\tlN^t}$.
\begin{proof}
It suffices to prove the result for $t=1$ and then iterate.  Note that
$E_{[-1,j]}v_j=v_{j+1}E_{[1,j+1]}$.  Simply observe
\begin{align*}
E_{[-1,j]}v_j
&= V_{j+1,1}V_{j+2,2}V_{j+3,3}\cdots V_{2j,j}V_{2j+1,j+1} 
    E_jE_{j-1}\cdots E_1 \\
&= V_{j+1,2}E_1V_{j+2,3}E_2V_{j+3,4}E_3\cdots V_{2j,j+1}E_jV_{2j+1,j+2} 
    E_{j+1}E_jE_{j-1}\cdots E_1 \\
&= V_{j+1,2}V_{j+2,3}V_{j+3,4}\cdots V_{2j,j+1}V_{2j+1,j+2} 
    E_1E_2 \cdots E_jE_{j+1}E_jE_{j-1}\cdots E_1 \\
&= V_{j+1,2} E_{[1,j+1]} E_1 \\
&= v_{j+1} E_{[1,j+1]} .
\end{align*}
Now 
\begin{align*}
\hspace{50pt} & \hspace{-50pt}
\br{\pi_j(R) \tN \id} \hatr{x_1 \tN \cdots \tN x_{j+2}} \\
&= \delta^{j+1} E^{M_{2j+1}}_{M_j}\br{R x_1v_1x_2v_2\cdots v_jx_{j+1}
    E_{[-1,j]}} v_{j+1}x_{j+2} \\
&= \delta^{j+2} E_{M_j}\br{E^{M_{2j+1}}_{M_{j+1}}\br{R x_1v_1x_2v_2\cdots 
     v_jx_{j+1}E_{[-1,j]}}e_{j+1}} e_{j+1} v_j x_{j+2} \\
&= \delta^{j+2} \tau E^{M_{2j+1}}_{M_{j+1}}\br{R x_1v_1x_2v_2\cdots 
     v_jx_{j+1}E_{[-1,j]}} e_{j+1} v_j x_{j+2} \\
&= \delta^j E^{M_{2j+1}}_{M_{j+1}}\br{R x_1v_1x_2v_2\cdots 
     v_jx_{j+1}E_{[-1,j]}v_jx_{j+2}} \\
&= \delta^j E^{M_{2j+1}}_{M_{j+1}}\br{R x_1v_1x_2v_2\cdots 
     v_jx_{j+1} v_{j+1}x_{j+2}E_{[1,j+1]}} \\
&= \pi_{j+1}(R) \hatr{x_1 \tN \cdots \tN x_{j+2}} ,
\end{align*}
where we use the Pull-down Lemma (Lemma~\ref{lemma: fin pull-down}) to
obtain the fourth line and the initial observation that
$E_{[-1,j]}v_j=v_{j+1}E_{[1,j+1]}$ to obtain the sixth line.
\end{proof}
\end{prop}

\begin{cor}
$\pi^k_j$ could be alternatively defined as follows.  For $R\in M_k$,
 $j\leq k \leq 2j+1$, define $\pi^k_j(R)\in \B(\Ltwo(M_j))$ by
\[
\pi^k_j(R) \tN \br{\id_{\Ltwo(M)}}^{\tlN^{k-j-1}}
= \pi^k_{k-1}(R)
\]
where $\pi^k_{k-1}$ is the defining representation of $M_k$ on
$\Ltwo(M_{k-1})$.  This is the definition coming from the multi-step
basic construction as described in Enock and Nest~\cite{EnockNest1996}
(see Prop~\ref{prop: additional info basic constr} for details in the
infinite index $\IIone$ case).
\begin{proof}
$R\in M_{2j+1}$, so $\pi_{k-1}(R) =\pi_{j+(k-j-1)}(R) =\pi_j(R) \tN
\br{\id}^{\tlN^{k-j-1}}$.
\end{proof}
\end{cor}


\subsection{Modular theory}
\label{modular theory}

The construction of planar algebras will not require all of the
modular theory that we discuss here and it therefore possible to skip
the next two sections and go immediately to section~\ref{sect: PA of
fin index subfactor}.  However, the following theorem is necessary to
know that the modular operators on the finite dimensional relative
commutants that arise in section~\ref{sect: PA of fin index subfactor}
are in fact the restrictions of the modular operators on the spaces
$\Ltwo(M_k,\varphi)$.

Modular theory involves a number of unbounded operators.  However, the
relative commutants are invariant under the operators of the modular
theory and finite dimensionality implies that the restrictions of
these operators to the relative commutants are bounded.  We collect
these and other technical results below.

\begin{thm}
\label{basic modular results}
Let $S_k$ be the operator $\hatt{x}\mapsto \widehat{x^*}$ on
$\widehat{M_k}
\subset \Ltwo(M_k,\varphi)$, $\Delta_k=S_k^*\overline{S_k}$ the modular
operator.  Then
\begin{description}
\item{(i)} $\widehat{N' \cap M_k}$ is in the domain of $S_k$, $S_k^*$ and 
$\Delta_k^t$ for all $t \in \mathbb{C}$.  $\widehat{N' \cap M_k}$ is
invariant under all of these operators, and also invariant under
$J_k$.
\item{(ii)}
$\Delta_k(\hatt{x})=\hatr{\sum_b bE_N(xb^*)}$ and
$\Delta_k^{-1}(\hatt{x})=\hatr{\sum_b E_N(bx)b^*}$ where $B$ is any
basis for $M$ over $N$.
\item{(iii)} $S_{k+1}$, $S_{k+1}^*$, $J_{k+1}$ and $\Delta_{k+1}^t$
($t \in \mathbb{C}$) restrict to $S_k$, $S_k^*$, $J_k$ and
$\Delta_k^t$ respectively on $\widehat{N' \cap M_k}$, and hence will
be denoted $S$, $S^*$, $J$ and $\Delta^t$.
\item{(iv)} $S_k \pi_k(N' \cap M_{2k+1}) S_k = \pi_k(N' \cap M_{2k+1})$.
More precisely, for $z\in N' \cap M_{2k+1}$, $S_k \pi_k(z^*) S_k$ has
dense domain in $\Ltwo(M_k,\varphi)$ and extends to a bounded operator
on $\Ltwo(M_k,\varphi)$ given by $\pi_k(R_k(z))$ where $R_k(z) \in N'
\cap M_{2k+1}$ is given by
\[
R_k(z)=\tau^{-k-1} \sum_{\ob \in \OB} E_{M_k}(e \ob z)e \ob^*
\]
where $e=e_{[-1,k]}$ implements the conditional expectation
$E^{M_k}_N$ and $\OB$ is any basis for $M_k$ over $N$.

Consequently $S_k\cdot S_k$, $S_k^* \cdot S_k^*$ and $J_k
\cdot J_k$ are all conjugate-linear automorphisms of the vector space
$\pi_k(N' \cap M_{2k+1})$ and $\Delta_k^{-1/2} \cdot
\Delta_k^{1/2}$, $\Delta_k^{1/2} \cdot \Delta_k^{-1/2}$ are linear
automorphisms of $\pi_k(N' \cap M_{2k+1})$.


\item{(v)} 
In addition to (iv), $S_k\cdot S_k$, $S_k^* \cdot S_k^*$,
$\Delta_k^{-1/2} \cdot \Delta_k^{1/2}$, $\Delta_k^{1/2} \cdot
\Delta_k^{-1/2}$ and $J_k \cdot J_k$ all map $\pi_k(N' \cap M_k)$ 
onto $\pi_k(M_k' \cap M_{2k+1})$.

\end{description}
\begin{proof}
\begin{description}
\item
\item{(i)} Let $x \in N' \cap M_k$.  Then $S_k(\hatt{x})= 
\hatt{x^*} \in \widehat{N' \cap M_k}$.

\noindent For all $y \in M_k$
\begin{align*}
\ip{\hatt{x},\hatt{y^*}}
&= \varphi(yx) 
 = \sum_b \varphi(bE_N(b^*y)x) 
 = \sum_b \varphi(bxE_N(b^*y)) \\
&= \sum_b \varphi(E_N(bx)b^*y) 
 = \sum_b \ip{\hatt{y},\hatr{bE_N(x^*b^*)}} .
\end{align*}
Hence $\hatt{x} \in D(S_k^*)$ and $S_k^*(\hatt{x})=\hatr{\sum_b
bE_N(x^*b^*)}$ which is thus independent of the basis $B$.  Given any
$v\in \mathcal{U}(N)$, $vB$ is also a basis for $M_k$ over $N$ and
hence
\[
v\left(\sum_b bE_N(x^*b^*)\right)v^* = \sum_b (vb)E_N(x^*(b^*v^*))
		   = \sum_b bE_N(x^*b^*)
\]
so that $\sum_b bE_N(x^*b^*)\in N' \cap M_k$.  Consequently $\hatt{x}$
is in the domain of $\Delta_k=\overline{S_k}S_k^*$ and
$\Delta_k(\hatt{x})\in \widehat{N' \cap M_k}$.


The fact that $\widehat{N' \cap M_k}$ is in the domain of $\Delta_k^t$
and invariant under it is a basic exercise in spectral theory.  There
exists a measure space $(X,\mu)$, a positive function $F$ on $X$ and a
unitary operator $v:\Ltwo(X,\mu) \rightarrow \Ltwo(M_k,\varphi)$ such
that $v^* \Delta_k v=M_F$, the operator of multiplication by $F$.  As
$v^*(\widehat{N' \cap M_k})$ is finite dimensional and invariant under
$M_F$, it has a basis $\{ f_1, \ldots , f_n\}$ of eigenvectors of
$M_F$.  Thus there exist $\lambda_i > 0$ such that $Ff_i=\lambda_i
f_i$ and hence $F=\lambda_i$ a.e. on the support of $f_i$.  This
implies that $F^t=\lambda_i^t$ a.e. on the support of $f_i$ and hence
$v^*(\widehat{N' \cap M_k})=\rmspan\{ f_1, \ldots ,f_n \}$ is in the
domain of $\br{M_F}^t$ and invariant under it.

Lastly, $J_k=\overline{S_k} \Delta_k^{-1/2}$ so $N' \cap M_k$ is
invariant under $J_k$.

\item{(ii)} From above
\begin{eqnarray*}
\Delta_k(\hatt{x})	& = & S_k S_k^* \hatt{x}
			  = \hatr{\sum_b b E_N(xb^*)} \\
\Delta_k^{-1}(\hatt{x})	& = & S_k^* S_k \hatt{x}
			  = \hatr{\sum_b E_N(bx)b^*} .
\end{eqnarray*}

\item{(iii)}  Let $x \in N' \cap M_k$.  Obviously $S_k \hatt{x} =
\widehat{x^*} = S_{k+1} \hatt{x}$.  To see that $S_k^* \hatt{x} =
S_{k+1}^* \hatt{x}$ observe that for all $y \in M_{k+1}$
\begin{align*}
\ip{S_k^*\hatt{x},\hatt{y}}
= \ip{S_k^*\hatt{x},\widehat{E_{M_k}(y)}}
&=\ip{\widehat{E_{M_k}(y)^*},\hatt{x}} 
= \ip{\widehat{E_{M_k}(y^*)},\hatt{x}} \\
&=\ip{\widehat{y^*},\hatt{x}}
= \ip{S_{k+1}^*\hatt{x},\hatt{y}} .
\end{align*}
From the above $\Delta_{k+1}\hatt{x}=\Delta_k \hatt{x}$.  The spectral
theory argument in (i) implies that $\Delta_{k+1}^t \hatt{x} =
\Delta_k^t \hatt{x}$.  Finally $J_{k+1}\hatt{x}= S_{k+1}
\Delta_{k+1}^{-1/2} \hatt{x}=S_k \Delta_k^{-1/2} \hatt{x} =J_k
\hatt{x}$.

\item{(iv)} It suffices to prove this in the case $k=0$.  Let $b$ be a
basis for $M$ over $N$.  Since $Be_1$ is a basis for $M_1$ over $M$,
there exist $x_b \in M$ with $z=\sum_b be_1x_b$.  Then for all $y \in
M$ we have
\begin{equation}
\label{S_0 . S_0}
S_0 z^* S_0 \hatt{y} = \sum_b S_0 x_b^* \widehat{E_N(b^*y^*)}
		     = \sum_b \hatr{E_N(yb)x_b} ,
\end{equation}
which demonstrates that $S_0 z^* S_0$ has domain $\widehat{M}$ which
is of course dense in $\Ltwo(M,\varphi)$.

Next observe that
\begin{align*}
\sum_b E_M(e_1 b z) e_1 b^* \hatt{y}
&= \sum_b E_M(e_1 b z) \hatt{E_N(b^* y)} 
= \sum_b \hatr{E_M(e_1 b z E_N(b^* y))} \\
&= \sum_b \hatr{E_M(e_1 b E_N(b^* y) z)} 
= \hatt{E_M(e_1 y z)} \\
&= \sum_b \hatr{E_M(e_1 y b e_1 x_b)} 
= \sum_b \hatr{E_M(E_N(yb) e_1 x_b)} \\
&= \tau \sum_b \hatr{E_N(yb) x_b}
\end{align*}
so that the bounded operator $R_0(z)=\tau^{-1} \sum_b E_M(e_1 b z) e_1 b^*$
agrees with $S_0 z^* S_0$ on its domain.  $R_0(z) \in N'$
because~\ref{S_0 . S_0} is clearly $N$-linear in $y$.  Since
$S_0^2=1$, $z \mapsto S_0 z^* S_0$ is an automorphism of $N'
\cap M_1$.

Taking adjoints, using the fact that $J_0 (N' \cap M_1) J_0 = N' \cap
M_1$ in addition to $\overline{S_0}=J_0 \Delta_0^{1/2}=
\Delta_0^{-1/2} J_0$ and $S_0^*=J_0 \Delta_0^{-1/2}= \Delta_0^{1/2} J_0$
the second part of (iv) is reasonably clear.  However, a little care
must be taken with domains.

For $\xi \in \mathcal{D}(S_0^*)$ and $y \in \mathcal{D}(S_0)=
\widehat{M}$ we have
\[
\ip{S_0 \hatt{y}, z S_0^* \xi}
= \ip{\xi, S_0 z^* S_0 \hatt{y}}
=\ip{\xi, R_0(z) \hatt{y}}
=\ip{R_0(z)^* \xi, \hatt{y}}
\]
so that $z S_0^* \xi \in \mathcal{D}(S_0^*)$ and $S_0^* z S_0^* \xi =
R_0(z)^* \xi$.  Hence $\mathcal{D}(S_0^* z S_0^*)= \mathcal{D}(S_0^*)$
and $S_0^* z S_0^*$ extends to the bounded operator $R_0(z)^*$, which
is in $N' \cap M_1$.  Since $(S_0^*)^2=1$, $z \mapsto S_0^* z S_0^*$
is a conjugate-linear automorphism of $N' \cap M_1$.

For $\xi \in \mathcal{D}(\Delta_0^{1/2}) =\mathcal{D}(\overline{S_0})$
and $\eta \in \mathcal{D}(\Delta_0^{-1/2}) =\mathcal{D}(S_0^*)$ we
have
\begin{align*}
\ip{\Delta_0^{-1/2} \eta, z \Delta_0^{1/2} \xi}
&= \ip{J_0 S_0^* \eta, z J_0 \overline{S_0} \xi}
 = \ip{J_0 z J_0 \overline{S_0} \xi, S_0^* \eta} \\
&= \ip{\overline{S_0} \xi, (J_0 z^* J_0) S_0^* \eta}
 = \ip{S_0^* (J_0 z^* J_0) S_0^* \eta, \xi} \\
&= \ip{R_0(J_0 z J_0)^* \eta, \xi}
 = \ip{\eta, R_0(J_0 z J_0) \xi}
\end{align*}
so that $z \Delta_0^{1/2} \xi \in \mathcal{D}(\Delta_0^{1/2})$ and
$\Delta_0^{-1/2} z \Delta_0^{1/2} \xi = R_0(J_0 z J_0) \xi$.  \\ Hence
$\mathcal{D}(\Delta_0^{-1/2} z \Delta_0^{1/2})=
\mathcal{D}(\Delta_0^{1/2})$ and $\Delta_0^{-1/2} z \Delta_0^{1/2}$
extends to the bounded operator $R_0(J_0 z J_0)$, which is in $N' \cap
M_1$.

A similar argument shows that $\mathcal{D}(\Delta_0^{1/2} z
\Delta_0^{-1/2})= \mathcal{D}(\Delta_0^{-1/2})$ and $\Delta_0^{1/2} z
\Delta_0^{-1/2}$ extends to the bounded operator $R_0(J_0 z^* J_0)^*$,
which is in $N' \cap M_1$.

In detail, 
\begin{align*}
\ip{\Delta_0^{1/2} \xi, z \Delta_0^{-1/2} \eta}
&= \ip{J_0 \overline{S_0} \xi, z J_0 S_0^* \eta}
 = \ip{J_0 z J_0 S_0^* \eta, \overline{S_0} \xi} \\
&= \ip{\xi, S_0^* (J_0 z J_0) S_0^* \eta}
 = \ip{\xi, R_0(J_0 z^* J_0)^* \eta}
\end{align*}
so that $z \Delta_0^{-1/2} \eta \in \mathcal{D}(\Delta_0^{1/2})$ and
$\Delta_0^{1/2} z \Delta_0^{-1/2} \eta = R_0(J_0 z^* J_0)^* \eta$.

Finally, since the two maps $z \mapsto \Delta_0^{-1/2} z
\Delta_0^{1/2}$ and $z \mapsto \Delta_0^{1/2} z \Delta_0^{-1/2}$ are
inverse to each other, they are linear automorphisms of $N' \cap M_1$.

\item{(v)} First note that $J_0 (N' \cap M) J_0 = M' \cap M_1$ and so 
these two spaces have the same dimension.  Next let $x \in N' \cap M$.
Since $S_0 x S_0$ is right multiplication by $x^*$, which commutes
with the (left) action of $M$ on $\Ltwo(M,\varphi)$, we have $S_0 x
S_0 \in M' \cap M_1$.  This map is injective and hence, by a dimension
count, also surjective.  The other maps can all be expressed in terms
of these two and adjoints, so also map $N' \cap M$ onto $M' \cap M_1$.
\end{description}
\end{proof}
\end{thm}

\begin{remark}
We will use $S_k(x)$, $S_k^*(x)$ and $\Delta_k^t(x)$ to denote
$\Lambda_k^{-1}(S_k(\Lambda_k(x)))$, \\ 
$\Lambda_k^{-1}(S_k^*(\Lambda_k(x)))$ and
$\Lambda_k^{-1}(\Delta_k^t(\Lambda_k(x)))$ respectively.  Note that
this is quite different to the operator product, for example $S_k x$
(more precisely denoted $S_k \pi_k(x)$).
\end{remark}


\subsection{The rotation operator}
\label{section: rotation}

Exactly as in Jones~\cite{Jones1999}~4.1.12 we define a rotation
operator on the relative commutants.  We prove that the rotation is
quasi-periodic, namely $\rho_k^{k+1}=\Delta_k^{-1}$.  We will obtain
this result again as a consequence of our work on planar algebras, but
we include a self-contained proof in this section.

Note that in the extremal
$\IIone$ case of~\cite{Jones1999} $E_{M'}$ below can be taken to be
the trace preserving conditional expectation onto $M' \cap M_{k+1}$,
but in the non-extremal type $\mathrm{II}$ case one must use the
commutant trace preserving expectation.  In the full generality
presented here such choices do not occur and the correct path is
clearer.

\begin{definition}[Rotation]
\label{define rotation}
On $N' \cap M_k$ define an operator $\rho_k$ by
\[
\rho_k(x)=\delta^2 E_{M_k}(v_{k+1}E_{M'}(x v_{k+1})) .
\]
Given a basis $B$ for $M$ over $N$, define $r^B_k$ on $M_k$ by
\[
r^B_k(x)=\sum_{b\in B} E_{M_k}\br{v_{k+1}bx v_{k+1}b^*} ,
\]
where $x=\theta\br{x_1 \tlN x_2 \tlN \cdots \tlN
x_{k+1}}$.  By Prop~\ref{E_M' via basis}, $\rho_k(x)=r^B_k(x)$ for
$x\in N' \cap M_k$.
\end{definition}

\begin{lemma}
\label{rotn in terms of basis}
For $x=\theta\left(x_1 \tN x_2 \tN \cdots \tN
x_{k+1}\right)$ and $B$ any basis for $M$ over $N$ we have
\[
r^B_k(x)=\sum_{b \in B}
	\theta\left(E_N(bx_1)x_2\tN x_3 \tN \cdots \tN
	  x_{k+1} \tN b^*\right) .
\]
\begin{proof}
The proof is exactly the same as Jones~\cite{Jones1999}~4.1.14:
\[
x v_{k+1}
	=(x_1 v_1 x_2 v_2 \cdots v_k x_{k+1}) v_{k+1}
	=\theta\left(x_1 \tN \cdots x_{k+1} \tN 1\right)
	=x_1 v_{k+1}^* x_2 v_k^* \cdots x_{k+1} v_1^* 1
\]
so that
\begin{eqnarray*}
r^B_k(x)
	& = & \sumb E_{M_k}(v_{k+1} bx_1 v_{k+1}^* x_2 v_k^* \cdots
		x_{k+1} v_1^* b^*) \\
	& = & \sumb E_{M_k}(\delta E_{k+1} E_N(bx_1) x_2 v_k^* \cdots
		x_{k+1}v_1^* b^*) \\
	& = & \sumb E_N(bx_1) x_2 v_k^* \cdots x_{k+1} v_1^* b^* \\
	& = & \sumb \theta\left(E_N(bx_1)x_2\tN x_3 \tN 
		\cdots \tN x_{k+1} \tN b^*\right)
\end{eqnarray*}
\end{proof}
\end{lemma}


\begin{thm}
\label{rho quasi-periodic}
The rotation is quasi-periodic: $\rho_k^{k+1}=\Delta_k^{-1}$.

\begin{proof}
\begin{description}
\item
\item{(i)} First note that if $x=\theta(x_1\tlN x_2 \tlN \cdots \tlN
x_{k+1})$ and $y=\theta(y_{k+1}\tlN y_k \tlN \cdots \tlN y_1)$, then
$E_N(xy)=E_N(x_1 E_N(x_2 \cdots E_N(x_{k+1} y_{k+1}) \cdots y_2)y_1)$.

Proceed by induction.  The result is true for $k=0$.  Suppose it is
true for $k-1$.  Let 
$\tilde{x}=\theta\left(x_1\tlN x_2 \tlN \cdots \tlN x_k\right)$ and let 
$\tilde{y}=\theta\left(y_k\tlN y_{k-1} \tlN\cdots\tlN y_1\right)$.
Then $x=\tilde{x} v_k x_{k+1}$ and $y=y_{k+1} v_k^* \tilde{y}$.  
Note that $v_k x_{k+1} y_{k+1} v_k^* = \delta E_k E_N(x_{k+1}
y_{k+1})$ and hence $E_{M_{k-1}}(v_k x_{k+1} y_{k+1}
v_k^*)=E_N(x_{k+1} y_{k+1})$.  Thus
\begin{eqnarray*}
E_N(xy)
 & = & E^{M_{k-1}}_N(\tilde{x} E_{M_{k-1}}(v_k x_{k+1} y_{k+1} v_k^*) \tilde{y}) \\
 & = & E^{M_{k-1}}_N(\tilde{x} E_N(x_{k+1} y_{k+1}) \tilde{y}) \\
 & = & E_N(x_1 E_N(x_2 \cdots E_N(x_{k+1} y_{k+1}) \cdots y_2)y_1)
\end{eqnarray*} 

\item{(ii)}
Let $B$ be a basis for $M$ over $N$.  By
Lemma~\ref{basislemma} $B_k=\{\theta(b_{i_1}\tlN \cdots
\tlN b_{i_{k+1}}) : i_j \in I \}$ is a basis for $M_k$ over $N$.
Let $x\in N'\cap M_k$.  Then we can write $x$ as 
\[
x = \sum_{c\in B_k} c n_c = \sum_{i_1,\ldots ,i_{k+1}} \theta\left(b_{i_1}\tN \cdots \tN b_{i_{k+1}}\right) n_{i_1,\ldots ,i_{k+1}} 
\]
where $n_c$ (or $n_{i_1,\ldots ,i_{k+1}}$) are in $N$.
Hence
\begin{eqnarray*}
\rho_k(x) & = & \sum_{i_1,\ldots ,i_{k+1}, j_1} 
          \theta\left(E(b_{j_1}b_{i_1})b_{i_2}\tN \cdots \tN
	  b_{i_{k+1}} n_{i_1,\ldots ,i_{k+1}} \tN b_{j_1}^*\right) 
	   \\
\left(\rho_k\right)^{k+1}(x) & = & 
        \sum_{\stackrel{i_1,\ldots ,i_{k+1}}{j_1,\ldots ,j_{k+1}}}
	\theta\left(E(b_{j_{k+1}}E(b_{j_k} \cdots E(b_{j_2}E(b_{j_1}
	b_{i_1})b_{i_2}) \cdots b_{i_k})b_{i_{k+1}})
	n_{i_1,\ldots ,i_{k+1}}\right) \cdot  \\
&   &   \;\;\;\;\;\;\;\;\;\;\;\;\;\;\;\;\;\;\;\;\;\;\;\;\;\;\;\;\;\;\;\;\;\;
	\cdot \theta\left(b_{j_1}^* \tN \cdots \tN b_{j_{k+1}}^*
	\right) \\
& = &   \sum_{\stackrel{i_1,\ldots ,i_{k+1}}{j_1,\ldots ,j_{k+1}}}
	E_N\left(\theta\left(b_{j_{k+1}}\tN \cdots b_{j_1}\right)
	    \theta\left(b_{i_1}\tN \cdots \tN b_{i_{k+1}}\right)n_{i_1,\ldots ,i_{k+1}}\right) \cdot  \\
&   &   \;\;\;\;\;\;\;\;\;\;\;\;\;\;\;\;\;\;\;\;\;\;\;\;\;\;\;\;\;\;\;\;\;\;
	\cdot \theta\left(b_{j_1}^* \tN \cdots \tN b_{j_{k+1}}^*\right) 
	 \\
& = &   \sum_{c\in B_k} E_N(c x) c^* \\
& = &  \Delta_k^{-1}(x) ,
\end{eqnarray*}
where the last equality comes from Theorem~\ref{basic modular results}.
\end{description}
\end{proof}
\end{thm}


\section{The Planar Algebra of a Finite Index Subfactor}
\label{sect: PA of fin index subfactor}

We are now in a position to define the planar algebra associated to a
finite index subfactor.  By Theorem~\ref{rho quasi-periodic} the
rotation is not quite periodic in the general case and this requires
us to change the axioms of a planar algebra slightly.  We will first
define a {\em rigid \cstar-planar algebra} in which boxes cannot be
rotated, which then gives rise to an induced {\em modular extension}
in which boxes can be rotated, but this changes the action of the
tangle in a specified way.

We show that the standard invariant of a finite index subfactor forms
a rigid \cstar-planar-algebra.  Along the way we will see that the
(quasi-)periodicity of the rotation is a trivial consequence of the
rigid planar algebra structure.

We show that any rigid \cstar-planar algebra gives rise to a spherical
\cstar-planar algebra.  As a corollary we see that for any finite
index subfactor there exists an extremal $\IIone$ subfactor with the
same (algebraic) standard invariant, which recovers a result
originally due to Izumi.

Finally we consider the inverse construction from a spherical
\cstar-planar algebra with some additional data to a rigid 
\cstar-planar algebra.  Lifting this construction to the subfactor 
level we show that any rigid \cstar-planar algebra arises from a
finite index subfactor.

These results justify the focus on $\IIone$ subfactors rather than
more general inclusions, at least as far as the study of the standard
invariant of finite index subfactors is concerned.  We see all
possible standard invariants of finite index subfactors by considering
$\IIone$ subfactors and if our interest lies only in the algebraic
structure without reference to Jones projections and conditional
expectations we need only consider $\IIone$ extremal subfactors with
the trace-preserving conditional expectation.


\subsection{Rigid planar algebras}
\label{section: rigid PA}

In essence the only difference between a planar algebra and a rigid
planar algebra is that we use rigid planar isotopy classes of tangles
in place of (full) planar isotopy classes of tangles.  A rigid planar
isotopy is one under which the internal discs undergo no rotation.
The set of rigid planar isotopy classes of tangles forms the rigid
planar operad $\PAr$.  A rigid planar algebra is an algebra over
$\PAr$ in the sense of May~\cite{May1997a}, that is to say a set of
vector spaces $V_k^+, V_k^-$ ($k\geq 0$) and a morphism of colored
operads from $\PAr$ to $\Hom$, the operad of linear maps between
tensor products of the $V_k^{\pm}$'s.  We describe these ideas in
detail below.

\begin{definition}[(Rigid) planar $k$-tangle]
A {\em planar $k$-tangle} is defined exactly as in
Jones~\cite{Jones1999}, except that we require the boundary points of
the discs to be evenly spaced, the strings to meet the discs normally
and the distinguished boundary segments can be either white or black.
Thus the definition below follows Jones almost verbatim.  For
simplicity when drawing tangles we will often draw the boundary points
so that they are not evenly spaced.  In general diagrams will be drawn
with the distinguished boundary segment on the left.  When this is not
the case we will denote the distinguished segment with a star.

A {\em planar $k$-tangle} will consist of the unit disc $D(=D_0)$ in
$\mathbb{C}$ together with a finite (possibly empty) set of disjoint
subdiscs $D_1,D_2,\ldots,D_n$ in the interior of $D$.  Each disc
$D_i$, $i\geq 0$, will have an even number $2k_i\geq 0$ of evenly
spaced marked points on its boundary (with $k=k_0$).  Inside $D$ there
is also a finite set of disjoint smoothly embedded curves called {\em
strings} which are either closed curves or whose boundaries are marked
points of the $D_i$'s.  Each marked point is the boundary point of
some string, which meets the boundary of the corresponding disc
normally.  The strings all lie in the complement of the interiors
$D_i^{\text{o}}$ of the $D_i$, $i>0$.  The connected components of the
complement of the strings in $D^{\text{o}}\backslash \bigcup_{i=1}^n D_i$ are
called regions and are shaded black and white so that regions whose
closures meet have different shadings.  The shading is part of the
data of the tangle, as is the choice, at every $D_i$, $i\geq 0$, of a
region whose closure meets that disc, or equivalently: a {\em
distinguished arc} between consecutive marked points (the whole
boundary of the disc if there are no marked points); or a {\em
distinguished point} which we will call $p_i$ at the midpoint of the
arc.  Define $\sigma_i$, the {\em sign} of $D_i$, to be $+$ if the
distinguished region is white and $-$ if it black.

A {\em rigid planar $k$-tangle} is a planar $k$-tangle such that for
every point $p_i$, $i\geq 0$, the phase relative to the center $x_i$
of $D_i$ is the same.  In other words the angle between the line from
$p_i$ to $x_i$ and a ray in the positive $x$-direction emanating from
$x_i$ is the same for all $i$.

We will call a (rigid) planar $k$-tangle, with $\sigma_0$ the sign of
$D$, a {\em (rigid) planar $(\sigma_0,k)$-tangle}.

Composition of tangles is defined as follows.  Given a planar
$k$-tangle $T$, a $k'$-tangle $S$, and an internal disc $D_i$ of $T$
with $k_i=k'$ and $\sigma_i(T)=\sigma_0(S)$ we define the $k$-tangle
$T \circ_i S$ by radially scaling $S$, then rotating and translating
it so that its boundary, together with the distinguished point,
coincides with that of $D_i$.  The boundary of $D_i$ is then removed
to obtain the tangle $T \circ_i S$.

\noindent Note: We will often use $0$ and $1$ in place or the signs $+$ and $-$
respectively.  Sometimes we will drop the sign completely in the case
of $\sigma=+$ and use a $\wtilde{\phantom{A}}$ to denote a disc with
$\sigma=-$.
\end{definition}

\begin{definition}[(Rigid) planar operad] 
A {\em planar isotopy} of a tangle $T$ is an orientation preserving
diffeomorphism of $T$, preserving the boundary of $D$.  A {\em rigid
planar isotopy} is a planar isotopy such that the restriction to
$D_i$, $i>0$, only scales and translates the disc $D_i$.

The {\em (full) planar operad} $\PA$ is the set of all planar isotopy
classes of planar $k$-tangles, $k$ being arbitrary.  The {\em rigid
planar operad} $\PAr$ is the set of all rigid planar isotopy classes
of planar $k$-tangles.  Clearly composition of tangles passes to the
planar operad and to the rigid planar operad.
\end{definition}

\begin{remark}
In a rigid planar tangle the phase of the external distinguished point
$p_0$ is not important, only the relative position of the internal
distinguished points.  Thus we could require that $p_0$ is the
leftmost point of $D$ (relative to some underlying orthogonal
coordinate system for the plane) and the condition could then be
stated as: for all $i>0$, $p_i$ is the leftmost point of $D_i$.

The use of discs in the definition of tangles and planar operads is
convenient, but we could contract the internal discs to points and
formulate the definitions in these terms, or use boxes (with $p_i$ on
the left edge) in place of discs which would lead to another
reformulation.  We will use these descriptions interchangeably.
\end{remark}

\begin{definition}[Rigid planar algebra]
A rigid planar algebra $(Z,V)$ is an algebra over $\PAr$.  That is to
say we have a disjoint union $V$ of vector spaces $V_k^+, V_k^-$
(sometimes denoted $V_k^0,V_k^1$ respectively), $k\geq 0$, and a
morphism of colored operads from $Z:\PAr \rightarrow \Hom(V)$, the
operad of linear maps between tensor products of the $V_k^{\pm}$'s.
In other words, for every (equivalence class of) rigid planar
$k$-tangle $T$ in $\PAr$ there is a linear map $Z(T):\tensor_{i=1}^n
V_{k_i}^{\sigma_i} \rightarrow V_k^{\sigma_0}$ (which is thus
unchanged by rigid planar isotopy).  The map $Z$ satisfies $Z(T\circ_i
S)=Z(T) \circ_i Z(S)$.

For any tangle $T$ with no internal discs, the empty tensor product is
just $\mathbb{C}$ and $Z(T):\mathbb{C} \rightarrow V_k^{\pm}$ and is
thus just multiplication by some element of $V_k^{\pm}$.  We will also
use $Z(T)$ to denote this element.

We also require finite dimensionality ($\dim(V_k^{\pm})<\infty$ for
all $k$), $\dim(V_0^{\pm})=1$ and
$Z\br{\includegraphics[-1pt,5pt][20pt,20pt]{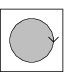}},
Z\br{\includegraphics[-1pt,5pt][20pt,20pt]{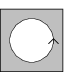}}
\neq 0$.  (note that these tangles are a white disc with one closed
string and a black disc with one closed string respectively - we have
drawn the discs as boxes to make this clear).  Finally, we require
that $Z(\text{annular tangle with $2k$ radial strings})$ be the
identity (for each $k$ there are two such tangles, either $\sigma_0=0$
or $\sigma_0=1$).
\end{definition}

\begin{remarks}
Without the last condition in the definition we would only know that
$Z(\text{annular tangle with $2k$ radial strings})$ is some idempotent
element $e_k^{\pm}$ in $\End(V_k^{\pm})$.  However, by inserting this
annular tangle around the inside of any internal $k$-disc and the
outside of any $k$-tangle it is clear that the planar algebra would
never ``see'' $(1-e_k^{\pm})V_k^{\pm}$.  Thus we may as well replace
$V_k^{\pm}$ with $e_k^{\pm}V_k^{\pm}$ and then the annular radial
tangles act as the identity.

$V_k^{\pm}$ becomes an algebra under the multiplication given by the
tangle below.  The thick black string represents $k$ regular strings.
\begin{figure}[htbp]
\begin{center}
\psfrag{D_1}{$D_1$}
\psfrag{D_2}{$D_2$}
\psfrag{*}{$*$}
\includegraphics{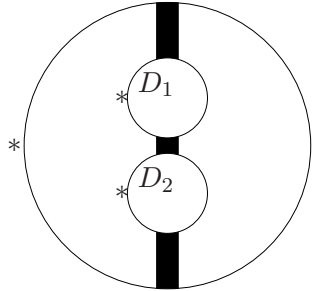}
\caption{The multiplication tangle}
\label{figure: multiplication}
\end{center}
\end{figure}
\vspace{-20pt}

\noindent For example, for $k=3$ we have the following tangles (for
$V_3^+$ and $V_3^-$ respectively).  The marked points have been evenly
spaced in this example as they should be, but for simplicity we will
usually draw the marked points with uneven spacing.
\begin{figure}[htbp]
\begin{center}
\psfrag{D_1}{$D_1$}
\psfrag{D_2}{$D_2$}
\psfrag{*}{$*$}
\includegraphics{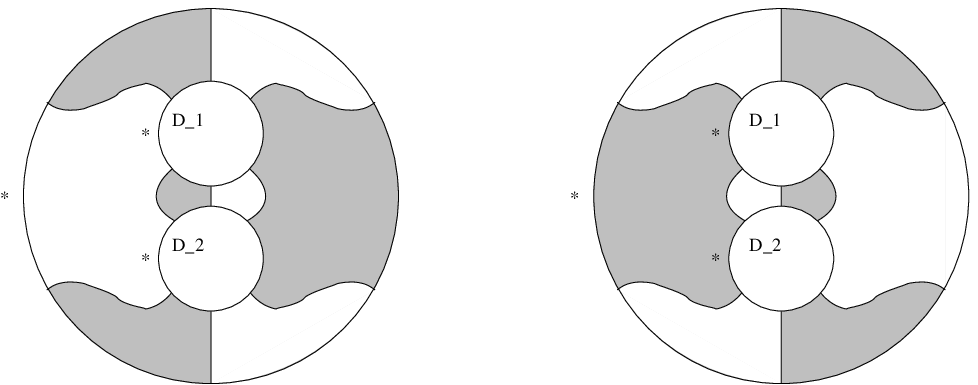}
\end{center}
\end{figure}
\vspace{-20pt}

\noindent The multiplicative identity in $V_k^{\pm}$ is just the image
of the tangle below (more precisely $Z$ of this tangle applied to $1
\in \mathbb{C}$, the empty tensor product).
\begin{figure}[htbp]
\begin{center}
\psfrag{*}{$*$}
\includegraphics{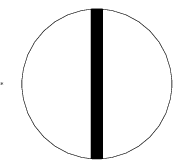}
\end{center}
\end{figure}
\vspace{-20pt}

\noindent Note that $V_0^{\pm}$ is now a 1-dimensional algebra with identity
given by $Z(\text{empty disc})$ (the picture above with no strings),
hence $V_0^{\pm}$ may be identified with $\mathbb{C}$ with
$Z(\text{empty disc})=1$.

With $V_0^{\pm}$ thus identified with $\mathbb{C}$, 
$Z\br{\includegraphics[-1pt,5pt][20pt,20pt]{figures/loop_black.eps}}
=\delta_1$ and 
$Z\br{\includegraphics[-1pt,5pt][20pt,20pt]{figures/loop_white.eps}}
=\delta_2$ for some $\delta_1,\delta_2 \in
\mathbb{C}\backslash\{ 0\}$.
Unless stated otherwise, we will require $\delta_1=\delta_2\defeq
\delta$ (one can always ``rescale'' to achieve this).

The fact that $\delta\neq 0$ implies that
$\iota_k^{\pm}:V_k^{\pm}\rightarrow V_{k+1}^{\pm}$, defined by adding
a string on the right, is injective.  Similarly
$\gamma_k^{\pm}:V_k^{\pm} \rightarrow V_{k+1}{\mp}$ by adding a string
on the left is injective.
\end{remarks}

\begin{definition}[Rigid planar *-algebra]
For a (rigid) $k$-tangle $T$ define $T^*$ to be the tangle obtained by
reflecting $T$ in any line in the plane (well-defined up to rigid
planar isotopy).
A {\em rigid planar *-algebra} is a rigid planar algebra $(Z,V)$
equipped with a conjugate linear involution $*$ on each $V_k^{\pm}$
such that $\sqbr{Z(T)(v_1 \tensor\cdots\tensor v_n)}^*=Z(T^*)(v_1^*
\tensor\cdots\tensor v_n^*)$.
\end{definition}

\begin{definition}[Rigid \cstar-planar algebra]
A {\em rigid \cstar-planar algebra} is a rigid \mbox{*-planar} algebra
$(Z,V)$ such that the map $\Phi=\Phi_k^{\pm}:V_k^{\pm}\rightarrow
\mathbb{C}$ given below is positive definite on $V_k^{\pm}$,
i.e. $\Phi(x^*x)>0$ for $0\neq x\in V_k^{\pm}$.
\begin{figure}[htbp]
\begin{center}
\psfrag{D_1}{}
\psfrag{*}{\small $*$}
\psfrag{Phi=Z( )}{$\Phi=Z
  \begin{pmatrix} &&&&&\\ &&&&&\\ &&&&&\\ &&&&&\end{pmatrix}$}
\includegraphics{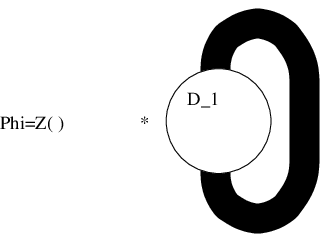}
\end{center}
\end{figure}
\vspace{-15pt}

\noindent As usual the thick string represents $k$ strings and we have
suppressed the outer $0$-disc.
\end{definition}

\begin{remark}
Let $(Z,V)$ be a rigid *-planar algebra.  Note that if $\Phi_k^+$ is
positive definite then so too is $\Phi_k^-$ by adding a string on the
left and applying $\Phi_{k+1}^+$.  Hence it suffices to check that
$\Phi_k^+$ is positive definite.  

Every $V_k^{\pm}$ is semi-simple since for any nonzero ideal $N$ take
nonzero $x\in N$, then $x^*xx^*x=(x^*x)^*(x^*x) \in N^2$ and is
nonzero by positive definiteness of $\Phi$.  Hence every $V_k^{\pm}$
is a direct sum of matrix algebras over $\mathbb{C}$ or, equivalently,
a finite dimensional \cstar-algebra.
\end{remark}

\begin{definition}
Define $\Phi'_k:V_k^{\pm}\rightarrow \mathbb{C}$ by
\begin{figure}[hbtp]
\begin{center}
\psfrag{D_1}{}
\psfrag{*}{\small $*$}
\psfrag{Phi'=Z( )}{$\Phi'=Z 
  \begin{pmatrix} &&&&&\\ &&&&&\\ &&&&&\\ &&&&&\end{pmatrix}$}
\includegraphics{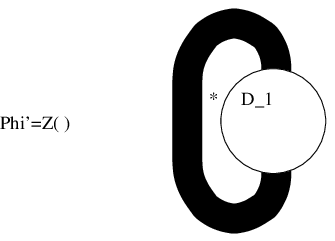}
\end{center}
\end{figure}

\newpage

\noindent Define $\cj_k^{(r)}=Z(TJ_k^{(r)}):V_k^{(r)}\rightarrow
V_k^{(r+k)}$ where $TJ_k^{(r)}$ is the following tangle ($r=0,1$)
\begin{figure}[htbp]
\begin{center}
\psfrag{D_1}{}
\psfrag{*}{\small $*$}
\psfrag{TJ_k=}{$TJ_k=$}
\includegraphics{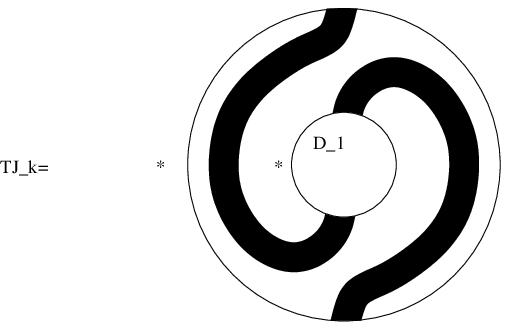}
\end{center}
\end{figure}
\vspace{-35pt}
\end{definition}

\begin{remarks}
Note that $\cj_k^{(r)}$ is invertible, with inverse
$Z\br{\br{TJ_k^{(r+k)}}^*}$.  $\cj$ has the following properties.  Let
$x\in V_k^{(r)}$.  Then, with indices suppressed,
\begin{itemize}
\item $\cj(xy)=\cj(y)\cj(x)$;
\item $\br{\cj(x)}^*=\cj^{-1}(x^*)$;
\item $\Phi(\cj(x))=\Phi(\cj^{-1}(x))=\Phi'(x)$;
\item $\Phi'(\cj(x))=\Phi'(\cj^{-1}(x))=\Phi(x)$;
\item $\Phi(xy)=\Phi'(\cj(x)\cj^{-1}(y))$, $\Phi'(xy)=\Phi(\cj^{-1}(x)\cj(y))$
because
\begin{figure}[htbp]
\begin{center}
\psfrag{x}{$x$}
\psfrag{y}{$y$}
\psfrag{Phi(j^{-1}(x)j(y))=}{$\hspace{-20mm}\Phi\br{\cj^{-1}(x)\cj(y)}=$}
\psfrag{=}{$=$}
\psfrag{=Phi'(xy)}{$=\Phi'(xy)$}
\includegraphics{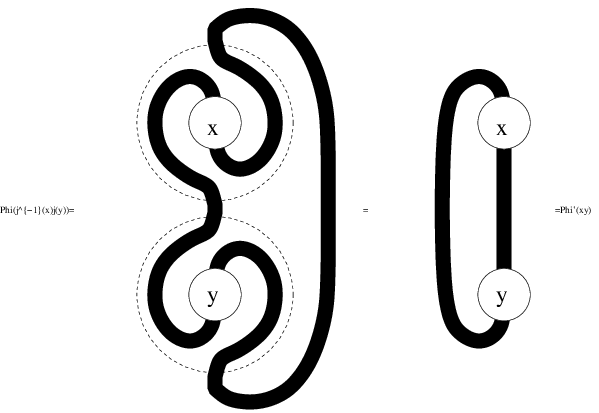}
\end{center}
\end{figure}
\end{itemize}
\end{remarks}

\begin{lemma}
Let $(Z,V)$ be a rigid \cstar-planar algebra.  Then
$\Phi'=\Phi'_k:V_k^{\pm}\rightarrow \mathbb{C}$ is positive definite.
\begin{proof}
For nonzero $x\in V_k^{\pm}$ let $y=\cj(x)\neq 0$, then
\[
0<\Phi(y^*y)=\Phi(\br{\cj(x)}^*\cj(x))=\Phi(\cj^{-1}(x^*)\cj(x))=\Phi'(x^*x) .
\]
\end{proof}
\end{lemma}

\begin{remarks}
The requirement that an annular tangle with $2k$ radial strings act as
the identity in the definition of a rigid planar algebra is not
necessary for a rigid \cstar-planar algebra.  Let
$\alpha=Z(\text{annular tangle with $2k$ radial strings})$.  Then
$\Phi(x\alpha(y))=\Phi(xy)$ for all $x\in V_k^{\pm}$ and $\alpha(x)=x$
by the positive definiteness of $\Phi$.

$\varphi_k=\delta^{-k}\Phi_k$ and $\varphi'_k=\delta^{-k}\Phi'_k$.
are normalized so that $\varphi_k(1)=\varphi'_k(1)=1$ and both
$\varphi$ and $\varphi'$ are compatible with the inclusions
$\iota:V_k^{\pm}\rightarrow V_{k+1}^{\pm}$ and $\gamma:V_k^{\pm}
\rightarrow V_{k+1}^{\mp}$.
\end{remarks}

\begin{notation}
For $r=0,1$ let $V_{j,k}^{(r)}=\gamma_{k-1}^{(r+1)}\gamma_{k-2}^{(r)}
\cdots\gamma_{k-j+1}^{(r+j+1)}\gamma_{k-j}^{(r+j)}\br{V_{k-j}^{(r+j)}}$
denote the subspace of $V_k^{(r)}$ obtained by adding $j$ strings to
the left of $V_{k-j}^{(r+j)}$ (where all upper indices are computed
mod $2$).  In other words the image of the map defined by
\begin{figure}[htbp]
\begin{center}
\psfrag{j}{$\underbrace{\hspace{10mm}}_{j}$}
\psfrag{k-j}{$\underbrace{\hspace{8mm}}_{k-j}$}
\includegraphics{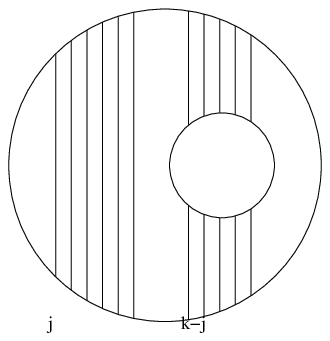}
\end{center}
\end{figure}
\vspace{-15pt}

\noindent Call this map the {\em shift by $j$}.

Let $w_{j,k}^{\pm}$ denote the Radon-Nikodym derivative of
$\br{\varphi'}_k^{\pm}$ with respect to $\varphi_k^{\pm}$ on
$V_{i,k}^{\pm}$.  Thus $w_{j,k}^{\pm}$ is the unique element in
$V_{j,k}^{\pm}$ such that $\varphi'(x)=\varphi(w_{j,k}^{\pm} x)$ for
all $x \in V_{j,k}^{\pm}$.  $w_{j,k}^{\pm}$ exists because
$V_{j,k}^{\pm}$ is finite dimensional and $\varphi$ is positive
definite.  By the positive definiteness of $\varphi'$, $w_{j,k}^{\pm}$
is positive and invertible.

Let $w_k^{\pm}=w^{\pm}_{k-1,k}$ and $z_k^{\pm}=w_{0,k}^{\pm}$.  As
mentioned earlier we will sometimes suppress the $+$ index and use
$\wtilde{\phantom{A}}$ in place of $-$.  Set $w=w_1=w^+_1$,
$\wtilde{w}=\wtilde{w}_1=w^-_1$.  Then $w_{2r+1}$ is just $w$ with
$2r$ strings to the left and $w_{2r}$ is $\wtilde{w}$ with $2r-1$
strings to the left.  i.e. $w_{2r+1}=\br{\gamma^-\gamma^+}^r(w)$ and
$w_{2r}=\br{\gamma^-\gamma^+}^{r-1}\gamma^-(\wtilde{w})$.
\end{notation}

\begin{lemma}
\label{lemma: properties of RN derivs}
\begin{description}
\item
\item
\item{1.} $\wtilde{w}=\br{\cj_1}^{-1}(w^{-1})=$
\vspace{-62pt}
\begin{figure}[htbp]
\begin{center}
\hspace{-98pt}
\psfrag{w^{-1}}{$w\!^{-1}$}
\includegraphics{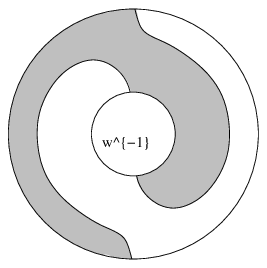}
\end{center}
\end{figure}

\newpage
\item{2.} $z_k=w_1w_2\cdots w_k$ and is thus given by
\begin{figure}[htbp]
\begin{center}
\psfrag{w}{$w$}
\psfrag{wt}{$\wtilde{w}$}
\psfrag{w^{-1}}{$w^{-1}$}
\psfrag{=}{$=$}
\psfrag{...}{$\ldots$}
\psfrag{z_k}{$z_k$}
\includegraphics{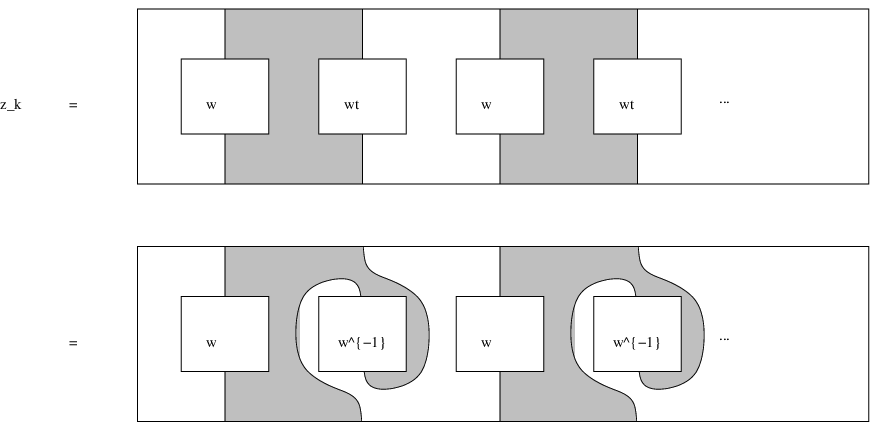}
\end{center}
\end{figure}
\vspace{-30pt}
\end{description}

\begin{proof}

\begin{description}
\item
\item{1.} For $x\in V_1^-$,
$
\wtilde{\varphi}_1(\cj_1^{-1}(w^{-1})x)
=\varphi'_1(w^{-1}\cj_1^{-1}(x))
=\varphi_1(\cj_1^{-1}(x))
=\wtilde{\varphi}'_1(x)
$
and hence $\wtilde{w}=\br{\cj_1}^{-1}(w^{-1})$.  

\vspace{8pt}
\item{2.} Let $y=w_1w_2\cdots w_k=$
\vspace{-42pt}
\begin{figure}[htbp]
\begin{center}
\hspace{-40pt}
\psfrag{w}{$w$}
\psfrag{wt}{$\wtilde{w}$}
\psfrag{...}{$\ldots$}
\psfrag{y=w1..wk=}{}  
\includegraphics{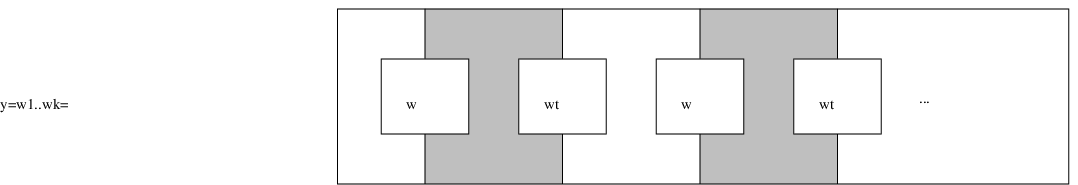}
\end{center}
\end{figure}
\vspace{-20pt}

Then, for $x\in V_k^+$,

\vspace{60pt}

\noindent $\Phi(yx)=$
\vspace{-82pt}
\begin{figure}[h]
\begin{center}
\psfrag{w}{$w$}
\psfrag{wt}{$\wtilde{w}$}
\psfrag{x}{$x$}
\psfrag{...}{$\ldots$}
\psfrag{=}{$=$}
\includegraphics{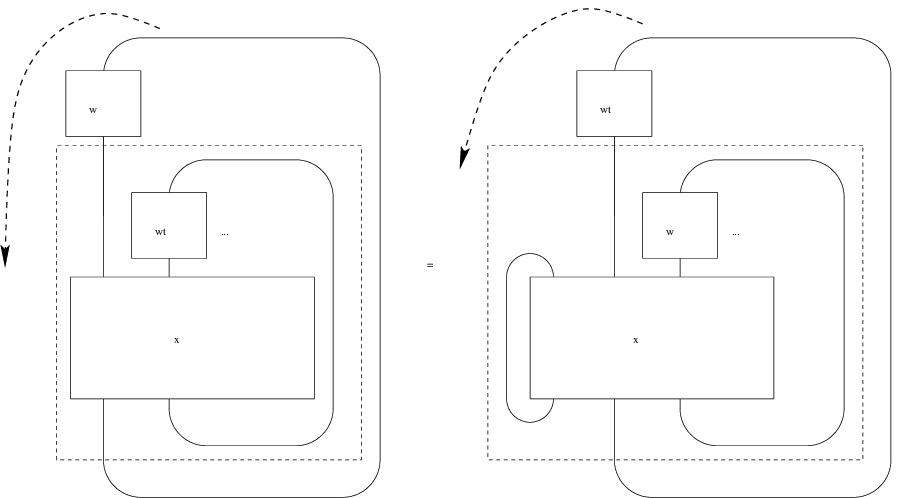}
\end{center}
\end{figure}
\vspace{-20pt}
\begin{figure}[htbp]
\begin{center}
\psfrag{w}{$w$}
\psfrag{x}{$x$}
\psfrag{...}{$\ldots$}
\psfrag{..}{$.\;.$}
\psfrag{=}{$=$}
\psfrag{= ... =}{$= \ldots =$}
\psfrag{=Phi'(x)}{$=\Phi'(x)$}
\includegraphics{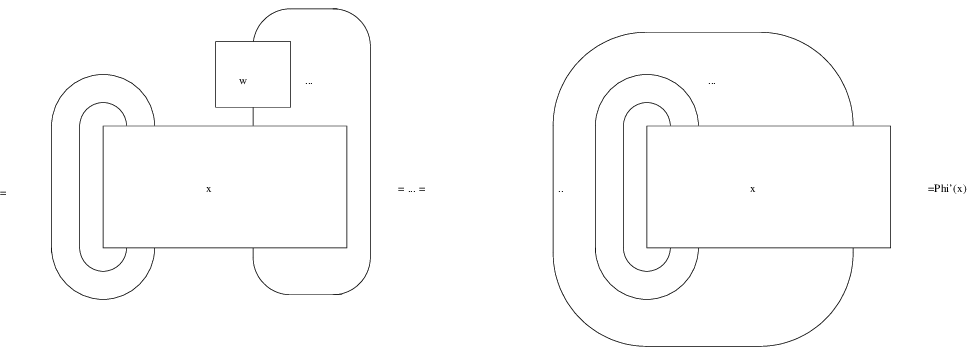}
\end{center}
\end{figure}
\vspace{-20pt}
\end{description}
\end{proof}
\end{lemma}

\begin{definition}
Define $\Delta_k^{\pm}:V_k^{\pm} \rightarrow V_k^{\pm}$ by
$\Delta=\cj^2$.  More precisely
$\Delta_k^{(r)}=\cj_k^{(r+k)}\circ\cj_k^{(r)}$.
\vspace{-20pt}
\begin{figure}[htbp]
\begin{center}
\psfrag{D(x)=}{$\Delta(x)=$}
\psfrag{x}{$x$}
\includegraphics{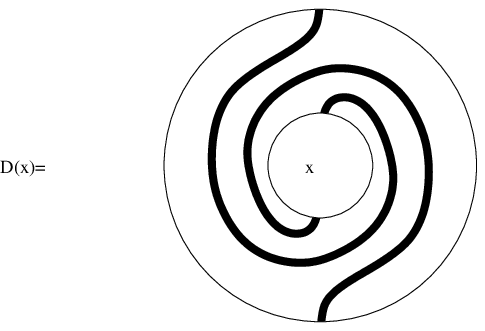}
\end{center}
\end{figure}
\vspace{-20pt}
\end{definition}

\begin{lemma}
\label{lemma: PA phi, phi' and Delta}
For all $x,y \in V_k^{\pm}$,
\begin{eqnarray*}
\varphi(xy) & = & \varphi(y \Delta(x)) = \varphi(\Delta^{-1}(y)x) , \\
\varphi'(xy) & = & \varphi'(y \Delta^{-1}(x)) = \varphi'(\Delta(y)x) .
\end{eqnarray*}
\begin{proof}
First
\[
\varphi(xy)
=\varphi'(\cj(x)\cj^{-1}(y))
=\varphi(\cj(\cj(x)\cj^{-1}(y)))
=\varphi(y\cj^2(x))
\]
and because $\varphi\circ\cj^2=\varphi'\circ\cj=\varphi$ we also have
$\varphi(xy)=\varphi(\cj^{-2}(y)x)$.  Secondly
\[
\varphi'(xy)
=\varphi(\cj^{-1}(x)\cj(y))
=\varphi(y\cj^2(x))
=\varphi'(\cj^{-2}(y)x) .
\]
\end{proof}
\end{lemma}

\begin{remark}
This result shows that $\Delta$ is in fact the modular operator for
$(V_k^{\pm},\varphi)$ and $\Delta^{-1}$ the modular operator for
$(V_k^{\pm},\varphi')$.
\end{remark}

\begin{cor}
$\Delta_{(+,k)}$ is a positive definite operator on
$\Ltwo(V_k^+,\varphi)$ and for $x \in V_{j,k}^+$ and $t \in
\mathbb{R}$, $\Delta_{(+,k)}^t \hatt{x} =
\hatr{w_{j,k}^{-t/2}xw_{j,k}^{t/2}}$ .
\begin{proof}
$\Delta_k^+$ is positive definite because, for nonzero $x\in V_k^+$,
\[
\ip{\Delta(x),x}
=\varphi\br{x^*\cj^2(x)}
=\varphi'\br{\cj(x)\cj^{-1}(x^*)}
=\varphi'\br{\cj(x)\br{\cj(x)}^*}
> 0 .
\]
Let $W=w_{j,k}$.  The result is true for $t=2$ because for all $x,y\in
V_{j,k}^+$
\[
\varphi(Wxy)
=\varphi'(xy)
=\varphi'(y\Delta^{-1}(x))
=\varphi(Wy\Delta^{-1}(x))
=\varphi(\Delta^{-2}(x)Wy) ,
\]
using Lemma~\ref{lemma: PA phi, phi' and Delta} twice.  Hence
$Wx=\Delta^{-2}(x)W$ so that $\Delta^{-2}(x)=WxW^{-1}$ and also \\
$\Delta^2(x)=W^{-1}xW$.

To prove the result for general $t$ we will first show that $x \mapsto
W^{-s} x W^s$ is a positive operator on $\Ltwo(V_{j,k}^+,\varphi)$ for
$s\in \mathbb{R}$.  Note that $\Delta^2
\hatt{W^s}=W^{-1}(W^s)W=W^s$ so that $\hatt{W^s}$ is an eigenvector of
eigenvalue $1$ for $\Delta^2$ and hence also for all $\Delta^r$ since
$\Delta$ is positive.
Therefore $W^s$ is in the centralizer of $\varphi$ on $V_{k}^+$
since, for all $x \in V_{k}^+$,
\[
\varphi(W^sx)
=\varphi(x\Delta(W^s))
=\varphi(xW^s) .
\]
This implies that for $s \in \mathbb{R}$, $A_s:x\mapsto W^{-s}xW^s$ is
a positive operator on $\Ltwo(V_{j,k}^+,\varphi)$:
\[
\ip{W^{-s}xW^s,x}
	=\varphi(x^*W^{-s}xW^s)
	=\varphi(W^{s/2}x^*W^{-s}xW^{s/2})
	\geq 0 .
\]
Thus $\{A_s\}$ is a continuous one-parameter family of positive
operators on $\Ltwo(V_{j,k}^+,\varphi)$, and $\Delta^2=A_1$.  A simple
spectral theory argument implies that $\Delta^r=A_{r/2}$ for all
$r\in\mathbb{R}$ and so $\Delta^t(x)=W^{-t/2}xW^{t/2}$.
\end{proof}
\end{cor}

\begin{cor}
\label{cor: PA rho periodic}
If $\varphi$ is tracial then $\Delta=\id$, $\varphi'$ is tracial and
the rotation operator $\rho_k:V_k^{\pm}\rightarrow V_k^{\pm}$ defined
below is periodic.
\begin{figure}[htbp]
\begin{center}
\psfrag{..}{$.\; .$}
\psfrag{rho=Z()}{$\rho_k=Z
  \begin{pmatrix} &&&&&&&&\\ &&&&&&&&\\ &&&&&&&&\\ &&&&&&&&\end{pmatrix}$}
\psfrag{k-2}{$\underbrace{\hspace{8mm}}_{k-2}$}
\includegraphics{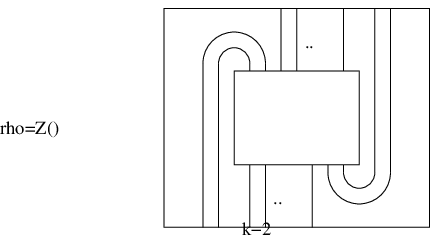}
\end{center}
\end{figure}
\vspace{-10pt}
\end{cor}

\subsection{Subfactors give rigid planar algebras} 
\label{section: subfactor to PA}

In this section we extend Jones'~\cite{Jones1999} result that the
standard invariant of an extremal finite index $\IIone$ subfactor has
a spherical \cstar-planar algebra structure.

\begin{thm}
\label{thm: subfactor to RPA}
Let $(N,M,E)$ be a finite index subfactor.  Let $V_k^+=N'\cap M_{k-1}$
and $V_k^-=M'\cap M_k$.  Then $V$ has a rigid \cstar-planar algebra
structure $Z^{(N,M,E)}$ satisfying:
\begin{description}
\item{(1)} $\delta_1=\delta_2=\delta\defeq \Ind(E)^{1/2}$.
\item{(2)} The inclusion maps $\iota_k^{\pm}:V_k^{\pm}\rightarrow
V_{k+1}^{\pm}$ are the usual inclusion maps $N'\cap M_{k-1}
\hookrightarrow N'\cap M_k$ and $M'\cap M_k \hookrightarrow M'\cap
M_{k+1}$.
\item{(3)} The map $\gamma_k^-:V_k^-\rightarrow V_{k+1}^+$ is the
inclusion map $M'\cap M_k \hookrightarrow N'\cap M_k$.  The map
$\gamma_k^+:V_k^+\rightarrow V_{k+1}^-$ is the shift map $sh:N'\cap
M_{k-1} \rightarrow M_1'\cap M_{k+1}$ defined by 
$R\in \End_{N-*}\br{\Ltwo(M_j)} \mapsto \id \tlN R \in 
\End_{M_1-*}\br{\Ltwo(M_{j+1})}$.

\item{(4)} The multiplication given by the standard multiplication
tangles agrees with the multiplication on the algebras $N' \cap M_k$.

\item{(5)} 
\vspace{-22pt}
\begin{figure}[htbp]
\hspace{24pt}
\psfrag{..}{$\cdots$}
\psfrag{k}{\scriptsize $k$}
\psfrag{k+1}{\scriptsize $k+1$}
\psfrag{=de_k=E_k}{$=\delta e_k=E_k$}
\includegraphics{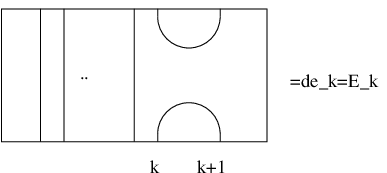}
\end{figure}
\vspace{-10pt}

\item{(6)}
\vspace{-22pt}
\begin{figure}[htbp]
\hspace{24pt}
\psfrag{..}{$.\;.$}
\psfrag{k}{\hspace{-4pt}$\underbrace{\hspace{24pt}}_{\scriptsize k}$}
\psfrag{=dE_M}{$=\delta E_{M_{k-1}}:N'\cap M_k \rightarrow N'\cap M_{k-1}$}
\includegraphics{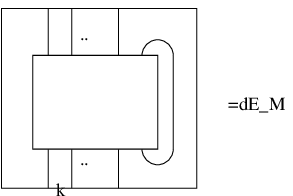}
\end{figure}

\item{(7)}
\vspace{-22pt}
\begin{figure}[htbp]
\hspace{30pt}
\psfrag{..}{$.\;.$}
\psfrag{j+1}{\hspace{-10pt}
  $\underbrace{\hspace{28pt}}_{\scriptsize j+1}$}
\psfrag{j1}{\hspace{-4pt}
  $\underbrace{\hspace{8pt}}_{\scriptsize j+1}$}
\psfrag{=dE_M'}{$=\delta E_{M'}$ \ \ and in general }
\psfrag{=dE_M_j'}{$=\delta^{j+1} E^{N'}_{M_j'}$}
\includegraphics{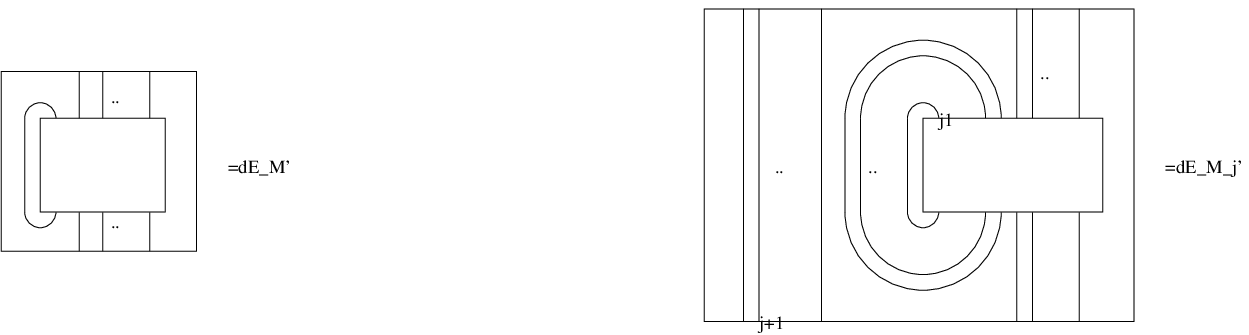}
\end{figure}
\end{description}
The rigid \cstar-planar algebra structure on $V$ is uniquely
determined by properties (4), (5), (6) and (7) (for $j=1$).
\end{thm}

\begin{remark}
As a consequence of (6) and (7) we see that $\varphi$ and $\varphi'$
coming from the planar algebra structure agree with those already
defined on the standard invariant.
\end{remark}

\Proof {\it of Theorem~\ref{thm: subfactor to RPA}.}

Given a $k$-tangle $T$ with internal $k_i$-discs $D_1, \ldots, D_n$
and elements $v_i\in V_{k_i}^{\sigma_i}$ we need to define
$Z(T)\br{v_1 \tensor \cdots \tensor v_n} \in V_k^{\sigma_0}$ and show
that the element we define is independent of rigid isotopy of $T$.
Finally we need to show that $Z$ is a morphism of colored operads.

Our line of proof will be very close to that of
Jones~\cite{Jones1999}, Section~4.2, although our construction will
differ slightly.  The advantages of this construction are an explicit
description of the action on $\Ltwo(M_i)$ and the fact that the
construction will also apply to the case of bimodule homomorphisms
$\Hom_{*-*}(\Ltwo(M_j),\Ltwo(M_k))$ rather than just endomorphisms
$\End_{*-*}(\Ltwo(M_i))$.

Consider first the case of a $(+,k)$-tangle $T$.  We will define
$Z(T)\br{v_1 \tensor \cdots \tensor v_n} \in
\End_{N-M_{-t}}(\Ltwo(M_r))\isom N'\cap M_{2r+t}$, where $k=2r+t$, $t=0,1$.

We will say that a tangle is in {\em standard form} if it is a
vertical concatenation of basic tangles of the following three types:
\vspace{10pt}
\begin{figure}[htbp]
\begin{center}
\psfrag{..}{$.\;.$}
\psfrag{TI=}{$TI_{k,p}=$}
\psfrag{TE=}{$TE_{k,p}=$}
\psfrag{TS=}{$TS_{k,l,p}=$}
\psfrag{k}{\scriptsize $k$}
\psfrag{k+1}{\scriptsize $k+1$}
\psfrag{p1}{$\underbrace{\hspace{102pt}}_{p}$}
\psfrag{p2}{$\overbrace{\hspace{102pt}}^{p}$}
\psfrag{p3}{$\overbrace{\hspace{120pt}}^{p}$}
\psfrag{k1}{$\underbrace{\hspace{25pt}}_{k}$}
\psfrag{l}{$\underbrace{\hspace{25pt}}_{l}$}
\includegraphics{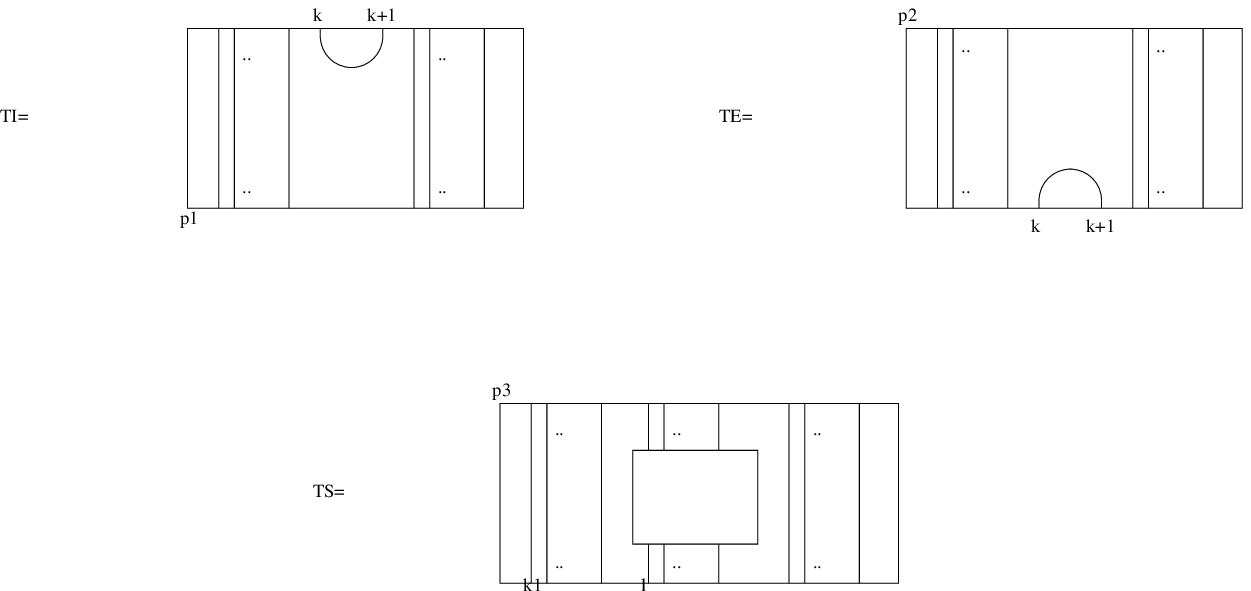}
\end{center}
\end{figure}

\noindent Any rigid planar $(+,k)$-tangle can be put in standard form
by rigid planar isotopy.

Note that the basic tangles have a different number of marked points
at the top and bottom of the box.  A tangle with sign $\sigma$, $j$
marked points at the bottom and $k$ marked points at the top will be
called a $(\sigma,j,k)$-tangle.

It may seem that we are introducing more general types of tangles
here, but we are really just used a variety of ``multiplication''
tangles rather than just the standard one in Figure~\ref{figure:
multiplication}.  A $(\sigma,j,k)$-tangle can be considered as a
$(\sigma,\frac{j+k}{2})$-tangle together with additional data given by
$j$.  The additional data determines the type of multiplication to be
used.  For example multiplication of a $(5,1)$-tangle and a
$(3,5)$-tangle is given by
\vspace{10pt}
\begin{figure}[htbp]
\begin{center}
\includegraphics{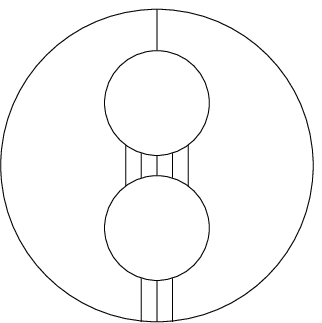}
\end{center}
\end{figure}
\vspace{-10pt}

In order to define $Z(T)$ for a general rigid tangle we will first
define it for the basic tangles.  Given a basic
$(+,2j+1+t,2k+1+t)$-tangle ($t=0,1$), $Z(T)$ will be an element of
$\Hom_{N-M_{-t}}(\Ltwo(M_j),\Ltwo(M_k))$.  We first need a preparatory
lemma.

\begin{lemma}
\label{lemma: alpha, beta, tilde(beta)}
Define $\alpha:\Ltwo(M)\rightarrow \Ltwo(M_1)$ to be $\delta$ times
the inclusion map.  Identifying $\Ltwo(M_1)$ with $\Ltwo(M)\tlN
\Ltwo(M)$ using $u_1$ from Proposition~\ref{prop: L^2 tensor and alg
tensor} allows us to write $\alpha$ as $\alpha(x)=\sum_b xb\tlN b^*
=\sum_b b \tlN b^*x$.  

\noindent Define maps
$\beta,\wtilde{\beta}:\Ltwo(M)\rightarrow \Ltwo(M_1)$ by
$\beta:x\mapsto x\tlN 1=xE_1$, $\wtilde{\beta}:x\mapsto 1\tlN x=E_1x$
(note that $\id \tlN \wtilde{\beta} = \beta \tlN \id$).  Then $\alpha,
\beta, \wtilde{\beta}$ are all continuous,
$\alpha\in\Hom_{M-M}\br{\Ltwo(M),\Ltwo(M_1)}$,
$\beta\in\Hom_{M-N}\br{\Ltwo(M),\Ltwo(M_1)}$,
$\wtilde{\beta}\in\Hom_{N-M}\br{\Ltwo(M),\Ltwo(M_1)}$ and
\begin{align*}
\alpha^*=\delta e_2
  &: x\tN y \mapsto xy \\
\beta^*
  &: x\tN y \mapsto xE(y) \\
\wtilde{\beta}^*
  &: x\tN y \mapsto E(x)y
\end{align*}
\begin{proof}
As we noted above, $\alpha$ is just $\delta$ times the inclusion map,
so $\alpha$ is continuous and $\alpha^*=\delta e_2$.  $\delta e_2
\hatr{xE_1y}=\delta^2 E_M(xe_1 y)=xy$.  
$\alpha\in\Hom_{M-M}\br{\Ltwo(M),\Ltwo(M_1)}$ because inclusion
$\Ltwo(M)\hookrightarrow \Ltwo(M_1)$ preserves left multiplication by
$M$ and also the right action of $M$ (since $J_1|_{\Ltwo(M)}=J_0$).

$\beta$ is continuous because $\varphi(E_1x^*xE_1)
=\delta\varphi(E_N(x^*x)E_1)=\varphi(E_N(x^*x))=\varphi(x^*x)$.
$\beta^*(x\tlN y)=xE(y)$ because
$
\ip{x \tlN y,z \tlN 1}
=\varphi\br{E(z^*x)y} =\varphi\br{z^*xE(y)} =\ip{xE(y),z} .
$
$\beta$ is clearly left $M$-linear.  To show right $N$-linearity
observe that
\[
\hatt{x}\cdot n \tN \hatt{1}
 = \br{J_0 n^* J_0 \hatt{x}} \tN \hatt{1}
 = J_1 \br{\hatt{1} \tN n^* J_0 \hatt{x}}
 = J_1 \br{\hatt{n^*} \tN J_0 \hatt{x}}
 = J_1 n^* J_1 \br{\hatt{x} \tN \hatt{1}} .
\]
$\wtilde{\beta}=J_1\beta J_0$ and hence $\wtilde{\beta}$ is continuous
and $N-M$ linear.  $\wtilde{\beta}^*(x\tlN y)=E(x)y$ because
\[
\ip{x \tN y,1 \tN z}
=\varphi\br{z^*E(x)y} =\ip{E(x)y,z} .
\]
\end{proof}
\end{lemma}

We can now define the elements of
$\Hom_{*-*}\br{\Ltwo(M_j),\Ltwo(M_k)}$ associated to each basic tangle.
\begin{align*}
Z\br{TI_{2i,p}} 
&= \delta^{-1/2} \br{\id}^{\tlN^{i-1}} \tN \alpha \tN
  \id \tN \id \cdots \tN \id \\
Z\br{TI_{2i-1,p}}
&= \delta^{1/2} \begin{cases}
    \wtilde{\beta}^* \tN \id \tN \id \tN \cdots \tN \id
     & i=1, \\
    \br{\id}^{\tlN^{i-2}} \tN \beta^* \tN \id \tN \cdots \tN \id
     & i\geq 2
   \end{cases} \\
Z\br{TE_{k,p}}
&= Z\br{TI_{k,p}}^* \\
Z\br{TS_{2i+t,k,p}}(R)
&= \br{\id}^{\tlN^i} \tN R \tN \id \tN \cdots \tN \id
\end{align*}

\noindent Letting $x=x_1 \tN x_2 \tN \cdots \tN x_k$ ($p=2k$ or
$2k+1$) we can see these maps explicitly as
\begin{align*}
Z\br{TI_{2i,p}}: x 
&\mapsto \delta^{-1/2} x_1 \tN \cdots \tN x_{i-1} \tN
  x_ib \tN b^* \tN x_{i+1} \tN \cdots \tN x_k \\
&= \delta^{-1/2} x_1 \tN \cdots \tN x_{i-1} \tN
  b \tN b^*x_i \tN x_{i+1} \tN \cdots \tN x_k \\
Z\br{TI_{2i-1,p}}: x
&\mapsto \begin{cases}
    \delta^{1/2} 1 \tN x_1 \tN x_2 \tN \cdots \tN x_k 
     & i=1, \\
    \delta^{1/2} x_1 \tN \cdots x_{i-1} \tN 1 \tN x_i 
     \tN \cdots \tN x_k 
     & i\geq 2
   \end{cases} \\
Z\br{TE_{2i,p}}: x 
&\mapsto \delta^{-1/2} x_1 \tN \cdots \tN x_{i-1} \tN
  x_ix_{i+1} \tN x_{i+2} \tN \cdots \tN x_k \\
Z\br{TE_{2i-1,p}}: x
&\mapsto \begin{cases}
    \delta^{1/2} E(x_1)x_2 \tN x_3 \tN \cdots \tN x_k 
     & i=1, \\
    \delta^{1/2} x_1 \tN \cdots x_{i-2} \tN x_{i-1}E(x_i) 
     \tN x_{i+1} \tN \cdots \tN x_k 
     & i\geq 2
   \end{cases} \\
Z\br{TS_{2i,2j-t,p}}(R):x
&\mapsto x_1 \tN \cdots x_i \tN
    \sqbr{\pi_{j-1}(R)\br{x_{i+1} \tN \cdots \tN x_{i+j}}} 
    \tN x_{i+j+1} \tN \cdots \tN x_k \\
Z\br{TS_{2i+1,2j+t,p}}(R):x
&\mapsto x_1 \tN \cdots x_i \tN
    \sqbr{\pi_{j}(R)\br{x_{i+1} \tN \cdots \tN x_{i+j+1}}} 
    \tN x_{i+j+2} \tN \cdots \tN x_k 
\end{align*}
where $t=0$ or $1$.

Finally let us check that a basic tangle with $2k-t$ lower strings and
$2k-t$ or $2k-t\pm 2$ upper strings defines an element of
$\End_{N-M_{-t}}\br{\Ltwo(M_k)}$ or
$\Hom_{N-M_{-t}}\br{\Ltwo(M_k),\Ltwo(M_{k\pm 1})}$ respectively
(recall that we are currently only considering the case of tangles
with $\sigma_0=+$).  All of these maps are left-$N$-linear and all are
right $M$-linear except for three cases.
$Z\br{TE_{2i-1}}=\br{\id}^{\tlN^{i-2}} \tlN \beta$ and
$Z\br{TI_{2i-1}}=\br{\id}^{\tlN^{i-2}} \tlN \beta^*$ are both
$N$-linear only on the right, but this is all that is required since
both have an odd number of strings.  $Z\br{TS_{2i+t,2j+1-t}}(R)=
\br{\id}^{\tlN^i} \tlN R$ is also $N$-linear only on the right, but
again has an odd number of strings.

Now, given a rigid planar $(+,2r+t)$-tangle $T$ (where $t=0,1$) we can
define $Z(T)\br{v_1 \tensor \cdots \tensor v_n} \in
\End_{N-M_{-t}}(\Ltwo(M_r))\isom N'\cap M_{2r+t}$ by putting $T$ in 
standard form by rigid planar isotopy and then composing the maps we
get from the basic tangles in the standard form.  We need to show that
the element of $\End_{N-M_{-t}}(\Ltwo(M_r))$ that we obtain is
independent of the particular standard picture we chose.  Then the map
$Z$ will not only be well-defined on a tangle $T$, but invariant under
rigid isotopy as any two isotopic tangles can be put in the same
standard form.

At this point our argument starts to closely resemble that of Jones.
If we have two standard pictures that are equivalent by a rigid planar
isotopy $h$ then, by contracting the internal discs to points and
putting the isotopy in general position, we see that the isotopy
results in a finite sequence of changes from one standard form to
another.  These changes happen in one of two ways
\begin{description}
\item{(i)} The $y$-coordinate of some string has a point of inflection
and the picture, before and after, looks locally like one of the
following.
\begin{figure}[htbp]
\begin{center}
\psfrag{..}{$.\;.$}
\psfrag{1}{(1)}
\psfrag{2}{(2)}
\psfrag{3}{(3)}
\psfrag{4}{(4)}
\psfrag{2i}{$\underbrace{\hspace{20pt}}_{2i}$}
\psfrag{2i+1}{$\underbrace{\hspace{20pt}}_{2i+1}$}
\includegraphics{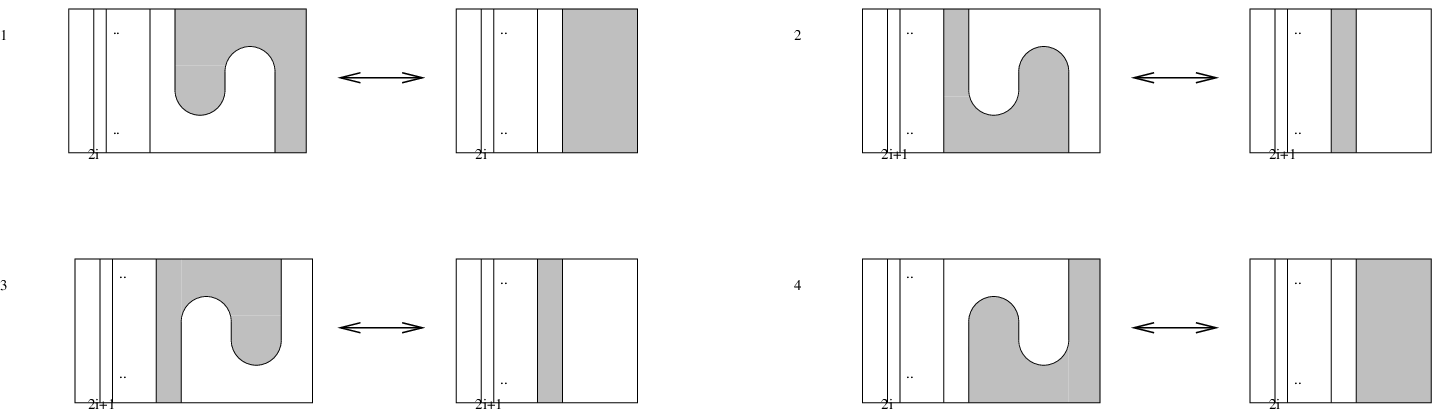}
\end{center}
\end{figure}

\item{(ii)} The $y$-coordinates of two internal discs, or maxima/minima
of strings (which we will refer to as caps/cups respectively),
coincide and change order while the $x$-coordinates remain distinct.
\end{description}

\noindent In case (i) note that $\alpha^*\beta=\alpha^*\wtilde{\beta}= \beta^*\alpha
=\wtilde{\beta}^*\alpha=\id$.  Then
\begin{description}
\item{(1)} This is either $\br{\alpha^*\wtilde{\beta}} \tlN 
(\id)^{\tlN^j}=\id$ ($i=0$), or ($i\geq 1$),
\[
(\id)^{\tlN^{i-1}}\tN \sqbr{\br{\id \tN \alpha^*}
 \br{\beta \tN \id}} \tN (\id)^{\tlN^j}
=(\id)^{\tlN^i}\tN \alpha^*\wtilde{\beta} \tN 
 (\id)^{\tlN^j}
=\id .
\]
\item{(2)} $(\id)^{\tlN^i} \tN \beta^*\alpha \tN
(\id)^{\tlN^j}=\id$.
\item{(3)} $(\id)^{\tlN^i} \tN \alpha^*\beta \tN
(\id)^{\tlN^j}=\id$.
\item{(4)} If $i=0$ then $\br{\wtilde{\beta}^*\alpha} \tlN 
(\id)^{\tlN^j}=\id$.  For $i\geq 1$,
\[
(\id)^{\tlN^{i-1}}\tN \sqbr{\br{\beta^* \tN \id}
 \br{\id \tN \alpha}} \tN (\id)^{\tlN^j}
=(\id)^{\tlN^i}\tN \wtilde{\beta}^*\alpha \tN 
 (\id)^{\tlN^j}
=\id .
\]
\end{description}

\noindent In case (ii) we have two cups, two caps or a cap and a cup
passing each other.  It is trivial to see that the two compositions of
maps from the basic tangles are the same in cases when the caps/cups
are separated and affect different parts of the tensor product.
Otherwise one uses associativity of the tensor product and the fact
that $\sum xb \tlN b^*=\sum b \tlN b^*x$.

In the case of one box and one cap/cup the maps again affect different
parts of the tensor product except for the following two cases (and
their adjoints).
\begin{figure}[htbp]
\begin{center}
\psfrag{..}{$..$}
\psfrag{1}{(i)}
\psfrag{2}{(ii)}
\psfrag{R}{$R$}
\psfrag{2i+1}{$\underbrace{\hspace{36pt}}_{2i+1}$}
\psfrag{2i+1 too}{\hspace{-2pt}$\underbrace{\hspace{18pt}}_{2i+1}$}
\includegraphics{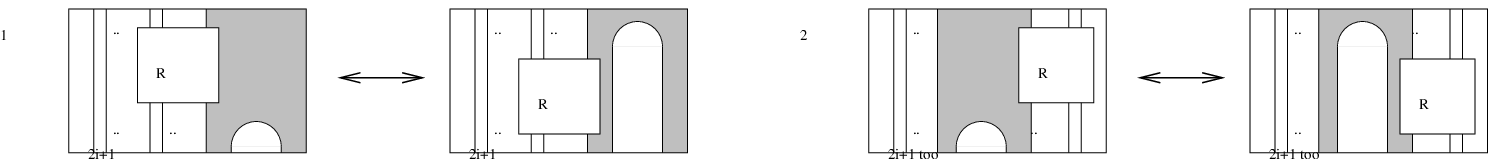}
\end{center}
\end{figure}

\noindent In case (i) we may assume that there is at most one string
to the left of the box since any other strings involve parts of the
tensor product where the tangles act as the identity.  Then the first
map is $R(x_1 \tlN \cdots \tlN x_{i-1} \tlN x_ix_{i+1})$ while the
second is $R(x_1 \tlN \cdots \tlN x_i)x_{i+1}$.  The two
expressions are the same since $R\in \End_{-M}$.  In
(ii) we can similarly assume that $i=0$.  Then $R(x_1x_2 \tlN x_3
\tlN \cdots \tlN x_j) =x_1R(x_2 \tlN x_3 \tlN \cdots \tlN x_j)$ because
$R\in \End_{M-}$.

The last case to consider is that of two boxes.  The only case which
does not involve distinct parts of the tensor product is the following.
\begin{figure}[htbp]
\begin{center}
\psfrag{..}{$\cdots$}
\psfrag{R}{$R$}
\psfrag{S}{$S$}
\includegraphics{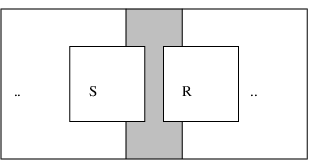}
\end{center}
\end{figure}
\vspace{-20pt}

\noindent Using the previous cases we see
that the following pictures represent the same linear maps.
\begin{figure}[htbp]
\begin{center}
\psfrag{..}{$..$}
\psfrag{R}{$R$}
\psfrag{S}{$S$}
\includegraphics{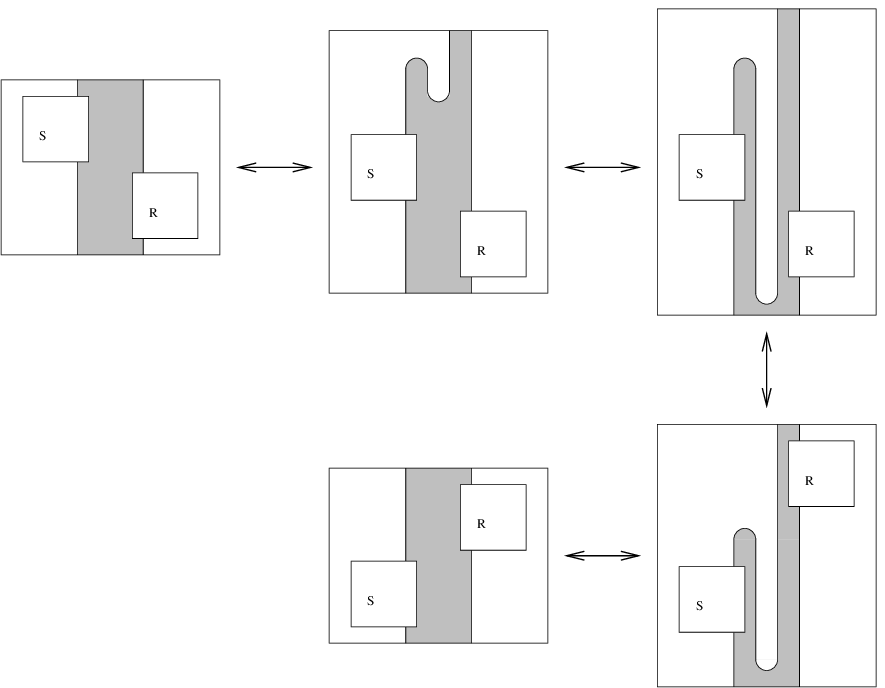}
\end{center}
\end{figure}

This completes our argument that $Z(T)$ does not depend on the choice
of standard picture for $T$ and is thus also rigid isotopy invariant.
For $(-,k)$-tangles just add a string to the left and hence obtain an
element of $N'\cap M_i$.  Choose a standard picture for $T$ which
leaves this string straight up.  Then the set of basic tangles which
make up $T$ cannot include $TI_{1,*}$, $TE_{1,*}$ or $TS_{0,*,*}$.
But these are the only basic tangles defining maps that are not
$M-$linear.  Hence we have an element of $M'\cap M_i$.

$Z$ is thus a map from $\PAr$ to $\Hom(V)$.  It is clearly an operad
morphism as we have defined $Z$ by composition of linear maps.  The
*-planar algebra property is also obvious from the way we have defined
$Z$.  Thus $(Z,V)$ is a rigid *-planar algebra.  It remains to prove
that properties (1) through (6) are satisfied.

\begin{description}
\item{(1)}  $\delta_1=\delta$ because
$Z\br{\includegraphics[-1pt,5pt][20pt,20pt]{figures/loop_black.eps}}\hatt{n}
=Z\br{\includegraphics[-1pt,5pt][20pt,20pt]{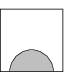}}
  \br{\delta^{1/2} n}
=\delta E_N(n)
=\delta n$
($n \in N$).

Similarly $\delta_2=\delta$ because 
$Z\br{\includegraphics[-1pt,5pt][20pt,20pt]{figures/loop_white.eps}}
\in \End_{M-M}\br{\Ltwo(M)}$ is given by
\[
Z\br{\includegraphics[-1pt,5pt][20pt,20pt]{figures/loop_white.eps}}\hatt{x}
=Z\br{\includegraphics[-1pt,5pt][24pt,20pt]{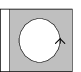}}\hatt{x}
=Z\br{\includegraphics[-1pt,5pt][24pt,20pt]{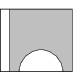}}
  \br{\delta^{-1/2}\sum_b xb \tlN b^*}
=\delta^{-1} \sum_b xbb^*
=\delta x .
\]

\item{(2)}  $\iota_{2k+1}^{+}:N'\cap M_{2k} \rightarrow N'\cap
M_{2k+1}$ is the usual inclusion map by construction (both ($2k+1$)- and
($2k+2$)-boxes are defined by their action on $\Ltwo(M_k)$).

$\iota_{2k+2}^{+}:N'\cap M_{2k+1} \rightarrow N'\cap M_{2k+2}$
involves an additional tensoring with the identity on the right, but
by Prop~\ref{prop: EN and PP multi-steps are the same} this just
changes the representation $\pi_k$ on $\Ltwo(M_k)$ to $\pi_{k+1}$ on
$\Ltwo(M_{k+1})$.  Hence $\iota_{2k+2}^{+}$ is also the usual
inclusion.

The result for $\iota_k^-$ is then immediate.

\item{(3)} $\gamma_k^-:V_k^-\rightarrow V_{k+1}^+$ is the
inclusion map $M'\cap M_k \hookrightarrow N'\cap M_k$ by construction.

$\gamma_k^+$ maps $R\in N' \cap M_{k-1}$ to $\id \tlN R$.  By
Prop~\ref{prop: EN and PP multi-steps are the same}
$\pi_{j+1}(M_1)=\pi_0(M_1) \tlN \br{\id}^{\tlN^j}$ and hence $\id \tlN
R \in M_1' \cap M_{k+1}$.

\item{(4)} By construction stacking boxes is multiplication of the
corresponding linear maps.

\item{(5)}  We want to show that $Z(TI_{k,k+1})Z(TE_{k,k+1})=E_k$.  By (2)
it suffices to show that $Z(TI_{k,2k})Z(TE_{k,2k})=\pi_{k-1}(E_k)$.
Let $x=x_1 \tlN \cdots \tlN x_k$.  Note that by definition
\[
Z(TE_{k,2k})\hatt{x}
=\begin{cases}
   \delta^{1/2} x_1 \tN \cdots \tN x_{r-1} E(x_r) \tN \cdots \tN x_k
   &k=2r-1 \\
   \delta^{-1/2} x_1 \tN \cdots \tN x_r x_{r+1} \tN \cdots \tN x_k
   &k=2r
  \end{cases}
\] 
so that $Z(TE_{k,2k})$ is simply $\delta^{1/2} E_{M_{k-2}}$, using
Proposition~\ref{prop: properties of tensorN}.  The map $Z(TI_{k,2k})$ is the
adjoint of $Z(TE_{k,2k})$, which is just $\delta^{1/2}$ times the
inclusion map.  Hence $Z(TI_{k,2k})Z(TE_{k,2k})=\delta
\pi_{k-1}(e_k)=\pi_{k-1}(E_k)$.

\item{(6)}  Using (5) we have
\begin{figure}[htbp]
\hspace{40pt}
\psfrag{x}{$x$}
\psfrag{y}{$E_{M_{k-1}}(x)$}
\psfrag{..}{$.\;.$}
\psfrag{k}{\hspace{-4pt}$\underbrace{\hspace{24pt}}_{\scriptsize k}$}
\psfrag{k2}{$\underbrace{\hspace{44pt}}_{\scriptsize k}$}
\psfrag{=ExE}{$=E_{k+1}xE_{k+1}
  =\delta E_{k+1}E_{M_{k-1}}(x)
  =\delta 
$}
\includegraphics{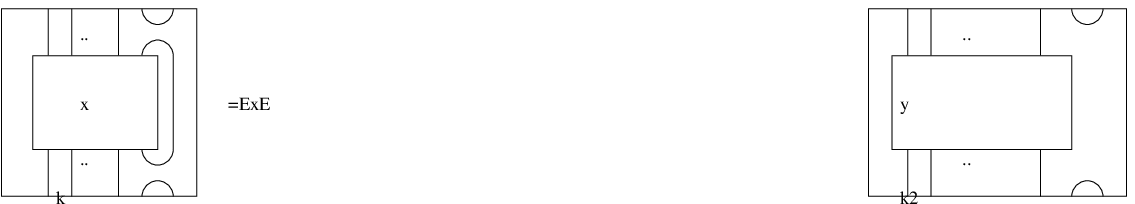}
\end{figure}

\noindent We then apply a tangle to both sides to close up the caps
and cups, which simply yields $\delta^2$.  After dividing by $\delta^2$
we obtain
\begin{figure}[htbp]
\hspace{150pt}
\psfrag{x}{$x$}
\psfrag{..}{$.\;.$}
\psfrag{k}{\hspace{-4pt}$\underbrace{\hspace{24pt}}_{\scriptsize k}$}
\psfrag{=EM}{$=\delta E_{M_{k-1}}(x)$}
\includegraphics{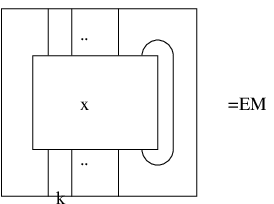}
\end{figure}

\item{(7)} Let $R\in N'\cap M_{2k}$ or $N' \cap M_{2k+1}$.  Then
(adding a string on the right if necessary)
\begin{figure}[htbp]
\hspace{70pt}
\psfrag{R}{$R$}
\psfrag{..}{$.\;.$}
\psfrag{k}{$\underbrace{\hspace{50pt}}_{\scriptsize 2k+2}$}
\psfrag{Z}{$Z
  \begin{pmatrix} & \hspace{60pt} \\ & \\ & \end{pmatrix}
  =Z\br{TE_{2,2k+2}} \br{\id \tN R} Z\br{TI_{2,2k+2}}$}
\includegraphics{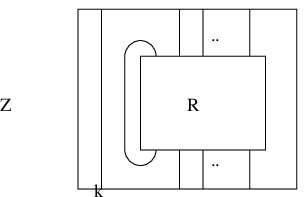}
\end{figure}

\noindent Calling this element $S$ and letting $x=x_1 \tlN \cdots \tlN
x_{k+1}$ we have 
\begin{align*}
S \hatt{x} 
&= Z\br{TE_{2,2k+2}} \br{\id \tN R} 
  \br{\delta^{-1/2} \sum_b b \tN b^*x_1 \tN x_2 \tN \cdots \tN x_{k+1}}\\
&= Z\br{TE_{2,2k+2}} 
  \br{\delta^{-1/2} \sum_b b \tN R\br{b^*x_1 \tN x_2 \tN \cdots \tN x_{k+1}}}\\
&= \delta^{-1} \sum_b b R\br{b^*x_1 \tN x_2 \tN \cdots \tN x_{k+1}} \\
&= \delta^{-1} \sum_b b R b^* \hatt{x} 
 = \delta E_{M'}(R) \hatt{x}
\end{align*}
which proves that $S=\delta E_{M'}(R)$ as required.

An immediate consequence is that $E_2 \br{\id_{\Ltwo(M)} \tlN R}
E_2=E_2 \br{\id \tlN \delta E_{M'}(R)}$, as we see below.
\begin{figure}[htbp]
\begin{center}
\psfrag{R}{$R$}
\psfrag{=}{$=$}
\psfrag{..}{$.\;.$}
\includegraphics{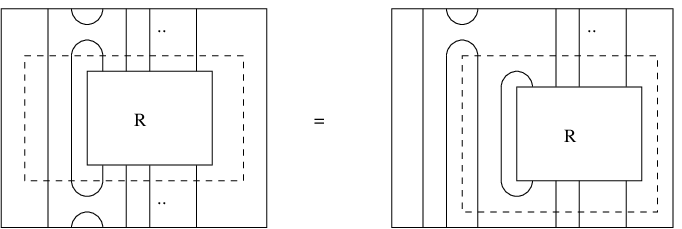}
\end{center}
\end{figure}
\vspace{-20pt}

Applying this to $N \subset M_j$ we obtain
\begin{equation}
\label{eq: e2Re2}
E_{[j,2j+1]}
\br{\id_{\Ltwo(M_j)} \tN R} E_{[j,2j+1]}=E_{[j,2j+1]}
\br{\id_{\Ltwo(M_j)} \tN \delta^{j+1} E_{M_j'}(R)}
\end{equation}
where $E_{[j,2j+1]}$ is one of the (scaled) multi-step Jones
projections from Theorem~\ref{multi-step bc thm}.  Multiplying out the
expression defining $E_{[j,2j+1]}$ in terms of $E_i$'s we obtain
\begin{figure}[htbp]
\begin{center}
\psfrag{..}{$.\;.$}
\psfrag{j+1}{$\underbrace{\hspace{30pt}}_{j+1}$}
\psfrag{E=}{$E_{[j,2j+1]}=$}
\includegraphics{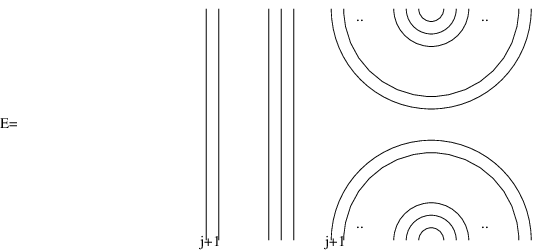}
\end{center}
\end{figure}
\vspace{-5pt}

\noindent Writing a thick string for $j+1$ regular strings, equation
(\ref{eq: e2Re2}) yields
\begin{figure}[htbp]
\begin{center}
\psfrag{R}{$R$}
\psfrag{EMR}{$E_{M_j'}(R)$}
\psfrag{=}{$=\delta^{j+1}$}
\psfrag{..}{$.\;.$}
\includegraphics{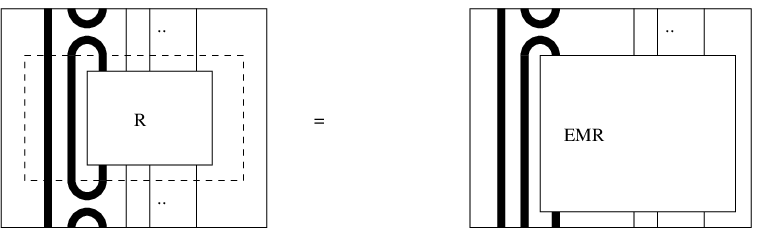}
\end{center}
\end{figure}
\vspace{-20pt}

\noindent and hence
\begin{figure}[htbp]
\begin{center}
\psfrag{R}{$R$}
\psfrag{EMR}{$E_{M_j'}(R)$}
\psfrag{=}{$=\delta^{j+1}$}
\psfrag{..}{$.\;.$}
\includegraphics{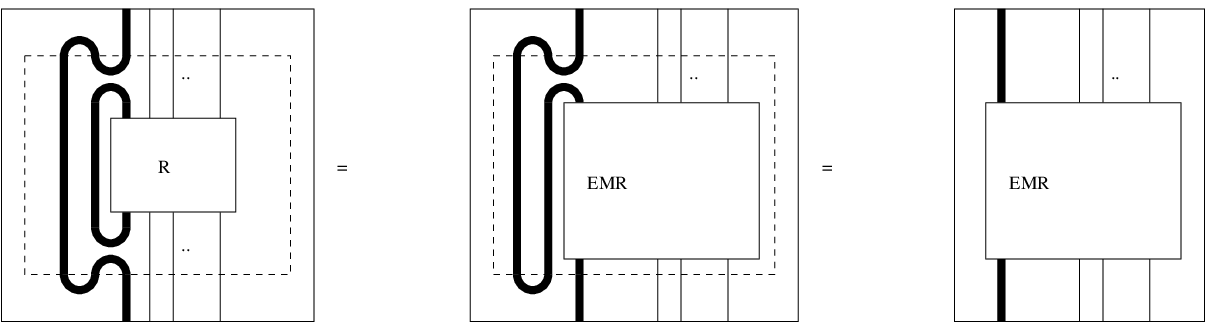}
\end{center}
\end{figure}
\vspace{-20pt}

\noindent which completes the proof of (7).
\end{description}

\noindent The proof of uniqueness is exactly the same as that in
Jones~\cite{Jones1999} 4.2.1.  The only difference is that our
property (4) is not required by Jones because it is built into his
axioms.

\proofend

\begin{remark}
The method of proof in (7) yields more general results.  In fact any
equation involving diagrams that holds for a general finite index
subfactor will also hold with thick strings (i.e. multiple strings) in
place of regular strings.  Let us loosely describe the procedure.  By
introducing extra caps and cups at top and bottom and by inserting
additional closed loops we may assume that any tangle has a standard
form made up of shifted boxes, $TS$, and Temperly-Lieb tangles.  Note
that the process of inserting closed loops can be accommodated simply
by dividing by $\delta$ and addition of caps/cups at top and bottom
can be inverted by applying an annular tangle to turn these into
closed loops.

We can thus write an equivalent equation in terms of shifts
(i.e. tensoring with the identity) and Jones projections, with no
reference to tangles.  Applying this to $N \subset M_j$ we can then
convert back to an equation in terms of tangles, but a shift by 2
becomes a shift by 2j and the Jones projections are 
\begin{figure}[htbp]
\begin{center}
\psfrag{..}{$.\;.$}
\psfrag{j+1}{$\underbrace{\hspace{30pt}}_{j+1}$}
\psfrag{kj+1}{$\underbrace{\hspace{30pt}}_{k(j+1)}$}
\includegraphics{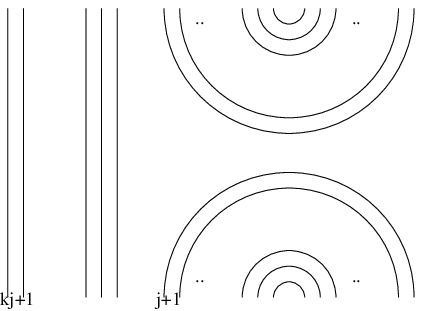}
\end{center}
\end{figure}
\vspace{-5pt}

\noindent We thus obtain the result for thick strings.
\end{remark}

\begin{prop}
\label{prop: cJ is S.S}
For $R \in N' \cap M_{2k+1}=V^+_{2k+2}$,
\[
\cj_{2k+2}(R) = S_k^* R^* S_k^* , \hspace{40pt}
\cj^{-1}_{2k+2}(R) = S_k R^* S_k .
\]
where $S_k$ on $\Ltwo(M_k)$ is the (in general unbounded) operator
defined by $\hatt{x}\mapsto \hatt{x^*}$.
\begin{proof}
For $k=0$
\[
\cj_2^{-1}(R)\hatt{x}
=Z\psfrag{R}{$R$}\br{\includegraphics[-1pt,14pt][37pt,37pt]{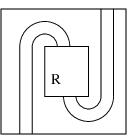}}\hatt{x}
=Z\br{\includegraphics[-1pt,14pt][37pt,37pt]{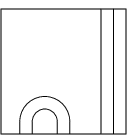}}
  \sqbr{\br{\sum_b x \tN R(b)} \tN b^*}
=\sum_b E_N\br{xR(b)} b^* .
\]
Now for any $R=ce_1 d\in M_1$,
\[
\sum_b E_N\br{xR(b)} b^*
=\sum_b E_N\br{xcE_N\br{db}} b^*
=E_N(xc)d
=\br{R^*\hatt{x^*}}^*
=S_0R^*S_0 \hatt{x} .
\]
Hence $\cj_2^{-1}(R)=S_0R^*S_0$.  Since $S^2=1$ we also have
$\cj_2(R)=S_0^*R^*S_0^*$.  The general result follows, as described in
the previous remark, by amplifying to the $j$-string case.
\end{proof}
\end{prop}

\begin{remark}
This extends the result of Bisch and Jones~\cite{BischJones2000} Prop
3.3 that $\cj(R)=JR^*J$ in the finite index extremal $\IIone$ case.
\end{remark}

\begin{cor}
The shift map is also given by
\[
sh(R)=\cj^{-1}_{2k+4}\cj_{2k+2}(R)
     = S_{k+1}\br{S_k R S_k} S_{k+1} , \hspace{40pt} R \in N' \cap M_{2k+1} .
\]
\begin{proof}
The first equality follows from
\begin{figure}[htbp]
\begin{center}
\psfrag{R}{$R$}
\psfrag{=}{$=$}
\includegraphics{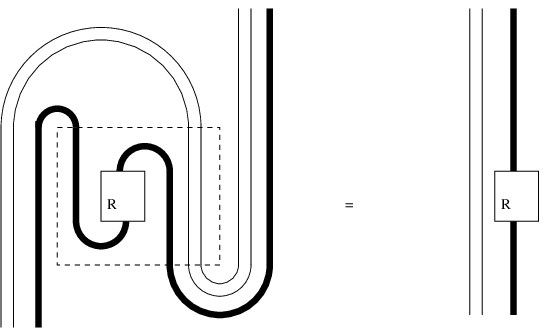}
\end{center}
\end{figure}
\vspace{-5pt}

\noindent and the second is then immediate by Prop~\ref{prop: cJ is S.S}.
\end{proof}
\end{cor}


Finally we show that our two definitions of the rotation agree.

\begin{lemma}
The rotation defined in Corollary~\ref{cor: PA rho periodic} agrees,
in the case of a rigid \cstar-planar algebra coming from a finite
index subfactor, with that in Definition~\ref{define rotation}.
\begin{proof}
Simply draw the appropriate tangle for $\delta^2
E_{M_k}(v_{k+1}E_{M'}(x v_{k+1}))$.  First note that $v_k$ is given by
\vspace{-10pt}
\begin{figure}[htbp]
\hspace{170pt}
\psfrag{..}{$\cdots$}
\includegraphics{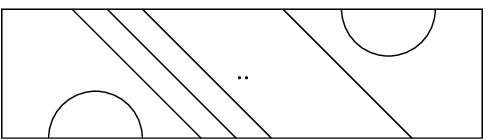}
\end{figure}
\vspace{-20pt}

\noindent and hence
\begin{figure}[htbp]
\hspace{40pt}
\psfrag{d=}{$\delta^2 E_{M_k}(v_{k+1}E_{M'}(x v_{k+1}))=$}
\psfrag{x}{$x$}
\includegraphics{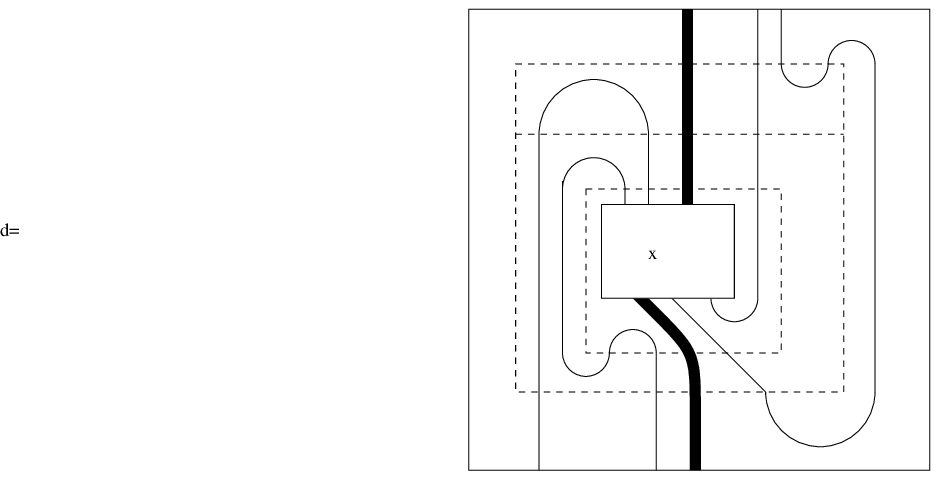}
\end{figure}

\end{proof}
\end{lemma}

\subsection{Modular structure} 
\label{section: modular structure}

Let $(Z^r,V)$ be a rigid \cstar-planar algebra.  We will extend $Z$ to
the full planar operad $\PA$.  Let $T$ be a $k$-tangle with internal
discs $D_1,\ldots,D_n$.  Take any planar isotopy $h_t$ taking $T$ to a
rigid planar $k$-tangle $T_0$.  Let $\theta_i$ be the amount of
rotation of $D_i$ under $h_t$.  Define
\[
Z(T)\br{v_1 \tN \cdots \tN v_n}
=Z^r(T_0) \br{\Delta_{(\sigma_1,k_1)}^{-\theta_1/2\pi}(v_1) \tN \cdots 
  \tN \Delta_{(\sigma_n,k_n)}^{-\theta_n/2\pi}(v_n)} .
\]

\begin{prop}
$Z(T)$ is independent of the choice of $h$ and $T_0$.
\begin{proof}
We use the same argument as in Jones~\cite{Jones1999}.  We can
surround each internal disc $D_i$ with a larger disc $\wtilde{D}_i$
such that for all $t\in [0,1]$ the images of these larger discs under
$h_t$ are disjoint and do not intersect any closed strings.  Let $x_i$
be the center of $D_i$.  Then there exists $r>0$ such that
$D(h_t(x_i),r)\subset h_t(D_i)^{\text{o}}$ for all $t$.  Define a rigid planar
isotopy $g$ in three steps
\begin{description}
\item{1.} Radially shrink $D_i$ to $D(x_i,r)$ while keeping 
$\mathbb{R}^2\backslash \wtilde{D}_i^{\text{o}}$ fixed.
\item{2.} Take $h_t$ on $\mathbb{R}^2\backslash \wtilde{D}_i^{\text{o}}$ 
and translation by $h_t(x_i)-x_i$ on $D(x_i,r)$ and interpolate in
any way in between.
\item{3.} Radially expand $D(x_i,r)+h_1(x_i)-x_i$ to $h_1(D_i)$ 
while keeping $\mathbb{R}^2\backslash h_1(\wtilde{D}_i)^{\text{o}}$ fixed.
\end{description}
Then $Z^r\br{g^{-1}(T_0)}=Z^r(T_0)$, so following $h$ with $g^{-1}$ we
may assume without loss of generality that $h$ is the identity outside
$\wtilde{D}_i^{\text{o}}$.

Suppose we have another such isotopy $\overline{h}$ taking $T$ to a
rigid planar tangle $\overline{T}_0$.  Then the rotation of $D_i$
under $\overline{h}$ is $\overline{\theta}_i=\theta_i+2\pi l_i$ for
some $l_i \in
\mathbb{Z}$.

$h^{-1}$ followed by $\overline{h}$ is a planar isotopy taking $T_0$
to $\overline{T}_0$ that is the identity outside $\wtilde{D}_i^{\text{o}}$ and
that rotates $D_i$ by $2\pi l_i$.  The mapping class group of
diffeomorphisms of the annulus that are the identity on the boundary is
generated by a single Dehn twist of $2\pi$.  Hence the difference
between $\overline{T}_0$ and $T_0$ is $\Delta_{(\sigma_i,k_i)}^{l_i}$
inside $D_i$, so
\[
\overline{T}_0
=\br{\br{T_0 \circ_1 \Delta_{(\sigma_1,k_1)}^{l_1}} \circ_2
  \Delta_{(\sigma_2,k_2)}^{l_2}}\cdots \circ_n
  \Delta_{(\sigma_n,k_n)}^{l_n} ,
\]
\[
Z^r(\overline{T}_0)
=\br{Z^r(T_0)} \circ \br{\Delta_{(\sigma_1,k_1)}^{l_1} \tN \cdots \tN
  \Delta_{(\sigma_n,k_n)}^{l_n}}
\]
and hence
\begin{align*}
Z^r(\overline{T}_0)\br{\!\Delta_{(\sigma_1,k_1)}^{-\overline{\theta}_1/2\pi}
  (v_1) \tN \cdots \tN
  \Delta_{(\sigma_n,k_n)}^{-\overline{\theta}_n/2\pi}(v_n)\!}
&=Z^r(T_0)\br{\!\Delta_{(\sigma_1,k_1)}^{l_1-\overline{\theta}_1/2\pi}(v_1)
  \tN \cdots \tN
  \Delta_{(\sigma_n,k_n)}^{l_n-\overline{\theta}_n/2\pi}(v_n)\!} \\ 
&=Z^r(T_0)\br{\Delta_{(\sigma_1,k_1)}^{-\theta_1/2\pi}(v_1)
  \tN \cdots \tN
  \Delta_{(\sigma_n,k_n)}^{-\theta_n/2\pi}(v_n)} .
\end{align*}
\end{proof}
\end{prop}

\begin{remark}
Now that $Z$ is well-defined, considering {\em any} two tangles $T$
and $T_0$ connected by a planar isotopy that rotates the internal disc
$D_i$ by $\theta_i$, we have
\[
Z(T)\br{v_1 \tN \cdots \tN v_n}
=Z(T_0)\br{\Delta_{(\sigma_1,k_1)}^{-\theta_1/2\pi}(v_1) \tN \cdots 
  \tN \Delta_{(\sigma_n,k_n)}^{-\theta_n/2\pi}(v_n)} .
\]

The extended map $Z$ is a mapping the set of {\em all} planar
$k$-tangles, modulo rigid planar isotopy, to $\Hom(V)$.  We still need
to establish that $Z$ is an operad morphism.
\end{remark}

\begin{prop}
For $r \in \mathbb{R}$,
\[
\br{\Delta_{(\sigma_0,k)}}^r Z(T) \br{v_1 \tN \cdots \tN v_n}
=Z(T)\br{\br{\Delta_{(\sigma_1,k_1)}}^r v_1 \tN \cdots \tN
 \br{\Delta_{(\sigma_n,k_n)}}^r v_n } .
\]
\begin{proof}
Since
$\br{\Delta_{(\sigma_i,k_i)}}^r
\br{\Delta_{(\sigma_i,k_i)}}^{-\theta_i/2\pi} v_i
=\br{\Delta_{(\sigma_i,k_i)}}^{-\theta_i/2\pi}
\br{\Delta_{(\sigma_i,k_i)}}^r v_i$, 
it suffices to prove the result for $T\in\PAr$.
Recall that $\br{\Delta_{(\sigma_i,k_i)}}^r =
\br{z_{(\sigma_i,k_i)}}^{-r/2} \br{\spdot}
\br{z_{(\sigma_i,k_i)}}^{r/2}$ and hence it suffices to check the
result for basic tangles, with $\Delta^r$ replaced by
$\br{z_{(\sigma,k)}}^{-r/2} \br{\spdot} \br{z_{(\sigma,j)}}^{r/2}$ for a 
$(\sigma,j,k)$-tangle.
Assume that $\sigma=+$.  Recall that $z_k=w_1 \cdots w_k$ and that the
$w_i$'s commute, so $z_k^s=w_1^s \cdots w_k^s$.  Let $s=r/2$.  Then

\begin{figure}[htbp]
\begin{center}
\psfrag{zTz=}{$z_{p+2}^{-s}TI_{2i-1,p}z_p^s=$}
\psfrag{..}{$\cdots$}
\psfrag{w-}{\scriptsize $\,w^{-s}$}
\psfrag{wt-}{\scriptsize $\,\wtilde{w}^{-s}$}
\psfrag{w+}{\scriptsize $w^{s}$}
\psfrag{wt+}{\scriptsize $\wtilde{w}^{s}$}
\includegraphics{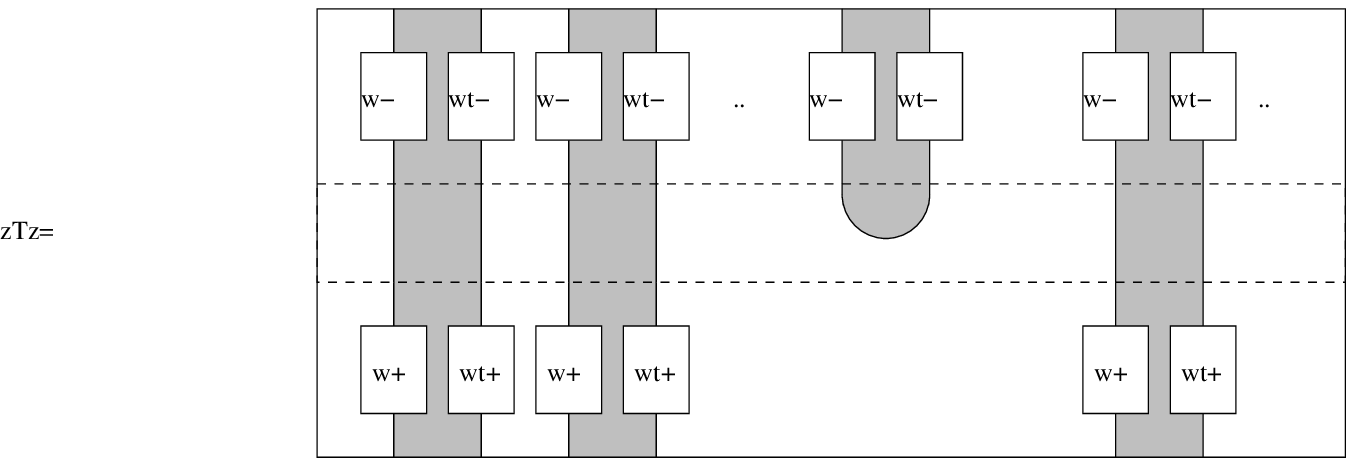}
\end{center}
\end{figure}
\vspace{-20pt}

\begin{figure}[htbp]
\begin{center}
\psfrag{zTz=}{$\phantom{z_{p+2}^{-s}TI_{2i-1,p}z_p^s}=$}
\psfrag{..}{$\cdots$}
\psfrag{w-}{\scriptsize $\,w^{-s}$}
\psfrag{w+}{\scriptsize $w^{s}$}
\includegraphics{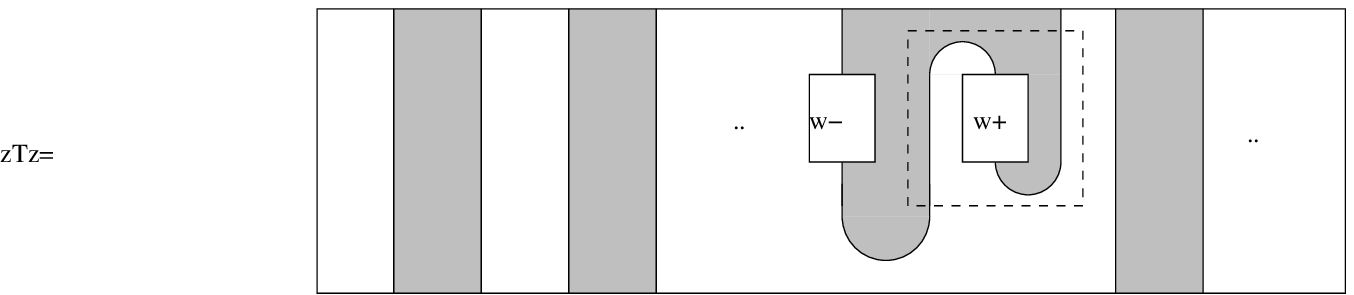}
\end{center}
\end{figure}
\vspace{-20pt}

\begin{figure}[htbp]
\begin{center}
\psfrag{zTz=}{$\phantom{z_{p+2}^{-s}TI_{2i-1,p}z_p^s}=$}
\psfrag{..}{$\cdots$}
\includegraphics{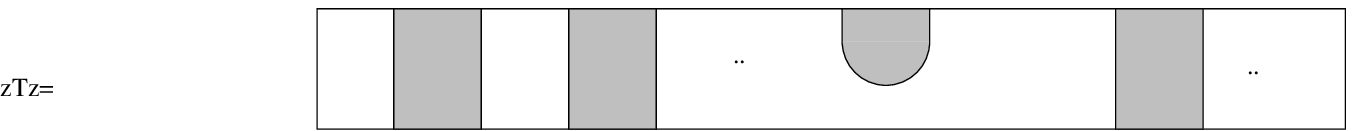}
\end{center}
\end{figure}

\noindent so $z_{p+2}^s TI_{2i-1,p} z_p^s=TI_{2i-1,p}$ and similarly
$z_{p+2}^s TI_{2i,p} z_p^s=TI_{2i,p}$.  Taking adjoints gives the
result for $TE_{k,p}$.  For $TS_{k,l,p}$ the result is trivial.  The
result for $(-,k)$-tangles follows by symmetry.
\end{proof}
\end{prop}

\begin{cor}
$Z$ is an operad morphism: $Z\br{S \circ_i T}=Z(S) \circ_i Z(T)$.
\end{cor}

\begin{definition}
We call the map $Z$ {\em the modular extension of the rigid planar
algebra $(Z^r,V)$}.
\end{definition}

\begin{remark}
The modular extension $(Z,V)$ is not a planar algebra because tangles
that are planar isotopy equivalent need not define the same linear
maps, but the difference comes down to rotations of the interior discs
and is thus controlled by the modular operators $\Delta$.  In turn all
the modular operators are controlled by $w$, the Radon-Nikodym
derivative of $\varphi$ with respect to $\varphi'$ on $V_1^+$.
\end{remark}

\begin{cor}
If $\varphi$ is tracial (for example in the case of
$(Z^{(N,M,E)},V^{(N,M,E)})$ for a finite index $\IIone$ subfactor with
trace-preserving conditional expectation $E$) then $\Delta=1$,
$\varphi'$ is tracial and $Z$ is invariant under full planar isotopy
so that we have a true planar algebra structure on the standard
invariant.
\end{cor}

\subsection{From rigid to spherical and back again}
\label{section: rigid to spherical}

From a rigid \cstar-planar algebra we have seen how to form the
modular extension where we are allowed to rotate the boxes in a tangle
but must pay a price in terms of $\Delta$ for doing so.  We will now
see how to modify $Z$ to obtain a spherical \cstar-planar algebra.

Let $(Z^r,V)$ be a rigid \cstar-planar algebra, $(Z,V)$ its modular
extension.  Given a planar tangle $T$ and a string $s:[0,1]\rightarrow
\mathbb{R}^2$ in $T$ define $\Theta(s)$ to be the {\em total angle}
along $s$, with $s$ parameterized so that $s$ bounds a black region on
its right.  $\Theta(s)$ may be computed as $\Theta(s)=\int_0^1
\frac{\dee \phi}{\dee t} \dee t$ where $\phi(t)$ is the angle at
$s(t)$, or as $\Theta(s)=\int_0^L \kappa \dee l$ where $L$ is the
length of the curve, $\kappa$ the curvature and $\dee l$ the length
element.  For example the string below has $\Theta(s)=\pi$.
\vspace{-30pt}
\begin{figure}[htbp]
\begin{center}
\includegraphics{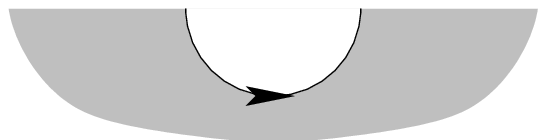}
\end{center}
\end{figure}
\vspace{-10pt}

\noindent For a planar tangle $T$ presented using boxes, define a {\em
spherically averaged tangle} $\mu(T)$ to be a tangle obtained from $T$
as follows: on every string $s$ insert a $(+,1)$-box containing
$w^{-\Theta(s)/4\pi}$.  Then define $Z^{\sph}:\PA \rightarrow \Hom(V)$
by $Z^{\sph}(T)=Z(\mu(T))$.

Note that $\mu(T)$ is not unique, but $T\mapsto Z(\mu(T))$ is
well-defined because the position of the inserted $(+,1)$-box does not
affect $Z(\mu(T))$.  To see this observe that moving the box along
the string may change the angle of the box by some amount $\theta$,
but
\[
\Delta_{+,1}^{-\theta/2\pi}\br{w^{-\Theta(s)/4\pi}}
=w^{-\theta/2\pi} w^{-\Theta(s)/4\pi} w^{\theta/2\pi}
= w^{-\Theta(s)/4\pi} .
\]

\begin{remark}
Note that it is important that we use boxes rather than discs at this
point.  If we defined $Z^{\sph}$ as above, but used discs rather than
boxes, we would introduce a multitude of complications.  For example
the multiplication tangle would be different under $Z^{\sph}$ and
would not be invariant under the standard embeddings of $V_i$ in
$V_{i+1}$.  Of course we could formulate a definition in terms of
discs, the essence of which would be using not $\Theta(s)$ but the
difference between $\Theta(s)$ and what it ``ought to be'' for a
string joining those two points, but the book-keeping is much cleaner
with boxes.
\end{remark}

\begin{thm}
\label{thm: construct spherical}
$(Z^{\sph},V)$ is a spherical \cstar-planar algebra with
$\delta^{\sph}=\lambda \delta$ where $\lambda=\varphi(w^{1/2})$.
\begin{proof}
$Z^{\sph}$ is invariant under rigid planar isotopy since the total
angle along a string does not change under rigid planar isotopy.  If
the initial angle is $\theta_{\mathrm{init}}$
and the final angle is $\theta_{\mathrm{fin}}$ then the total angle is
$\theta_{\mathrm{fin}}-\theta_{\mathrm{init}}+2\pi l$ for some
$l\in\mathbb{Z}$.  Under a rigid isotopy $\theta_{\mathrm{init}}$ and
$\theta_{\mathrm{fin}}$ are constant.  Since the total angle must vary
continuously under isotopy, the total angle must also be constant.

$Z^{\sph}$ is an operad morphism because, when we compose tangles $T$
and $S$ to form $T\circ_i S$, the total angle is additive for strings
from $T$ and $S$ that meet at $\curlyd D_i$ to form a single string in
$T\circ_i S$.  Multiplying powers of $w$ is of course also additive in
the exponent.  It is worth noting at this point that we are using the
fact that strings meet discs normally.

$(Z^{\sph},V)$ is in fact a rigid planar *-algebra.  Reflection
changes orientation and in doing so reverses the direction in which the
string is parameterized, but each of these two changes multiplies the
total angle by $-1$, so there no net change to the total angle.

Let $x \in V_1^+$.  Then
\newpage
\begin{figure}[htbp]
\hspace{44pt}
\psfrag{P=}{$\Phi^{\sph}(x)=Z^{\sph}
  \begin{pmatrix} &&&&&\\ &&&&&\\ &&&&&\end{pmatrix}
  =Z\begin{pmatrix} &&&&&\\ &&&&&\\ &&&&&\\ &&&&&\\ &&&&&\\ &&&&&\end{pmatrix}
  =\Phi\br{w^{1/2}x}$,}
\psfrag{x}{$\,x$}
\psfrag{w}{$w^{\frac{1}{2}}$}
\includegraphics{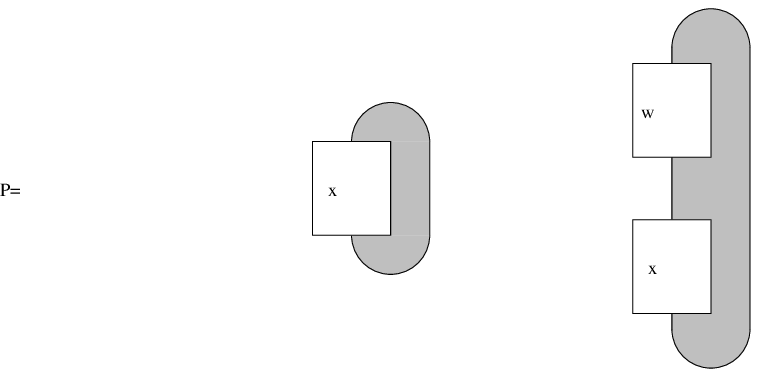}
\end{figure}


\begin{figure}[htbp]
\hspace{50pt}\psfrag{P=}{
\hspace{-16pt}$\hspace{-5pt}\br{\Phi^{\sph}}'(x)=Z^{\sph}
  \begin{pmatrix} &&&&&\\ &&&&&\\ &&&&&\end{pmatrix}
  =Z\begin{pmatrix} &&&&&\\ &&&&&\\ &&&&&\\ &&&&&\\ &&&&&\\ &&&&&\end{pmatrix}
  = \Phi'(w^{-1/2}x)
 = \Phi(w^{1/2}x)$.}
\psfrag{x}{$\,x$}
\psfrag{w}{$\!\!w^{-\frac{1}{2}}$}
\includegraphics{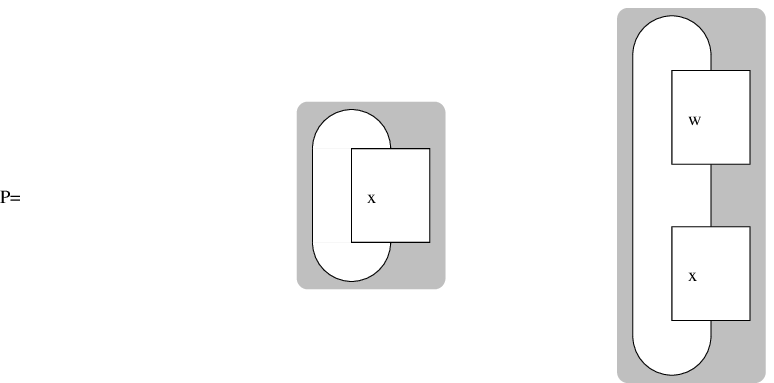}
\end{figure}

Thus $\Phi^{\sph}=\br{\Phi^{\sph}}'$ and so
$\delta_1^{\sph}=\Phi^{\sph}(1)=\br{\Phi^{\sph}}'(1)=\delta_2^{\sph}$ In
addition all Radon-Nikodym derivatives are $1$, $\Delta=1$ and
$\varphi^{\sph}$ is tracial.  $\Trace=\Phi^{\sph}$ is positive
definite because
\begin{figure}[htbp]
\hspace{55pt}
\psfrag{Tr=}{$\Trace(x^*x)=Z
  \begin{pmatrix} &&&&&\\ &&&&&\\ &&&&&\\ &&&&&\\ &&&&&\end{pmatrix}
  =Z\begin{pmatrix} &&&&&\\ &&&&&\\ &&&&&\\ &&&&&\\ &&&&&\\ &&&&&\end{pmatrix}
  =\Phi\br{w^{1/4}x^*xw^{1/4}}$,}
\psfrag{x}{$\!x^*x$}
\psfrag{w}{$w^{\frac{1}{2}}$}
\psfrag{w1}{$w^{\frac{1}{4}}$}
\includegraphics{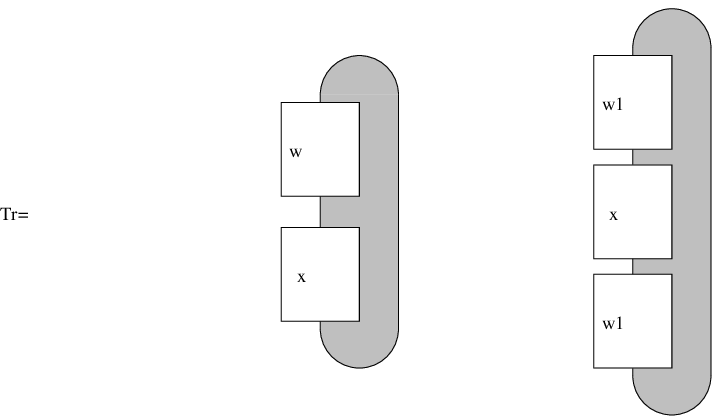}
\end{figure}

\noindent where we have used the fact that $w^r$ can move along
strings without changing $Z$.  Note that $\Trace(x^*x)=0$ iff
$xw^{1/4}=0$ iff $x=0$.  Hence $(Z^{\sph},V)$ is a spherical
\cstar-planar algebra.

Finally note that $\delta^{\sph}=\Phi^{\sph}(1)=\Phi(w^{1/2})=\delta
\varphi(w^{1/2}) =\lambda \delta$.
\end{proof}
\end{thm}

\begin{cor}
\label{cor: alg std inv may as well be IIone}
Let $(N,M,E)$ be a finite index subfactor.  Then there exists an
extremal $\IIone$ subfactor $\wtilde{N}\subset \wtilde{M}$ with a
lattice of higher relative commutants that is algebraically isomorphic
to that of $(N,M,E)$ (although of course the conditional expectations
and Jones projections may differ).
\begin{proof}
Let $(Z,V)$ be the rigid planar algebra constructed on the standard
invariant of $(N,M,E)$.  Then apply Popa~\cite{Popa1995} in the form
of Jones~\cite{Jones1999} Theorem 4.3.1 to $(Z^{\sph},V)$ to obtain
the extremal $\IIone$ subfactor $\wtilde{N}\subset \wtilde{M}$.
\end{proof}
\end{cor}

\begin{remark}
It is well known that, algebraically (i.e. without reference to the
conditional expectations), all standard invariants of finite index
subfactors can be realized using $\IIone$ subfactors.  One can tensor
$N \subset M$ with a $\III_1$ factor (this leaves the standard
invariant the same) to obtain a $\III_1$ subfactor.  Taking the
crossed product with the modular group one obtains a finite index
inclusion of $\II_{\infty}$ factors that splits as a $\II_1$ subfactor
tensored with a $\mathrm{I}_{\infty}$ factor.  Izumi~\cite{Izumi1993}
shows that this type $\II$ inclusion has the same principal graph, and
hence the same algebraic standard invariant, as the original $\III_1$
inclusion.  We thank Dietmar Bisch for bringing this to our attention.
\end{remark}

\begin{remarks}
The idea of the construction of $Z^{\sph}$ is to insert powers of $w$
to cancel the effects of rotations (implemented by $\Delta$).  One
might ask what would happen to a (non-spherical) \cstar-planar algebra
under this procedure, since $\Delta=1$ and we {\em already} have full
planar isotopy invariance.  Certainly we change the action of tangles,
but the only additional isotopy invariance we gain is the ability to
move strings past the point at infinity.  Let us examine the effects
of isotopy on $Z^{\sph}$ in this case.

Suppose we rotate a box labeled with $x$ by an angle $\theta$.
The change in total angle for the upper strings is the negative of the
corresponding change for the lower strings.
\begin{figure}[htbp]
\begin{center}
\includegraphics{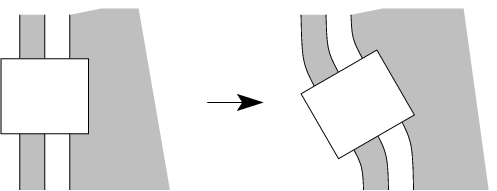}
\end{center}
\end{figure}
\vspace{-10pt}

\noindent Hence there is no change in $x$ because $z_k^{-\theta/4\pi}
x z_k^{\theta/4\pi}=x$ (since $z_k$ is central in $V_k$).

Considering strings alone, rigid planar isotopies to not change the total
angle.  Only spherical isotopy moving a string past the point at
infinity can change the total angle.  For example
\vspace{-10pt}
\begin{figure}[htbp]
\begin{center}
\psfrag{2p}{$2\pi$}
\psfrag{-2p}{$-2\pi$}
\includegraphics{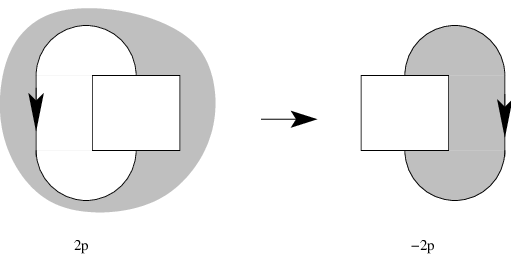}
\end{center}
\end{figure}
\vspace{-10pt}

\noindent The total angle changes by $-4\pi$, so the difference in
$Z^{\sph}$ is the insertion of an additional $w^1$ which we know is
the necessary correction to change $\varphi'$ to $\varphi$.
\end{remarks}

We can reverse the construction of $(Z^{\sph},V,w)$ from $(Z,V)$ as we
see below.  Note that in $(Z^{\sph},V,w)$ we have
$\trace(w^{1/2})=\lambda^{-1}\varphi'(w^{-1/2}w^{1/2}) =\lambda^{-1}$
and $\trace(w^{-1/2})=\lambda^{-1}\varphi(w^{1/2}w^{-1/2})
=\lambda^{-1}$.

\begin{thm}
\label{thm: spherical to rigid}
Let $(Z,V)$ be a spherical \cstar-planar algebra with modulus
$\delta^{\sph}$ and let $w$ be a positive, invertible element of
$V_1^+$ with $\trace(w^{1/2})=\trace(w^{-1/2})=\lambda^{-1}$.

For a planar tangle $T$ define the $\nu(T)$ to be the tangle obtained
from $T$ as follows: on every string $s$ insert a $(+,1)$-box
containing $w^{\Theta(s)/4\pi}$.  Define $Z^{\modular}:\PAr\rightarrow
\Hom(V)$ by $Z^{\modular}(T)=Z(\nu(T))$.  Let $Z^{\rig}$ be the
restriction of $Z^{\modular}$ to rigid planar tangles.

Then $(Z^{\rig},V)$ is a rigid \cstar-planar algebra with
$w^{\rig}=w$ and $\delta=\lambda \delta^{\sph}$.

\begin{proof}
Note that the position of the inserted $(+,1)$-box on a string does
not affect $Z$ because $(Z,V)$ is spherical.  As in the proof of
Theorem~\ref{thm: construct spherical} we have: (i) rigid planar
isotopy does not change the total angle along a string; (ii) the total
angle is additive under composition; (iii) the total angle in
invariant under reflection.  Hence $(Z^{\rig},V)$ is a rigid *-planar
algebra.

Let $\Trace=\Phi^{(Z,V)}$ and $\trace=\varphi^{(Z,V)}$.  As in
Theorem~\ref{thm: construct spherical}
$\Phi^{\rig}(x)=\Trace(w^{-1/2}x)$ which is positive definite and
$\delta_1^{\rig}=\Phi_{(+,1)}^{\rig}(1)=\Trace_1(w^{-1/2}x)=\delta^{\sph}
\trace(w^{-1/2})=\lambda^{-1} \delta^{\sph}$.  Similarly 
$\delta_2^{\rig}=\delta^{\sph}\trace(w^{1/2})=\lambda^{-1} \delta^{\sph}$.

Finally, $\br{\Phi^{\rig}}'(x)=\Trace(w^{1/2}x)=\Trace(w^{-1/2}(wx))
=\Phi^{\rig}(wx)$ so $w^{\rig}=w$.
\end{proof}
\end{thm}

\begin{remarks}
The condition $\trace(w^{1/2})=\trace(w^{-1/2})$ is only used to make
sure that $\delta^{\rig}_1=\delta^{\rig}_2$.  Given any positive
invertible $w$ in $V_1^+$ we can scale $w$ by
$\trace(w^{-1/2})/\trace(w^{1/2})$ to obtain $\overline{w}$ with
$\trace(\overline{w}^{1/2})=\trace(\overline{w}^{-1/2})$.

The two constructions $(Z,V)\mapsto (Z^{\sph},V,w)$ and
$(Z,V,w)\mapsto (Z^{\rig},V)$ are obviously inverse to each other.
\end{remarks}

\subsection{Planar algebras give subfactors}
\label{section: PA to subfactors}

Consider the two main results of the previous section.  Starting with
a rigid \cstar-planar algebra $(Z^r,V)$ we can construct the
associated spherical \cstar-planar algebra $(Z^{\sph},V,w)$ with
distinguished element $w$ and then construct the rigid \cstar-planar
algebra $\br{\br{Z^{\sph}}^{\rig},V}$.  This two part construction
simply reproduces $(Z^r,V)$.  By Popa's standard lattice
result~\cite{Popa1995}, applied in Jones~\cite{Jones1999} Theorem
4.3.1, there exists a subfactor with standard invariant
$(Z^{\sph},V)$.  If we can ``lift'' the second part of the planar
algebra construction to the subfactor level then we will have a
subfactor with $(Z^{(N,M,E)},V^{(N,M,E)})\isom(Z^r,V)$.

\begin{center}
\begin{tabular}{ccccccc}
\vspace{3pt} \\
 & & $(N,M,E_{\trace-\text{preserving}})$
 & $\overset{\text{?}}{\longrightarrow}$
 & $(N,M,?)$ & & \\
\vspace{3pt} \\
& & $\downarrow$ & & $\downarrow$ & & \\
\vspace{3pt} \\
$(Z^r,V)$ & $\longrightarrow$
 & $(Z^{\sph},V,w)$
 & $\longrightarrow$
 & $\br{(Z^{\sph}}^{\rig},V)$ & $\isom$ & $(Z^r,V)$ \\
\vspace{3pt} 
\end{tabular}
\end{center}

\begin{thm}
Let $(Z,V)$ be a rigid \cstar-planar algebra.  Then there exists a
finite index $\IIone$ subfactor $(N,M,E)$ such that $(Z^{(N,M,E)},V)
\isom (Z,V)$.  In other words, there exists an isomorphism $\Psi:V
\rightarrow V^{(N,M,E)}$ such that for all $T\in \PAr$,
\[
\br{Z^{(N,M,E)}(T)}\circ\br{\Psi_{(\sigma_1,k_1)} \tensor \cdots \tensor
  \Psi_{(\sigma_n,k_n)}}
= \Psi_{(\sigma_0,k)}\circ Z(T) .
\]
\end{thm}

\Proof
As noted in the discussion preceding the statement of the theorem, we
can assume that we have an extremal, finite-index $\IIone$ subfactor
$N\subset M$ giving rise to the associated spherical \cstar-planar
algebra $(Z,V)=\br{Z^{(N,M,E_N)},V^{(N,M,E_N)}}$, where $E_N$ is the
trace-preserving conditional expectation.  We also have a positive,
invertible element $w\in N'\cap M=V_1^+$ satisfying
$\trace(w^{1/2})=\trace(w^{-1/2})=\lambda^{-1}$ for some $\lambda>0$.
We want to show that the rigid planar algebra $(Z^r,V)$ constructed
using Theorem~\ref{thm: spherical to rigid} can be realized as
$\br{Z^{(N,M,\ov{E})},V^{(N,M,\ov{E})}}$ for some new conditional
expectation $\ov{E}:M\rightarrow N$.

Let $w_i$ denote the Radon-Nikodym derivatives in $(Z^r,V)$.  From
Theorem~\ref{thm: spherical to rigid} $w_1=w$.  Recall from
Lemma~\ref{lemma: properties of RN derivs} that
\vspace{10pt}
\begin{figure}[htbp]
\begin{center}
\psfrag{wt}{\scriptsize $\wtilde{w}$}
\psfrag{w-}{\scriptsize $w^{-1}$}
\psfrag{w2=}{$w_2=Z^r
  \begin{pmatrix} &&&&&&&\\ &&&&&&&\\ &&&&&&&\end{pmatrix}
  =Z \begin{pmatrix} &&&&&&&\\ &&&&&&&\\ &&&&&&&\end{pmatrix}
  =Z \begin{pmatrix} &&&&\\ &&&&\\ &&&&\end{pmatrix}
$}
\includegraphics{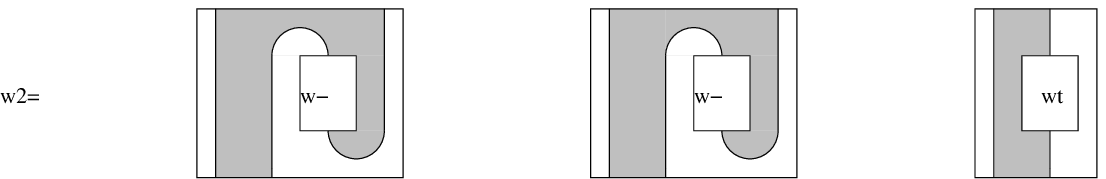}
\end{center}
\end{figure}

\noindent where we have used the fact that the extra $w^{-1/4}$ and
$w^{1/4}$ terms, involved in going from $Z^r$ to $Z$, will cancel.  

\newpage
\noindent In general
\begin{figure}[htbp]
\begin{center}
\psfrag{..}{$\cdots$}
\psfrag{wt}{\scriptsize $\wtilde{w}$}
\psfrag{w}{\scriptsize $w$}
\psfrag{2i}{$\underbrace{\hspace{68pt}}_{2i}$}
\psfrag{2i+1}{$\underbrace{\hspace{68pt}}_{2i+1}$}
\psfrag{w2=}{$w_{2i}=Z
  \begin{pmatrix} & \hspace{88pt} \\ & \\ &
  \end{pmatrix}$}
\psfrag{w3=}{$w_{2i+1}=Z
  \begin{pmatrix} & \hspace{88pt} \\ & \\ &
  \end{pmatrix}$}
\includegraphics{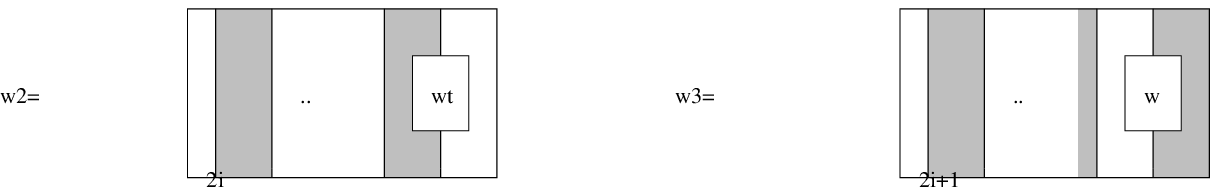}
\end{center}
\end{figure}

\noindent By Proposition~\ref{prop: cJ is S.S},
for $R\in N' \cap M_{2k+1}$
\begin{figure}[htbp]
\hspace{64pt}
\psfrag{R}{$R$}
\psfrag{oR}{\hspace{-1pt}\rotatebox{180}{$R$}}
\psfrag{J=}{$\cj_{2k+2}(R)=Z
  \begin{pmatrix} & \hspace{60pt} \\ & \\ & \end{pmatrix}
  =Z \begin{pmatrix} & \hspace{42pt} \\ & \\ & \end{pmatrix}
  = J_k R^* J_k
$}
\psfrag{*}{$*$}
\includegraphics{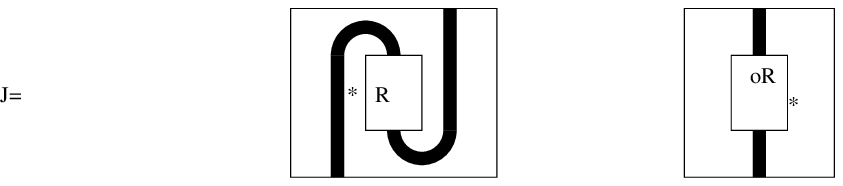}
\end{figure}
\vspace{-5pt}

\noindent where the thick string as usual represents $2k+2$ regular
strings.  So, with $R=w^{-1}$ or $\wtilde{w}^{-1}$ and suppressing
mention of $Z$,
\begin{figure}[htbp]
\psfrag{..}{$\cdots$}
\psfrag{R}{$R$}
\psfrag{k}{\scriptsize $k$}
\psfrag{2k}{\scriptsize $2k$}
\psfrag{Jw-J=}{$J_{k-1}w_k^{-1}J_{k-1}=\cj
  \begin{pmatrix} & \hspace{118pt} \\ & \\ & \end{pmatrix}$}
\includegraphics{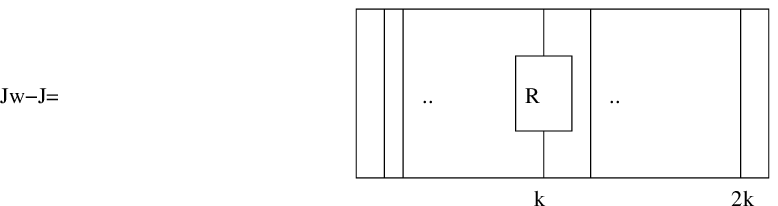}
\end{figure}

\begin{figure}[htbp]
\psfrag{..}{$\cdots$}
\psfrag{oR}{\hspace{-1pt}\rotatebox{180}{$R$}}
\psfrag{k+1}{\scriptsize $k+1$}
\psfrag{2k}{\scriptsize $2k$}
\psfrag{Jw-J=}{$\phantom{J_{k-1}w_k^{-1}J_{k-1}}=$}
\includegraphics{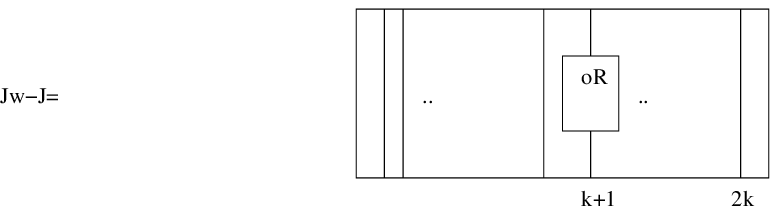}
\end{figure}

\vspace{-82pt}
\begin{figure}[htbp]
\hspace{145pt}
\psfrag{..}{$\cdots$}
\psfrag{R}{$R$}
\psfrag{k+1}{\scriptsize $k+1$}
\psfrag{2k}{\scriptsize $2k$}
\psfrag{Jw-J=}{\hspace{15pt}$\phantom{J_{k-1}w_k^{-1}J_{k-1}}=
  \hspace{140pt}=w_{k+1}$}
\includegraphics{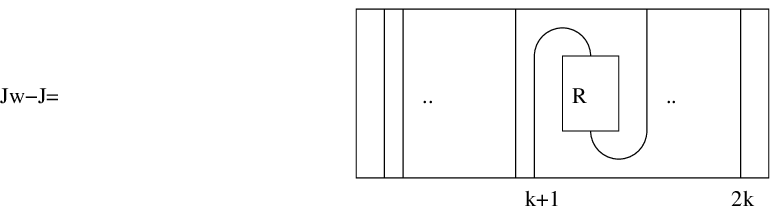}
\end{figure}

\begin{lemma}
\label{lemma: going between the basic constructions}
Let $N\subset M$ be a $\IIone$ subfactor with index
$[M:N]=\br{\delta^{\sph}}^2$ and let $w_1=w\in N'\cap M$ be a
positive invertible element with $\trace\br{w^{1/2}}
=\trace\br{w^{-1/2}}=\lambda^{-1}$, for some $\lambda>0$.  Let
$E_N$ denote the trace-preserving conditional expectation from $M$
onto $N$ and define $\ov{E}:M\rightarrow N$ by
\[
\ov{E}(x)=\lambda E_N\br{w^{-1/2}x} .
\]
Then $\ov{E}$ is a conditional expectation with
$\Ind(\ov{E})=\delta^2$, where $\delta=\lambda^{-1} \delta^{\sph}$.

Take the state $\ov{\varphi}=\trace$ on $N$ and extend it to $M$
by $\ov{\varphi}\circ\ov{E}=\lambda \trace\br{w^{-1/2}\spdot}$.
For $x\in M$ write $x$ for the standard action on $\Ltwo(M,\trace)$ by
left multiplication.  Write $\pi(x)$ for the action on
$\Ltwo(M,\ov{\varphi})$ by left multiplication.  Let
$\ov{e}_1$ be the projection in $\B(\Ltwo(M,\ov{\varphi}))$
given by $\ov{E}$.  Let $\ov{M}_1$ denote the result of the
basic construction applied to $(N,M,\ov{E})$,
i.e. $\ov{M}_1=\ov{J}_0 \pi(N)' \ov{J}_0
=\pi(M)\ov{e}_1\pi(M)$.

\noindent Then there exists a unitary operator $U:\Ltwo(M,\trace)\rightarrow
\Ltwo(M,\ov{\varphi})$ such that
\begin{description}
\item{(1)} $\pi(x)=UxU^*$.
\item{(2)} $U^*\ov{M}_1 U=M_1$.
\item{(3)} $U^*\ov{E}^{\ov{M}_1}_M\br{U\spdot U^*}U
=\lambda E_M\br{w_2^{-1/2}\spdot}$, where
$w_2=J_{\Ltwo(M,\trace)} w^{-1} J_{\Ltwo(M,\trace)}$.
\end{description}
\begin{proof}
Note that for any $y\in N'\cap M$, $E_N(xy)=E_N(yx)$ because for all
$n\in N$
\[
\ip{xy,n}=\ip{x,ny^*}=\ip{x,y^*n}=\ip{yx,n} .
\]
Thus $\ov{E}(x)=\lambda E_N\br{w^{-1/2}x}=\lambda
E_N\br{w^{-1/4}xw^{-1/4}}$ is positive.  $\ov{E}$ is $N$-linear on the
left and the right, and $\ov{E}(1)=\lambda
E_N(w^{-1/2})=\lambda\trace(w^{-1/2}) =1$.  Hence $\ov{E}$ is a
conditional expectation.

Take $\trace$ as the state on $N$.  Then the state $\ov{\varphi}$
on $M$ is
\[
\ov{\varphi}(x)
=\trace\br{E(x)}
=\lambda \trace\br{E_N(w^{-1/2}x)}
=\lambda \trace(w^{-1/2}x) .
\]


Define the unitary operator $U:\Ltwo(M,\trace)\rightarrow
\Ltwo(M,\ov{\varphi})$ by defining a bijection from $M$ to $M$ by
$x\mapsto \lambda^{-1/2}xw^{1/4}$ and noting that this map is
isometric
\begin{align*}
||\lambda^{-1/2}xw^{1/4}||_{\Ltwo(M,\ov{\varphi})}^2
= \lambda^{-1} \ov{\varphi}\br{w^{1/4}x^*xw^{1/4}}
&=\trace\br{w^{-1/2}w^{1/4}x^*xw^{1/4}} \\
&=\trace(x^*x)
=||x||_{\Ltwo(M,\trace)}^2 .
\end{align*}
Observe that $U^*=U^{-1}:x \mapsto \lambda^{1/2} x w^{-1/4}$ and
$UxU^*=\pi(x)$ for $x\in M$.  Also,
\begin{align*}
U^* \ov{e}_1 U:x
&\mapsto \lambda^{1/2} E\br{\lambda^{-1/2}xw^{1/4}} w^{-1/4} \\
&= w^{-1/4} E\br{xw^{1/4}} \\
&= \lambda w^{-1/4} E_N\br{w^{-1/2}xw^{1/4}} \\
&= \lambda w^{-1/4} E_N\br{w^{-1/4}x} ,
\end{align*}
so $U^* \ov{e}_1 U=\lambda w^{-1/4} e_1 w^{-1/4}$.  Thus, if $\{ b \}$
be a basis for $M$ over $N$ with respect to $E_N$, then
$\ov{b}=\lambda^{-1/2}bw^{1/4}$ is a basis for $M$ over $N$ with
respect to $E$.  Now let us compute the index of $\ov{E}$.  Using
Lemma~\ref{lemma: Kosaki 1998 3.4}
\begin{align*}
\Ind(\ov{E})
 = \ov{E}^{-1}(1)
&= \ov{E}^{-1}\br{\sum \ov{b}\ov{e}_1\ov{b}^*}
 = \sum \ov{E}^{-1}\br{\theta^{\trace}\br{\hatt{\ov{b}},
    \hatt{\ov{b}}}} \\
&= \sum \theta^{\trace\circ\ov{E}}\br{\hatt{\ov{b}},
    \hatt{\ov{b}}}
 = \sum \theta^{\ov{\varphi}}\br{\hatt{\ov{b}},\hatt{\ov{b}}} \\
&= \sum \ov{b}\ov{b}^*
 = \lambda^{-1} \sum bw^{1/2}b^* \\
&= \lambda^{-1} E^{-1}(w^{1/2})
 = \lambda^{-1} \br{\delta^{\sph}}^2 \trace\br{w^{1/2}} \\
&= \lambda^{-2} \br{\delta^{\sph}}^2
 = \delta^2 .
\end{align*}

\noindent Now $\ov{M}_1=\pi(M)\ov{e}_1\pi(M)$ so
\[
U^*\ov{M}_1U
=\pi(M)U^*\ov{e}_1 U\pi(M)
=M\br{\lambda w^{-1/4}e_1 w^{-1/4}}M
=M e_1 M
=M_1 .
\]

\noindent To show (3) note that for $a,b \in M$,
$\ov{E}_M(a\ov{e}_1b)=\frac{1}{\delta^2} ab$ and
\begin{align*}
\lambda E_M\br{w_2^{-1/2}\br{U^*a\ov{e}_1bU}}
&= \lambda^2 E_M\br{w_2^{-1/2}aw^{-1/4}e_1}w^{-1/4}b \\
&= \lambda^2 E_M\br{aw^{-1/4}e_1w_2^{-1/2}}w^{-1/4}b \\
&= \lambda^2 aw^{-1/4}E_M\br{e_1w_2^{-1/2}}w^{-1/4}b
\end{align*}
where we have used the fact that $w_2\in M'\cap M_1$ to write
$E_M\br{w^r \spdot}=E_M\br{\spdot w^r}$. 

\noindent Now
\begin{figure}[htbp]
\psfrag{w}{$w^{\frac{1}{2}}$}
\psfrag{EM=}{$E_M\br{e_1w_2^{-1/2}}=\frac{1}{\br{\delta^{\sph}}^2}
  Z \begin{pmatrix} & \hspace{96pt} \\ & \\ & \\ & \\ & \\ & \\ & \\ &
  \end{pmatrix}$}
\psfrag{=dw}{$=\frac{1}{\br{\delta^{\sph}}^2}
  Z \begin{pmatrix} & \hspace{41pt} \\ & \\ &
  \end{pmatrix}
$}
\includegraphics{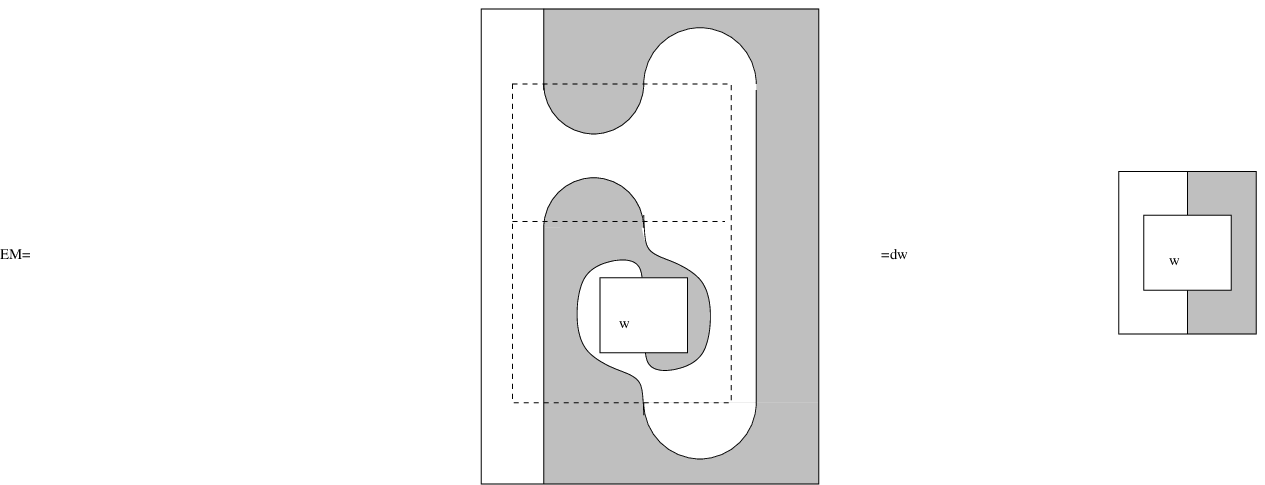}
\end{figure}

\noindent So
\[
\lambda E_M\br{w_2^{-1/2}\br{U^*a\ov{e}_1bU}}
=\br{\frac{\lambda}{\delta^{\sph}}}^2 aw^{-1/4}w^{1/2}w^{-1/4}b
=\frac{1}{\delta^2} ab .
\]
\end{proof}
\end{lemma}

Lemma~\ref{lemma: going between the basic constructions} says that we
have an isomorphism of the (short) towers which fixes $M$ and is given by
\begin{center}
\begin{tabular}{cccccc}
\vspace{3pt} \\
$N$ & $\underset{\ov{E}}{\subset}$ & $M$ & 
  $\underset{\ov{E}_M}{\subset}$ & $\ov{M}_1$ & \\
\vspace{3pt} \\
& & $\downarrow$ & & & \vspace{10pt} $\isom$ via $U^*\spdot U$ \\
\vspace{3pt} \\
$N$ & $\underset{\ov{E}}{\subset}$ & $M$ & 
  $\underset{\br{E_M}\ov{\phantom{A}}}{\subset}$ & $M_1$ & \\
\vspace{3pt} 
\end{tabular}
\end{center}
$\br{E_M}\ov{\phantom{A}}$ is obtained from $M \subset M_1$ and
$w_2=J_0 w_1^{-1} J_0$ as in the first part of the Lemma.
i.e. $\br{E_M}\ov{\phantom{A}}(x)=\lambda E_M\br{w_2^{-1/2}x}$.  We
thus have the following immediate corollary.

\begin{cor}
\label{cor: overM to M}
There exists an isomorphism of towers $\bigE:\{ \ov{M}_i
\}_{i\geq -1} \rightarrow \{ M_i \}_{i \geq -1}$ such that
\begin{description}
\item{1.} $\bigE|_{\ov{M}_i}:\ov{M}_i\rightarrow M_i$ is 
a *-algebra isomorphism.
\item{2.} $\bigE\circ \ov{E}_{\ov{M}_i} \circ \bigE^{-1} = \lambda
E_{M_i}\br{w_{i+1}^{-1/2}\spdot}$ where $w_1=w$,
$w_{i+1}=J_{i-1}w_i^{-1}J_{i-1}$.
\end{description}
\begin{proof}
By Lemma~\ref{lemma: going between the basic constructions} we have
isomorphism up to $i=1$ given by $\bigE_1=U^*\spdot U$.  Suppose we have
such an isomorphism $\bigE_i$ at level $i$.  Recalling that the
basic construction is independent of the state, take the state
$\phi=\trace\circ\bigE_i$ on $\ov{M}_{i-1}$.  Then
\begin{center}
\begin{tabular}{cccc}
\vspace{3pt} \\
$(\ov{M}_{i-1},\phi)$ 
  & $\underset{\ov{E}_{\ov{M}_{i-1}}}{\subset}$ 
  & $\ov{M}_i$ & \\
\vspace{3pt} \\
& $\downarrow$ & & \vspace{10pt} $\isom$ via $\bigE_i$ \\
\vspace{3pt} \\
$(M_{i-1},\trace)$ 
  & $\underset{\br{E_{M_{i-1}}}\ov{\phantom{A}}}{\subset}$ 
  & $M_i$ & \\
\vspace{3pt} 
\end{tabular}
\end{center}
So $\bigE_i$ extends to an isomorphism $\bigE_i^{\mathrm{ext}}$
of the basic constructions.  Combining this with Lemma~\ref{lemma: going
between the basic constructions} we have
\begin{center}
\begin{tabular}{cccccc}
\vspace{3pt} \\
$(\ov{M}_{i-1},\phi)$ 
  & $\underset{\ov{E}_{\ov{M}_{i-1}}}{\subset}$ 
  & $\ov{M}_i$ 
  & $\underset{\ov{E}_{\ov{M}_i}}{\subset}$ 
  & $\ov{M}_{i+1}$ & \\
\vspace{3pt} \\
& & $\downarrow$ & & & \vspace{10pt} $\isom$ via $\bigE_i^{\mathrm{ext}}$\\
\vspace{3pt} \\
$(M_{i-1},\trace)$ 
  & $\underset{\br{E_{M_{i-1}}}\ov{\phantom{A}}}{\subset}$ 
  & $M_i$ 
  & $\underset{E^B_{M_i}}{\subset}$ 
  & $B$ & \\
\vspace{3pt} \\
& & $\downarrow$ & & & \vspace{10pt} $\isom$ via $U_{i+1}^*\spdot U_{i+1}$\\
\vspace{3pt} \\
$M_{i-1}$ 
  & $\underset{\br{E_{M_{i-1}}}\ov{\phantom{A}}}{\subset}$ 
  & $M_i$ 
  & $\underset{\br{E_{M_i}}\ov{\phantom{A}}}{\subset}$ 
  & $M_{i+1}$ & \\
\vspace{3pt} 
\end{tabular}
\end{center}
So we let
$\bigE_{i+1}=\mathrm{Ad}(U_{i+1}^*)\circ\bigE_i^{\mathrm{ext}}$
and note that $\bigE_{i+1}$ restricts to $\bigE_i$ on
$\ov{M}_i$.  $\bigE$ is just the direct limit of the
$\bigE_i$'s.
\end{proof}
\end{cor}

Restrict $\bigE$ to $N' \cap \ov{M}_i$ to obtain an isomorphism
$\Psi:V^{(N,M,\ov{E})}\rightarrow V^{(N,M,E_N)}=V$.
Let $\overline{Z}$ be defined by 
\[
\overline{Z}(T)
=\Psi\circ\br{Z^{(N,M,E)}(T)}\circ\br{\Psi^{-1}_{(\sigma_1,k_1)} 
  \tensor \cdots \tensor \Psi^{-1}_{(\sigma_n,k_n)}} .
\]
We want to show that $\overline{Z}$ is just the rigid \cstar-planar
algebra $Z^r$.  By the uniqueness part of Theorem~\ref{thm: subfactor
to RPA} it suffices to check that properties (4) through (7) are true.
Property (4) is obvious because $\Psi$ is an algebra isomorphism and
the multiplication tangles are unchanged by the constructions of
Section~\ref{section: rigid to spherical}.

Before proving the other properties, consider the Radon-Nikodym
derivatives.  On $N'\cap M$, $\ov{\varphi}=\ov{E}=\lambda
\trace\br{w^{-1/2}\spdot}$.  Meanwhile
\begin{align}
\ov{E}'(x)
=\frac{1}{\delta^2}\sum \ov{b}x\ov{b}^*
& =\frac{\lambda^{-1}}{\delta^2}\sum bw^{1/4}xw^{1/4}b^* \nonumber \\
& =\br{\frac{\delta^{\sph}}{\delta}}^2\lambda^{-1} E'\br{w^{1/4}xw^{1/4}}
=\lambda E'\br{w^{1/2}x} \label{eq: ovE' to E'}
\end{align}
so that $\ov{\varphi}'(x)=\lambda\trace'\br{w^{1/2}x}
=\lambda\trace\br{w^{1/2}x}$ and hence $w^{(N,M,\ov{E})}=w$.

Similarly, using part 2 of Corollary~\ref{cor: overM to M} and the
construction of a basis in Lemma~\ref{lemma: going between the basic
constructions}, we obtain $\Psi\br{w_k^{(N,M,\ov{E})}}=w_k$.  Thus the
Radon-Nikodym derivatives for $\ov{Z}$ are the same as those for
$Z^r$.  Now
\begin{description}
\item{(5)} Consider $k$ odd (for $k$ even just use $\wtilde{w}$ in
place of $w$ in the second diagram).
\begin{figure}[htbp]
\hspace{30pt}
\psfrag{k}{\scriptsize $k$}
\psfrag{k+1}{\scriptsize $k+1$}
\psfrag{w}{$w^{\frac{1}{4}}$}
\psfrag{..}{$\cdots$}
\psfrag{Z=}{$\br{\ov{Z}}^{\sph}
  \begin{pmatrix} & \hspace{120pt} \\ & \\ & \\ & \end{pmatrix}
  =\ov{Z}\begin{pmatrix} & \hspace{120pt} \\ & \\ & \\ & \end{pmatrix}$}
\includegraphics{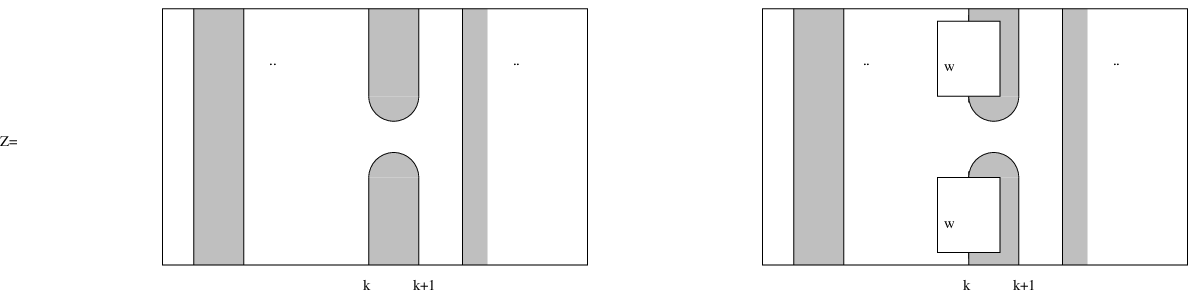}
\begin{align*}
\hspace{6pt} &= \delta w_k^{1/4} \Psi\br{\ov{e}_k} w_k^{1/4} 
 = \delta \lambda w_k^{1/4} \br{w_k^{-1/4} e_k w_k^{-1/4}} w_k^{1/4} \\
&= \delta^{\sph} e_k
\end{align*}
\end{figure}
\vspace{-30pt}

\item{(6)} Consider $k$ even (for $k$ odd use $\wtilde{w}$ in
place of $w$ in the second diagram).
\begin{figure}[htbp]
\hspace{30pt}
\psfrag{k}{$\underbrace{\hspace{30pt}}_{k}$}
\psfrag{x}{$x$}
\psfrag{w}{$w^{\frac{1}{2}}$}
\psfrag{..}{$.\;.$}
\psfrag{Z}{$\br{\ov{Z}}^{\sph}
  \begin{pmatrix} & \hspace{76pt} \\ & \\ & \\ & \end{pmatrix}
  =\ov{Z}\begin{pmatrix} & \hspace{76pt} \\ & \\ & \\ & \\ & \end{pmatrix}$}
\includegraphics{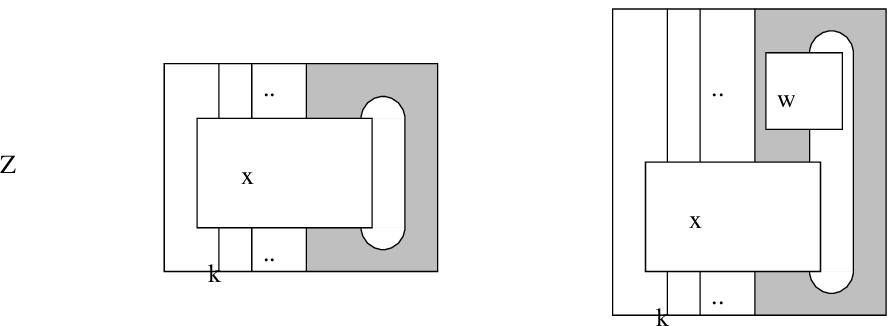}
\vspace{13pt}
\begin{align*}
\hspace{16pt} &= \delta \Psi\br{\ov{E}_{\ov{M}_{k-1}}\br{\Psi^{-1}
  \br{w_k^{1/2}x}}} 
 = \delta \lambda E_{M_{k-1}}\br{w_k^{-1/2} \br{w_k^{1/2}x}} \\
&= \delta^{\sph} E_{M_{k-1}}(x)
\end{align*}
\end{figure}
\vspace{-25pt}

\noindent using Corollary~\ref{cor: overM to M} part 2.

\item{(7)}
\begin{figure}[htbp]
\vspace{-20pt}
\hspace{30pt}
\psfrag{x}{$x$}
\psfrag{w}{$\!w^{-\frac{1}{2}}$}
\psfrag{..}{$.\;.$}
\psfrag{Z}{$\br{\ov{Z}}^{\sph}
  \begin{pmatrix} & \hspace{76pt} \\ & \\ & \\ & \end{pmatrix}
  =\ov{Z}\begin{pmatrix} & \hspace{76pt} \\ & \\ & \\ & \\ & \end{pmatrix}$}
\includegraphics{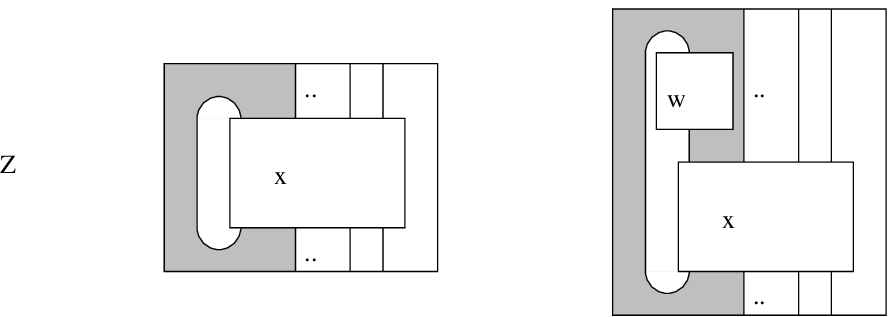}
\vspace{5pt}
\begin{align*}
\hspace{-3pt} &= \delta \Psi\br{\ov{E}'\br{\Psi^{-1}
  \br{w_{-1/2}x}}} 
 = \delta \lambda E'\br{w^{1/2} \br{w^{-1/2}x}} \\
&= \delta^{\sph} E'(x)
\end{align*}
\end{figure}
\vspace{-25pt}

\noindent using equation~(\ref{eq: ovE' to E'}) in the third line.
Hence $\overline{Z}=Z^r$ and $\Psi$ is an isomorphism of $(Z^r,V)$
with $(Z^{(N,M,\ov{E})},V^{(N,M,\ov{E})})$.
\end{description}
\proofend


\section{Finite index $\IIone$ subfactors}
\label{sect: finite index IIone}

Here we consider the case of a finite index $\IIone$ subfactor $N
\subset M$ with the unique trace preserving conditional expectation
$E_N$.  Recall,

\begin{definition}[Extremal]
A finite index $\IIone$ subfactor $N \subset M$ is called {\em
extremal} if the unique traces $\trace'$ and $\trace$ on $N'$ and $M$
respectively coincide on $N' \cap M$ (where $N'$ is calculated on any
Hilbert space on which $M$ acts with finite $M$-dimension).
\end{definition}

\begin{remark}
Note that in particular an irreducible finite index $\IIone$ subfactor
($N' \cap M = {\mathbb C}$) is extremal.  In~\cite{PimsnerPopa1986}
and~\cite{PimsnerPopa1988} Pimsner and Popa show that if $N \subset M$
is extremal then $\trace'=\trace$ on all $N' \cap M_i$.
\end{remark}

We will show that the two rotations defined in Huang~\cite{Huang2000}
are the same if and only if the subfactor is extremal, and illuminate
the connection between these two rotations in the general case.  Our
approach will yield a new proof of the periodicity of the rotation and
provide the correct formulation to generalize to the infinite index
$\IIone$ case in Chapter~\ref{chapter: infinite index}.

Motivated by the relationship between the two rotations, we produce a
two-parameter family of rotations.  We conclude with a collection of
results on general finite index $\IIone$ subfactors including a new
proof of some of Pimsner and Popa's characterizations of extremality
in~\cite{PimsnerPopa1986}.


\subsection{Huang's two rotations for a nonextremal $\IIone$ subfactor}
\label{section: two rotations}

In~\cite{Huang2000} Huang defines two rotations on the standard
invariant of a finite index $\IIone$ subfactor, shows that each is
periodic and conjectures that the two are equal.  One is the rotation
$\rho_k$ defined in~\ref{define rotation}.  The other is defined as
follows:
\begin{definition}
\label{def: tilde(rho)}
Define $\wtilde{\rho}_k:M_k \rightarrow N' \cap M_k$ by defining its
action on basic tensors $x=x_1 \tensorN x_2 \tensorN \cdots \tensorN
x_{k+1}$ as
\begin{equation}
\label{central rotn}
\wtilde{\rho}_k(x)
= P_c \br{x_2 \tensorN x_3 \tensorN \cdots \tensorN x_{k+1} \tensorN x_1},
\end{equation}
where $P_c$ is the orthogonal projection onto the $N$-central vectors
in $\Ltwo(M_k,\trace)$, which is just the finite dimensional subspace
$N' \cap M_k$.  $P_c$ above could thus be replaced by $E_{N' \cap
M_k}$, the unique trace-preserving conditional expectation from
$(M_k,\trace)$ to $(N' \cap M_k,\trace)$.
\end{definition}

\begin{remark}
In the extremal case Jones~\cite{Jones1999} shows that
$\rho=\wtilde{\rho}$ and periodicity of $\rho$ follows from that of
$\wtilde{\rho}$.  As we have already mentioned, in the nonextremal
case Huang~\cite{Huang2000} shows that $\rho$ and $\wtilde{\rho}$ are
periodic.
\end{remark}

We begin by formulating another equivalent definition of the
rotation $\rho$.

\begin{lemma}
For $x \in N' \cap M_k$,
\[
\hatr{\rho_k(x)}
=\sum_{b \in B} R_{b^*} \br{ L_{b^*}}^* \hatt{x} ,
\]
where $L_b, R_b:\Ltwo(M_{k-1},\trace) \rightarrow \Ltwo(M_k,\trace)$
by $L_b \xi=\hatt{b} \tensorlN \xi$ and $R_b \xi = \xi \tensorlN
\hatt{b}$.
\begin{proof}
From \ref{notation: left tensor}, $\br{L_b}^* \br{\hatt{c} \tensorlN
\eta} = \br{L_b}^* L_c \eta=E_N(b^*c)\eta$ (for $\eta\in \Ltwo(M_{k-1})$).
Hence
\begin{align*}
\rho_k(x)
&=\sum_{b \in B}
  \sum_i E_N\br{bx^{(i)}_1}x^{(i)}_2\tensorN x^{(i)}_3 \tensorN 
  \cdots \tensorN x^{(i)}_{k+1} \tensorN b^* \\
&=\sum_{b \in B}
  R_{b^*} \br{ L_{b^*}}^* x .
\end{align*}
\end{proof}
\end{lemma}

\begin{prop}
\label{prop: rho inner product}
For all $x \in N' \cap M_k$ and for all $y_i \in M$,
\[
\ip{\rho_k(x), y_1 \tensorN y_2 \tensorN \cdots \tensorN y_{k+1}}
=\ip{x, y_{k+1} \tensorN y_1 \tensorN \cdots \tensorN y_k}
\]
where the inner product is that on $\Ltwo(M_k,\trace)$.
\begin{proof}
\begin{align*}
\hspace{20pt} & \hspace{-20pt}
\ip{\rho_k(x),y}
 =\sum_b \ip{x, L_{b^*} \br{ R_{b^*}}^* y} 
 =\sum_b \ip{x, b^* \tensorN y_1 \tensorN\cdots\tensorN y_k E_N(y_{k+1}b)} \\
&=\sum_b \ip{xE_N(y_{k+1}b)^*,b^* \tensorN y_1 \tensorN\cdots\tensorN y_k} 
 =\sum_b \ip{E_N(y_{k+1}b)^*x,b^* \tensorN y_1 \tensorN\cdots\tensorN y_k} \\
&=\sum_b \ip{x,E_N(y_{k+1}b)b^* \tensorN y_1 \tensorN \cdots \tensorN y_k} 
 = \ip{x,y_{k+1} \tensorN y_1 \tensorN \cdots \tensorN y_k} .
\end{align*}
\end{proof}
\end{prop}

\begin{cor}
$\rho_k$ is periodic, $\br{\rho_k}^{k+1}=\id$.  On $\Ltwo(N' \cap
M_k,\trace)$ $\wtilde{\rho}_k=\br{\rho_k^{-1}}^*$ (and hence
$\wtilde{\rho}_k$ is also periodic, $\br{\wtilde{\rho}_k}^{k+1}=\id$).
\end{cor}

\begin{cor}
\label{cor: rho to tilde rho via RN}
$\tilde{\rho}_0=\rho_0=\id$ and for $k \geq 1$ and $y \in N' \cap
M_k$, $\wtilde{\rho}_k(y)=\rho_k(y)z_1=\rho_k(z_1^{-1}y)$ (recall
$z_1$ is the Radon-Nikodym derivative of $\trace'$ with respect to
$\trace$ on $N' \cap M_1$, so $\trace'(x)=\trace(z_1x)$ for all $x\in
N'\cap M_1$).
\begin{proof}
Prop~\ref{prop: rho inner product} implies that
$(\rho_k(x^*))^*=\rho_k^{-1}(x)$ for $x \in N' \cap M_k$.  To see this
note that for $y= y_1 \tensorlN \cdots \tensorlN y_{k+1}$,
\begin{align*}
\ip{\rho_k(x^*)^*,y}
&=\ip{y_{k+1}^* \tensorN y_k^* \tensorN \cdots \tensorN y_1^*,
  \rho_k(x^*)} 
 =\ip{y_1^* \tensorN y_{k+1}^* \tensorN \cdots \tensorN y_2^*,
  x^*} \\
&=\ip{x,
  y_2 \tensorN \cdots \tensorN y_{k+1} \tensorN y_1 }
 =\ip{\rho_k^{-1}(x),y} .
\end{align*}
Let $y \in N' \cap M_k$.  Then, writing a thick string to represent
$k-1$ regular strings,
\begin{align*}
&\delta^{k+1}\trace(x^* \wtilde{\rho}_k(y))
= \delta^{k+1}\ip{\wtilde{\rho}_k(y),x}
= \delta^{k+1}\ip{y,\rho_k^{-1}(x)}
= \delta^{k+1}\ip{y,(\rho_k(x^*))^*} \\
& \phantom{\delta^{k+1}}= \delta^{k+1}\trace(\rho_k(x^*)y) \\
\end{align*}
\vspace{-10mm}
\begin{figure}[htbp]
\begin{center}
\psfrag{y}{$y$}
\psfrag{x^*}{$x^*$}
\psfrag{=}{$=$}
\psfrag{z_1}{\scriptsize $z_1$}
\psfrag{z_1^-1}{\scriptsize $z_1^{-1}$}
\includegraphics{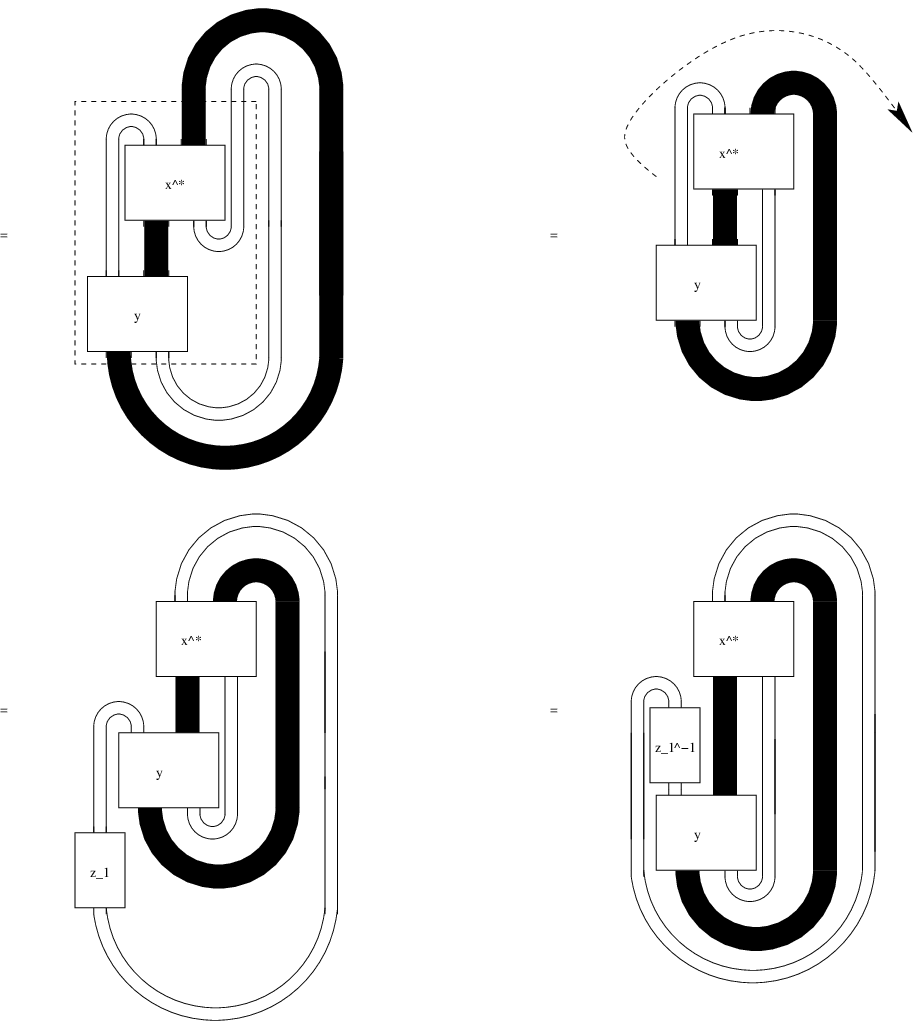}
\end{center}
\end{figure}
\vspace{-5mm}
\begin{align*}
&\hspace{-10mm} = \delta^{k+1} \trace(x^* \rho_k(y) z_1) 
&=\delta^{k+1} \trace(x^* \rho_k(z_1^{-1} y))
\end{align*}
Hence $\wtilde{\rho}_k(y)=\rho_k(y) z_1=\rho_k(z_1^{-1} y)$.
\end{proof}
\end{cor}



\subsection{A two-parameter family of rotations}
\label{section: 2-param family}

In fact we can extend the relationship between $\rho$ and
$\tilde{\rho}$ to define a two-parameter family of rotations (periodic
automorphisms of the linear space $N' \cap M_k$ of period $k+1$).

\begin{definition}
For $r,s \in \mathbb{R}$ define $\rho_k^{(r,s)}:N' \cap M_k
\rightarrow N' \cap M_k$ by
\[
\rho^{(r,s)}_k(x)=w_{k-2,k}^r \rho_k(x) w_{-1,1}^s
\]
where (recall) $w_{i,j}$ is the Radon-Nikodym derivative of $\trace'$
with respect to $\trace$ on $\M{i}{j}$.
\end{definition}

\begin{prop}
$\left(\rho^{(r,s)}_k \right)^{k+1}=\id$.
\end{prop}

\Proof .

\noindent First note that for $k$ odd (so an even number of strings),
\vspace{10pt}
\begin{figure}[htbp]
\hspace{70pt}
\psfrag{x}{$x$}
\psfrag{=}{$\hspace{40pt}=$}
\psfrag{a}{\scriptsize $w^{r}$}
\psfrag{b}{\scriptsize $\wtilde{w}^{r}$}
\psfrag{c}{\scriptsize $w^{s}$}
\psfrag{d}{\scriptsize $\wtilde{w}^{s}$}
\psfrag{e}{\scriptsize $\!w^{-s}$}
\psfrag{f}{\scriptsize $\!\wtilde{w}^{-s}$}
\psfrag{g}{\scriptsize $\!w^{-r}$}
\psfrag{h}{\scriptsize $\!\wtilde{w}^{-r}$}
\psfrag{..}{$\cdots$}
\psfrag{r=}{$\rho_k^{(r,s)}(x)=$}
\includegraphics{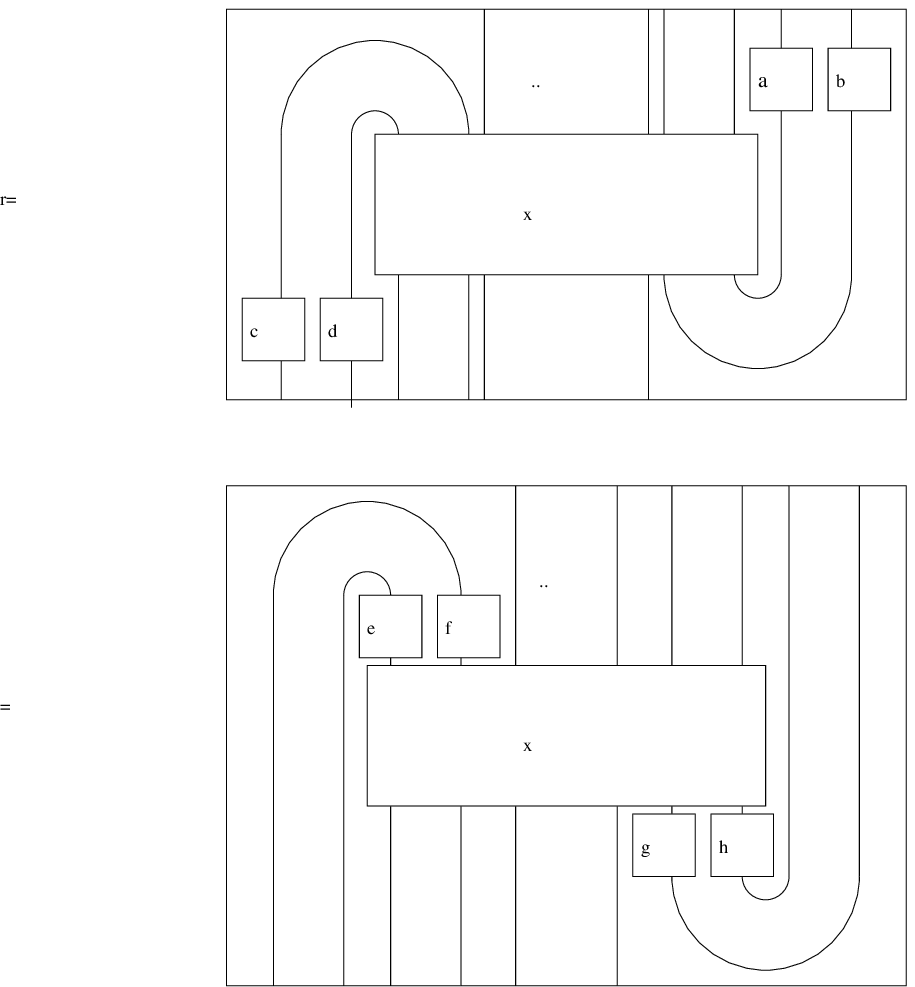}
\end{figure}
\vspace{10pt}

\noindent For $k$ even just switch $w^r$ and $\wtilde{w}^r$ (resp.
$w^{-r}$ and $\wtilde{w}^{-r}$).

We could describe $\rho_k^{(r,s)}$ by saying that every time we pull a
string down on the left-hand side we put $w^s$ or $\wtilde{w}^s$
(whichever makes sense) on the end of the string away from the central
box, or alternatively $\wtilde{w}^{-s}$ or $w^{-s}$ at the end of the
string near the box.  On the right-hand side we put $w^r$ or
$\wtilde{w}^r$ on the end of the string away from the box, or
alternatively $\wtilde{w}^{-r}$ or $w^{-r}$ at the end of the string
near the box.

$\br{\rho_k^{(r,s)}}^{k+1}$ will pull every string through a full
counter-clockwise rotation back to its starting point and every string
will pick up two boxes with powers of $w$.  We obtain the following
\begin{figure}[htbp]
\hspace{70pt}
\psfrag{x}{$x$}
\psfrag{a}{\scriptsize $w^{r}$}
\psfrag{b}{\scriptsize $\wtilde{w}^{r}$}
\psfrag{c}{\scriptsize $w^{s}$}
\psfrag{d}{\scriptsize $\wtilde{w}^{s}$}
\psfrag{e}{\scriptsize $\!w^{-s}$}
\psfrag{f}{\scriptsize $\!\wtilde{w}^{-s}$}
\psfrag{g}{\scriptsize $\!w^{-r}$}
\psfrag{h}{\scriptsize $\!\wtilde{w}^{-r}$}
\psfrag{..}{$\cdots$}
\psfrag{r=}{$\br{\rho_k^{(r,s)}}^{k+1}(x)=$}
\includegraphics{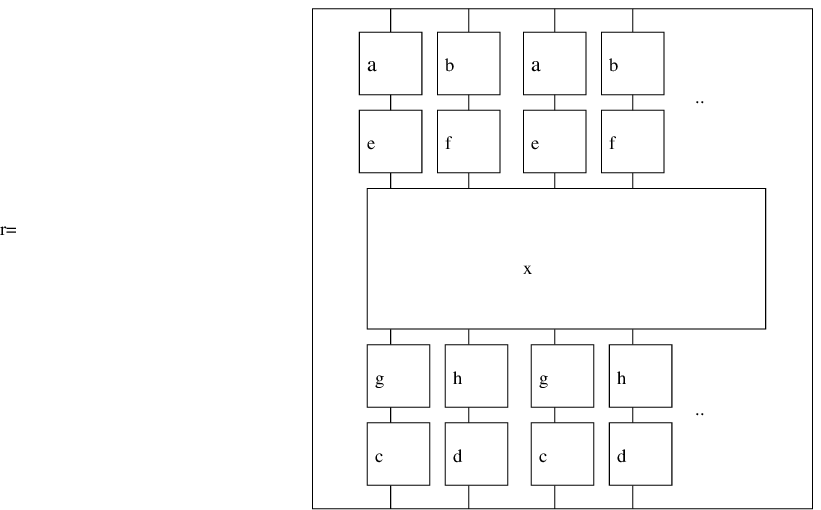}
\vspace{0pt}
\begin{align*}
\hspace{-30pt} &=z_k^{r-s} x z_k^{s-r}=x 
\end{align*}
\end{figure}
\vspace{-10pt}

\noindent where the last line uses the fact that $z_k$ is central.


\subsection{Additional results on finite index $\IIone$ subfactors}
\label{section: additional results IIone}

\begin{lemma}
$E_{M' \cap M_1}(e_1)=\tau w_1$ where $E_{M' \cap M_1}$ is the unique
$\trace$-preserving conditional expectation onto $M' \cap M_1$ and
$w_1$ is the Radon-Nikodym derivative of $\trace'$ with respect to
$\trace$ on $M' \cap M_1$.
\begin{proof}
For $x \in M' \cap M_1$
\begin{align*}
\trace(e_1 x)
&= \overline{\trace'(J_0 e_1 x J_0)} 
 = \overline{\trace'(e_1 J_0 x J_0)} \\
&= \overline{\trace'(e_1 E_N(J_0 x J_0))}
 = \overline{\trace'(e_1 \trace(J_0 x J_0))} \\
&= \overline{\trace'(e_1)}\trace'(x) 
 = \tau \trace(w_1 x)
\end{align*}
\end{proof}
\end{lemma}

\begin{prop}
\label{characterizing extremal}
Let $N \subset M$ be a finite index $\IIone$ subfactor.  Then the
following are equivalent: 
\begin{description}
\item{(i)} $N \subset M$ is extremal ($\trace=\trace'$ on $N' \cap M$)
\item{(ii)} $\trace=\trace'$ on all $N' \cap M_k$
\item{(iii)} $E_{M' \cap M_1}(e_1)=\tau$ where $E_{M' \cap M_1}$ is the 
unique $\trace$-preserving conditional expectation onto $M' \cap M_1$.
\item{(iv)} $\rho_k=\tilde{\rho}_k$ for all $k \geq 0$
\item{(v)} $\rho_k=\tilde{\rho}_k$ for some $k \geq 1$
\end{description}
\begin{proof}
The equivalence of (i) and (ii) is proved in
Pimsner-Popa~\cite{PimsnerPopa1986}, but follows easily from our
knowledge of Radon-Nikodym derivatives.  $N \subset M$ is extremal iff
$\trace=\trace'$ on $N' \cap M$, iff $z_0=w_0=1$, iff $z_k=1$ for all
$k \geq 0$ (by Lemma~\ref{lemma: properties of RN derivs}), iff
$\trace=\trace'$ on all $N' \cap M_k$.  Equivalence with (iii) follows
from the preceding corollary.

If $N \subset M$ is extremal then $z_1=1$ and $\rho=\tilde{\rho}$.  If
$\rho_k=\tilde{\rho}_k=\rho_k(\cdot) z_1$ then, since $\rho_k : N'
\cap M_k \rightarrow N' \cap M_k$ is periodic and hence surjective,
$z_1=1$ and hence $N \subset M$ is extremal.
\end{proof}
\end{prop}

\begin{prop}
\label{other Jones proj}
For $j \geq 1$ let $\tilde{e}_i=w_i^{-1/2}e_iw_i^{-1/2}$.  Then $\{
\tilde{e}_i \}$ are also Jones projections and $E_{M' \cap
M_1}(\tilde{e}_1)=\tau$ ($\trace$-preserving conditional expectation).
\begin{proof}
We need to check that $\tilde{e}_i^2=\tilde{e}_i$, that $\tilde{e}_i
\tilde{e}_{i\pm 1} \tilde{e}_i = \tau \tilde{e}_i$ and $[e_i,e_j]=0$
for $|i-j|>1$.  With the tools we have developed this could be done
simply by drawing the appropriate tangles, but can be obtained with a
little more insight as we see below.

Applying the construction of the spherical \cstar-planar algebra in
Section~\ref{section: rigid to spherical} we obtain new Jones
projections $\ov{e}_j=\lambda^{-1}w_{j-1}^{1/4} e_j
w_{j-1}^{1/4}=\lambda^{-1}w_{j}^{-1/4} e_j w_{j}^{-1/4}$, with
$\delta$ changed to $\delta^{\sph}=\delta \lambda$, where
$\lambda=\trace(w^{1/2})$.  For example
\begin{figure}[htbp]
\psfrag{w}{\scriptsize $w^{\frac{1}{4}}$}
\psfrag{wt-}{\scriptsize $\,\wtilde{w}^{-\frac{1}{4}}$}
\psfrag{e=}{$\ov{E}_1=Z
  \begin{pmatrix} & \hspace{20pt} \\ & \\ & \end{pmatrix}
  =Z^{\sph} \begin{pmatrix}\hspace{40pt} \\ & \\ & \\ & \\ & \end{pmatrix}
  =Z^{\sph} \begin{pmatrix}\hspace{63pt} \\ &\\ &\\ &\\ &\\ &\\ &\end{pmatrix}
  =Z^{\sph} \begin{pmatrix}\hspace{40pt} \\ & \\ & \\ & \\ & \end{pmatrix}
$}
\includegraphics{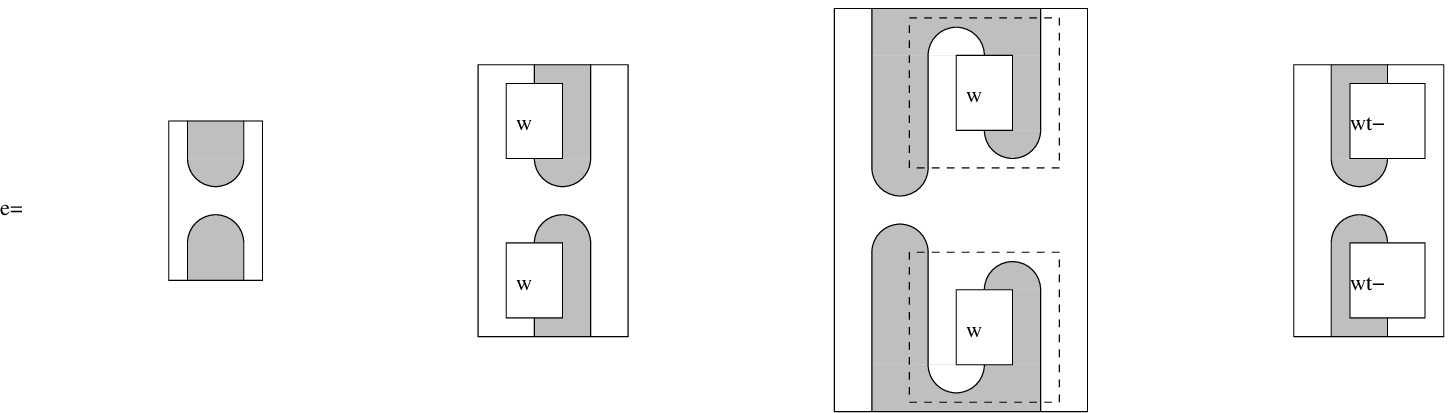}
\end{figure}

\noindent so $\ov{E}_1=w^{1/4}E_1 w^{1/4}=\wtilde{w}^{-1/4} E_1
\wtilde{w}^{-1/4}$, which we could also write as
$\ov{E}_1=w_1^{1/4}E_1 w_1^{1/4}=w_2^{-1/4} E_1 w_2^{-1/4}$.

We have produced a spherical \cstar-planar algebra $(Z^{\sph},V,w)$
with distinguished element $w$.  Take instead the triple
$(Z^{\sph},V,w^{-1})$ and apply the construction of a rigid
\cstar-planar algebra in Theorem~\ref{thm: spherical to rigid} to
obtain $(\wtilde{Z},V)$ (in fact an ordinary \cstar-planar algebra).
Then, using the fact that
$\trace\br{(w^{-1})^{1/2}}=\trace\br{w^{-1/2}}=\trace\br{w^{1/2}}
=\lambda$,
\[
\wtilde{e}_1
=\lambda w^{1/4} \br{\lambda^{-1}w^{1/4} e_1 w^{1/4}}w^{1/4}
=w^{1/2} e_1 w^{1/2}
\]
and the Jones projections $\wtilde{e}_i$ are those given above.  Also
note that $\wtilde{\delta}=\lambda^{-1}\delta^{\sph}=\delta$ so that
$\wtilde{\tau}=\tau$.

Finally, for all $x \in M' \cap M_1$
\begin{align*}
\trace(\tilde{e}_1 x)
&= \trace(e_1 w_1^{-1/2} x w_1^{-1/2}) 
 = \trace(e_1 x w_1^{-1}) \\
&= \overline{\trace'(e_1 J_0 x J_0 J_0 w_1^{-1} J_0)} 
 = \overline{\trace'(e_1 J_0 x J_0 w_0)} \\
&= \overline{\trace'(e_1 E_N(J_0 x J_0 w_0))} 
 = \overline{\trace'(e_1 \trace(J_0 x J_0 w_0))} \\
&= \overline{\trace'(e_1) \trace'(J_0 x J_0)} 
 = \tau \trace(x) .
\end{align*}
and hence $E_{M' \cap M_1}(\tilde{e}_1)=\tau$.
\end{proof}
\end{prop}

\begin{remark}
By changing $w$ to $w^{-1}$ in the triple $(Z^{\sph},V,w)$ and then
constructing $(\wtilde{Z},V)$ we have essentially switched the roles
of $\trace$ and $\trace'$.

This is more than just an interesting trick.  We can do the same thing
for any \cstar-planar algebra.  Because the trace preserving
conditional expectation onto $V_{1,2}$ of $\tilde{e}_1$ is a scalar (a
property that $e_1$ does not possess) one can show that, even in the
nonextremal case, the horizontal limit algebras constructed from a
$\lambda$-lattice in Popa~\cite{Popa1995} form a tunnel with index
$\tau^{-1}$ and the proof in~\cite{Popa1995} is valid in general, with
a small number of modifications.  

These sort of ideas appear in a small part of~\cite{Popa2002} where
Popa refines his axiomatization of the standard invariant of a finite
index $\IIone$ subfactor in terms of $\lambda$-lattices.  See Lemma
1.6 where the element $a'$ is nothing other than the Radon-Nikodym
derivative $w$ and the conditional expectations constructed are simply
those in $(\wtilde{Z},V)$.
\end{remark}


\chapter{Infinite Index Subfactors of Type $\II$}
\label{chapter: infinite index}

Most the literature on subfactors is concerned with finite index
subfactors, particularly $\IIone$ extremal subfactors.  The study of
infinite index subfactors really began with Herman and
Ocneanu~\cite{HermanOcneanu1989}.  The results that they announced
were proved and expanded upon by Enock and Nest~\cite{EnockNest1996},
where the basic results for infinite index subfactors are laid down,
although the main purpose of their paper is to characterize the
subfactors arising as cross-products by Kac algebras of discrete or
compact type.


We begin in Section~\ref{sect: inf index background} with some
background material on Hilbert-module bases before giving a summary of
results from Enock and Nest~\cite{EnockNest1996} on the basic
construction.

In Section~\ref{sect: IIone case} we exploit the additional structure
present for an infinite index inclusion of $\IIone$ factors to develop
computational tools based on the $k$-fold relative tensor product of
$M$ that sits densely in $M_{k-1}$.

After defining extremality and showing that our definition possesses
the usual properties in Section~\ref{sect: extremality}, followed by a
brief diversion into $N$-central vectors in Section~\ref{sect:
N-central vectors}, we are ready for the main results.  Motivated by
the finite index case we define the rotation operators on the
$N$-central vectors in $\Ltwo(M_k)$ in Section~\ref{sect: rotations}
and in Section~\ref{sect: rotations and extremality} show that the
rotations exist iff the subfactor is approximately extremal.

Cross products by outer actions of an infinite discrete groups are
extremal and provide the simplest examples.  The restriction to the
$\Ltwo$-spaces avoids the sort of pathologies that we see in
Section~\ref{type III rel comm}.  However, as the example of Izumi,
Longo and Popa~\cite{IzumiLongoPopa1998} shows, there exist
irreducible subfactors which are not approximately extremal.  Future
work involves defining a rotation on a certain subspace of $N' \cap
M_k$ for any infinite index $\IIone$ subfactor.




\section{Background on Infinite Index Subfactors}
\label{sect: inf index background}

We will make heavy use of the material on operator-valued weights,
Hilbert $A$-modules and relative tensor products described in
Sections~\ref{sect: op-val weights} and~\ref{sect: Hilb mod and rel
tensor prod}.  The reader is advised to reacquaint themself with those
sections before proceeding.
In this section we add some additional results about bases before
introducing the basic construction, where we mostly follow Enock and
Nest~\cite{EnockNest1996}.  They consider arbitrary inclusions of
factors equipped with normal faithful semifinite (n.f.s.) weights.  We
will generally stick to inclusions $P \subset Q$ of arbitrary type
$\II$ factors, with traces $\Trace_P$ and $\Trace_Q$ respectively.

\subsection{Bases}

\begin{definition}
\label{def: basis}
An {\em $\leftidx{_A}{\h}{}$-basis} is a set $\{ \xi_i \} \subset
D(\leftidx{_A}{\h}{})$ such that
\[
\sum_i R(\xi_i) R(\xi_i)^* = 1_{\h} .
\]
An {\em $\h_A$-basis} is $\{ \xi_i \}
\subset D(\h_A)$ such that
\[
\sum_i L(\xi_i) L(\xi_i)^* = 1_{\h} .
\]
A $\h_A$- (resp. $\leftidx{_A}{\h}{}$-) basis is called {\em
orthogonal} if $L(\xi_i) L(\xi_i)^*$ (resp. $R(\xi_i) R(\xi_i)^*$) are
pairwise orthogonal projections.  Equivalently:
$L(\xi_i)^*L(\xi_j)=\delta_{i,j}p_i$
(resp. $R(\xi_i)^*R(\xi_j)=\delta_{i,j}p_i$) for some projections $p_i
\in A$.  If $p_i=1$ for all $i$ we say that the basis is {\em
orthonormal}.
\end{definition}

\begin{remark}
The existence of an $\leftidx{_A}{\h}{}$-basis or an
$\h_A$-basis is proved by Connes, Prop~3(c)
of~\cite{Connes1980}.  Given a conjugate-linear isometric involution
$J$ on $\h$, $JAJ$ is isomorphic to $A^{\text{op}}$, so that
$\h$ has a natural right-Hilbert-$A$-module structure, which
we will denote $\h_A$.  In this case, if $\{ \xi_i \}$ is an
$\leftidx{_A}{\h}{}$- (resp. $\h_A$-) basis then $\{
J \xi_i \}$ is an $\h_A$-
(resp. $\leftidx{_A}{\h}{}$-) basis.
\end{remark}

\begin{lemma}
\label{lem:Tr_Q' by basis}
Let $Q$ be a type $\II$ factor represented on a Hilbert space $\h$.
Let $\xi\in D(\leftidx{_Q}{\h}{})$ and let $\{ \xi_i \}$ be a
$\leftidx{_Q}{\h}{}$-basis.  Then
\begin{description}
\item{(i)} $\Trace'_{Q'\cap\B(\Ltwo(Q))}\br{R(\xi)^*R(\xi)}
\defeq \Trace_Q\br{J_QR(\xi)^*R(\xi)J_Q}
=||\xi||^2$.
\item{(ii)} $\sum_i \ip{R(\xi)R(\xi)^*\xi_i,\xi_i}=||\xi||^2$.
\item{(iii)} For $x\in\br{Q'\cap\B(\h)}_+$, $\sum_i \ip{x\xi_i,\xi_i}$
is independent of the basis used and, up to scaling,
\[
\Trace'_Q = \sum_i \ip{\spdot\xi_i,\xi_i} .
\]
\end{description}

\begin{proof}
Let $\psi'=\Trace'_{Q'\cap\B(\Ltwo(Q))}$.
\begin{description}
\item{(i)} This is simply Lemma~4 of Connes~\cite{Connes1980}.  In the
type $\II$ case the proof is particularly simple: Take projections
$p_i \in \gothn_{\Trace_Q}$ with $p_i \nearrow 1$.  Then
$\Trace_Q=\lim_i <\spdot \hatt{p_i},\hatt{p_i}>$ and hence
\[
\Trace_Q(J_Q R(\xi)^*R(\xi) J_Q)
 = \lim_i ||R(\xi)\hatt{p_i}||^2
 = \lim_i ||p_i \xi||^2
 = ||\xi||^2 .
\]
\item{(ii)} Let $x\in \br{Q'}_+$.  Using (i),
\begin{align}
\ip{x\xi_i,\xi_i}
=||x^{1/2}\xi_i||^2
&=\psi'\br{R(x^{1/2}\xi_i)^*R(x^{1/2}\xi_i)} \nonumber \\
&=\psi'\br{R(\xi_i)^*xR(\xi_i)} . \label{eq: 1}
\end{align}
Hence
\begin{align*}
\sum_i \ip{R(\xi)R(\xi)^*\xi_i,\xi_i}
&= \sum_i \psi'\br{R(\xi_i)^*R(\xi)R(\xi)^*R(\xi_i)} \\
&= \sum_i \psi'\br{R(\xi)^*R(\xi_i)R(\xi_i)^*R(\xi)} \\
&= \psi'\br{R(\xi)^*R(\xi)}
 = ||\xi||^2 .
\end{align*}
\item{(iii)} Given a basis $\Xi=\{ \xi_i \}$ define a normal weight
$\phi'_\Xi$ on $\br{Q'\cap\B(\h)}_+$ by $\phi'_{\Xi}=\sum_i
\ip{\spdot\xi_i,\xi_i}$.  For bases $\Xi$ and $\wtilde{\Xi}$ and
$x\in Q'\cap\B(\h)$, use (\ref{eq: 1}) to obtain
\begin{align*}
\phi'_\Xi(x^*x)
&=\sum_i \psi'\br{R(\xi_i)^*x^* 1_{\h} x R(\xi_i)} \\
&=\sum_i \sum_j \psi'((R(\xi_i)^*x^*R(\wtilde{xi}_j))
		      (R(\wtilde{\xi}_j)^*x R(\xi_i))) \\
&=\sum_i \sum_j \psi'(R(\wtilde{\xi}_j)^*x R(\xi_i)
		      R(\xi_i)^* x^*R(\wtilde{\xi}_j))
 &\text{since $\psi'$ is tracial} \\
&=\sum_j \sum_i \psi'(R(\wtilde{\xi}_j)^*x R(\xi_i)
		      R(\xi_i)^* x^*R(\wtilde{\xi}_j)) \\
&=\sum_j \psi'(R(\wtilde{\xi}_j)^*xx^* R(\wtilde{\xi}_j)) \\
&=\phi'_{\wtilde{\Xi}}(xx^*)
\end{align*}
Hence $\phi'_\Xi$ is tracial (taking $\wtilde{\Xi}=\Xi$) and
$\psi'_\Xi=\psi'_{\wtilde{\Xi}}$.  By (ii) $\psi'_\Xi \neq \infty$ and
so by uniqueness of the trace on a type $\II$ factor (up to scaling in
the $\II_\infty$ case), $\psi'_\Xi=\Trace_{Q'}$.
\end{description}
\end{proof}
\end{lemma}

\begin{definition}
\label{def: cannonical trace}
Let $Q$ be a type $\II$ factor with trace $\Trace_Q$.  Given
$\leftidx{_Q}{\h}{}$ there is a canonical choice of scaling for the
trace on $Q'\cap\B(\h)$ given by
\[
\Trace_{Q'\cap\B(\h)} = \sum_i \ip{\spdot\xi_i,\xi_i}
\]
where $\{ \xi_i \}$ is any $\leftidx{_Q}{\h}{}$-basis.  [Note that if
$Q'$ is a $\IIone$ factor this may not be the normalized trace on $Q'$].
\end{definition}

\begin{cor}
\label{cor: Tr_Q' via basis extended}
For all $x \in \plushat{(Q' \cap \BofH)}$ we have
\[
\Trace_{Q'}(x)=\sum_i x(\omega_{\xi_i}) ,
\]
where $\Trace_{Q'}$ now denotes the extension of the trace on $Q'$ to
$\plushat{(Q')}$.
\begin{proof}
Let $x \in \plushat{(Q' \cap \BofH)}$.  Take $x_k \in (Q' \cap
\BofH)_+$ with $x_k \nearrow x$ (Prop~\ref{prop: Haagerup ext pos
part}).  Then
\begin{align*}
\Trace_{Q'}(x)
&= \lim_k \Trace_{Q'}(x_k) & \text{by Prop~\ref{prop: Haagerup 1.10}} \\
&= \lim_k \sum_i <x_k \xi_i, \xi_i> & \text{by Lemma~\ref{lem:Tr_Q' by basis}}\\
&= \sum_i \lim_k <x_k \xi_i, \xi_i> & \text{since $x_k$ is increasing} \\
&= \sum_i \lim_k x_k(\omega_{\xi_i}) \\
&= \sum_i x(\omega_{\xi_i}) & \text{by definition.}
\end{align*}
\end{proof}
\end{cor}


\subsection{The basic construction}

\begin{definition}
Let $P \subset Q$ be an inclusion of type $\II$ factors.  The basic
construction applied to $P \subset Q$ is $P \subset Q \subset Q_1$,
where
\[
Q_1=J_Q P' J_Q .
\]
\end{definition}
Enock and Nest~\cite{EnockNest1996} show in 2.3 that an alternative
description of the basis construction is
\[
Q_1=\{L(\xi)L(\eta)^* : \xi, \eta \in D(\Ltwo(Q)_P)\}'' .
\]

\begin{prop}[Enock and Nest~\cite{EnockNest1996} 10.6, 10.7]
\label{prop: Enock Nest 10.6-10.7}
Let $T:Q_+\rightarrow \hatt{P}_+$ be the unique trace-preserving
operator-valued weight.  Then
\begin{itemize}
\item $\gothn_T \cap \gothn_{\Trace_Q}$ is weakly dense in $Q$ and
also dense in $\Ltwo(Q)$.
\item For $x \in \gothn_T$ there is a
bounded operator $\Lambda_T(x) \in \Hom_{-P}(\Ltwo(P),\Ltwo(Q))$
defined by
\begin{align*}
&\Lambda_T(x)\hatt{a}=\hatt{xa} 
&\text{for all } a \in \gothn_{\Trace_P}
\end{align*}
(in addition $xa \in \gothn_T \cap \gothn_{\Trace_Q}$).
\item The adjoint of $\Lambda_T(x)$ satisfies
\begin{align*}
&\Lambda_T(x)^*\hatt{z}=\hatt{T(x^*z)} 
&\text{for all } z \in \gothn_T \cap \gothn_{\Trace_Q} .
\end{align*}
\item For $x,y \in \gothn_T$,
\[
\Lambda_T(x)^*\Lambda_T(y)=T(x^*y) .
\]
\item $Q_1=\{ \Lambda_T(x)\Lambda_T(y)^* : x,y \in \gothn_T \}''$
\item The n.f.s trace-preserving operator valued weight $T_Q:\br{ Q_1
}_+ \rightarrow \plushat{Q}$ satisfies
\[
T_Q(\Lambda_T(x)\Lambda_T(y)^*)=xy^* .
\]
\end{itemize}
\end{prop}

\begin{notation}
Herman and Ocneanu~\cite{HermanOcneanu1989} use the notation $x
\tensor_P y^*$ for the operator $\Lambda_T(x) \Lambda_T(y)^*$ because
$\Lambda_T(x) \Lambda_T(y)^*$ is clearly $P$-middle-linear and, for $z
\in \gothn_T \cap \gothn_{\Trace_Q}$,
\[
(x \tensor_P y^*) \hatt{z} = \hatt{xT(y^*z)} ,
\]
so that in the finite index $\IIone$ case $x \tensor_P y^*$ is simply
$xe_1y^*$ which is $x \tensorlP y^*$ under the isomorphism $Q_1 \isom
Q \tensorlP Q$.
\end{notation}

The basic construction can be iterated to obtain a tower
\[
P \overset{T_P}{\subset} Q \overset{T_Q}{\subset} Q_1 
  \overset{T_{Q_1}}{\subset} Q_2 \cdots 
\]
and many results from finite index carry over to the general case.

\begin{prop}
\label{prop: additional info basic constr}
\begin{description}
\item
\item{$\bullet$} $\Ltwo(Q_1) \isom \Ltwo(Q) \tensor_P \Ltwo(Q)$ via
$L(\xi)L(\eta)^* \mapsto \xi \tensor_P J_Q \eta$ .  Note that for $x,y
\in \gothn_T \cap \gothn_{\Trace_Q}$ this agrees with Herman and
Ocneanu's $\tensor_P$ notation.
\item{$\bullet$} Consequently there is a unitary operator
$\theta_k:\tensor_P^{k+1} \Ltwo(Q) \rightarrow \Ltwo(Q_k)$ such that
\begin{align*}
\theta_k^* J_k \theta_k \br{ \xi_1 \tensorP \cdots \tensorP \xi_{k+1} }
&= J_0\xi_{k+1} \tensorP J_0\xi_k \tensorP \cdots \tensorP J_0\xi_1 \\
\theta_k^* x \theta_k
&= \pi_{k-1}(x) \tensorP \id ,
\end{align*}
where $x \in Q_k$ acts by left multiplication on $\Ltwo(Q_k)$ and
$\pi_{k-1}(x)$ is the (defining) representation of $Q_k$ on
$\Ltwo(Q_{k-1})$.
\item{$\bullet$} {\bf Multi-step basic construction:} $P \subset Q_i \subset Q_{2i+1}$
is also a basic construction.  In more detail, represent $Q_{2i+1}$ on
$\Ltwo(Q_i)\tensor_P \Ltwo(Q_{i-1})$ using $u=(\theta_i \tensor_P
\theta_{i-1})\theta_{2i}^*:\Ltwo(Q_{2i})\rightarrow \Ltwo(Q_i)\tensor_P
\Ltwo(Q_{i-1})$, then
\[
u Q_{2i+1} u^*
= J_i P' J_i \tensorP \id .
\]
By using $Q_j$ in place of $P$ we obtain a representation of $Q_k$
($i+1\leq k\leq 2i+1$) on $\Ltwo(Q_i)$.  Denote this representation
$\pi^k_i$ or simply $\pi_i$ if $k$ is clear.
\item{$\bullet$} {\bf Shifts:} By the above multi-step basic construction
$j_i\defeq J_i (\spdot)^* J_i$ gives anti-isomorphisms $j_i:N' \cap
M_{2i+1} \rightarrow N' \cap M_{2i+1}$, $j_i:M' \cap M_{2i+1}
\rightarrow N' \cap M_{2i}$, $j_i:N' \cap M_{2i} \rightarrow M' \cap
M_{2i+1}$.  Hence the {\em shift} $sh_i=j_{i+1}j_i$ gives isomorphisms
$sh_i:N' \cap M_{2i+1} \rightarrow M_1' \cap M_{2i+3}$, $sh_i:N' \cap
M_{2i} \rightarrow M_1' \cap M_{2i+2}$, $sh_i:M' \cap M_{2i+1}
\rightarrow M_2' \cap M_{2i+3}$, $sh_i:M' \cap M_{2i} \rightarrow M_2'
\cap M_{2i+2}$.
\end{description}
\end{prop}

\begin{remark}
Prop~\ref{prop: additional info basic constr} (i) is Theorem 3.8 of
Enock-Nest~\cite{EnockNest1996}, which is a reformulation of 3.1 of
Sauvageot~\cite{Sauvageot1983}.  Their general result includes the
spatial derivative, in this case
$\dee\Trace_{Q_1}/\dee\Trace_{P^{\text{op}}}$, but
$\Trace_{Q_1}\br{L(\xi)L(\xi)^*}=||\xi||^2$ by Lemma~\ref{lem:Tr_Q' by
basis} and hence $\dee\Trace_{Q_1}/\dee\Trace_{P^{\text{op}}}=1$.

We remark once again that for a finite index inclusion of $\IIone$
factors the canonical trace on $M_i$ is not the normalized trace.  If
one wishes to use normalized traces then Prop~\ref{prop: additional
info basic constr} and many other results here will need to be
modified with appropriate constants.
\end{remark}

\newpage

\section{The $\IIone$ Case}
\label{sect: IIone case}

We now consider the special case of an inclusion $N \subset M$ of
$\IIone$ factors with $[M:N]=\infty$.  We will reserve $N \subset M$
for inclusions of type $\IIone$ factors and use tend to use $P \subset
Q$ for general type $\II$ inclusions.  For $N \subset M$ additional
structure is provided by the existence of some of the Jones
projections and the embedding of $M$ in $\Ltwo(M)$.  We prove some
technical lemmas leading up to existence of the odd Jones projections.
We then give explicit definitions of the isomorphisms $\theta_k$ from
Prop~\ref{prop: additional info basic constr} and establish a number
of useful properties of these maps.

We conclude this section by constructing a basis for $M$ over $N$,
which will also allow us to construct bases for $M_j$ over $M_k$.

\subsection{Odd Jones projections and conditional expectations}

\begin{lemma}
\label{lemma: Trace = norm}
Let $P \subset Q$ be an inclusion of type $\II$ factors and let $P
\subset Q \subset Q_1$ be the basic construction.  Then:
\begin{description}
\item{(i)} For $\eta \in D(\Ltwo(Q)_P)$
\[
\Trace_{Q_1}\br{ L(\eta) L(\eta)^* }
= \Trace_P \br{ L(\eta)^* L(\eta) }
= || \eta ||^2 .
\]
\item{(ii)} Hence, for $\xi$ also in $D(\Ltwo(Q)_P)$,
$L(\eta)L(\xi)^*$ is trace-class and
\[
\Trace_{Q_1}\br{L(\eta)L(\xi)^*} = < \eta, \xi > .
\]
\item{(iii)} For $x \in \gothn_T$
\[
\Trace_{Q_1}\br{ x \tensorP x^* } = \Trace_Q(xx^*) .
\]
\end{description}
\begin{proof}
\begin{description}
\item{(i)} This is just Lemma~\ref{lem:Tr_Q' by basis} with $P$ in
place of $Q$ and a right-module in place of a left-module.
\item{(ii)} This is the usual polarization trick and the fact that the
product of two Hilbert-Schmidt operators is trace-class.
\item{(iii)} Let $\{ \xi_i \}$ be a $\Ltwo(Q)_P$-basis and $p_j$ a
sequence of projections in $P$ increasing to $1$.  Then
\begin{align*}
\Trace_{Q_1}\br{ \Lambda_T(x) \Lambda_T(x)^* }
&= \sum_i \ip{ \Lambda_T(x)\Lambda_T(x)^* \xi_i , \xi_i }_{\Ltwo(Q)} \\
&= \sum_i \lim_j \ip{ L(\xi)^* \Lambda_T(x)\Lambda_T(x)^* L(\xi_i) \hatt{p_j}
       , \hatt{p_j} }_{\Ltwo(P)} \\
&= \sum_i \Trace_P\br{\Lambda_T(x)^* L(\xi_i)L(\xi)^* \Lambda_T(x) } \\
&= \Trace_P\br{\Lambda_T(x)^*\Lambda_T(x) } \\
&= \Trace_P(T(x^*x))
= \Trace_Q(xx^*) .
\end{align*}
\end{description}
\end{proof}
\end{lemma}

\begin{lemma}
Let $P \subset Q$ be type $\II$ factors, $P \subset Q \subset Q_1
\subset Q_2 \subset \cdots$ the tower.  Suppose
$\Trace_{Q_1}(Q_+)=\{0,\infty\}$.  Then
$\Trace_{Q_3}\br{\br{Q_2}_+}=\{0,\infty\}$.
\begin{proof}
Take a $\Ltwo(Q)_P$-basis $\{ \xi_i \}$.  Then $\{ \xi_i \tensor_P
\xi_j \}$ is an $\Ltwo(Q_1)_P$-basis (Theorem 3.15
of~\cite{EnockNest1996}).  Let $x\geq 0$ be an element of $Q_2 =
\mathrm{End}_{-Q}\br{ \Ltwo(Q_1) }$.  Then
\begin{align*}
\Trace_{Q_3}(x)
&= \sum_{i,j} \ip{ x \xi_i \tensor_P \xi_j , \xi_i \tensor_P \xi_j 
       }_{\Ltwo(Q_1)} \\
&= \sum_{i,j} \ip{ \br{ L_{\xi_i}}^* x L_{\xi_i} \xi_j , \xi_j
       }_{\Ltwo(Q)} \\
&= \sum_i \Trace_{Q_1}\br{ \br{ L_{\xi_i}}^* x L_{\xi_i}} 
\in \{ 0,\infty \} ,
\end{align*}
where we have used the fact that $\br{ L_{\xi_i}}^* x
L_{\xi_i} \in \mathrm{End}_{-Q}\br{ \Ltwo(Q)} = Q$, so each
term in the sum is either $0$ or $\infty$.
\end{proof}
\end{lemma}

\begin{lemma}
\label{lemma: TraceP=TraceQ}
Let $P \subset Q$ be an inclusion of type $\II$ factors with
$\left. \Trace_{Q} \right|_{P}=\Trace_P$.  Let $T:Q_+ \rightarrow
\plushat{P}$ be the trace-preserving operator-valued weight.  Let $e$
denote orthogonal projection from $\Ltwo(Q)$ onto $\Ltwo(P)$.  Then
\begin{description}
\item{(i)} $T$ is in fact a conditional expectation, which we will
denote $E$.
\item{(ii)} $e \hatt{x} = \hatt{E(x)}$ for $x \in \gothn_{\Trace_Q}$.
\item{(iii)} $exe = E(x)e$ for $x \in Q$.
\item{(iv)} $e \in P' \cap Q_1$.
\item{(v)} For $x,y\in \gothn_T$, $x \tensorlP y*=xey^*$. 
\end{description}
\begin{proof}
\begin{description}
\item{(i)} $\Trace_Q(a1a^*)=\Trace_P(a1a^*)$ for all $a \in P$.  Hence
$T_P(1)=1$ so that $T_P$ is a conditional expectation.
\item{(ii)} Let $x \in \gothn_{\Trace_Q}$.  Then for all $y \in
\gothn_{\Trace_P}$:
\[
\ip{ \hatt{E(x)},\hatt{y}}
=\Trace_P(E(x)y^*)
=\Trace_P(E(xy^*))
=\Trace_Q(xy^*)
=\ip{\hatt{x},\hatt{y}}
=\ip{e\hatt{x},\hatt{y}} .
\]
To justify the third equality note that $xy^*$ is trace-class and
hence a linear combination of positive trace-class elements.  For each
of these positive elements $a$, $\Trace_P(E(a))=\Trace_Q(a)$ and hence
$\Trace_P(E(xy^*)) =\Trace_Q(xy^*)$.
\item{(iii)} Let $a \in \gothn_{\Trace_Q}$.  Then
$exe \hatt{a}
=ex \hatt{E(a)}
=\hatt{E(xE(a))}
=E(x)\hatt{E(a)}
=E(x)e \hatt{a}$.
\item{(iv)} For $x \in \gothn_{\Trace_Q}$, $a,b\in P$,
$e \hatt{(axb)}
=\hatt{E(axb)}
=\hatt{aE(x)b}
=a (e\hatt{x}) \cdot b$.
Hence $e \in End_{P-P}(\Ltwo(Q))=P' \cap Q_1$.
\item{(v)} For $z\in \gothn_T\cap\gothn_{\Trace_Q}$, $\br{x\tensorlP
y^*}\hatt{z}=\br{xT(y^*z)}\hatt{\phantom{A}}=xey^*\hatt{z}$.
\end{description}
\end{proof}
\end{lemma}

\begin{prop}
\label{prop: e_is}
\begin{description}
\item{}
\item{(i)} $\left. \Trace_{2i} \right|_{M_{2i-1}}=\Trace_{2i-1}$ so
that $T_{M_{2i-1}}:(M_{2i})_+ \rightarrow \plushat{(M_{2i-1})}$ is in
fact a conditional expectation, which we will denote $E_{M_{2i-1}}$.

As $\left. \Trace_{2i} \right|_{M_{2i-1}}=\Trace_{2i-1}$, let
$e_{2i+1}$ denote the orthogonal projection from $\Ltwo(M_{2i})$ onto
the subspace $\Ltwo(M_{2i-1})$.  Then $e_{2i+1} \in M_{2i-1}' \cap
M_{2i+1}$ and $E_{M_{2i-1}}$ is implemented by $e_{2i+1}$, i.e.
\begin{align*}
&e_{2i+1} x e_{2i+1} = E_{M_{2i-1}}(x) e_{2i+1}
&x \in M_{2i} .
\end{align*}

\item{(ii)} $\left. \Trace_{2i+1} \right|_{M_{2i}}\neq \Trace_{2i}$.
In fact $\Trace_{M_{2i+1}}((M_{2i})_+)=\{0,\infty\}$.  Hence
$T_{M_{2i}}:(M_{2i+1})_+ \rightarrow \plushat{(M_{2i})}$ is not a
conditional expectation.

\item{(iii)} $T_{M_{2i}}\br{ e_{2i+1}}=1$ and
\[
e_{2i+1} \tensorSE{M_{2i}} e_{2i+1} = e_{2i+1} .
\]

\item{(iv)} $\Trace_{2i+1}\br{e_1e_3\cdots e_{2i+1}}=1$.
\end{description}
\begin{proof}
\begin{description}
\item{(i)} By Lemma~\ref{lemma: TraceP=TraceQ} we simply need to show
that $\left. \Trace_{2i} \right|_{M_{2i-1}}=\Trace_{2i-1}$ for $i \geq
0$.  $\trace_M=\trace_N$ so the result is true for $i=0$.  Suppose
this is true for some $i$.  Let $P=M_{2i-1}$, $Q=M_{2i}$ so that
$Q_1=M_{2i+1}$ and $Q_2=M_{2i+2}$.  Let $\{ \xi_i \}$ be an
$\Ltwo(Q)_P$-basis.  Take projections $q_k \in Q$ with $\sum_k q_k =1$
and $\Trace_Q(q_k)<\infty$.  Then $\{ \hatt{q_k} \}$ is an
$\Ltwo(Q)_Q$-basis and $\{ \xi_i \tensor_P \hatt{q_k} \}$ is an
$\Ltwo(Q_1)_Q$-basis (easily checked since $L(\xi_i \tensor_P
\hatt{q_k})=L_{\xi_i} q_k$).  Now, for $x \in (M_{2i+1})_+$,
\begin{align*}
\Trace_{2i+2}(x)
&= \sum_{i,k} \ip{ x \xi_i \tensor_P \hatt{q_k}
       ,\xi_i \tensor_P \hatt{q_k} }_{\Ltwo(Q_1)} 
 = \sum_{i,k} \ip{ (x \xi_i) \tensor_P \hatt{q_k}
       ,\xi_i \tensor_P \hatt{q_k} }_{\Ltwo(Q_1)} \\
&= \sum_{i,k} \ip{ (x \xi_i) \ip{\hatt{q_k},\hatt{q_k}}_P
       ,\xi_i }_{\Ltwo(Q)} 
 = \sum_{i,k} \ip{ (x \xi_i) E_P\br{q_k}
       ,\xi_i }_{\Ltwo(Q)} \\
&= \sum_{i} \ip{ (x \xi_i) E_P(1)
       ,\xi_i }_{\Ltwo(Q)} 
 = \sum_{i} \ip{ x \xi_i
       ,\xi_i }_{\Ltwo(Q)} \\
&= \Trace_{2i+1}(x) .
\end{align*}

\item{(ii)} Note that $M_1$ is a $\II_{\infty}$ factor, hence
$\Trace_{M_1}(1)=\infty$.  For any nonzero projection in $M$ there is
a finite set of similar projections in $M$ with sum dominating $1$ and
hence infinite trace in $M_1$.  Thus $\Trace_{M_1}(p)=\infty$ for all
nonzero projections $p$ in $M$.  By spectral theory this implies that
$\Trace_{M_1}(x)=\infty$ for all nonzero $x$ in $M_+$.  Using the
preceding lemma we see that $\Trace_{M_{2i+1}}\br{\br{
M_{2i}}_+}=\{ 0,\infty \}$.  Hence for nonzero $a \in
M_{2i}$
\[
\infty
=\Trace_{2i+1}(a^*a)
=\Trace_{2i}(T_{M_{2i}}(a^*a))
=\Trace_{2i}(a^*T_{M_{2i}}(1)a)) ,
\]
so that $T_{M_{2i}}(1) \notin \mathbb{C}1$ and so $T_{M_{2i}}$ is not a
conditional expectation and not a multiple of a conditional
expectation.

\item{(iii)} Note that $\gothn_{T_{2i-1}}=\gothn_{E_{2i-1}}=M_{2i}$.
By Lemma~\ref{lemma: Trace = norm}
$\Trace_{2i+1}(xe_{2i+1}x^*)=\Trace_{2i}(xx^*)$ for all $x \in
M_{2i}$.  Thus $T_{M_{2i}}(e_{2i+1})=1$.

To show that $e_{2i+1} \tensor_{M_{2i}} e_{2i+1} = e_{2i+1}$ let
$z=ae_{2i+1}b$ where $a,b \in \gothn_{\Trace_{2i}}$ (and hence $z \in
\gothn_{\Trace_{2i+1}}\cap\gothn_{T_{Q_{2i}}}$ because
$z^*z=b^*E_{M_{2i-1}}(a^*a)e_{2i+1}b$).  Then
\begin{align*}
\br{e_{2i+1} \tensorSE{M_{2i}} e_{2i+1}} \hatt{z}
&= \sqbr{e_{2i+1} T_{M_{2i}}\br{ e_{2i+1} z }}\hatt{\phantom{A}} 
 = \sqbr{e_{2i+1} T_{M_{2i}}\br{ E_{M_{2i-1}}(a)e_{2i+1}b}}\hatt{\phantom{A}} \\
&= \sqbr{e_{2i+1} E_{M_{2i-1}}(a) b}\hatt{\phantom{A}} 
 = e_{2i+1} \hatt{z} .
\end{align*}
Since the span of such elements $z$ is dense in $\Ltwo(M_{2i+1})$ we
have shown $e_{2i+1} \tensor_{M_{2i}} e_{2i+1} = e_{2i+1}$.
\item{(iv)} The proof is by induction.
$\Trace_1(e_1)=\trace(T(e_1))=\trace(1)=1$ and 
\begin{align*}
\Trace_{2i+1}\br{e_1\cdots e_{2i+1}}
&=\Trace_{2i}\br{T_{M_{2i}}(e_1\cdots e_{2i+1})} 
 =\Trace_{2i}\br{e_1\cdots e_{2i-1}T_{M_{2i}}(e_{2i+1})} \\
&=\Trace_{2i}\br{e_1\cdots e_{2i-1}} 
 =\Trace_{2i-1}\br{e_1\cdots e_{2i-1}} .
\end{align*}
\end{description}
\end{proof}
\end{prop}


\subsection{Properties of the isomorphisms $\theta_k$}

In this section we explicitly define the isomorphisms $\theta_k:
\tensor_N^{k+1} \Ltwo(M) \rightarrow \Ltwo(M_k)$.

\begin{definition}
Given $P \subset Q$ an inclusion of type $\II$ factors.  For $r \geq
1$ define $J=J_{Q,P,r}:\tensor_P^r \Ltwo(Q) \rightarrow \tensor_P^r
\Ltwo(Q)$ by
\[
J\br{ \xi_1 \tensorSE{P} \xi_2 \tensorSE{P} 
        \cdots \tensorSE{P} \xi_r }
= J_Q\xi_r \tensorSE{P} \cdots \tensorSE{P} J_Q\xi_1 .
\]
$J$ is a conjugate-linear isometry onto $\tensor_P^r \Ltwo(Q)$.
\end{definition}

\begin{notation}
Let $v_{k+1}:\Ltwo(Q_k) \tensor_{Q_{k-1}} \Ltwo(Q_k) \rightarrow
\Ltwo(Q_{k+1})$ denote the isomorphism defined by
\begin{align*}
&\xi \tensorSE{Q_{k-1}} J_k \eta \mapsto L(\xi)L(\eta)^*
&\xi,\eta \in D(\Ltwo(Q_k)_{Q_{k-1}}) .
\end{align*}
Let $\iota_k:\Ltwo(Q_k) \tensor_{Q_k} \Ltwo(Q_k) \rightarrow
\Ltwo(Q_k)$ denote the isomorphism defined by
\begin{align*}
& \hatt{x} \tensorSE{Q_k} \hatt{y} \mapsto \hatt{xy}
&x,y \in \gothn_{\Trace_k} .
\end{align*}
\end{notation}

Note that both of these maps are $Q_k$-$Q_k$ bimodule maps and both
{\em preserve $J$}, i.e.
\begin{align*}
J_{k+1} v_{k+1}
&= v_{k+1} J_{Q_k,Q_{k-1},2} \\
J_k \iota_k
&= \iota_k J_{Q_k,Q_k,2}
\end{align*}

\begin{definition}
Define a $Q_k$-$Q_k$ bimodule isomorphism
\[
\psi_{k,r+1}: \tensorSE{Q_{k-1}}^{r+1} \Ltwo(Q_k)\rightarrow
              \tensorSE{Q_k}^r \Ltwo(Q_{k+1})
\]
by
\[
\psi_{k,r+1} = \br{\tensor_{Q_k}^r v_{k+1} } \circ
               \br{ \id_k \tensorSE{Q_{k-1}}
                  \br{ \tensorSE{Q_{k-1}}^{r-1} \iota_k^* }
                  \tensorSE{Q_{k-1}} \id_k }
\]
i.e.
\begin{gather*}
\Ltwo(Q_k) \tensorS{Q_{k-1}} \Ltwo(Q_k) \tensorS{Q_{k-1}} \cdots
  \tensorS{Q_{k-1}} \Ltwo(Q_k) \tensorS{Q_{k-1}} \Ltwo(Q_k) \\
\downarrow \\
\Ltwo(Q_k) \! \tensorS{Q_{k-1}} \! 
  \br{\! \Ltwo(Q_k)\tensorS{Q_k}\Ltwo(Q_k)\!}
  \tensorS{Q_{k-1}} 
  \! \cdots \!
  \tensorS{Q_{k-1}} \!
  \br{\!\Ltwo(Q_k)\tensorS{Q_k}\Ltwo(Q_k)\!}
  \tensorS{Q_{k-1}} \! \Ltwo(Q_k) \\
|| \\
\br{\Ltwo(Q_k)\tensorS{Q_{k-1}}\Ltwo(Q_k)} \tensorS{Q_k} \cdots
  \tensorS{Q_k} \br{\Ltwo(Q_k)\tensorS{Q_{k-1}}\Ltwo(Q_k)} \\
\downarrow \\
\Ltwo(Q_{k+1}) \tensorS{Q_k} \cdots \tensorS{Q_k} \Ltwo(Q_{k+1})
\end{gather*}
\end{definition}

\begin{definition}
Define a $Q$-$Q$ bimodule isomorphism 
$\theta_r :\tensorlP^{r+1} \Ltwo(Q)\rightarrow \Ltwo(Q_r)$ by
\[
\theta_r = \psi_{r-1,2} \circ \psi_{r-2,3} \cdots \psi_{0,r+1} .
\]
i.e.
\begin{gather*}
\Ltwo(Q) \tensorS{P} \Ltwo(Q) \tensorS{P} \Ltwo(Q) \tensorS{P} \cdots 
 \tensorS{P} \Ltwo(Q) \tensorS{P} \Ltwo(Q) \tensorS{P} \Ltwo(Q) \\
\downarrow \\
\Ltwo(Q_1) \tensorS{Q} \Ltwo(Q_1) \tensorS{Q} \cdots \tensorS{Q} \Ltwo(Q_1) 
  \tensorS{Q} \Ltwo(Q_1) \\
\downarrow \\
\vdots \\
\Ltwo(Q_{r-2}) \tensorS{Q_{r-3}}\Ltwo(Q_{r-2})\tensorS{Q_{r-3}} \Ltwo(Q_{r-2}) \\
\downarrow \\
\Ltwo(Q_{r-1}) \tensorS{Q_{r-2}} \Ltwo(Q_{r-1}) \\
\downarrow \\
\Ltwo(Q_r)
\end{gather*}
In general define $\theta_r^i: \tensor_{Q_{i-1}}^{r+1}\Ltwo(Q_i)
\rightarrow \Ltwo(Q_{i+r})$ by
\[
\theta_r^i = \psi_{i+r-1,2} \circ \psi_{i+r-2,3} \cdots \psi_{i,r+1} .
\]
\end{definition}

Note that the maps $\psi$ and hence $\theta$ also preserve $J$,
because $V$ and $\iota$ do.  Hence
\[
J_{r} \theta_r\br{ \xi_1 \tensorP \cdots \tensorP \xi_{r+1} }
= \theta_r\br{ J_0\xi_{r+1} \tensorP \cdots \tensorP J_0\xi_1 } .
\]
Also note that $\theta_{k+1}=\theta^0_{k+1}=\theta^1_k\circ
\psi_{0,k+2}$.

\begin{lemma}
\label{lemma: TC=HS^2}
Let $P \subset Q$ be an inclusion of type $\II$ factors.  Let
$T=T^Q_P$.  Any element $x \in \gothm_T \cap \gothm_{\Trace_Q}$ can be
written as $x=y^*z$ for some $y,z\in \gothn_T \cap \gothn_{\Trace_Q}$.
\begin{proof}
Consider the polar decomposition $x=u|x|$ and let $y=u|x|^{1/2}$,
$z=|x|^{1/2}$.  Then $z^*z=x$ and $y^*y=|x|^{1/2}u^*u|x|^{1/2}\leq
|x|$.  Hence $z^*z, y^*y \in \gothm_T \cap \gothm_{\Trace_Q}$ and
$z, y \in \gothn_T \cap \gothn_{\Trace_Q}$.
\end{proof}
\end{lemma}

\begin{lemma}
\label{lemma: HS tensor HS in TC}
Let $P \subset Q$ be an inclusion of type $\II$ factors.  Let
$T=T^Q_P$.  Let $a,b\in \gothn_T \cap \gothn_{\Trace_Q}$.  Then
\[
a \tensorSE{P} b^* \in \gothm_{T^{Q_1}_Q} \cap \gothm_{\Trace_{Q_1}} .
\]
\begin{proof}
Note that $L(\hat{a})=\Lambda_T(a)$, $L(\hat{b})=\Lambda_T(b)$, so $a
\tensorlP b^* = L(\hat{a})L(\hat{b})^* \in \gothm_{\Trace_{Q_1}}$ by
Lemma~\ref{lemma: Trace = norm}~(ii).  In addition $a \tensorlP b^* =
\Lambda_T(a)\Lambda_T(b)^* \in \gothm_{T^{Q_1}_Q}$ by Prop~10.7
of~\cite{EnockNest1996}.
\end{proof}
\end{lemma}

\begin{remark}
Note $\gothm_T \cap \gothm_{\Trace_Q} \subset \gothn_T \cap \gothn_T^*
\cap \gothn_{\Trace_Q}$.  Hence for $x \in \gothm_T \cap
\gothm_{\Trace_Q}$, $\hat{x}\in D(\Ltwo(Q)_P) \cap
D(\leftidx{_P}{\Ltwo(Q)}{})$ and $\Lambda_T(x)=L(\hat{x})$, $J_Q
\Lambda_T(x^*) J_P=R(\hat{x})$.
\end{remark}

\begin{prop}
\label{prop: tensorP}
Let $P \subset Q$ be an inclusion of type $\II$ factors.  Let
$T=T^Q_P$ and let $\wtilde{Q}=\gothm_T \cap \gothm_{\Trace_Q}$.  Let
$a_i, b_i \in \wtilde{Q}$, $i=1,2,\ldots\,$.  Let $s_j \in P$ be
defined inductively by $s_0=1$, $s_{j+1}=T(a_{j+1}^*s_jb_{j+1})$,
$j\geq 0$.  Then
\begin{description}
\item{(i)}
\[
\theta_k\br{ \hatt{a_1} \tensorP \cdots \tensorP \hatt{a_{k+1}} }
\in \hatr{\gothm_{T^{Q_k}_{Q_{k-1}}} \cap \gothm_{\Trace_{Q_k}}}
\]
and hence defines an element of $Q_k$ which we will denote $a_1
\tensorlP \cdots \tensorlP a_{k+1}$.  In the case $k=1$ this is the
same element as that represented in the Herman-Ocneanu notation by
$a_1 \tensorlP a_2$.
\item{(ii)}
\[
\br{ a_1 \tensorP a_2 \tensorP \cdots \tensorP a_{k+1} }^*
= a_{k+1}^* \tensorP a_k^* \tensorP \cdots \tensorP a_1^* .
\]
\item{(iii)}
\begin{multline*}
\br{a_1 \tensorP \cdots \tensorP a_{k+1} }^* \tensorSE{Q_{k-1}}
  \br{b_1 \tensorP\cdots\tensorP b_{k+1} } \\
= \begin{cases}
   a_{k+1}^* \tensorP \cdots \tensorP a_{r+1}^* \tensorP a_r^* s_{r-1} b_r 
     \tensorP \cdots \tensorP b_{k+1}
   &k=2r-1\\
   a_{k+1}^* \tensorP\cdots\tensorP a_{r+2}^* \tensorP a_{r+1}^* s_r \tensorP 
     b_{r+1} \tensorP \cdots \tensorP b_{k+1}
   &k=2r\\
 \end{cases}
\end{multline*} 
\item{(iv)}
\begin{multline*}
T^{Q_k}_{Q_{k-1}}\br{ a_1 \tensorP \cdots \tensorP a_{k+1} } \\
= \begin{cases}
   a_1 \tensorP \cdots \tensorP a_r a_{r+1} \tensorP \cdots \tensorP a_{k+1}
   &k=2r-1\\
   a_1 \tensorP \cdots \tensorP a_r T(a_{r+1}) \tensorP \cdots \tensorP a_{k+1}
   &k=2r\\
  \end{cases}
\end{multline*}
\item{(v)}
\begin{multline*}
\br{ a_1 \tensorP \cdots \tensorP a_{k+1} }^*
\br{ b_1 \tensorP \cdots \tensorP b_{k+1} } \\
= \begin{cases}
   a_{k+1}^* \tensorP \cdots \tensorP a_{r+1}^* s_r 
     \tensorP b_{r+1} \tensorP \cdots \tensorP b_{k+1}
   &k=2r-1\\
   a_{k+1}^* \tensorP \cdots \tensorP a_{r+1}^*
     s_r b_{r+1} \tensorP \cdots \tensorP b_{k+1}
   &k=2r\\
  \end{cases}
\end{multline*}
\item{(vi)} 
\begin{multline*}
\pi^k_i\br{a_1 \tensorP \cdots \tensorP a_{k+1}}^*
  \hatr{b_1 \tensorP \cdots \tensorP b_{i+1}} \\
=\begin{cases}
\hatr{a_{k+1}^* \tensorP \cdots \tensorP a_{r+1}^*s_r \tensorP b_{r+1}
    \tensorP \cdots \tensorP b_{i+1}}
& k=2r-1 \\
\hatr{a_{k+1}^* \tensorP \cdots \tensorP a_{r+1}^*s_rb_{r+1}
    \tensorP \cdots \tensorP b_{i+1}}
& k=2r
\end{cases}
\end{multline*}
\end{description}
\begin{proof}
By Lemma~\ref{lemma: TC=HS^2} there exist $w_i^*, x_i, y_i^*, z_i \in
\gothn_T \cap \gothn_{\Trace_Q}$ such that \\ $a_i=w_ix_i, b_i=y_iz_i$.  Let
$A_1=a_1 \tensorlP w_2, B_1=b_1 \tensorlP y_2$ and for $i\geq 2$ let
$A_i=x_i \tensorlP w_{i+1}, B_i=z_i \tensorlP y_{i+1}$.  For $j\geq 1$
let $\overline{A}_j=x_j \tensorlP a_{j+1}, \overline{B}_j=z_j
\tensorlP y_{j+1}$.  By Lemma~\ref{lemma: HS tensor HS in TC} $A_i,
B_i, \overline{A}_i, \overline{B}_i \in \gothm_{T^{Q_1}_Q} \cap
\gothm_{\Trace_{Q_1}}$.

\noindent Let $S_0=1$, $S_{j+1}=T^{Q_1}_Q(A_{j+1}^*S_jB_{j+1})$,
$j\geq 0$.  We claim that $S_j = w_{j+1}^* s_j y_{j+1}$ for $j \geq
1$.  For $j=1$, use the fact that
$\Lambda_T(a_1)^*\Lambda_T(b_1)=T(a_1^*b_1)$ from Prop~\ref{prop:
Enock Nest 10.6-10.7} to obtain
\[
A_1^* B_1
= \br{w_2^* \tensorP a_1^* }\br{b_1 \tensorP y_2 }
= w_2^* T(a_1^*b_1) \tensorP y_2
= w_2^* s_1 \tensorP y_2 ,
\]
Hence $S_1 = T_Q(A_1^* B_1) = w_2^* s_1 y_2$.  In general, if the
result holds for $j\geq 1$ then
\begin{align*}
S_{j+1}
&= T_Q(A_{j+1}^* S_j B_{j+1}) \\
&= T_Q\br{\br{w_{j+2}^* \tensorP x_{j+1}^*}w_{j+1}^*s_jy_{j+1} 
          \br{z_{j+1} \tensorP y_{j+2}}} \\
&= T_Q\br{w_{j+2}^* T(x_{j+1}^* w_{j+1}^*s_jy_{j+1}z_{j+1})
          \tensorP y_{j+2}} \\
&= w_{j+2}^* T(a_{j+1}^*s_jb_{j+1}) y_{j+2} \\
&= w_{j+2}^* s_{j+1} y_{j+2} .
\end{align*}

\begin{description}
\item{(i)} The result is trivially true for $k=0$.  For $k=1$, 
$\theta_1=\psi_{0,2}=v_1$ and we have, from the definition of $v_1$,
\[
\theta_1\br{ \hatt{a_1} \tensorP \hatt{a_2}}
=\hatr{L(\hatt{a_1})L(J_Q\hatt{a_2})^*}
=\hatr{\Lambda_T(a_1)\Lambda_T(a_2^*)^*} .
\]
$\Lambda_T(a_1)\Lambda_T(a_2^*)^*$ is in $\gothm_{T^{Q_1}_Q} \cap
\gothm_{\Trace_{Q_1}}$ by Lemma~\ref{lemma: HS tensor HS in TC} and is
simply $a_1 \tensorlP a_2$ in Herman-Ocneanu notation.

Assume that (i) is true for some $k\geq 1$.  Then
\begin{align*}
\begin{split}
&\theta_{k+1}\br{ \hatt{a_1} \tensorP \hatt{a_2} \tensorP \cdots
   \tensorP \hatt{a_{k+2}} } \\ 
&=\theta_{k+1}\br{ \hatt{a_1} \tensorP \hatt{w_2x_2} \tensorP \cdots
   \tensorP \hatt{w_{k+1}x_{k+1}} \tensorP \hatt{a_{k+2}}} \\ 
&=\theta^1_k\circ \psi_{0,k+2}\br{ \hatt{a_1} \tensorP \hatt{w_2x_2} 
   \tensorP \cdots \tensorP \hatt{w_{k+1}x_{k+1}} \tensorP \hatt{a_{k+2}}} \\ 
&=\theta^1_k\biggl( \hatr{a_1 \tensorP w_2} \tensorSE{Q}\hatr{x_2 \tensorP w_3}
   \tensorSE{Q} \cdots \\
& \phantom{=\theta^1_k\biggl( \hatr{a_1 \tensorP w_2}}
  \cdots \tensorSE{Q} \hatr{x_k \tensorP w_{k+1}} \tensorSE{Q}
   \hatr{x_{k+1} \tensorP a_{k+2}}\biggr) \\
&=\theta^1_k\br{ \hatt{A_1} \tensorSE{Q} \hatt{A_2} \tensorSE{Q} \cdots 
   \tensorSE{Q} \hatt{\overline{A}_{k+1}} } \\
&\in \hatt{\gothm_{T^{Q_{k+1}}_{Q_k}} \cap \gothm_{\Trace_{Q_{k+1}}}}
\end{split}
\end{align*}
by assumption (and the fact that if $C=Q$ and $D=Q_1$ then the tower
obtained from $C\subset D$ satisfies $D_j=Q_{j+1}$).

Note that in addition we have shown that
\[
a_1 \tensorP a_2 \tensorP \cdots \tensorP a_{k+2}
= A_1 \tensorSE{Q} A_2 \tensorSE{Q} \cdots \tensorSE{Q} A_{k}
  \tensorSE{Q} \overline{A}_{k+1} .
\]

\item{(ii)} Simply observe that
\[
J_k \theta_k\br{ \hatt{a_1} \tensorP \cdots \tensorP \hatt{a_{k+1}} }
= \theta_k\br{ J_0\hatt{a_{k+1}} \tensorP \cdots \tensorP J_0\hatt{a_1} 
           } 
= \theta_k\br{ \hatt{a_{k+1}^*} \tensorP \cdots \tensorP \hatt{a_1^*} 
           }
\]

\item{(iii)} For $k=0$, $a_1^* s_0 \tensorlP b_1 = a_1^* \tensorlP
b_1$.  For $k=1$, from the proof of (i),
\begin{align*}
a_2^* \tensorP a_1^* s_0 b_1 \tensorP b_2
&= a_2^* \tensorP a_1^* b_1 \tensorP b_2 \\
&= \br{a_2^* \tensorP a_1^*} \tensorSE{Q} \br{b_1 \tensorP b_2} \\
&= \br{a_1 \tensorP a_2}^* \tensorSE{Q} \br{b_1 \tensorP b_2} .
\end{align*}

Assume the result is true for some $k\geq 1$.  Note that
\begin{align*}
A_r^* S_{r-1} B_r
&=\br{w_{r+1}^*\tensorP x_r^*} w_r^*s_{r-1}y_r\br{z_r\tensorP y_{r+1}} \\
&=w_{r+1}^*T(x_r^* w_r^*s_{r-1}y_rz_r) \tensorP y_{r+1} \\
&=w_{r+1}^*s_r \tensorP y_{r+1} \\
A_{r+1}^* S_r
&= w_{r+2}^* \tensorP x_{r+1}^* w_{r+1}^*s_r y_{r+1} \\
&= w_{r+2}^* \tensorP a_{r+1}^*s_r y_{r+1} .
\end{align*}
Hence, with the first of the two cases always denoting $k=2r-1$ and
the second $k=2r$,
\begin{align*}
&\br{a_1 \tensorP \cdots \tensorP a_{k+2} }^* \tensorSE{Q_k}
  \br{b_1 \tensorP \cdots \tensorP b_{k+2} } \\
&=\br{A_1 \tensorSE{Q} \cdots \tensorSE{Q} A_k \tensorSE{Q} 
  \overline{A}_{k+1}}^* \tensorSE{Q_k} \br{B_1 \tensorSE{Q} \cdots 
  \tensorSE{Q} B_k \tensorSE{Q} \overline{B}_{k+1}} \\
&=\begin{cases}
 \overline{A}_{k+1}^* \tensorQ A_k \tensorQ \cdots \tensorQ A_r^* S_{r-1} 
  B_r \tensorQ \cdots \tensorQ B_k \tensorQ \overline{B}_{k+1} \\
 \overline{A}_{k+1}^* \tensorQ A_k \tensorQ \cdots \tensorQ A_{r+1}^* S_r 
  \tensorQ B_{r+1} \tensorQ \cdots \tensorQ B_k \tensorQ \overline{B}_{k+1}
 \end{cases} \\
&= \begin{cases}
 a_{k+2}^* \tensorP \cdots \tensorP a_{r+2}^* \tensorP 
  x_{r+1}^*w_{r+1}^*s_r \tensorP y_{r+1}z_{r+1} \tensorP b_{r+2} \tensorP 
  \cdots \tensorP b_{k+2} \\
 a_{k+2}^* \tensorP \cdots \tensorP a_{r+3}^* \tensorP x_{r+2}^*w_{r+2}^*
  \tensorP a_{r+1}^*s_r y_{r+1}z_{r+1} \tensorP b_{r+2} \tensorP 
  \cdots \tensorP b_{k+2} 
 \end{cases} \\
&= \begin{cases}
 a_{k+2}^* \tensorP a_{k+1}^* \tensorP \cdots \tensorP
  a_{r+1}^*s_r \tensorP b_{r+1} \tensorP \cdots
  \tensorP b_{k+2} \\
 a_{k+2}^* \tensorP a_{k+1}^* \tensorP \cdots \tensorP
  a_{r+1}^*s_r b_{r+1} \tensorP b_{r+2} \tensorP \cdots
  \tensorP b_{k+2}
 \end{cases}
\end{align*}
Hence the result is true for $k+1$ and the general result follows by
induction.

\item{(iv)} For $k=0$, $T^{Q_0}_{Q_1}(a_1)=T(a_1)$.  Assume the result
holds for some $k\geq 0$.  Note that
\begin{align*}
A_rA_{r+1}
&= x_r T(w_{r+1}x_{r+1}) \tensorP w_{r+2}
 = x_r T(a_{r+1}) \tensorP w_{r+2} \\
A_rT^{Q_1}_Q(A_{r+1}) ,
&= x_r \tensorP w_{r+1}x_{r+1}w_{r+2}
 = x_r \tensorP a_{r+1}w_{r+2} .
\end{align*}
Hence, with the first of the two cases once again denoting $k=2r-1$
and the second $k=2r$,
\begin{align*}
&T^{Q_{k+1}}_{Q_k}\br{a_1 \tensorP \cdots \tensorP a_{k+2}} \\
&= T^{Q_{k+1}}_{Q_k}\br{A_1 \tensorQ \cdots \tensorQ A_k \tensorQ 
  \overline{A}_{k+1}} \\
&= \begin{cases}
 A_1 \tensorQ \cdots \tensorQ A_rA_{r+1} \tensorQ \cdots \tensorQ A_k 
  \tensorQ \overline{A}_{k+1} \\
 A_1 \tensorQ \cdots \tensorQ A_rT^{Q_1}_Q(A_{r+1}) \tensorQ \cdots 
  \tensorQ A_k \tensorQ \overline{A}_{k+1}
 \end{cases} \\
&= \begin{cases}
 a_1 \tensorP \cdots \tensorP a_{r-1} \tensorP w_rx_r T(a_{r+1}) \tensorP 
  w_{r+2}x_{r+2} \tensorP a_{r+3} \tensorP \cdots \tensorP a_{k+2} \\
 a_1 \tensorP \cdots \tensorP a_{r-1} \tensorP w_rx_r \tensorP 
  a_{r+1}w_{r+2}x_{r+2} \tensorP a_{r+3} \tensorP \cdots a_{k+2}
 \end{cases} \\
&= \begin{cases}
 a_1 \tensorP \cdots \tensorP a_{r-1} \tensorP a_rT(a_{r+1}) \tensorP
  a_{r+2} \tensorP \cdots \tensorP a_{k+2} \\
 a_1 \tensorP \cdots \tensorP a_r \tensorP a_{r+1} a_{r+2} \tensorP
  a_{r+3} \tensorP \cdots \tensorP a_{k+2}
 \end{cases}
\end{align*}

\item{(v)} We use the fact that $T_Q(x\tensorP y)=xy$ combined with
(iii) and (iv).  Once again, let the first of the two cases denote
$k=2r-1$ and the second $k=2r$.
\begin{align*}
&\br{a_1 \tensorP \cdots \tensorP a_{k+1}}^*
 \br{b_1 \tensorP \cdots \tensorP b_{k+1}} \\
&= T^{Q_{k+1}}_{Q_k}\br{\br{a_1 \tensorP \cdots \tensorP a_{k+1}}^*
 \tensorSE{Q_{k-1}} \br{b_1 \tensorP \cdots \tensorP b_{k+1}}} \\
&= \begin{cases}
 T^{Q_{k+1}}_{Q_k}\br{a_{k+1}^* \tensorP \cdots a_{r+1}^* \tensorP a_r^* 
  s_{r-1} b_r \tensorP \cdots \tensorP b_{k+1}} \\
 T^{Q_{k+1}}_{Q_k}\br{a_{k+1}^* \tensorP \cdots a_{r+1}^* \tensorP 
  a_{r+1}^* s_r \tensorP b_{r+1} \tensorP \cdots \tensorP b_{k+1}}
 \end{cases} \\
&= \begin{cases}
 a_{k+1}^* \tensorP \cdots \tensorP a_{r+1}^* T(a_r^* s_{r-1} b_r) 
  \tensorP b_{r+1} \tensorP \cdots \tensorP b_{k+1} \\
 a_{k+1}^* \tensorP \cdots \tensorP a_{r+1}^* s_r b_{r+1} 
  \tensorP b_{r+2} \tensorP \cdots \tensorP b_{k+1}
 \end{cases} \\
&= \begin{cases}
 a_{k+1}^* \tensorP \cdots \tensorP a_{r+1}^* s_r
  \tensorP b_{r+1} \tensorP \cdots \tensorP b_{k+1} \\
 a_{k+1}^* \tensorP \cdots \tensorP a_{r+1}^* s_r b_{r+1} 
  \tensorP b_{r+2} \tensorP \cdots \tensorP b_{k+1}
 \end{cases}
\end{align*}

\item{(vi)} Let $j=2i-k$.  We will induct on $j$.  For $j=-1$
($k=2i+1$, $i \geq 0$), let $u=\br{\theta_i \tensorlP
\theta_{i-1}}\theta_{2i}^*$.  From Prop~\ref{prop: additional info
basic constr}
\begin{align*}
&\sqbr{\pi^{2i+1}_i\br{a_1 \tensorP \cdots \tensorP a_{2i+2}}^*
       \hatr{b_1 \tensorP \cdots \tensorP b_{i+1}}}
 \tensorP \hatr{c_1 \tensorP \cdots \tensorP c_i} \\
&=u \br{a_1 \tensorP \cdots \tensorP a_{2i+2}}^*
  \hatr{b_1 \tensorP \cdots \tensorP b_{i+1} \tensorP c_1 \tensorP \cdots 
    \tensorP c_i} \\
&=u \hatr{a_{2i+2}^* \tensorP \cdots \tensorP a_{i+2}^*s_{i+1} \tensorP 
    c_1 \tensorP \cdots \tensorP c_i} \\
&=\hatr{a_{2i+2}^* \tensorP \cdots \tensorP a_{i+2}^*s_{i+1}} \tensorP
   \hatr{c_1 \tensorP \cdots \tensorP c_i} .
\end{align*}
Suppose (vi) holds for some $\bar{\jmath}$ (and all $i>\bar{\jmath}$).
Let $j=\bar{\jmath}+1$, $i>j$ and $k=2i-j$.  Set $\bar{\imath}=i-1$,
$\bar{k}=2\bar{\imath}-\bar{\jmath}=k-1$.  Let $\bar{P}=Q$,
$\bar{Q}=Q_1$ and $\bar{\pi}$ denote the representations
constructed from Prop ~\ref{prop: additional info basic constr}
applied to $\bar{P} \subset \bar{Q}$.  Then
\begin{align*}
&\pi^k_i\br{a_1 \tensorP \cdots \tensorP a_{k+1}}^*
       \hatr{b_1 \tensorP \cdots \tensorP b_{i+1}} \\
&=\bar{\pi}^{\bar{k}}_{\bar{\imath}}
  \br{A_1 \tensorQ \cdots \tensorQ A_{\bar{k}} \tensorQ 
      \bar{A}_{\bar{k}+1}}^*
  \hatr{B_1 \tensorQ \cdots \tensorQ B_{\bar{\imath}} B_{\bar{\imath}+1}} \\
&=\begin{cases}
 \hatr{\bar{A}_{\bar{k}+1}^* \tensorQ \cdots \tensorQ 
  A_{\bar{r}+1}^*S_{\bar{r}} \tensorQ B_{\bar{r}+1} \tensorQ 
  \cdots \tensorQ \bar{B}_{\bar{\imath}+1}}
 & \bar{k}=2\bar{r}-1 \\
 \hatr{\bar{A}_{\bar{k}+1}^* \tensorQ \cdots \tensorQ 
 A_{\bar{r}+1}^*S_{\bar{r}} B_{\bar{r}+1} \tensorQ 
 \cdots \tensorQ \bar{B}_{\bar{\imath}+1}}
 & \bar{k}=2\bar{r}
 \end{cases} \\
&=\begin{cases}
 \hatr{a_{k+1}^* \tensorP \cdots \tensorP 
  a_{\bar{r}+1}^*s_{\bar{r}} b_{\bar{r}+1} \tensorP 
  \cdots \tensorP b_{i+1}}
 & k=2\bar{r} \\
 \hatr{a_{k+1}^* \tensorP \cdots \tensorP 
 a_{\bar{r}+2}^*s_{\bar{r}+1} \tensorP b_{\bar{r}+2} \tensorP 
 \cdots \tensorP b_{i+1}}
 & k=2\bar{r}+1
 \end{cases} \\
&=\begin{cases}
 \hatr{a_{k+1}^* \tensorP \cdots \tensorP 
  a_{r+1}^*s_r b_{r+1} \tensorP 
  \cdots \tensorP b_{i+1}}
 & k=2r \\
 \hatr{a_{k+1}^* \tensorP \cdots \tensorP 
 a_{r+1}^*s_r \tensorP b_{r+1} \tensorP 
 \cdots \tensorP b_{i+1}}
 & k=2r-1
 \end{cases}
\end{align*}
\end{description}
\end{proof}
\end{prop}

\begin{notation}
Let $\wtilde{Q}_k=\rmspan\{ a_1\tensorlP \cdots \tensorlP a_{k+1}:
a_i\in\wtilde{Q}\}$ (not to be confused with
$\br{Q_k}\wtilde{\phantom{A}}=\gothm_{T^{Q_k}_{Q_{k-1}}} \cap
\gothn_{\Trace_k}$).
\end{notation}

\begin{cor}
\label{cor: tensor 1 = prod e_odd}
Let $N \subset M$ be an inclusion of $\IIone$ factors.  Then
\begin{description}
\item{(i)} 
\[
\underbrace{1 \tensorN 1 \tensorN 1 \cdots \tensorN 1}_{2r \text{ or } 2r+1}
= e_1e_3\cdots e_{2r-1} .
\]
\item{(ii)} For $a_1, a_2, \ldots, a_{2r} \in M$
\[
a_1 \tensorN \cdots \tensorN a_r \tensorN 1 \tensorN a_{r+1} \tensorN \cdots
\tensorN a_{2r}
= a_1 \tensorN \cdots \tensorN a_{2r} .
\]
\item{(iii)} $\wtilde{M}_k$ is dense in $\Ltwo(M_k)$.
\end{description}
\begin{proof}
\begin{description}
\item{(i)} First note that $1=1$, $1 \tensorN 1=1e_11=e_1$.  Suppose that
\[
\underbrace{1 \tensorN 1 \tensorN 1 \cdots \tensorN 1}_{2r}
= e_1e_3\cdots e_{2r-1} .
\]
Then, by Prop~\ref{prop: tensorP} (iii),
\begin{align*}
\underbrace{1 \tensorN 1 \tensorN 1 \cdots \tensorN 1}_{2r+1}
&= \br{\underbrace{1 \tensorN 1 \tensorN 1 \cdots \tensorN 1}_{2r}}
 \tensorSE{M_{2r-2}}
 \br{\underbrace{1 \tensorN 1 \tensorN 1 \cdots \tensorN 1}_{2r}} \\
&= \br{e_1 e_3 \cdots e_{2r-1}} \tensorSE{M_{2r-2}}
  \br{e_1 e_3 \cdots e_{2r-1}} \\
&= \br{e_1 e_3 \cdots e_{2r-3}} e_{2r-1} \tensorSE{M_{2r-2}} e_{2r-1}
  \br{e_1 e_3 \cdots e_{2r-3}} \\
&= \br{e_1 e_3 \cdots e_{2r-3}} e_{2r-1}
  \br{e_1 e_3 \cdots e_{2r-3}} \\
&= e_1 e_3 \cdots e_{2r-1} ,
\end{align*}
using Prop~\ref{prop: e_is} (iii).  Also
\begin{align*}
\underbrace{1 \tensorN 1 \tensorN 1 \cdots \tensorN 1}_{2r+2}
&= \br{\underbrace{1 \tensorN 1 \tensorN 1 \cdots \tensorN 1}_{2r+1}}
 \tensorSE{M_{2r-1}}
 \br{\underbrace{1 \tensorN 1 \tensorN 1 \cdots \tensorN 1}_{2r+1}} \\
&= \br{e_1 e_3 \cdots e_{2r-1}} e_{2r+1}
  \br{e_1 e_3 \cdots e_{2r-1}} \\
&= e_1 e_3 \cdots e_{2r-1} e_{2r+1} .
\end{align*}
The result now follows by induction.

\item{(ii)} Let $X$ denote $a_1 \tensorlN \cdots \tensorlN a_r$ and
let $Y$ denote $a_{r+1} \tensorlN \cdots \tensorlN a_{2r}$.  Using
part (i) and Prop~\ref{prop: tensorP} (iii) and then (v),
\begin{align*}
& a_1 \tensorN \cdots \tensorN a_r \tensorN 1 \tensorN a_{r+1} \tensorN 
 \cdots \tensorN a_{2r} \\
&= \br{X \tensorN \underbrace{1 \tensorN \cdots \tensorN 1}_{r}} 
 \tensorSE{M_{2r-2}}
 \br{\underbrace{1 \tensorN \cdots \tensorN 1}_{r} \tensorN Y} \\
&= \sqbr{\!\br{\! X \tensorN \underbrace{1 \tensorN \cdots \tensorN 1}_{r}} 
    \!\!\br{\underbrace{1 \tensorN \cdots \tensorN 1}_{2r}}\!}
 \tensorSE{M_{2r-2}}
   \sqbr{\!\br{\underbrace{1 \tensorN \cdots \tensorN 1}_{2r}}
    \!\!\br{\underbrace{1 \tensorN \cdots \tensorN 1}_{r} \tensorN Y\!}\!} \\
&= \br{X \tensorN \underbrace{1 \tensorN \cdots \tensorN 1}_{r}} 
 \br{\underbrace{1 \tensorN \cdots \tensorN 1}_{2r+1}}
 \br{\underbrace{1 \tensorN \cdots \tensorN 1}_{r} \tensorN Y} \\
&= \br{X \tensorN \underbrace{1 \tensorN \cdots \tensorN 1}_{r}} 
 \br{\underbrace{1 \tensorN \cdots \tensorN 1}_{2r}}
 \br{\underbrace{1 \tensorN \cdots \tensorN 1}_{r} \tensorN Y} \\
&= \br{X \tensorN \underbrace{1 \tensorN \cdots \tensorN 1}_{r}} 
 \br{\underbrace{1 \tensorN \cdots \tensorN 1}_{r} \tensorN Y} \\
&= a_1 \tensorN \cdots \tensorN a_{2r} .
\end{align*}
\item{(iii)} $\wtilde{M}=M$ is dense in $\Ltwo(M)$.  The general result
follows from the following claim:
\begin{claim}
Suppose $\wtilde{\h}\subset D(\h_B)$ is dense in
$\leftidx{_A}{\h}{_B}$ and $\wtilde{\mathcal{K}} \subset
D(\leftidx{_B}{\mathcal{K}}{})$ is dense in
$\leftidx{_B}{\mathcal{K}}{_C}$.  Then $\rmspan\{ \xi \tensor_B \eta :
\xi \in \wtilde{\h}, \eta \in \wtilde{\mathcal{K}} \}$ is dense
in $\h \tensor_B \mathcal{K}$.
\begin{proof}
Given $\xi \in D(\h_B)$ and $\eta \in
D(\leftidx{_B}{\mathcal{K}}{})$ take $\eta_n \in \wtilde{\mathcal{K}}$
with $\eta_n \rightarrow \eta$.  Then
\begin{align*}
||\xi \tensorSE{B} \eta - \xi \tensor \eta_n||^2
&= ||\xi \tensorSE{B} \eta||^2 + ||\xi \tensorSE{B} \eta_n||^2
 -2\mathrm{Re}\ip{\xi \tensorSE{B} \eta,\xi \tensorSE{B} \eta_n} \\
&= \ip{\ip{\xi,\xi}_B \eta, \eta} + \ip{\ip{\xi,\xi}_B \eta_n, \eta_n} 
 -2\mathrm{Re}\ip{\ip{\xi,\xi}_B \eta_n, \eta_n} \\
&\rightarrow 0 ,
\end{align*}
so we can approximate $\xi \tensor_B \eta$ by $\xi \tensor \eta_0$ for
some $\eta_0 \in \wtilde{\mathcal{K}}$.  Similarly we can approximate
$\xi \tensor_B \eta_0$ by $\xi_0 \tensor_B \eta_0$ for some $\xi_0 \in
\wtilde{\h}$.
\end{proof}
\end{claim}
The result now follows by induction because $\wtilde{M}_k \subset
\gothm_{T^{M_k}_N} \cap \gothm_{\Trace_k}$ by iterating
Prop~\ref{prop: tensorP} (iv) and $\gothm_{T^{M_k}_N} \cap
\gothm_{\Trace_k} \subset D(\Ltwo(M_k)_N)$.
\end{description}
\end{proof}
\end{cor}


\subsection{Bases revisited}
\label{sect: IIone case - bases} 

\begin{definition}
Let $P \subset Q$ be an inclusion of type $\II$ factors.  A
$Q_P$-basis (also called a basis for $Q$ over $P$) is $\{ b \} \subset
\gothn_T \cap \gothn_{\Trace_Q}$ such that
\[
\sum_b b \tensorP b^* =1 .
\]
A $\leftidx{_P}{Q}{}$-basis is $\{ b \} \subset \gothn_T^* \cap
\gothn_{\Trace_Q}$ such that
\[
\sum_b b^* \tensorP b =1 .
\]
\end{definition}

\begin{remark}
Note that a $Q_P$-basis is a special case of a $\Ltwo(Q)_P$-basis.  We
define orthogonal and orthonormal $Q_P$- and $\leftidx{_P}{Q}{}$-bases
as in~\ref{def: basis}.

\noindent If $\{ b \}$ is a $Q_P$- (resp. $\leftidx{_P}{Q}{}$-) basis, then $\{
b^* \}$ is a $\leftidx{_P}{Q}{}$- (resp. $Q_P$-) basis.
\end{remark}

We will show that for a $\IIone$ inclusion $N \subset M$ there exists
an orthonormal basis for $M$ over $N$ and orthogonal bases for all
$M_k$ over $M_j$.  We will then relate bases for $Q$ over $P$ to the
commutant operator-valued weight.

\subsubsection{Existence}

We begin with the infinite index version of the so-called {\em
pull-down lemma}.

\begin{lemma}[Pull-down lemma]
\label{pull-down lemma}
Let $z \in \gothn_T^*$.  Then $z e_1=T(ze_1)e_1$.
\begin{proof}
Since $ze_1, T(ze_1)e_1 \in \gothn_{Tr_1}$ equality is
proved by taking inner products against $ae_1b$, where $a,b \in M$:
\begin{align*}
\Trace_1(ze_1 ae_1b)
&=\Trace_1(ze_1 E_N(a)b)
=\trace(T(ze_1)E_N(a)b) \\
&=\Trace_1(T(ze_1)e_1E_N(a)b)
=\Trace_1(T(ze_1)e_1ae_1b) .
\end{align*}
\end{proof}
\end{lemma}

\begin{prop}
\label{prop: basis!}
Let $N \subset M$ be a $\IIone$ subfactor of infinite index.  Then
there exists an orthogonal $M_N$-basis.
\begin{proof}
Let $T=T^{M_1}_M$.  Take an $\Ltwo(M)_N$-basis $\{ \xi_\beta \}$ with
$L(\xi_\beta)L(\xi_\beta)^*$ pairwise orthogonal projections (by 2.2
of Enock-Nest~\cite{EnockNest1996}).  In particular
\[
\Ltwo(M) = \bigoplus_\beta L(\xi_\beta)\Ltwo(N)
\]

For $\xi \in D(\Ltwo(M)_N)$, $T(L(\xi)L(\xi)^*)$ has finite trace
and so has spectral decomposition
\[
T(L(\xi)L(\xi)^*)=\int_0^{\infty} \lambda \dee q_{\lambda} .
\]
Let $f_1=q_1$, $f_{n+1}=q_{n+1}-q_n$.  Let $z_n=f_n L(\xi) e_1 \in
M_1$.  Then $z_n \in \gothn_T^*$ since
\[
T(z_nz_n^*)
=T(f_n L(\xi) e_1 L(\xi)^* f_n)
=T(f_n L(\xi) L(\xi)^* f_n)
=f_nT(L(\xi) L(\xi)^*) f_n
\]
which has norm less than $n$.  By the Pull-down Lemma
(Lemma~\ref{pull-down lemma}) $z_n=b_n e_1$ for some $b_n \in M$.
$L(\xi)e_1=\sum_n f_n L(\xi)e_1 =\sum b_ne_1$ and so
\[
L(\xi)\Ltwo(N) = \bigoplus_n b_n \Ltwo(N) .
\]

Applying this procedure to each $\xi_\beta$ we obtain a decomposition
of $\Ltwo(M)$ as
\[
\Ltwo(M) = \bigoplus_\alpha b_\alpha \Ltwo(N)
\]
where $b_\alpha \in M$.

For $b \in M$ let $p$ be the orthogonal projection onto $b\Ltwo(N)$.
Consider the spectral decomposition of $b e_1 b^*$
\[
b e_1 b^* = \int_0^{K} \lambda \dee q_\lambda
\]
Let $p_n=\CHI_{(K/(n+1),K/n]}(b e_1 b^*)=q_{K/n}-q_{K/(n+1)}$.  Note
that $\CHI_{\{0\}}(be_1b^*)$ is orthogonal projection onto
$\mathrm{Ker}(be_1b^*)=\mathrm{Ker}(e_1b^*)=(\mathrm{Range}(be_1))^{\perp}$
so that $\sum p_n=1-\CHI_{\{0\}}(be_1b^*)=p$.

Let
\[
z_n = \int_{K/(n+1)}^{K/n} \lambda^{-1/2} \dee q_\lambda ,
\]
so that $z_n be_1b^* z_n^*=p_n$.  Note that $z_n \leq ((n+1)/K)^{3/2}
b e_1 b^*$ so $z_n \in \gothm_T \subset \gothn_T^*$.  By the Pull-down
Lemma there exists $b_n\in M$ with $z_n be_1=b_n e_1$.  Then
\[
\sum_n b_ne_1b_n^* = \sum_n z_n be_1b^* z_n^* = \sum_n p_n = p .
\]

Applying this procedure to each $b_\alpha$ we obtain a collection of
elements $b_\iota$ in $M$ with
\[
\sum_\iota b_\iota e_1 b_\iota^* =1 
\]
and $b_\iota e_1 b_\iota^*$ pairwise orthogonal projections.
\end{proof}
\end{prop}

\begin{prop}
\label{prop: ON basis!}
Let $N \subset M$ be a $\IIone$ subfactor of infinite index.  Then
there exists an orthonormal $M_N$-basis.
\begin{proof}
Consider $N \subset M \subset M_1$, $T=T_M$.  By Prop~\ref{prop:
basis!} there exists a sequence of projections $p_i \in \gothm_T$ with
$\sum p_i =1$.  Observe that $\sum \Trace_1(p_i)=\Trace_1(1)=\infty$.
By adding finite sets of projections and taking subprojections,
operations under which $\gothm_T$ is closed, we may assume that
$\Trace_1(p_i)=1$ for all $i$.  Thus $p_i$ is equivalent in $M_1$ to
$e_1$, so there exist $v_i \in M_1$ with $v_iv_i^*=p_i$,
$v_i^*v_i=e_1$.  $v_i^*=v_i^*p_i \in \gothn_T$, so $v_i \in
\gothn_T^*$ and by the Pull-down Lemma there exists $b_i\in M$ with
$v_i=v_i e_1= b_i e_1$.  Then $b_i e_1 b_i^*= v_i v_i^* = p_i$ so that
$\sum_i b_i e_1 b_i^* =1$.  In addition
\[
e_1E_N(b_i^*b_j)
=e_1b_i^*b_je_1
=v_i^*v_j
=\delta_{i,j}e_1 ,
\]
so, applying $T$, $E_N(b_i^*b_j)=\delta_{i,j}1$.
\end{proof}
\end{prop}

\begin{remark}
In the proof of Prop~\ref{prop: ON basis!} the only fact we use from
Prop~\ref{prop: basis!} is the existence of a set of projections $\{
p_i \} \subset \gothm_T$ with $\sum_i p_i=1$.  Herman and Ocneanu
claim the existence of such a partition of unity for {\em any}
semi-finite factors $N \subset M$ on a separable Hilbert space (Prop 2
of~\cite{HermanOcneanu1989}), however no proof is available for this
result.
\end{remark}

\begin{lemma}
\label{lem: (M_r)_N basis}
Let $B=\{ b_i \}$ be an $M_N$-basis.  Then $B_r=\{ b_{j_1} \tensorlN
\tensorlN \cdots b_{j_{r+1}} \}$ is an $(M_r)_N$ basis.  If $B$ is
orthonormal then so is $B_r$.
\begin{proof}
\begin{equation}
\label{eq: (M_r)_N basis}
L\br{b_{j_1} \tensorN \cdots \tensorN \ b_{j_{r+1}}}
= L_{b_{j_1}} \cdots L_{b_{j_r}} 
   L\br{b_{j_{r+1}}}
\end{equation}
so that
\begin{multline*}
\sum L\br{b_{j_1} \tensorN \cdots \tensorN b_{j_{r+1}}}
     L\br{b_{j_1} \tensorN \cdots \tensorN b_{j_{r+1}}}^* \\
= \sum L_{b_{j_1}} \cdots L_{b_{j_r}} 
     L\br{b_{j_{r+1}}} L\br{b_{j_{r+1}}}^* 
     L_{b_{j_r}}^* \cdots L_{b_{j_1}}^*
= 1 .
\end{multline*}
If $B$ is orthonormal then
$L(b_j)^*L(b_i)=L_{b_j}^*L_{b_i}=\delta_{i,j}1$, so from equation 
(\ref{eq: (M_r)_N basis}) $L(b)^*L(\tilde{b})=\delta_{b,\tilde{b}}1$
for $b, \tilde{b}\in B_r$.
\end{proof}
\end{lemma}

\begin{lemma}
Let $\{ b_i \}$ be an orthonormal $M_N$-basis.  Let $k=2r-1$ or $2r$.
Define $p_{j_1,\ldots, j_r} \in M_k$ by
\[
p_{j_1,\ldots, j_r} 
= \begin{cases}b_{j_1} \tensorN \cdots \tensorN b_{j_r} \tensorN b_{j_r}^* 
          \tensorN \cdots \tensorN b_{j_1}^*,&k=2r-1\\
      b_{j_1} \tensorN \cdots \tensorN b_{j_r} \tensorN 1 \tensorN b_{j_r}^* 
          \tensorN \cdots \tensorN b_{j_1}^*,&k=2r .
    \end{cases}
\]
Then $p_{j_1,\ldots, j_r}$ are orthogonal projections with sum $1$.
\begin{proof}
For $k=2r-1$ note that by Lemma~\ref{lem: (M_r)_N basis} $B_{r-1}=\{
b_{j_1} \tensorlN \tensorlN \cdots b_{j_r} \}$ is an orthonormal
basis and hence $\br{b_{j_1} \tensorlN \tensorlN \cdots b_{j_r}}
\tensorN \br{b_{j_1} \tensorlN \tensorlN \cdots b_{j_r}}^*$ are
orthogonal projections with sum $1$.

For $k=2r$ simply note that
\[
b_{j_1} \tensorN \cdots \tensorN b_{j_r} \tensorN 1 \tensorN b_{j_r}^* 
          \tensorN \cdots \tensorN b_{j_1}^*
=b_{j_1} \tensorN \cdots \tensorN b_{j_r} \tensorN b_{j_r}^* 
          \tensorN \cdots \tensorN b_{j_1}^*
\]
by Corollary~\ref{cor: tensor 1 = prod e_odd}.
\end{proof}
\end{lemma}

\begin{prop}
\label{prop: M_l_M_k-basis}
Let $\{ b_i \}$ be an orthonormal $M_N$-basis.  Then for $-1\leq k
\leq l$ an orthogonal $(M_l)_{M_k}$-basis is given by:
\[
b^{j_1,\ldots,j_r}_{i_1,\ldots,i_{l-k}}=\begin{cases}
  b_{i_1} \tensorN \cdots \tensorN b_{i_{l-k}} \tensorN
  b_{j_1} \tensorN \cdots \tensorN b_{j_r} \tensorN b_{j_r}^* \tensorN
  \cdots \tensorN b_{j_1}^* 
   & k=2r-1 \\
  b_{i_1} \tensorN \cdots \tensorN b_{i_{l-k}} \tensorN
  b_{j_1} \tensorN \cdots \tensorN b_{j_r} \tensorN 1\tensorN b_{j_r}^*
  \tensorN \cdots \tensorN b_{j_1}^* 
   & k=2r .
  \end{cases}
\]
\begin{proof}

Note that
\[
L\br{ b^{j_1,\ldots,j_r}_{i_1,\ldots,i_{l-k}} }
= \mathcal{L}_{b_{i_1}} \cdots \mathcal{L}_{b_{i_{l-k}}} p_{j_1,\ldots, j_r} .
\]
Hence
\begin{align*}
\sum L\br{ b^{j_1,\ldots,j_r}_{i_1,\ldots,i_{l-k}} }
     L\br{ b^{j_1,\ldots,j_r}_{i_1,\ldots,i_{l-k}} }^* 
&= \sum \mathcal{L}_{b_{i_1}} \cdots \mathcal{L}_{b_{i_{l-k}}} 
        p_{j_1,\ldots, j_r}
        \mathcal{L}_{b_{i_{l-k}}}^* \cdots \mathcal{L}_{b_{i_1}}^* \\
&= 1 
\end{align*}
and the terms in the sum are orthogonal projections.
\end{proof}
\end{prop}

\begin{cor}
\label{cor: traces consistent}
Let $\Trace'$ be the canonical trace on $M_k' \cap \B(\Ltwo(M_l))$
(definition~\ref{def: cannonical trace}).  Then
$\Trace'\br{J_l\br{\spdot}^*J_l}=\Trace_{2l-k}$.
\begin{proof}
Let $\Trace=\Trace'\br{J_l\br{\spdot}^*J_l}$.  Note that by uniqueness
of the trace up to scaling, $\Trace$ is a multiple of $\Trace_{2l-k}$.
Let $f_i=1 \tensorlN \cdots \tensorlN 1$ ($i+1$ terms).
$\Trace_{2l-k}\br{f_{2l-k}}=1$ so we simply need to show that
$\Trace\br{f_{2l-k}}=1$.

Let $B$ be an orthonormal $M_N$-basis and $B_{l,k}$ the resulting
$(M_l)_{M_k}$-basis from Prop~\ref{prop: M_l_M_k-basis}.  $B_{l,k}^*$
is a $\leftidx{_{M_k}}{\br{M_l}}{}$-basis, so by Lemma~\ref{lem:Tr_Q'
by basis} $\Trace'=\sum_{b\in rB_{l,k}} \ip{\spdot
\hatt{b^*},\hatt{b^*}}$ and hence $\Trace=\sum_{b\in B_{l,k}}
\ip{\spdot \hatt{b},\hatt{b}}$.  Take $r$ such that
$k=2r-1+t$ ($t=0$ or $1$).  Fix $i_1,\ldots,i_{l-k},j_1,\ldots,j_r$.
Let $c=b_{i_1} \tensorlN \cdots \tensorlN b_{i_{l-k}} \tensorlN
b_{j_1} \tensorlN \cdots \tensorlN b_{j_r}\in M_{l-r-t}$ and let
$b=b^{j_1,\ldots,j_r}_{i_1,\ldots,i_{l-k}}$.  Then, from 
Prop~\ref{prop: tensorP},
\[
\pi^{2l-k}_l\br{f_{2l-k}}\hatt{b}
=\underbrace{1 \tensorN \cdots \tensorN 1}_{l-r+1} s \tensorN b_{j_r}^* 
 \tensorN \cdots \tensorN b_{j_1}^*
\]
where $s=E_N\br{1 \cdots E_N\br{1 E_N\br{1b_{i_1}} b_{i_2}} \cdots
 b_{j_r}} =T^{M_{l-r-t}}_N\br{f_{l-r-t}c}$.  Thus
\begin{align*}
&\ip{\pi^{2l-k}_l\br{f_{2l-k}}\hatt{b},\hatt{b}}
 =\ip{f_{l-r}s,c}
 =\ip{f_{l-r-t},cs^*}
 & t=0, \\
&\ip{\pi^{2l-k}_l\br{f_{2l-k}}\hatt{b},\hatt{b}}
 =\ip{f_{l-r}s,c\tensorN 1}
 =\ip{f_{l-r},cs^* \tensorN 1}
 =\ip{f_{l-r-t},cs^*}
 & t=1.
\end{align*}
Finally
\begin{align*}
\Trace\br{f_{2l-k}}
=\sum_{b\in B_{l,k}} \ip{\pi^{2l-k}_l\br{f_{2l-k}}\hatt{b},\hatt{b}}
&=\sum_{c \in B_{l-r}} \ip{f_{l-r-t},cT^{M_{l-r-t}}_N\br{c^*f_{l-r-t}}} \\
&=\ip{f_{l-r-t},f_{l-r-t}}
=\Trace_{l-r-t}\br{f_{l-r-t}}
=1.
\end{align*}
\end{proof}
\end{cor}


\subsubsection{Commutant Operator-Valued Weight}

The main result of this section is:

\begin{thm}
\label{thm:T_Q' by basis}
Let $P \subset Q$ be type two factors represented on a Hilbert space
$\h$.  Suppose that there exists a $Q_P$-basis $B=\{ b \}$.  Then the
n.f.s trace-preserving operator-valued weight $T_{Q'}:\plushat{\br{P'}}
\rightarrow \plushat{Q'}$ satisfies
\begin{align*}
&T_{Q'}(x) = \sum_b bxb^*
&x \in \plushat{P'}
\end{align*}
\end{thm}

\begin{prop}
\label{prop: Phi(x) in Q'}
Let $P \subset Q$ be type two factors represented on a Hilbert space
$\h$.  Let $B=\{ b \}$ be a $Q_P$-basis.  Then
$\Phi_B(x)=\sum_b bxb^* \in \plushat{(\BofH)}$ is affiliated with $Q'$
and hence in $\plushat{(Q')}$.  In addition $\Phi_B$ is independent of
$B$ and will hence be denoted simply $\Phi$.
\begin{proof}
For $\xi \in \h$ define an unbounded operator $R(\xi): \Ltwo(Q)
\rightarrow \h$ with domain $\gothn_{\Trace_Q}$ by
$R(\xi)\hatt{a}=a\xi$.  For $\eta \in D(\leftidx{_Q}{\h}{})$ we see
that $\eta \in D(R(\xi)^*)$:
\begin{align*}
<R(\xi)\hatt{a},\eta>
=<a\xi,\eta>
=<\xi,a^*\eta>
=<\xi,R(\eta)\hatt{a^*}>
&=<R(\eta)^*\xi,\hatt{a^*}> \\
&=<\hatt{a},J_Q R(\eta)^* \xi> .
\end{align*}
Hence $D(R(\xi)^*)$ is dense so that $R(\xi)$ is pre-closed (Theorem
2.7.8 (ii) of~\cite{KadisonRingroseI}).

Let $x \in P' \cap \BofH$ and let $A=x^{1/2}\overline{R(\xi)}$.
Define $m \in \plushat{P'}$ by 
\[
m(\omega_{\eta})=\begin{cases}
||(A^*A)^{1/2}\eta||^2 \phantom{1234} & \eta \in D((A^*A)^{1/2}) \\
\infty & \text{otherwise} 
\end{cases}
\]
(note that $A^*A$ is a positive, self-adjoint operator on $\Ltwo(Q)$
and is affiliated with $P'$ (a simple computation)).

\noindent Consider the polar decomposition $A=v(A^*A)^{1/2}$.  Using
Corollary~\ref{cor: Tr_Q' via basis extended}
\begin{align*}
\Trace_{P'\cap \B(\Ltwo(Q))}(m)
&= \sum_i m(\omega_{\hatt{b_i^*}}) \\
&= \sum_i ||(A^*A)^{1/2} \hatt{b_i^*}||^2 \\
&= \sum_i ||A \hatt{b_i^*}||^2 \\
&= \sum_i ||x^{1/2} R(\xi) \hatt{b_i^*}||^2 \\
&= \sum_i ||x^{1/2} b_i^* \xi||^2 \\
&= \sum_i < b_i x b_i^* \xi, \xi> \\
&= (\Phi_B(x))(\omega_\xi) .
\end{align*}
Since any element of $\br{\B(\h)}_*^+$ is a sum $\sum \omega_{\xi_k}$,
$\Phi_B(x)$ is thus independent of $B$.  In particular, for $u \in
\mathcal{U}(Q)$, $\Phi_B(x)=\Phi_{uB}(x)=u\Phi_B(x)u^*$, so that
$\Phi_B(x)$ is affiliated with $Q'$ and hence, by Prop~\ref{prop:
Haagerup ext pos part}, $\Phi_B(x) \in \plushat{(Q')}$.
\end{proof}
\end{prop}

\Proof {\it of Theorem~\ref{thm:T_Q' by basis}.}

Observe that since $\{ b^* \}$ is a
$\leftidx{_P}{\Ltwo(Q)}{}$-basis, $\{b^* \xi_i \}$ is a
$\leftidx{_P}{\h}{}$-basis: simply note that
$R(b^*\xi_i)=R(\xi_i)R(b^*)$ is bounded and
\[
\sum_{i,b} R(b^*\xi_i)R(b^*\xi_i)^*
 =\sum_{i,b} R(\xi_i)R(b^*)R(b^*)^*R(\xi_i)^*
 =\sum_i R(\xi_i)R(\xi_i)^*
 =1_{\h}
\]
Hence, by Lemma~\ref{lem:Tr_Q' by basis} we have, for $y \in (Q')_+$
\begin{align*}
\Trace_{P'}(y^{1/2}xy^{1/2})
&= \sum_i \sum_b \ip{y^{1/2}xy^{1/2} b^* \xi_i, b^*\xi_i} \\
&= \sum_i \sum_b \ip{b y^{1/2}xy^{1/2} b^* \xi_i, \xi_i} \\
&= \sum_i \sum_b \ip{y^{1/2} bxb^* y^{1/2} \xi_i, \xi_i} \\
&= \Trace_{Q'}(y^{1/2} \Phi(x) y^{1/2})
\end{align*}
so that $T_{Q'}(x)=\Phi(x)$.

For $x\in \plushat{\br{P'}}$ take $x_k \in \br{P'}_+$ with $\sum_k x_k=x$.
Then
\[
T_{Q'}(x)
=\sum_k T_{Q'}(x_k)
=\sum_k \sum_b bx_kb^*
=\sum_b \sum_k bx_kb^*
=\sum_b bxb^* .
\]
\proofend

\begin{cor}
For $P \subset Q \subset R$ and $x \in \plushat{P' \cap R}$,
$T_{Q'}(x) \in \plushat{\br{Q' \cap R}}$ does not depend on the Hilbert
space on which we represent $R$.
\end{cor}

\newpage

\section{Extremality and Rotations}
\label{sect: extremality and rotations}


\subsection{Extremality}
\label{sect: extremality}

Note that on $N' \cap M$ the traces coming from $N'$ and from $M$ are
not equal or even comparable, $\pi(N)'$ being a type $\IIinf$ factor
for any representation $\pi$ of $M$, while $M$ is a $\IIone$ factor.
This phenomenon continues up through the tower on all $\M{2i-1}{2i}$.

In the finite index case irreducibility ($N' \cap M=\mathbb{C}$)
implies extremality.  The example constructed by Izumi, Longo and Popa
in~\cite{IzumiLongoPopa1998} shows that this is not true for infinite
index inclusions.  They construct an irreducible $\IIone$ subfactor of
infinite index where the two traces on $N' \cap M_1$ are not even
comparable.

All of this suggests that we should be looking only at $N' \cap
M_{2i+1}$ when defining extremality for general $\IIone$ inclusions.
In this section we give the first definition of extremality in the
infinite index case and show that this definition has as many of the
desired properties as we can expect.

\begin{definition}[Commutant Traces]
Let $r=0,1$.  On $N' \cap M_{2i+r}$ define a trace
$\Trace'_{2i+r}$ by
\[
\Trace'_{2i+r}(x)=\Trace_{2i+1}(J_i x^*J_i) .
\]
In general, on $M_j'\cap M_{2i+r}$ define a trace
$\Trace'_{j,2i+r}$ by
\[
\Trace'_{j,2i+r}(x)=\Trace_{2i-j}(J_i x^*J_i) .
\]
\end{definition}

\begin{remark}
Note that the traces defined above are simply those coming from
representing $M_j$ on $\Ltwo(M_i)$ and using the fact that $J_i M_j'
J_i = M_{2i-j}$ by the multi-step basic construction.  Note that by
Corollary~\ref{cor: traces consistent} $\Trace'_{2i+r}$ is the
canonical trace on $N' \cap \B(\Ltwo(M_i))$.

Finally note that $\Trace'_{2i+1}\br{e_1\cdots e_{2i+1}}=1$.  This is
a consequence of the fact that $\Trace_{2i+1}\br{e_1\cdots
e_{2i+1}}=1$ and $J_i e_1\cdots e_{2i+1} J_i=e_1\cdots e_{2i+1}$:
\begin{align*}
&J_i \pi^{2i+1}_i\br{e_1\cdots e_{2i+1}} J_i 
 \hatr{a_1 \tensorN \cdots \tensorN a_{i+1}} \\
&= J_i \pi^{2i+1}_i\br{1 \tensorN \cdots \tensorN 1}
 \hatr{a_{i+1}^* \tensorN \cdots \tensorN a_1^*} \\
&= J_i \hatr{1 \tensorN \cdots \tensorN 1 
 E_N\br{\cdots E_N\br{E_N\br{a_{i+1}^*}a_i^*}\cdots a_1^*}} \\
&= J_i \hatr{1 \tensorN \cdots \tensorN 1 
 E_N\br{a_{i+1}^*}E_N\br{a_i^*}\cdots E_N\br{a_1^*}} \\
&= \hatr{1 \tensorN \cdots \tensorN 1 
 E_N\br{a_1}\cdots E_N\br{a_{1+1}}} \\
&=\pi^{2i+1}_i\br{e_1\cdots e_{2i+1}}
 \hatr{a_1 \tensorN \cdots \tensorN a_{i+1}}
\end{align*}
Hence we could define $\Trace'_{2i+r}$ as the restriction of the
unique trace on $N'\cap\mathcal{B}(\Ltwo(M_i))$ scaled so that
$\Trace'_{2i+1}(e_1e_3\cdots e_{2i+1})=1$.
\end{remark}

\begin{definition}[Extremality]
Let $N \subset M$ be an inclusion of $\IIone$ factors.  $N \subset M$
is {\em extremal} if $\Trace'_1=\Trace_1$ on $N' \cap M_1$.  $N
\subset M$ is {\em approximately extremal} if $\Trace'_1$ and
$\Trace_1$ are equivalent on $N' \cap M_1$ (i.e. there exists $C>0$
such that $C^{-1} \Trace_1 \leq \Trace'_1 \leq C \Trace_1$ on $N' \cap
M_1$).
\end{definition}

\begin{remark}
We will abuse notation a little by writing
$\Trace'_{2j+1}=\Trace_{2j+1}$ when $\Trace'_{2j+1}=\Trace_{2j+1}$ on
$N' \cap M_{2j+1}$ and similarly $\Trace'_{2j+1}\sim \Trace_{2j+1}$
when $\Trace'_{2j+1}$ and $\Trace_{2j+1}$ are equivalent on $N' \cap
M_{2j+1}$.
\end{remark}

\begin{prop}
\label{prop: equiv defs of extremal}
\begin{description}
\item{(i)} If $N \subset M$ is extremal then
$\Trace'_{2i+1}=\Trace_{2i+1}$ for all $i \geq 0$.
\item{(ii)} If $N \subset M$ is approximately extremal then
$\Trace'_{2i+1}\sim \Trace_{2i+1}$ for all $i \geq 0$.
\item{(iii)} If $\Trace'_{2i+1}=\Trace_{2i+1}$ for some $i \geq 0$
then $N \subset M$ is extremal.
\item{(iv)} If $\Trace'_{2i+1}\sim \Trace_{2i+1}$ for some $i \geq 0$
then $N \subset M$ is approximately extremal.
\end{description}
\end{prop}

\begin{lemma}
\label{lem: Trace' consistent}
For $z\in \plushat{\br{N' \cap M_j}}$,
$T^{N'}_{M_j'}(z)=\Trace'_j(z)1$.
\begin{proof}
Assume $z\in\br{N' \cap M_j}_+$ (in general take $z_k\nearrow z$,
$z_k\in\br{N' \cap M_j}_+$).  $T^{N'}_{M_j'}(z)\in \plushat{M_j' \cap
M_j} = [0,\infty]1$, so $T^{N'}_{M_j'}(z)=\lambda 1$ for some
$\lambda\in[0,\infty]$.  Let $i$ be an integer such that $j=2i$ or
$2i-1$.  Let $\Trace'$ be the canonical trace on $M_j' \cap
\B(\Ltwo(M_i))$.  By Corollary~\ref{cor: traces consistent}
$\Trace'\br{J_i\br{\spdot}^*J_i}=\Trace_{2i-j}=\trace$ and hence
$\Trace'(1)=1$.  Thus
\[
\lambda
=\Trace'\br{T^{N'}_{M_j'}(z)}
=\Trace_{N' \cap \B(\Ltwo(M_i))}(z)
=\Trace'_j(z) .
\]
\end{proof}
\end{lemma}

\begin{lemma}
\label{lem: stepping up and down traces}
\begin{description}
\item
\item{(i)} $\Trace_{2i-1}(x)=\Trace_{2i+1}(xe_{2i+1})$ for $x\in
\plushat{\br{M_{2i-1}}}$.
\item{(ii)} $\Trace_{2i}(x)=\Trace_{2i+1}(e_{2i+1}xe_{2i+1})$ for
$x\in \plushat{\br{M_{2i}}}$.
\item{(iii)} $\Trace'_{2i-1}(x)=\Trace'_{2i+1}\br{xe_{2i+1}}$ for
$x\in \plushat{\br{N' \cap M_{2i-1}}}$.
\end{description}
\begin{proof}
\begin{description}
\item
\item{(i)} $\Trace_{2i+1}(\spdot e_{2i+1})$ is tracial on $M_{2i}$ and
hence a multiple of $\Trace_{2i-1}$.  The multiple is $1$ because
$\Trace_{2i-1}(e_1\cdots e_{2i-1})=1 =\Trace_{2i+1}(e_1\cdots
e_{2i-1}e_{2i+1})$.
\item{(ii)} $e_{2i+1}xe_{2i+1}=E_{M_{2i-1}}(x)e_{2i+1}$ so by (i)
\[
\Trace_{2i+1}\br{e_{2i+1}xe_{2i+1}}
=\Trace_{2i-1}\br{E_{M_{2i-1}}(x)}
=\Trace_{2i}\br{x} .
\]
\item{(iii)} First note that
$T^{M_{2i-1}'}_{M_{2i+1}'}\br{ e_{2i+1}}=1$: take a basis $\{ b_j \}$
for $M_{2i}$ over $M_{2i-1}$.  Then $\{ b_j e_{2i+1} b_k \}$ is a
basis for $M_{2i+1}$ over $M_{2i-1}$ and hence
\[
T^{M_{2i-1}'}_{M_{2i+1}'}\br{ e_{2i+1}}
= \sum_{j,k} b_j e_{2i+1} b_k e_{2i+1} b_k^* e_{2i+1} b_j^* 
= \sum_j b_j e_{2i+1} b_j^* = 1 .
\]

Let $x \in (N' \cap M_{2i-1})_+$.  Using Lemma~\ref{lem: Trace'
consistent}
\begin{align*}
\Trace'_{2i+1}\br{xe_{2i+1}} 1
&= T^{N'}_{M_{2i+1}'}\br{ xe_{2i+1} }  \\
&= T^{M_{2i-1}'}_{M_{2i+1}'}\br{ T^{N'}_{M_{2k-1}'}
     \br{ xe_{2i+1}}} \\
&= T^{M_{2i-1}'}_{M_{2i+1}'}\br{ T^{N'}_{M_{2i-1}'}(x) e_{2i+1}} \\
&= T^{M_{2i-1}'}_{M_{2i+1}'}\br{ \Trace'_{2i-1}(x) e_{2i+1}} \\
&= \Trace'_{2i-1}(x) T^{M_{2i-1}'}_{M_{2i+1}'}\br{ e_{2i+1}} \\
&= \Trace'_{2i-1}(x) 1 .
\end{align*}
\end{description}
\end{proof}
\end{lemma}

\Proof {\it of Prop~\ref{prop: equiv defs of extremal}.}

In the finite index case the proof is easily accomplished using planar
algebra machinery.  $\Trace_{2i+1}$ is just closing up a $2i+2$-box on
the right, $\Trace'_{2i+1}$ is closing up a $2i+2$-box on the left.
Extremality means that a $2$-box closed on the left is equal to the
same box closed on the right.  For a $2i+2$-box just move the strings
from right to left two at a time.  This is the approach that we will
take here, although of course we cannot use the planar algebra
machinery and must proceed algebraically.

\begin{description}
\item{(i)} Suppose $\Trace_{2i-1}=\Trace'_{2i-1}$.  We will show that
$\Trace_{2i+1}=\Trace'_{2i+1}$.  Let $z \in \br{N' \cap M_{2i+1}}_+$,
then,
\begin{align*}
&\Trace_{2i+1}(z) 1 \\
&=\Trace_{2i-1} \br{T_{M_{2i-1}}(z)} 1 \\
&=\Trace'_{2i-1} \br{T_{M_{2i-1}}(z)} 1 \\
&=T^{N'}_{M_{2i-1}'} \br{T_{M_{2i-1}}(z)}
  &\text{by Lemma~\ref{lem: Trace' consistent}} \\
&=\sum_b b T_{M_{2i-1}}(z) b^*
  &\text{where $\{ b \}$ is a basis for $M_{2i-1}$ over $N$} \\
&=\sum_b T_{M_{2i-1}}(bzb^*) \\
&=T_{M_{2i-1}}\br{T^{N'}_{M_{2i-1}'}(z)} 
  &\text{now represent everything on }\Ltwo(M_i) \\
&=j_i T^{N'}_{M_1'}\br{j_i\br{T^{N'}_{M_{2i-1}'}(z)}}
  &\text{because }j_i\br{M'_{2i-1}\cap M_{2i+1}}=M_1 \cap N' \\
&=\Trace'_1\br{j_i\br{T^{N'}_{M_{2i-1}'}(z)}}1 
  &\text{by Lemma~\ref{lem: Trace' consistent}} \\
&=\Trace_1\br{j_i\br{T^{N'}_{M_{2i-1}'}(z)}}1 \\
&=\Trace_{M_{2i-1}'\cap \B(\Ltwo(M_i))}\br{T^{N'}_{M_{2i-1}'}(z)}1 \\
&=\Trace_{N'\cap \B(\Ltwo(M_i))}(z)1 \\
&=\Trace'_{2i+1}(z)1
\end{align*}

\item{(iii)} Suppose $\Trace_{2i+1}=\Trace'_{2i+1}$, $i \geq 1$.
Using Lemma~\ref{lem: stepping up and down traces}, for $z \in \br{N'
\cap M_{2i-1}}_+$,
\[
\Trace_{2i-1}(z)
=\Trace_{2i+1}(ze_{2i+1})
=\Trace'_{2i+1}(ze_{2i+1})
=\Trace'_{2i-1}(z) .
\]

\item{(ii) and (iv)} Similar to (i) and (iii) with inequalities and
constants.

\end{description}
\proofend


\subsection{$N$-central vectors}
\label{sect: N-central vectors}

Here we will show that the set of $N$-central vectors in $\Ltwo(M_k)$
is $\overline{N' \cap \gothn_{\Trace_k}}$ (closure in $\Ltwo(M_k)$).
The proof is essentially an application of ideas from
Popa~\cite{Popa1981}.

\begin{definition}
$\xi\in\Ltwo(M_k)$ is an {\em $N$-central vector} if $n\xi=\xi n$ for
all $n\in N$.  The set of $N$-central vectors in $\Ltwo(M_k)$ is
denoted $N' \cap \Ltwo(M_k)$.
\end{definition}

\begin{lemma}
\label{2-norm lsc}
Let $M$ be a semifinite von Neumann algebras with n.f.s trace $\Trace$.
The 2-norm $||x||_2=[\Trace(x^*x)]^{1/2}$ is lower
semi-continuous with respect to the weak operator topology.
\begin{proof}
Suppose $\{ y_{\alpha} \}$ is a net in $M$ converging weakly to $y$.
Take a set of orthogonal projections $\{ p_k \} \subset
\gothn_{\Trace}$ with $\sum p_k =1$.  Then for $z\in \gothn_{\Trace}$
with $|| z ||_2\leq 1$ one has $|< y_{\alpha}
\hatt{p_k},\hatt{z}>|\leq || y_{\alpha} p_k||_2$ and hence
\[
|<y\hatt{p_k},\hatt{z}>|=\lim |< y_{\alpha} \hatt{p_k},\hatt{z}>|
	\leq \liminf_{\alpha} ||y_{\alpha} p_k ||_2
\]
so that
\[
||y p_k||_2
 =\sup\{|<y\hatt{p_k},\hatt{z}>| : z\in {\gothn_{\Trace}}, ||z||_2\leq 1 \}
 \leq \liminf ||y_{\alpha} p_k||_2 .	
\]
Finally, using Fatou's lemma for the second inequality,
\[
||y||_2^2=\sum_k ||y p_k||_2^2
	\leq \sum_k \liminf_{\alpha} ||y_{\alpha} p_k||_2^2
	\leq \liminf_{\alpha} \sum_k ||y_{\alpha} p_k||_2^2
	= \liminf_{\alpha} ||y_{\alpha}||_2^2 .
\]
\end{proof}
\end{lemma}

\begin{remark}
The fact that $||y||_2^2=\sum_k ||y p_k||_2^2$ even if $y\notin
\gothn_{\Trace}$ is established as follows.  Let
$q_N=\sum_{k=1}^N{p_k}$ and note that $q_N\nearrow 1$.  Hence
$yq_Ny^*\nearrow yy^*$ and so $Tr(q_Ny^*yq_N)= Tr(yq_N y^*)\nearrow
Tr(yy^*)=||y||_2^2$ (even though $q_Ny^*yq_N$ may {\em not} be
increasing).  Finally, $y^*yq_N\in \gothn_{\Trace}$ and
$Tr(q_Ny^*yq_N)=Tr(y^*yq_N)=\sum_{k=1}^N Tr(y^*yp_k)=\sum_{k=1}^N
||yp_k||_2^2$.
\end{remark}

The next lemma follows exactly the line of proof in
Popa~\cite{Popa1981} 2.3.

\begin{prop}
\label{convex averaging}
Let $N$ be a von Neumann subalgebra of a semifinite von Neumann
algebra $M$ equipped with a n.f.s. trace $\Trace$.  Let ${\mathcal
U}(N)$ denote the unitary group of $N$.  Let $x\in \gothn_{\Trace}$
and let $K_0=\{ \sum_{i=1}^n \lambda_i v_i x v_i^* : \lambda_i
\in[0,\infty], \sum \lambda_i=1, v_i\in {\mathcal U}(N) \}$.  Then
$N'\cap {\gothn_{\Trace}}\cap \overline{K_0}\neq \emptyset$, where
$\overline{K_0}$ denotes the weak closure of $K_0$.
\begin{proof}
Let $K=\overline{K_0}$.  Observe that for $v_i\in {\mathcal U}(N)$ and
$\lambda_i\in[0,\infty]$ such that $\sum\lambda_i=1$, we have
$||\sum_{i=1}^n \lambda_i v_i x v_i^* || \leq \sum_{i=1}^n \lambda_i
||v_ix v_i^* || = \sum_{i=1}^n \lambda_i ||x|| = ||x||$.  Hence
$K\subset B(0, ||x||)$ ($||y||=\sup\{ |<y\xi,\eta >|: \xi, \eta \in
({\mathcal H})_1\}$, so if $y_{\alpha}\rightarrow y$ and
$||y_{\alpha}||\leq r$ for all $\alpha$ then $||y||\leq r$).

Since $K$ is bounded and weakly closed it is weakly compact.  Let
$\omega=\inf_K {||y||_2}$ (note that $\omega<\infty$ since $x\in K$).
Take $y_n\in K$ such that $||y_n||_2\rightarrow \omega$.  Since $K$ is
weakly compact there exists a weakly convergent subsequence.  Hence we
may assume that $\{ y_n \}$ is weakly convergent.  Let $y=\lim y_n$.

Since $||\cdot ||_2$ is lower semi-continuous for the weak operator
topology we have $\omega\leq ||y||_2 \leq \liminf ||y_n||_2 = \omega$
and hence $||y||_2=\omega$.

Since $K$ is convex $y$ is the unique element of $K$ such that
$||y||_2=\omega$, for if $||y||_2=||z||_2=\omega$ then
$\frac{1}{2}(y+z)\in K$ so
\[
\omega^2\leq \frac{1}{4} ||y+z||^2=\frac{1}{2}\omega^2+\frac{1}{2}{\mathrm Re}
<y,z>
	\leq \frac{1}{2}\omega^2+\frac{1}{2}||y||_2||z||_2
	=\omega^2
\]
with equality iff $y=z$.

Finally, note that for all $v\in {\mathcal U}(N)$ we have $vyv^*\in K$
and $||vyv^*||_2=||y||_2$ so that $vyv^*=y$.  Hence $y\in N'\cap M$,
and since $||y||_2=\omega<\infty$, $y\in \gothn_{\Trace}$.
\end{proof}
\end{prop}

\begin{cor}
\label{cor: almost N'}
Let $N$ be a von Neumann subalgebra of a semifinite von Neumann
algebra $M$ equipped with a n.f.s. trace $\Trace$ and suppose that
$x\in M$ and $||uxu^*-x||_2\leq\delta$ for all $u\in {\mathcal U}(N)$.
Then there exists $y\in N'\cap M$ with $||x-y||_2\leq \delta$.
\begin{proof}
Let $v_i\in {\mathcal U}(N)$ and let $\lambda_i \in[0,1]$ such that
$\sum\lambda_i=1$.  Then
\[
||\sum \lambda_i v_i x v_i^*-x||_2
	=||\sum \lambda_i (v_i x v_i^*-x)||_2
	\leq \sum \lambda_i ||v_i x v_i^*-x||_2
	\leq \sum \lambda_i \delta
	= \delta
\]
So $||x-z||_2\leq\delta$ for all $z\in K_0$.  By the previous lemma
there exists $y\in N'\cap M$ and $\{y_n\}$ in $K_0$ with
$y_n\rightarrow y$ weakly.  Thus $x-y_n\rightarrow x-y$ weakly and by
the lower semi-continuity of $||\cdot ||_2$
\[
||x-y||_2\leq \liminf ||x-y_n||_2 \leq \delta .
\]
\end{proof}
\end{cor}

\begin{thm}
\label{L2=rel comm + commutators}
Let $N$ be a von Neumann subalgebra of a semifinite von Neumann
algebra $M$ equipped with a n.f.s. trace $\Trace$.  Then
\begin{description}
\item{(i)} $N' \cap \Ltwo(M)=\overline{N' \cap \gothn_{\Trace}}$.
\item{(ii)} $\br{N' \cap \Ltwo(M)}^\perp$ is the (span of) the
commutators in $\Ltwo(M)$.  Hence
\[
{\mathrm L}^2(M)=\overline{N'\cap {\gothn_{\Trace}}}\oplus
\overline{[\gothn_{\Trace},N]} .
\]
\end{description}
\begin{proof}
\begin{description}
\item{(i)} Clearly $N'\cap \gothn_{\Trace} \subseteq N'\cap\Ltwo(M)$.
Let $\xi\in N'\cap\Ltwo(M)$.  Take $\{ x_m \}$ in $\gothn_{\Trace}$
with $\hatt{x_m}\rightarrow \xi$ in $||\cdot||_2$.  Then for all $u\in
{\mathcal U}(N)$
\begin{multline*}
||ux_mu^*-x_m||_2
=||ux_m-x_mu||_2
=||u(x_m-\xi)-(x_m-\xi)u||_2 \\
 \leq ||u(x_m-\xi)||_2 +||(x_m-\xi)u||_2
 =2||x_m-\xi||_2
\end{multline*}

By Corollary~\ref{cor: almost N'} there are $y_m\in N'\cap
\gothn_{\Trace}$ with $||x_m-y_m||\leq 2||x_m-\xi||_2$ and thus
$||y_m-\xi||_2\leq 4 ||x_m-\xi||_2 \rightarrow 0$.  Thus $\xi\in
\overline{N'\cap \mathcal N}$, so ${N'\cap\Ltwo(M)}=\overline{N'\cap
\gothn_{\Trace}}$.
\item{(ii)} Simply note that 
$<\xi,n\hatt{m}-\hatt{m}n>=<n^* \xi-\xi n^*,\hatt{m}>$, so $\xi \in
N'\cap\Ltwo(M)$ iff $\xi \perp [\mathfrak{N},N]$.
\end{description}
\end{proof}
\end{thm}

\begin{cor}
\label{L2=rel comm + commutators, finite index}
If $N\subset M$ is a finite index $\IIone$ subfactor then the
$N$-central vectors in $\Ltwo(M_k)$ are precisely $N' \cap M_k$.
\begin{proof}
$\gothn_{\Trace}=M_k$ and $N' \cap M_k$ is finite dimensional so
$\overline{N' \cap M_k}=N' \cap M_k$.
\end{proof}
\end{cor}

\begin{cor}
\label{pcentral via strong limits of convex average}
Let $x\in \gothn_{\Trace}$ and let $K_0=\{ \sum_{i=1}^n \lambda_i v_i
x v_i^* : \lambda_i \in[0,\infty], \sum \lambda_i=1, v_i\in {\mathcal
U}(N) \}$.  Then $P_{N'\cap\Ltwo(M)}(x)$ is in the strong closure of
$K_0$.
\begin{proof}
We will show that the element $y \in N' \cap \gothn_{\Trace} \cap
\overline{K_0}$ whose existence is guaranteed by Prop~\ref{convex
averaging} is $P_{N'\cap\Ltwo(M)}(x)$.  Since $K_0$ is convex its weak
and strong closures coincide.

Note that for all $z \in N' \cap \gothn_{\Trace}$,
\[
\sum_{i=1}^n \ip{\lambda_i v_i x v_i^*, z}
=\sum_{i=1}^n \lambda_i \ip{x, v_i^* z v_i}
=\sum_{i=1}^n \lambda_i \ip{x,z}
=\ip{x,z} .
\]
Take $y_{\alpha}$ in $K_0$ such that $y_{\alpha} \rightarrow y$
weakly.  Then for all $z_1, z_2 \in N' \cap \gothn_{\Trace}$
\[
\ip{y, z_1 z_2^*}
=\ip{y \hatt{z_2},\hatt{z_1}}
=\lim_{\alpha} \ip{y_{\alpha} \hatt{z_2},\hatt{z_1}}
=\lim_{\alpha} \ip{y_{\alpha}, z_1 z_2^*}
=\ip{x,z_1 z_2^*} .
\]
Since $\{ z_1 z_2^* : z_i \in N' \cap \gothn_{\Trace}\}$ is dense in
$N'\cap\Ltwo(M)$, $y=P_{N'\cap\Ltwo(M)}(x)$.

The density of $\{ z_1 z_2^* : z_i \in N' \cap \gothn_{\Trace}\}$ in
$\overline{N' \cap \gothn_{\Trace}}$ can be seen as follows.  By
Takesaki~\cite{Takesaki1979}, Chapter V, Lemma 2.13, there exists a
maximal central projection $p$ of $N' \cap M$ such that $\Trace$ is
semifinite on $p(N' \cap M)$ and $\Trace=\infty$ on $[(1-p)(N' \cap
M)]_+\backslash \{ 0 \}$.  Let $A=p(N' \cap M)$.  Then
$\gothn_{\Trace_A}=N' \cap \gothn_{\Trace}$ and $\Ltwo(A)=\overline{N'
\cap \gothn_{\Trace}}$.  By Tomita-Takesaki theory applied to
$(A,\Trace)$ we see that $\{ z_1 z_2^* : z_i \in \gothn_{\Trace_A}\}$
is dense in $\Ltwo(A)$.
\end{proof}
\end{cor}


\subsection{Rotations}
\label{sect: rotations}

\begin{definition}
\begin{description}
\item
\item{1.} Let $B$ be a $M_N$-basis.  If $\sum_{b \in B} R_{b^*} \br{
L_{b^*}}^* \hatt{x}$ converges for all $x \in N' \cap
\gothn_{\Trace_k}$ and extends to a bounded operator from $N' \cap
\Ltwo(M_k)$ to $\Ltwo(M_k)$ we define $\rho^B_k$ to be this extension.
\item{2.} Let $P_c=P_{N'\cap\Ltwo(M_k)}$.  If
\[
\br{x_1 \tensorN x_2 \tensorN \cdots \tensorN x_{k+1}}
\mapsto P_c\br{x_2 \tensorN x_3 \tensorN \cdots 
  \tensorN x_{k+1} \tensorN x_1}
\]
extends to a bounded operator on all $\Ltwo(M_k)$ we define
$\wtilde{\rho}_k$ to be this extension.
\end{description}
\end{definition}

\begin{thm}
\label{thm: rho, tilde(rho)}
If $\rho^B_k$ and $\rho^{\overline{B}}_k$ both exist then
$\rho^B_k=\rho^{\overline{B}}_k$.

If $\wtilde{\rho}_k$ exists then for any orthonormal basis $B$,
$\rho^B_k$ exists.  If $\rho^B_k$ exists for some $B$ then
$\wtilde{\rho}_k$ exists.  In this case both operators map $N' \cap
\Ltwo(M_k)$ onto $N' \cap \Ltwo(M_k)$.  Restricted to $N' \cap
\Ltwo(M_k)$ both are periodic (with period $k+1$) and
$\br{\rho^B_k}^*=\br{\wtilde{\rho}_k}^{-1}$.
\begin{proof}
\begin{description}
\item{1.} If $\rho^B_k$ exists then, as in the finite index case of
Section~\ref{sect: finite index IIone}, for $x \in N' \cap
\gothn_{\Trace_k}$ and $y=y_1 \tensorN \cdots \tensorN y_{k+1}$,
\begin{align}
\ip{\rho^B_k(x),y}
&=\sum_b \ip{x, L_{b^*} \br{ R_{b^*}}^* y} 
  \nonumber \\
&=\sum_b \ip{x, b^* \tensorN y_1 \tensorN \cdots \tensorN y_k E_N(y_{k+1}b)} 
  \nonumber \\
&=\sum_b \ip{xE_N(y_{k+1}b)^*,b^* \tensorN y_1 \tensorN \cdots \tensorN y_k} 
  \nonumber \\
&=\sum_b \ip{E_N(y_{k+1}b)^*x,b^* \tensorN y_1 \tensorN \cdots \tensorN y_k} 
  \nonumber \\
&=\sum_b \ip{x,E_N(y_{k+1}b)b^* \tensorN y_1 \tensorN \cdots \tensorN y_k} 
  \nonumber \\
&=\ip{x,y_{k+1} \tensorN y_1 \tensorN \cdots \tensorN y_k} .
  \label{eq: rho inner prod}
\end{align}
Hence $\rho^B_k(x)$ is independent of the basis used.  Note that
for $u\in \mathcal{U}(N)$
\begin{align*}
\sum_{b \in B} R_{(ub)^*} \br{ L_{(ub)^*}}^* \hatt{x}
&= \sum_{b \in B} R_{b^*u^*} \br{ L_{b^*u^*}}^* \hatt{x} \\
&= \sum_{b \in B} R_{b^*u^*} u \br{ L_{b^*}}^* \hatt{x} \\
&= \sum_{b \in B} u R_{b^*u^*} \br{ L_{b^*}}^* \hatt{x} \\
&= \sum_{b \in B} u \br{R_{b^*} \br{ L_{b^*}}^* \hatt{x} }\cdot u^* \\
&= u (\rho_k(\hatt{x}))\cdot u^* .
\end{align*}
So if $\rho^B_k$ exists then $\rho^{uB}_k$ exists and
$u\rho^B_k(\hatt{x})u^*=\rho^{uB}_k(\hatt{x})=\rho^B_k(\hatt{x})$.
Thus $\rho^B_k(\hatt{x}) \in N' \cap \Ltwo(M_k)$.  From~(\ref{eq: rho
inner prod}), $\br{\rho^B_k}^{k+1}=\id$.

In addition $P_c\br{y_2 \tensorlN \cdots \tensorlN y_{k+1} \tensorlN
y_1 }=\br{\br{\rho^B_k}^k P_c}^*\hatt{y}$, so that $\wtilde{\rho}_k$
exists.  On $N' \cap \Ltwo(M_k)$, $\wtilde{\rho}_k
=\br{\br{\rho^B_k}^{-1}}^*$.  Hence $\br{\wtilde{\rho}_k}^{k+1}=\id$ and
$\br{\rho^B_k}^*=\wtilde{\rho}_k^{-1}$.

\item{2.} If $\wtilde{\rho}_k$ exists then let $\sigma_k=\br{ P_c J_k
\wtilde{\rho}_k J_k P_c }^*$.  Take an orthonormal basis $\{ b \}
= \{ b_i \}$.  Then, reversing the argument in~(\ref{eq: rho inner
prod}), for $\xi \in N' \cap \Ltwo(M_k)$,
\begin{equation}
\label{eq: weak conv?}
\ip{\sigma_k(\xi),y}
=\sum_b \ip{R_{b^*} \br{L_{b^*}}^*\xi, y} ,
\end{equation}
Let $\eta=\sigma_k(\xi)$.  Note that $\eta_i \defeq R_{b_i^*}\br{
L_{b_i^*} }^* \xi$ are pairwise orthogonal.  For $y=y_1 \tensorlN
\cdots y_k \tensorlN b_i^*$ the sum in (\ref{eq: weak conv?}) only has
one term, so the equality extends by continuity to all
$R_{b_i^*}\zeta$ ($\zeta \in \Ltwo(M_{k-1})$) and in particular to
$\eta_i$.  Thus $\ip{\eta,\eta_i} =\ip{\eta_i, \eta_i}$ and so
\[
\sum_{i=1}^K ||\eta_i||^2
=\ip{\eta, \sum_{i=1}^K \eta_i} \\
\leq ||\eta|| \left|\left|\sum_{i=1}^K \eta_i\right|\right| \\
= ||\eta|| \br{ \sum_{i=1}^K ||\eta_i||^2 }^{1/2}
\]
which yields
\[
\br{\sum_{i=1}^K ||\eta_i||^2}^{1/2} \leq ||\eta|| .
\]
Hence $\sum_{i=1}^K ||\eta_i||^2$ is bounded, so $\sum_{i=1}^K \eta_i$
converges and by~(\ref{eq: weak conv?}) converges to $\eta$.  Hence
$\rho^B_k$ exists.
\end{description}
\end{proof}
\end{thm}

\begin{remark}
In Corollary~\ref{cor: rotn for one basis implies rotn for all} we
will see that if $\rho^B_k$ exists for some basis $B$ then it exists
for all bases and we will then drop the reference to $B$ and use the
notation $\rho_k$.
\end{remark}


\subsection{Rotations and extremality}
\label{sect: rotations and extremality}

We prove that approximate extremality is necessary and sufficient for
the rotations to exist.

\begin{prop}
\label{prop: approx extremal implies rotation}
\begin{description}
\item{}
\item{(i)} If $N \subset M$ is approximately extremal then $\rho^B_k$
exists for all bases $B$ and for all $k \geq 0$.  In this case we will
use the notation $\rho_k$.
\item{(ii)} If $N \subset M$ is extremal then in addition to (i)
$\rho_k$ is a unitary operator on $N' \cap \Ltwo(M_k)$.  Hence
$\rho_k=\wtilde{\rho}_k$.
\end{description}
\end{prop}

\begin{lemma}
\label{lem: stuff about e_1}
\begin{description}
\item
\item{(i)} $L_1 \br{L_1}^*=e_1$.
\item{(ii)} $J_ie_1 J_i=e_{2i+1}$.
\item{(iii)} For $x\in \br{M' \cap M_k}_+$,
$e_1xe_1=T^{M'}_{M_1'}(x)e_1$.
\end{description}
\begin{proof}
\begin{description}
\item{(i)} With the usual notations established in the proof of
Prop~\ref{prop: tensorP},
\begin{align}
e_1 \hatr{a_1 \tensorN \cdots \tensorN a_{i+1}}
&= e_1 \hatr{A_1 \tensorM \cdots \tensorM \overline{A}_i} \nonumber \\
&= \hatr{\br{e_1A_1} \tensorM \cdots \tensorM \overline{A}_i} \nonumber \\
&= \hatr{\br{E_N(a_1)e_1w_1} \tensorM A_2 \tensorM \cdots 
   \tensorM \overline{A}_i} \nonumber \\
&= \hatr{E_N(a_1) \tensorN a_2 \tensorN \cdots \tensorN a_{i+1}} 
   \label{eq: e_1 action} \\
&= L_1 \hatr{E_N(a_1)a_2 \tensorN a_3 \tensorN \cdots \tensorN a_{i+1}} 
   \nonumber \\
&= L_1 \br{L_1^*L_{a_1}}\hatr{a_2 \tensorN a_3 \tensorN \cdots \tensorN 
   a_{i+1}} \nonumber \\
&= L_1L_1^* \hatr{a_1 \tensorN a_2 \tensorN \cdots \tensorN a_{i+1}} \nonumber .
\end{align}
\item{(ii)} From (\ref{eq: e_1 action}), $J_ie_1J_i\hatr{a_1 \tensorlN
\cdots \tensorlN a_{i+1}}=\hatr{a_1 \tensorlN \cdots \tensorlN a_i
\tensorlN E_N(a_{i+1})}$.  Let $\iota:\Ltwo(M_{2i-1})\rightarrow
\Ltwo(M_{2i})$ be the inclusion map.  From the definition of
$\pi^{2i+1}_i$,
\begin{align*}
&\sqbr{\pi^{2i+1}_i\br{e_{2i+1}}\hatr{a_1 \tensorN \cdots \tensorN a_{i+1}}} 
  \tensorN \hatr{b_1 \tensorN \cdots \tensorN b_i} \\
&=e_{2i+1} \hatr{a_1 \tensorN \cdots \tensorN a_{i+1}
  \tensorN b_1 \tensorN \cdots \tensorN b_i} \\
&=\iota \hatr{E_{M_{2i-1}}\br{a_1 \tensorN \cdots \tensorN a_{i+1}
  \tensorN b_1 \tensorN \cdots \tensorN b_i}} \\
&=\iota \hatr{a_1 \tensorN \cdots \tensorN a_iE_N(a_{i+1})
  \tensorN b_1 \tensorN \cdots \tensorN b_i} \\
&=\hatr{a_1 \tensorN \cdots \tensorN a_iE_N(a_{i+1}) \tensorN 1
  \tensorN b_1 \tensorN \cdots \tensorN b_i} \\
&=\sqbr{e_1\hatr{a_1 \tensorN \cdots \tensorN a_{i+1}}} 
  \tensorN \hatr{b_1 \tensorN \cdots \tensorN b_i} .
\end{align*}
\item{(iii)} For $x\in \br{M' \cap M_k}_+$,
\[
e_1xe_1
=J_ie_{2i+1} J_ixJ_i e_{2i+1}J_i
=J_i E_{M_{2i-1}}\br{J_ixJ_i}e_{2i+1} J_i
=T^{M'}_{M_1'}(x) e_1 .
\]
\end{description}
\end{proof}
\end{lemma}

\Proof {\it of Prop~\ref{prop: approx extremal implies rotation}.}
\begin{description}
\item{(i)} First consider $k=2l+1$.  Then, for $x \in N' \cap
\gothn_{\Trace_k}$,
\begin{align*}
\left|\left| \sum_{i=r}^s R_{b_i^*} \br{ L_{b_i^*} }^* x \right|\right|^2
&= \sum_{i=r}^s \sum_{j=r}^s
	\ip{ \br{R_{b_j^*}}^* R_{b_i^*} \br{L_{b_i^*}}^*x,
	  \br{L_{b_j^*}}^*x } \\
&= \sum_{i=r}^s \sum_{j=r}^s
	\ip{ \br{\br{L_{b_i^*}}^*x}\cdot E_N(b_i^*b_j),
	  \br{L_{b_j^*}}^*x } \\
&\leq \sum_{i=r}^s
	\ip{ \br{L_{b_i^*}}^*x,
	  \br{L_{b_j^*}}^*x } ,
\end{align*}
because $[E_N(b_i^* b_j)]_{i,j=r, \ldots, s} \in M_{s-r+1}(N)$ is
dominated by $1=\delta_{i,j}$ (basically because the infinite matrix
$[E_N(b_i^* b_j)]$ is a projection).  Hence
\begin{align*}
\left|\left| \sum_{i=r}^s R_{b_i^*} \br{ L_{b_i^*} }^* x \right|\right|^2
&\leq \sum_{i=r}^s ||\br{L_{b_i^*}}^*x ||^2 \\
&=    \sum_{i=r}^s ||\br{b_i^* L_1}^*x ||^2 \\
&=    \sum_{i=r}^s \ip{ L_1\br{L_1}^* b_i x, b_i x } \\
&=    \sum_{i=r}^s \ip{ e_1 b_i x, b_i x } 
      &\text{by Lemma~\ref{lem: stuff about e_1}} \\
&=    \sum_{i=r}^s \Trace_k\br{ e_1 b_i xx^* b_i^* e_1 } .
\end{align*}
$\sum_{i=1}^s b_i xx^* b_i^* \nearrow T^{N'}_{M'}(xx^*)$ and hence
$\sum_{i=1}^s e_1 b_i xx^* b_i^* e_1 \nearrow
E^{M'}_{M_1'}\br{T^{N'}_{M'}(xx^*)} e_1$, by Lemma~\ref{lem: stuff
about e_1}.  Thus,
\begin{align}
\sum_{i=1}^s \Trace_k\br{ e_1 b_i xx^* b_i^* e_1 }
&\nearrow \Trace_k\br{ e_1 T^{N'}_{M_1'}(xx^*) } \nonumber \\
&\leq C \Trace'_k\br{ e_1 T^{N'}_{M_1'}(xx^*) } \nonumber \\
&= C\Trace_{2l+1}\br{j_l\br{e_1 T^{N'}_{M_1'}(xx^*)}} \nonumber \\
&= C\Trace_{2l+1}\br{j_l\br{T^{N'}_{M_1'}(xx^*)}e_{2l+1}} 
   &\text{by Lemma~\ref{lem: stuff about e_1}} \nonumber \\
&= C\Trace_{2l-1}\br{j_l\br{T^{N'}_{M_1'}(xx^*)}} 
   &\text{by Lemma~\ref{lem: stepping up and down traces}} \nonumber \\
&= C\Trace_{M_1'\cap\B(\Ltwo(M_l))}\br{T^{N'}_{M_1'}(xx^*)}
   &\text{by Corollary~\ref{cor: traces consistent}} \nonumber \\
&= C\Trace_{N'\cap\B(\Ltwo(M_l))}(xx^*) \nonumber \\
&= C \Trace'_k(xx^*) \nonumber \\
&\leq C^2 \Trace_k(xx^*) \nonumber \\
&= C^2 ||x||_2^2 . \label{eq: rotn existence 1}
\end{align}
Hence $\sum_{i=r}^s \Trace_k\br{ e_1 b_i xx^* b_i^* e_1 } \rightarrow
0$ and so $\{ \sum_{i=1}^s R_{b_i^*} \br{ L_{b_i^*} }^* \hatt{x} \}$
is Cauchy and hence converges.  In addition $||\sum_{i=1}^\infty
R_{b_i^*} \br{ L_{b_i^*} }^* \hatt{x}|| \leq C||\hatt{x}||$, so that
$\rho^B_k$ exists and $||\rho^B_k|| \leq C$.

\vspace{4mm}

For $k=2l$ we begin as above. For $x \in N' \cap \gothn_{\Trace_k}$,
\[
\left|\left| \sum_{i=r}^s R_{b_i^*} \br{ L_{b_i^*} }^* x \right|\right|^2
\leq \sum_{i=r}^s \Trace_k\br{ e_1 b_i xx^* b_i^* e_1 } .
\]
By (\ref{eq: rotn existence 1}) and Lemma~\ref{lem: stepping up and
down traces},
\begin{align*}
\sum_{i=1}^s \Trace_k\br{ e_1 b_i xx^* b_i^* e_1 }
&= \sum_{i=1}^s \Trace_k\br{ e_1 b_i x e_{2l+1} x^* b_i^* e_1 } \\
&\leq C^2 \Trace_{2l+1}\br{ e_{2l+1}xx^*e_{2l+1} } \\
&= C^2 \Trace_{2l}(xx^*) ,
\end{align*}
and the remainder of the argument proceeds exactly as in the $k=2l+1$
case.

\item{(ii)} If $N \subset M$ is extremal then $C=1$ so that
$||\rho_k||\leq 1$.  As $\rho_k$ is periodic this implies that
$\rho_k$ is a unitary operator.
\end{description}
\proofend

In order to establish the converse result we connect $\rho_k$ to
$J\spdot J$.

\begin{prop}
\label{prop: rotn implies 1/2 turn = Jx^*J}
\begin{description}
\item{}
\item{(i)} If either $\rho^B_{2k-1}$ exists or $\rho^B_{2k}$ exists then
$\sigma: \Ltwo(M_{2k-1}) \rightarrow N' \cap \Ltwo(M_{2k-1})$ defined by
\[
y_1 \tensorN \cdots \tensorN y_{2k} \mapsto
 P_c\br{ y_{k+1} \tensorN y_{k+2} \tensorN \cdots \tensorN y_{2k}
           \tensorN y_1 \tensorN y_2 \tensorN \cdots \tensorN y_k } ,
\]
exists (extends to a bounded operator, also denoted $\sigma$).
\item{(ii)} In that case let $\mu=\br{ \sigma |_{N' \cap
\Ltwo(M_{2k-1})} }^*$.  Then $\mu\br{ N' \cap
\gothn_{\Trace_{2k-1}} } = N' \cap \gothn_{\Trace_{2k-1}}$ and
\begin{align*}
&\mu(x)= J_{k-1} x^* J_{k-1}
& \text{for } x \in N' \cap \gothn_{\Trace_{2k-1}}
\end{align*}
\end{description}
\begin{proof}

\begin{description}
\item
\item{(i)} If $\rho^B_{2k-1}$ exists then by (\ref{eq: rho inner
prod}) in the proof of Theorem~\ref{thm: rho, tilde(rho)},
$\sigma=\br{ \rho_{2k-1}^k \circ P_c }^*$.  If $\rho_{2k}$ exists
then, for $x \in N' \cap \gothn_{\Trace_{2k-1}}$ and $y=\sum_i
y_1^{(i)} \tensorlN \cdots \tensorlN y_{2k}^{(i)}$,
\begin{align*}
&\ip{x, \sum_i P_c \br{ y_{k+1}^{(i)} \tensorN \cdots \tensorN
             y_{2k}^{(i)} \tensorN y_1^{(i)} \tensorN \cdots \tensorN
             y_k^{(i)} } }_{\Ltwo(M_{2k-1})} \\
&=\ip{x,\sum_i y_{k+1}^{(i)} \tensorN \cdots \tensorN
             y_{2k}^{(i)} \tensorN y_1^{(i)} \tensorN \cdots \tensorN
             y_k^{(i)} }_{\Ltwo(M_{2k-1})} \\
&=\ip{x,\sum_i y_{k+1}^{(i)} \tensorN \cdots \tensorN
             y_{2k}^{(i)} \tensorN y_1^{(i)} \tensorN \cdots \tensorN
             y_k^{(i)} }_{\Ltwo(M_{2k})} \\
&=\ip{x,\sum_i y_{k+1}^{(i)} \tensorN \cdots \tensorN y_{2k}^{(i)} \tensorN 
             1 \tensorN y_1^{(i)} \tensorN \cdots \tensorN
             y_k^{(i)} }_{\Ltwo(M_{2k})} \\
&=\ip{\rho_{2k}^k(x), \sum_i 1 \tensorN y_1^{(i)} \tensorN \cdots \tensorN 
             y_{2k}^{(i)} } \\
&=\ip{\br{ L_1}^* \rho_{2k}^k(x), y} .
\end{align*}
Letting $z=\sum_i P_c \br{ y_{k+1}^{(i)} \tensorN \cdots \tensorN
             y_{2k}^{(i)} \tensorN y_1^{(i)} \tensorN \cdots \tensorN
             y_k^{(i)} }$,
\begin{align*}
|| z ||_2
&=\sup \{ |<x,z>| : x \in N' \cap \Ltwo(M_{2k-1}), ||x||_2=1 \} \\
&\leq \left|\left| \br{ L_1}^* \rho_{2k}^k \right|\right| ||y||_2 \\
&\leq \left|\left|\rho_{2k}\right|\right|^k ||y||_2 .
\end{align*}
Hence $\sigma$ exists.

\item{(ii)} We proceed in four steps:
\begin{description}
\item{(1)} For $y \in \gothn_{\Trace_{2k-1}}$ and for $\xi, \eta \in
\Ltwo(M_{k-1})$ satisfying either $\xi \in D(\Ltwo(M_{k-1})_N)$ or
$\eta \in D(\Ltwo(M_{k-1})_N)$, we have
\[
\ip{\hatt{y}, \xi \tensorN J_{k-1} \eta }_{\Ltwo(M_{2k-1})}
=\ip{y \eta, \xi }_{\Ltwo(M_{k-1})} .
\]
This is easily established by first taking both $\xi, \eta \in
D(\Ltwo(M_{k-1})_N)$, and the general result then follows by
continuity.  Using Lemma~\ref{lemma: Trace = norm},
\begin{align*}
\ip{ \hatt{y} , \xi \tensorN J_{k-1} \eta }_{\Ltwo(M_{2k-1})}
&= \ip{ y , L(\xi)L(\eta)^* }_{\Ltwo(M_{2k-1})} \\
&= \Trace_{2k-1}\br{y L(\eta) L(\xi)^* } \\
&= \Trace_{2k-1}\br{L(y\eta) L(\xi)^* } \\
&= \ip{y \eta, \xi} .
\end{align*}

\item{(2)} For $y \in N' \cap \gothn_{\Trace_{2k-1}}$ satisfying
$J_{k-1} y^* J_{k-1} \in \gothn_{\Trace_{2k-1}}$ we have
\[
\mu(\hatt{y}) = J_{k-1} y^* J_{k-1} .
\]
To prove this first note that for $\xi \in N' \cap \Ltwo(M_{2k-1})$
and $\eta, \zeta \in D(\Ltwo(M_{k-1})_N)$ or $\eta, \zeta \in
D(\leftidx{_N}{\Ltwo(M_{k-1})})$,
\[
\ip{ \mu(\xi) , \eta \tensorN \zeta }
= \ip{ \xi , \zeta \tensorN \eta } .
\]
The result is true for $\eta=\hatt{a}$, $\zeta=\hatt{b}$, where $a, b
\in \wtilde{M}_{k-1}$ and the general result follows by density.

Now, using the result from (1),
\begin{align*}
\ip{ J_{k-1} y^* J_{k-1} , \eta \tensorN \zeta }_{2k-1}
&= \ip{ J_{k-1} y^* J_{k-1} J_{k-1} \zeta , \eta }_{k-1} \\
&= \ip{ J_{k-1} \eta , y^* \zeta }_{k-1} \\
&= \ip{ y J_{k-1} \eta , \zeta }_{k-1} \\
&= \ip{ y , \zeta \tensorN \eta }_{2k-1} \\
&= \ip{ \mu(\hatt{y}) , \eta \tensorN \zeta }_{2k-1} .
\end{align*}
Hence $ J_{k-1} y^* J_{k-1} =\mu(\hatt{y})$.

\item{(3)} Let $p$ be the maximal central projection in $N' \cap
M_{2k-1}$ such that $\Trace_{2k-1}$ is semifinite on $p(N' \cap
M_{2k-1})$.  Let $A=p(N' \cap M_{2k-1})$.  We will show that if $x \in
N' \cap \gothn_{\Trace_{2k-1}}$ then $\xi=\mu(\hatt{x})$ is a bounded
vector in $A$ (i.e. in $D(\Ltwo(A)_A)$) and hence an element of $A$.

Let $a \in N' \cap \gothn_{\Trace_{2k-1}}$ and let $y=J_{k-1}a^*
J_{k-1} x$.  Note that $y$ satisfies the conditions of (2).  Now, with
all supremums taken over $\sum \eta_i \tensorlN \zeta_i$ such that
$\left|\left|\sum \eta_i \tensorlN \zeta_i \right|\right|=1$, we have:
\begin{align*}
|| \xi a ||
&= \sup \left| \ip{ \xi a , \sum \eta_i \tensorN \zeta_i } \right| \\
&= \sup \left| \sum \ip{ \xi , 
     J_{2k-1} a J_{2k-1} \eta_i \tensorN \zeta_i } \right| \\
&= \sup \left| \sum \ip{ \xi , 
     J_{2k-1} \br{a J_{k-1} \zeta_i \tensorN J_{k-1}\eta_i} } \right| \\
&= \sup \left| \sum \ip{ \xi , 
     J_{2k-1} \br{\br{a J_{k-1} \zeta_i} \tensorN J_{k-1}\eta_i} } \right| \\
&= \sup \left| \sum \ip{ \xi , 
     \eta_i \tensorN J_{k-1} a J_{k-1} \zeta_i } \right| \\
&= \sup \left| \sum \ip{ x , 
     J_{k-1} a J_{k-1} \zeta_i \tensorN \eta_i } \right| \\
&= \sup \left| \sum \ip{ J_{k-1} a^* J_{k-1} x , 
      \zeta_i \tensorN \eta_i } \right| \\
&= \sup \left| \sum \ip{ y , 
      \sigma\br{\eta_i \tensorN \zeta_i } } \right| \\
&= \sup \left| \sum \ip{ \mu(y) , 
      \eta_i \tensorN \zeta_i } \right| \\
&= \left|\left| \mu(y) \right|\right|_{\Ltwo(M_{2k-1})} \\
&= \left|\left| J_{k-1} y^* J_{k-1} \right|\right|_2 \\
&= \left|\left| a^* J_{k-1} x J_{k-1} \right|\right|_2 \\
&\leq ||a||_2 ||x|| .
\end{align*}
Hence $\xi$ is a bounded vector and thus an element of $A$.


\item{(4)} Finally, let $z \in N' \cap \gothn_{\Trace_{2k-1}}$ such
that $\mu(\hatt{x})=\hatt{z}$.  Then $z=J_{k-1} x^* J_{k-1}$ because
\begin{align*}
\ip{ J_{k-1} x^* J_{k-1} \eta, \zeta }
&= \ip{ x J_{k-1} \zeta , J_{k-1} \eta } \\
&= \ip{ x , J_{k-1} \eta \tensorN \zeta } \\
&= \ip{ \mu(x) , \zeta \tensorN J_{k-1} \eta } \\
&= \ip{ z , \zeta \tensorN J_{k-1} \eta } \\
&= <z \eta , \zeta > .
\end{align*}
\end{description}
\end{description}
\end{proof}
\end{prop}

\begin{prop}
\label{prop: rotn implies approx ext}
\begin{description}
\item{(i)} If $\rho^B_i$ exists then $N \subset M$ is approximately
extremal.
\item{(ii)} If $\rho^B_i$ exists and is a unitary operator then $N
\subset M$ is extremal.
\end{description}
\begin{proof}

\begin{description}
\item{(i)} Let $j$ be the largest odd number with $j \leq i$.  By
Prop~\ref{prop: rotn implies 1/2 turn = Jx^*J} there exists $\mu:N'
\cap \Ltwo(M_j) \rightarrow N' \cap \Ltwo(M_j)$ with $||\mu|| \leq
\left|\left| \rho_i \right|\right|^k$ (where $j=2k-1$) and
$\mu(x)=J_{k-1} x^* J_{k-1}$ for $x \in N' \cap
\gothn_{\Trace_j}$.

We first show that $\Trace'_j \leq ||\mu||^2 \Trace_j$.  Take $x \in
N' \cap M_j$.  If $\Trace_j(x^*x)=\infty$ then we are done.  Otherwise
$x \in \gothn_{\Trace_j}$ and
\begin{align*}
\Trace'_j(x^8x)
&= \Trace_j \br{ J_{k-1} x^* J_{k-1} J_{k-1} x J_{k-1} } \\
&= \Trace_j ( \mu(x) \mu(x)^*) \\
&\leq ||\mu||^2 ||x||_2^2 \\
&= ||\mu||^2 \Trace_j(x^*x) .
\end{align*}

Finally, 
$\Trace_j
=\Trace'_j\br{J_{k-1}\spdot J_{k-1}}
\leq ||\mu||^2 \Trace_j\br{J_{k-1}\spdot J_{k-1}}
= ||\mu||^2 \Trace'_j$.

\item{(ii)} If $\rho_i$ is unitary then $||\mu|| \leq 1$ so that
$\Trace_j \leq \Trace'_j \leq \Trace_j$ and hence $\Trace_j =
\Trace'_j$.
\end{description}
\end{proof}
\end{prop}

\begin{cor}
\label{cor: rotn for one basis implies rotn for all}
If there exists a basis $B$ such that $\rho^B_k$ exists then for any
basis $\overline{B}$, $\rho^{\overline{B}}_k$ exists and is
independent of the basis used.  Hence we will use $\rho_k$ to denote
$\rho^B_k$.
\begin{proof}

Suppose $\rho^B_k$ exists.  By Prop~\ref{prop: rotn implies approx
ext} the subfactor is approximately extremal.  By Prop~\ref{prop:
approx extremal implies rotation} $\rho^{\overline{B}}_k$ exists for
any basis $\overline{B}$.  By Theorem~\ref{thm: rho, tilde(rho)}
$\rho_k$ is independent of the basis used.
\end{proof}
\end{cor}

The results of Sections \ref{sect: extremality} through \ref{sect:
rotations and extremality} can be summarized as follows.

\begin{thm}
\begin{description}
\item{}
\item{(i)}The following are equivalent:
\begin{itemize}
\item $N \subset M$ is approximately extremal;
\item $\Trace'_{2i+1} \sim \Trace_{2i+1}$ for all $i \geq 0$;
\item $\Trace'_{2i+1} \sim \Trace_{2i+1}$ for some $i \geq 0$;
\item $\rho_k$ exists for all $k \geq 0$;
\item $\rho_k$ exists for some $k \geq 1$;
\item $\widetilde{\rho}_k$ exists for all $k \geq 0$;
\item $\widetilde{\rho}_k$ exists for some $k \geq 1$.
\end{itemize}
\item{(ii)} The following are equivalent:
\begin{itemize}
\item $N \subset M$ is extremal;
\item $\Trace'_{2i+1}=\Trace_{2i+1}$ for all $i \geq 0$;
\item $\Trace'_{2i+1}=\Trace_{2i+1}$ for some $i \geq 0$;
\item $\rho_k$ (or $\widetilde{\rho}_k$) exists for all $k \geq 0$ and
  is unitary
\item $\rho_k$ (or $\widetilde{\rho}_k$) exists for all $k \geq 0$ and
  $\rho=\widetilde{\rho}$;
\item $\rho_k$ (or $\widetilde{\rho}_k$) exists for some $k \geq 1$
  and is unitary
\item $\rho_k$ (or $\widetilde{\rho}_k$) exists for some $k \geq 1$
  and $\rho=\widetilde{\rho}$;
\end{itemize}
\end{description}
\end{thm}

\newpage


\section{A $\IIone$ Subfactor With Type $\III$ Component in a Relative Commutant}
\label{type III rel comm}

Here we construct an infinite index type $\IIone$ subfactor $N\subset
M$ such that $N'\cap M_1$ has a type $\rm III$ central summand.


\subsection{Outline}

Let $\Finfinity$ denote the free group on infinitely many generators.
We take a suitably chosen type $\rm III$ factor representation
$w:\Finfinity\rightarrow {\mathcal U}({\mathcal H}_0)$ on a separable
Hilbert space ${\mathcal H}_0$.  By tensoring this representation with
the trivial representation on $l^2(\mathbb N)$ we may assume that if
$w_{\gamma}$ is a Hilbert-Schmidt perturbation of the identity then
$w_{\gamma}$ must be the identity.

The corresponding Bogoliubov automorphisms of the canonical
anti-commutation relation algebra $A=CAR({\mathcal H}_0)$ provide an
action of $\Finfinity$ on $A$.  This action passes to an action
$\alpha$ on the hyperfinite $\IIone$ factor $R=\pi(A)''$ obtained via
the GNS representation $\pi$ on ${\mathcal H}={\rm L}^2(A,tr)$.

We conclude from Blattner's Theorem, characterizing inner Bogoliubov
automorphisms, that for each $\gamma\in\Finfinity$ either
$\alpha_{\gamma}$ is outer or $\alpha_{\gamma}=\rm id$.

Now we construct $M=R\rtimes_{\alpha} \Finfinity$ and $N={\mathbb
C}\rtimes_{\alpha} \Finfinity \cong vN(\Finfinity)$, the von Neumann
algebra of the left regular representation of $\Finfinity$.  We show
that $M$ is a $\IIone$ factor and that the basic construction for
$N\subset M$ yields $M_1={\mathcal B(\mathcal
H)}\rtimes_{\alpha}\Finfinity \cong {\mathcal B(\mathcal H)} \otimes
vN(\Finfinity)$, a ${\rm II}_{\infty}$ factor.

We show that by virtue of our choice of representation $w$ we
have $N'\cap M_1={\mathcal B(H)}^{\Finfinity}$ and finally we show
that ${\mathcal B(H)}^{\Finfinity}$, while not a factor, has a type
$\rm III$ central summand.

\subsection{Preliminary results}

We will use the following basic facts about von Neumann algebras:

\begin{lemma}
\label{dense unitaries}
Let $S$ be a von Neumann algebra with separable predual.  Then there
exists a countable set $\Lambda \subset {\mathcal U}(S)$ such that
$\Lambda''=S$.
\begin{proof}
Since $S_*$ has separable predual there exists a countable dense
subset $\{ \phi_n\}\subset (S_*)_1$.  Now the weak-$*$ topology on
$(S)_1$ (recall $S$ is the dual of $S_*$) is metrizable by
$d(x,y)=\sum_{n\in\mathbb{N}} 2^{-n}|\phi_n(x-y)|$, $x,y\in (S)_1$.
The $\sigma$-weak topology on $(S)_1$ coincides with the weak-$*$
topology we have shown that $(S)_1$, so that $(S)_1$ with the
$\sigma$-weak topology is not only compact, but also metrizable, and
hence separable.


Take $\{ x_i\}$ $\sigma$-weakly dense in $(S)_1$.  Write each $x_i$ as
a linear combination of four unitary operators and let $\Lambda$ be
the set of all these unitary operators.  Then $\Lambda''=S$.
\end{proof}
\end{lemma}

\begin{lemma}
\label{cutdown}
Let $Q$ be a von Neumann algebra, $p\in Q$ such that $pQ'$ is a
factor.  Let $z(p)$ denote the central support of $p$.  Then $z(p)Q'$
is also a factor.
\begin{proof}
Let $q=z(p)$.  If $qQ'$ is not a factor then there exist projections
$q_1, q_2\in Z(qQ')=qQ'\cap qQ=qZ(Q)$ such that $q_i\neq 0$ and
$q_1+q_2=q$.  Let $p_i=pq_i\in pQ'$ and note that $p_i=pq_ip\in pQp$
so $p_i\in pQ'\cap pQp=pQ' \cap (pQ')'={\mathbb C}p$.  Hence $p_i=0$
or $p_i=p$.  WLOG $p_1=p$ and $p_2=0$, but then $p\leq q_1$ so $q$ is
not the central support of $p$.  Hence $qQ'$ is a factor.
\end{proof}
\end{lemma}

\subsubsection{CAR algebra and Bogoliubov automorphisms}
\label{CAR and Bogoliubov basics}

Let $\mathcal H$ be a complex Hilbert space and let ${\mathcal
F}(\mathcal H)=\oplus_{n\geq 0}\wedge^n{\mathcal H}$ be the
anti-symmetric Fock space of $\mathcal H$.  The canonical
anticommutation relation algebra $A=CAR(\mathcal H)$ is the
\cstar-algebra generated by the creation and annihilation operators
$a(\xi), a^*(\xi)$ ($\xi\in\mathcal H$).  $A$ has a unique tracial
state $\trace$, namely the quasi-free state of covariance
$\frac{1}{2}$.  We have a representation $\pi$ of $A$ on
$\Ltwo(A,\trace)$ by left multiplication.  It is well known that
$\pi(A)''=R$ the hyperfinite $\IIone$ factor.

Each unitary operator $u\in{\mathcal U}(\mathcal H)$ gives rise to an
automorphism of the CAR algebra via $a(\xi)\rightarrow a(u\xi)$.  We
call this the Bogoliubov automorphism induced by $u$ and denote it
$Bog(u):A\rightarrow A$.

By uniqueness of the tracial state on $A$, any automorphism $\alpha$
of $A$ defines a unitary operator $W$ on $\Ltwo(A,\trace)$ by
$W\hatt{x}=\alpha(x)$ for $x\in A$.  This unitary operator implements
the automorphism $\alpha$ in the representation $\pi$:\\
$W\pi(x)W^*\hatt{y}=\hatr{\alpha(x\alpha^{-1}(y))}=\hatr{\alpha(x)y}
=\pi(\alpha(x))\hatt{y}$.  The automorphism $\alpha$ can thus be
extended to an automorphism $\wtilde{\alpha}$ of $R=\pi(A)''$ by
$\wtilde{\alpha}(x)=WxW^*$.  For a Bogoliubov automorphism $Bog(u)$ we
will also refer to $(Bog(u))\wtilde{\phantom{A}}$ as the Bogoliubov
automorphism induced by $u$ and denote it by $\alpha_u$.

A theorem of Blattner~\cite{Blattner1958} characterizes the inner
Bogoliubov automorphisms.  The theorem is usually stated in terms of
real Hilbert spaces, so we will briefly review the construction of the
CAR algebra via a real Hilbert space and the Clifford algebra over it.
For full details see for example de~la~Harpe and
Plymen~\cite{delaHarpePlymen1979}.

Let $E$ be a real Hilbert space, $Cl(E)$ the Clifford algebra of the
quadratic form $q(\xi)=<\xi|\xi>=||\xi||^2$ and $Cl(E)^{\mathbb
C}=Cl(E)\otimes_{\mathbb R}{\mathbb C}$ its complexification.
$Cl(E)^{\mathbb C}$ has a unique tracial ``state'' $\trace$.  The
\cstar-algebra generated by the left representation of $Cl(E)^{\mathbb
C}$ on its Hilbert space completion is called the CAR algebra over
$E$.  As before, the von Neumann algebra generated by $Cl(E)^{\mathbb
C}$ is the hyperfinite $\IIone$ factor $R$ and we have the following
theorem:

\begin{thm}[Blattner~\cite{Blattner1958}]
Let $v\in {\mathcal O}(E)$, the orthogonal group of $E$.  Then
$\alpha_v$ is inner iff either (i) the eigenspace corresponding to
eigenvalue $-1$ is even dimensional (or infinite) and $v$ is a
Hilbert-Schmidt perturbation of the identity; or (ii) the eigenspace
corresponding to eigenvalue $1$ is (finite) odd dimensional and $-v$
is a Hilbert-Schmidt perturbation of the identity.
\end{thm}

Given a complex Hilbert space $\mathcal H$ we may construct
$CAR(\mathcal H)$ as the CAR algebra over $E={\mathcal H}_{\mathbb R}$
using the Clifford algebra approach -- see~\cite{JonesdelaHarpe1995}
for details.  When considered as an operator on ${\mathcal H}_{\mathbb
R}$, a unitary operator $u\in{\mathcal U}(\mathcal H)$ is clearly
orthogonal.  The eigenspaces of $u$ are always even dimensional (if
$u\xi=\lambda \xi$ then $u(i\xi)=\lambda(i\xi)$) and if $u=1+x$ for
some $x\in{\mathcal B}({\mathcal H}_{\mathbb R})$ then $x$ must be
$\mathbb C$-linear since both $u$ and $1$ are $\mathbb C$-linear.  If
$x$ is a Hilbert-Schmidt operator on ${\mathcal H}_{\mathbb R}$ then
$x$ is also a Hilbert-Schmidt operator on $\mathcal H$ (if $\{\xi_n\}$
is an orthonormal basis for $\mathcal H$ then $\{\xi_n\}\cup\{ i\xi\}$
is an orthonormal basis for ${\mathcal H}_{\mathbb R}$ and
$||x||_{{\mathbb R}-HS}^2=\sum (||x\xi_n||^2+
||x(i\xi_n)||^2)=2||x||_{{\mathbb C}-HS}^2$).  Thus, in the case of
Bogoliubov automorphisms on the CAR algebra of a complex Hilbert space
$\mathcal H$ constructed via creation and annihilation operators on
anti-symmetric Fock space, Blattner's Theorem may be stated as:

\begin{thm}
\label{Blattner}
Let $\mathcal H$ be a complex Hilbert space and $u\in {\mathcal
U}({\mathcal H})$.  Then $\alpha_u$ is inner iff $u$ is a
Hilbert-Schmidt perturbation of the identity.
\end{thm}

\subsection{The construction}

Let $\Finfinity=<a_n>_{n=1}^\infty$ denote the free group on countably
many generators, with $\{ a_n \}_{n=1}^\infty$ a particular choice of
generators.  Fix a bijection $\phi:\{ a_n \}_{n=1}^{\infty}
\rightarrow \mathbb{N} \times \Finfinity$.  Let $p:\mathbb{N} \times
\Finfinity \rightarrow \Finfinity$ be projection onto the second
component, $p(n,\gamma)=\gamma$, and let $\varphi=p\circ\phi$.  Since
$\Finfinity$ is free, $\varphi$ extends to a homomorphism
$\overline{\varphi}:\Finfinity\rightarrow \Finfinity$.

Let $S$ be a type $\rm III$ factor acting on a separable Hilbert space
${\mathcal H}_0$.  By Lemma~\ref{dense unitaries} there exists a
countable subset $\{ v_n \}_{n=1}^\infty$ of ${\mathcal U}(S)$ which
is dense in the sense that it generates $S$ as a von Neumann algebra.

Define a homomorphism $\psi:\Finfinity\rightarrow\mathcal {U}(\h_0)$
by letting $\psi(a_n)=v_n$.  Then define a representation $w$ of
$\Finfinity$ on $\h_0$ by $w=\psi\circ\overline{\phi}$.  Note that
$w(\Finfinity)''=\{v_n\}''=S$.

Let $A=CAR({\mathcal H}_0)$, $R=\pi(A)''$ where $\pi$ is the GNS
representation of $A$ on $\h=\Ltwo(A,\trace)=\Ltwo(R,\trace)$.  As we
saw in section~\ref{CAR and Bogoliubov basics} each $w_\gamma$ induces
a unitary operator $W_\gamma$ on $\mathcal H$ and a Bogoliubov
automorphism of $R$ which we will denote $\alpha_\gamma
(=\alpha_{w_\gamma}={\rm Ad}(W_\gamma))$.  Thus we have a
representation $W: \Finfinity \rightarrow \mathcal H$ and an action
$\alpha_\gamma ={\rm Ad}(W_\gamma)$ on ${\mathcal B}(\mathcal H)$
which restricts to an action on $R$.

Replacing ${\mathcal H}_0$ and $w$ with ${\mathcal H}_0\otimes
l^2(\mathbb N)$ and $w\otimes 1$ respectively, we may assume that for
each $\gamma\in\Finfinity$ either $w_\gamma$ is the identity or
$w_\gamma$ is not a Hilbert-Schmidt perturbation of the identity.  By
Blattner's Theorem,~(\ref{Blattner}), $\alpha_\gamma|_R$ is either outer or
trivial.

Define:
\begin{eqnarray*}
N & = & {\mathbb C}\rtimes_{\alpha} \Finfinity \cong vN(\Finfinity) \\
M & = & R \rtimes_{\alpha} \Finfinity \\
M_1 & = & {\mathcal B}({\mathcal H})\rtimes_{\alpha} \Finfinity
	\cong {\mathcal B(\mathcal H)} \otimes vN(\Finfinity)
\end{eqnarray*}
We will show that $N \subset M$ has all the desired properties stated
in the introduction.

\begin{lemma}
$M$ is a $\IIone$ factor.
\begin{proof}
Let $x=\sum x_\gamma u_\gamma \in Z(M)$.  Then for all $y\in R$ we
have
\[
\sum_\gamma y x_\gamma u_\gamma 
=yx
=xy
=\sum_\gamma x_\gamma \alpha_\gamma(y) u_\gamma
\]
which yields $y x_\gamma=x_\gamma \alpha_\gamma(y)$ for all $y\in R$
and all $\gamma\in\Finfinity$.  Now, recall that either
$\alpha_\gamma$ is outer, in which case $x_\gamma=0$ since an outer
action on a factor is automatically free, or $\alpha_\gamma={\rm id}$
in which case $x_\gamma\in Z(R)={\mathbb C}$.  Thus every $x_\gamma$
is a scalar and so $x\in N$.  But since $x\in Z(M)$ we have $x\in N
\cap M' \subseteq N \cap N'={\mathbb C}$.  Thus $M$ is a factor.  $M$
is infinite dimensional and has trace $\trace\br{\sum_\gamma x_\gamma
u_\gamma}=\trace_R\br{x_e}$, so $M$ is a $\IIone$ factor.
\end{proof}
\end{lemma}

\begin{lemma}
$N\subset M \subset M_1$ is a basic construction.
\begin{proof}
In essence this is true because $\mathbb{C} \subset R \subset
\B(\h)=\B(\Ltwo(R,\trace))$ is a basic construction ($J_R
(\mathbb{C}1)'J_R=\B(\Ltwo(R))$).  Let $\overline{e}_1$ be the Jones
projection for $\mathbb{C} \subset R$ and note that this is just the
extension of $\trace_R$ to $\Ltwo(R)$.

Note that $\Ltwo(M)\isom \Ltwo(R) \tensor l^2(\Finfinity)$ via
$U:\hatr{\sum_\gamma x_\gamma u_\gamma} \mapsto \sum \hatt{x_\gamma}
\tensor \delta_\gamma$ and $U\Ltwo(N)=\mathbb{C} \tensor
l^2(\Finfinity)$, so $e_1=\overline{e}_1 \tensor 1$.

Let $\pi=U \spdot U^*$ be the representation of $M$ on $\Ltwo(R)
\tensor l^2(\Finfinity)=\h \tensor l^2(\Finfinity)$.  $\pi(x)=x
\tensor \id$ for $x\in R$ and $\pi\br{ u_\gamma}=w_\gamma \tensor
\lambda_\gamma$ where $\lambda$ is the left regular representation.
$M_1$ is given by
\[
M_1
=\br{\pi(M) \cup \{ e_1 \}}''
=\br{\br{\pi(R) \cup \{ e_1 \}}'' \cup \{ \pi(u_\gamma) \} }''
\]
and
\begin{align*}
\br{\pi(R) \cup \{ e_1 \}}''
&= \br{\br{R\tensor \mathbb{C}} \cup \{ \overline{e}_1 \tensor\id \}}'' \\
&= \br{R \cup \{ \overline{e}_1 \}}'' \tensor \mathbb{C} \\
&= \B(\Ltwo(R)) \tensor \mathbb{C} .
\end{align*}
Finally, for $x\in\B(\Ltwo(R))=\B(\h)$, $\pi(u_\gamma)(x\tensor
1)\pi(u_\gamma)^*=w_\gamma x w_\gamma^* \tensor 1=\alpha_\gamma(x)
\tensor 1$.  Thus $M_1= \B(\h) \rtimes_{\alpha} \Finfinity$ and every
element of $M_1$ can be written as $\sum_\gamma \br{x_\gamma\tensor 1}
\pi(u_\gamma)$.
\end{proof}
\end{lemma}

\begin{lemma}
$N'\cap M_1=\B(\h)^{\Finfinity}=W(\Finfinity)'$.
\begin{proof}
Suppose $x=\sum x_\gamma u_\gamma \in N'\cap M_1$, $x_\gamma\in
\B(\h)$.  Then for all $\rho \in \Finfinity$ we have
\[
\sum x_\gamma u_\gamma = x = u_\rho x u_\rho^{-1}
	=\sum \alpha_{\rho}(x_\gamma) u_{\rho\gamma\rho^{-1}}
	=\sum \alpha_{\rho}(x_{\rho^{-1}\gamma\rho}) u_\gamma
\]
which yields $\alpha_\rho(x_{\rho^{-1}\gamma\rho})=x_\gamma$ for all
$\gamma, \rho \in \Finfinity$.  We can rewrite this as
\begin{equation}
\label{twisted constant conjugacy classes}
x_{\rho^{-1}\gamma\rho}=\alpha_{\rho^{-1}}(x_\gamma) \hspace{2cm} 
	\forall \gamma, \rho \in \Finfinity
\end{equation}
In other words, the matrix entries of $x$ are constant on conjugacy
classes modulo a twist by the action of $\Finfinity$.

Now suppose that there exists $\gamma_0\neq e$ ($e$ the identity
element of $\Finfinity$) such that $x_{\gamma_0} \neq 0$, say
$\gamma_0=a_{m_1}^{\pm 1} \ldots a_{m_2}^{\pm 1}$ in reduced form.
Recall that $\phi:\{ a_n \}_{n=1}^{\infty} \rightarrow
\mathbb{N}\times\Finfinity$ is a bijection.  Choose infinitely many
distinct positive integers $\{ n_i\}_{i\in\mathbb N}$ such that
$n_i\notin\{ m_1,m_2\}$ and
$\varphi(a_{n_i})=\overline{\varphi}(\gamma_0)$.  Then the elements
$\gamma_i=a_{n_i}\gamma_0 a_{n_i}^{-1}$ are all distinct and
\[
w_{\gamma_i}
= \psi(\overline{\varphi}(a_{n_i}\gamma_0 a_{n_i}^{-1}))
= \psi(\overline{\varphi}(\gamma_0)\overline{\varphi}(\gamma_0)
   \overline{\varphi}(\gamma_0)^{-1})
= \psi(\overline{\varphi}(\gamma_0))
= w_{\gamma_0}
\]
Hence $\alpha_{\gamma_i}=\alpha_{\gamma_0}$.

Now take $\xi\in{\mathcal H}$ such that
$\alpha_{\gamma_0}^{-2}(x_{\gamma_0})\xi\neq 0$.  Consider $\xi\otimes
\CHI_e \in {\mathcal H}\otimes l^2(\Finfinity)$.  Noting that
$\gamma_i^{-1}\gamma_0\gamma_i$ are all distinct, we have:
\begin{eqnarray*}
x(\xi\otimes\CHI_e)
	& = & (\sum x_\gamma u_\gamma) (\xi\otimes\CHI_e) \\
	& = & \sum \alpha_{\gamma^{-1}}(x_\gamma)\xi \otimes \CHI_\gamma \\
\Rightarrow ||x(\xi\otimes\CHI_e)||^2
	& = & \sum ||\alpha_{\gamma^{-1}}(x_\gamma)\xi ||^2 \\
	& \geq & \sum_i ||\alpha_{\gamma_i^{-1}\gamma_0^{-1}\gamma_i}
		(x_{\gamma_i^{-1}\gamma_0\gamma_i})\xi ||^2 \\
	& = & \sum_i ||\alpha_{\gamma_0^{-1}}
		(\alpha_{\gamma_i^{-1}}(x_{\gamma_0}))\xi ||^2
		\hspace{20mm} 
		\mbox{(from \ref{twisted constant conjugacy classes})} \\
	& = & \sum_i ||\alpha_{\gamma_0}^{-2}(x_{\gamma_0})\xi||^2 \\
	& = & \infty
\end{eqnarray*}

We conclude that $x=x_e\in {\mathcal B}({\mathcal H})$.  Now
$[x,u_\rho]=0$ iff $\alpha_\rho(x)=x$, so $N'\cap M_1={\mathcal
B}({\mathcal H})^{\Finfinity}=W(\Finfinity)'$.
\end{proof}
\end{lemma}

\begin{lemma}
$N'\cap M_1=W(\Finfinity)'$ has a central summand which is a type $\rm
III$ factor.
\begin{proof}
Let $Q=W(\Finfinity)'$.  There is a canonical embedding $\iota:\h_0
\hookrightarrow \h=\Ltwo\br{CAR\br{\h_0}}$ given by
$\iota(\xi)=\sqrt{2}a(\xi)$ (one has
$\ip{\sqrt{2}a(\xi),\sqrt{2}a(\eta)}=2\trace(a^*(\xi)a(\eta))=\ip{\xi,\eta}$
because $\trace$ is the quasi-free state of covariance $\frac{1}{2}$).
The action of $\Finfinity$ commutes with $\iota$ so that
$\iota({\mathcal H}_0)$ is invariant under $W(\Finfinity)$.  If $p$
denotes the orthogonal projection onto $\iota({\mathcal H}_0)$ then
$p\in W(\Finfinity)'=Q$.

Let $q=z(p)$, the central support of $p$ in $Q$.  Then, using
Lemma~\ref{cutdown}, $qQ'$ is also a factor and $p\in (qQ')'=qQ$, so
$qQ'\cong pqQ'=pQ'=pW(\Finfinity)''={w_\gamma}''=S$ the type $\rm III$
factor we started with.  Hence the central summand $qQ=(qQ')'$ is also
a type $\rm III$ factor.
\end{proof}
\end{lemma}

\noindent In summary:

\begin{thm}
There exist infinite index $\IIone$ subfactors $N\subset M$ such that
upon performing the basic construction $N\subset M \subset M_1$, the
relative commutant $N'\cap M_1$ has a type $\rm III$ central summand.
\end{thm}






\bibliographystyle{plain}
\bibliography{vn}

\end{document}